\documentclass{article}
\usepackage[T1]{fontenc}
\usepackage{titlesec}
\usepackage[utf8]{inputenc}
\usepackage{amsmath}
\usepackage{amssymb}
\usepackage{amsfonts}
\usepackage{geometry}
\usepackage{xcolor}
\usepackage[hidelinks]{hyperref}
\newtheorem{theorem}{Theorem}
\newtheorem{lemma}{Lemma}
\newtheorem{prop}{Proposition}
\numberwithin{equation}{section}

\title{PhD}
\author{}
\date{January 2021}

\begin{document}

\begin{center}
    \textbf{On the Eigenstructure of Covariance Matrices with Divergent Spikes}
\end{center}

\begin{center}
     Simona Diaconu\footnote{Department of Mathematics, Stanford University, sdiaconu@stanford.edu}
\end{center}

\begin{abstract}
For a generalization of Johnstone's spiked model, a covariance matrix with eigenvalues all one but \(M\) of them, the number of features \(N\) comparable to the number of samples \(n: N=N(n), M=M(n), \gamma^{-1} \leq \frac{N}{n} \leq \gamma\) where \(\gamma \in (0,\infty),\) we obtain consistency rates in the form of CLTs for separated spikes tending to infinity fast enough whenever \(M\) grows slightly slower than \(n: \lim_{n \to \infty}{\frac{\sqrt{\log{n}}}{\log{\frac{n}{M(n)}}}}=0.\) Our results fill a gap in the existing literature in which the largest range covered for the number of spikes has been \(o(n^{1/6})\) and reveal a certain degree of flexibility for the centering in these CLTs inasmuch as it can be empirical, deterministic, or a sum of both. Furthermore, we derive consistency rates of their corresponding empirical eigenvectors to their true counterparts, which turn out to depend on the relative growth of these eigenvalues.
\end{abstract}

\tableofcontents

\section{Introduction}
\par
Covariance matrices arise naturally in a myriad of disciplines (finance, statistics, physics, etc.) and have been studied for a long time due to the tight connection between them and several techniques widely used across distinct fields such as principal component analysis, a method employed in modern applications such as dimension reduction for data visualization. The problem setup is as follows: consider i.i.d. samples \(\mathbf{x}_1,\mathbf{x}_2, \hspace{0.05cm} ... \hspace{0.05cm},\mathbf{x}_n \in \mathbb{R}^{N}\) (whose entries are often called \textit{features}: for instance, they can be measurements of certain characteristics or pixel intensities) drawn from the distribution of a random vector \(\mathbf{x} \in \mathbb{R}^N\) for which \(\mathbb{E}[\mathbf{x}]=\mathbf{0}, \mathbb{E}[\mathbf{x}\mathbf{x}^T]=\Sigma \in \mathbb{R}^{N \times N}.\) The most commonly posed questions are what can be inferred about the true covariance matrix, \(\Sigma,\) from its empirical counterpart, \(\mathbf{S}_n=\frac{1}{n}\mathbf{X}_n\mathbf{X}_n^T,\) where \(\mathbf{X}_n=[\mathbf{x}_1 \hspace{0.2cm} \mathbf{x}_2 \hspace{0.2cm} ... \hspace{0.2cm} \mathbf{x}_n] \in \mathbb{R}^{N \times n}\) is the matrix whose columns are the observations, and vice versa, what can be derived about \(\mathbf{S}_n\) from \(\Sigma.\) The former question is particularly relevant in practical contexts, as the ground truth of the assumed model, in this case \(\Sigma,\) can never be observed, while the second question has given rise to a theory of oracles (loosely speaking, oracles are estimators based on model parameters). Another common assumption is that \(\mathbf{X}_n=UD^{1/2}\mathbf{Z}_n,\) where \(U \in \mathbb{R}^N\) is an orthogonal matrix yielding a spectral decomposition of \(\Sigma=UDU^T,\) and \(\mathbf{Z}_n \in \mathbb{R}^{N \times n}\) a random matrix with i.i.d. entries of mean zero and variance one.
\par
Although for \(N\) fixed and \(n\) large, \(\mathbf{S}_n\) is a consistent estimator of \(\Sigma\) since \(\mathbf{S}_n \xrightarrow[]{a.s.} \Sigma\) as \(n \to \infty\), this is no longer the case when \(N\) grows at least at the same rate as \(n\) does. In particular, the eigenvalues of \(\mathbf{S}_n\) can differ significantly from their true counterparts, both globally and locally: for instance, if \(\Sigma=I_N, N=N(n),\lim_{n \to \infty}{\frac{N(n)}{n}}= \gamma \in (0,1),\) and the i.i.d. entries of \(\mathbf{Z}_n=(z_{ij})_{1 \leq i \leq N, 1 \leq j \leq n}\) have a distribution independent of \(n,\) then the empirical spectral distribution of \(\mathbf{S}_n\) converges almost surely to a Marchenko-Pastur law that is absolutely continuous with respect to the Lebesgue measure on the real line with p.d.f. 
\[p_\gamma(x)=\begin{cases} 
\frac{\sqrt{(b-x)(x-a)}}{2\pi \gamma x}, & x \in [a,b]\\
0, & x \not \in [a,b]
\end{cases},\] 
for \(a=(1-\sqrt{\gamma})^2, b=(1+\sqrt{\gamma})^2\) (see \cite{silverstein} and the seminal paper \cite{marchenkopastur}), while Bai and Yin~\cite{baiyin} proved that if additionally \(\mathbb{E}[z^4_{11}]<\infty,\) then almost surely the minimal and maximal eigenvalues of \(\mathbf{S}_n\) converge to \(a,b,\) respectively. In contrast, the empirical distribution of the ground truth \(\Sigma=I_N\) is the Dirac delta function with mass at one and its eigenvalues are evidently all equal to one. 
\par
The relationship between the number of features, \(N,\) and the number of samples, \(n,\) has turned out to be a driving factor that determines how the sample covariance matrix \(\mathbf{S}_n\) and its true complement \(\Sigma_N\) compare, and consequently three different asymptotic regimes have been analyzed:

1. \textit{classical}: \(N\) fixed, \(n \to \infty\) (Anderson~\cite{anderson}),

2. \textit{random matrix theory}: \(\frac{N}{n} \to \gamma \in (0,\infty),\)

3. \textit{high-dimensional low-sample size}: \(\frac{N}{n} \to \infty\) (Ahn et al.~\cite{ahn}, Hall et al.~\cite{hall}, Jung et al.~\cite{jung}).
\par
The eigenvalues of \(\Sigma\) are equally of uttermost importance, and two of the most extensively treated families are the \textit{isotropic} (i.e., \(\Sigma=I_N\)), and the \textit{spiked} models: at a high level, this latter class encompasses covariance matrices for which few eigenvalues (the \textit{spikes}) differ from the rest (the \textit{bulk}). Johnstone was the first to consider such a model in \cite{johnstone}: 
\[\Sigma=diag(l_1,l_2, \hspace{0.05cm} ... \hspace{0.05cm}, l_M,1,1,\hspace{0.05cm} ... \hspace{0.05cm}, 1) \in \mathbb{R}^{N \times N},\]
where the spikes \(l_1 \geq l_2 \geq ... \geq l_M>1\) are fixed, and subsequently several of its properties in the random matrix theory regime have been discovered. Particularly, for \(\gamma<1,\) Baik and Silverstein~\cite{baiksilv} showed that if an eigenvalue \(l\) satisfies \(l \geq 1,\) then its empirical counterpart \(\hat{l}\) tends almost surely to 
\[\phi(l)=\begin{cases} 
       (1+\sqrt{\gamma})^2, & l<1+\sqrt{\gamma} \\
       l+\frac{\gamma l}{l-1}, & l \geq 1+\sqrt{\gamma}
   \end{cases},\]
while Paul~\cite{paul} obtained the fluctuations underlying this convergence: namely, if \(l>1+\sqrt{\gamma}\) has multiplicity one, \(N/n=\gamma+o(n^{-1/2}),\) then
\[\sqrt{n}(\hat{l}-\phi(l)) \Rightarrow N(0,\sigma^2(l)), \hspace{0.4cm} \sigma^2(l)=2l^2(1-\frac{\gamma}{(l-1)^2}).\]
The former results assume \(\mathbf{Z}_n \in \mathbb{C}^{N \times n}\) has i.i.d. entries with distributions independent of \(n\)  and \(\mathbb{E}[z_{11}]=0, \mathbb{E}[|z_{11}|^2]=1, \mathbb{E}[|z_{11}|^4]<\infty,\) while the latter require the columns of \(\mathbf{Z}_n \in \mathbb{R}^{N \times n}\) (and consequently of \(\mathbf{X}_n\)) have a multivariate normal distribution. Moreover, Paul~\cite{paul} considered the eigenvectors \(p_l\) corresponding to such eigenvalues \(l:\) namely, for \(p_l^T=[p_{A,l}^T \hspace{0.2cm} p_{B,l}^T], p_{A,l} \in \mathbb{R}^M, p_{B,l} \in \mathbb{R}^{N-M},\) he showed that asymptotically \(\frac{p_{A,l}}{||p_{A,l}||}\) has a normal behavior, whereas \(\frac{p_{B,l}}{||p_{B,l}||}\) is distributed uniformly on the unit sphere. 
\par
A closely related model to Johnstone's that has been receiving more attention in recent years is one in which the spikes are allowed to vary with \(n,N,\) some of them increasing to infinity. Particularly, Shen et al.~\cite{shen} obtained eigenstructure consistency results of the form
\begin{equation}\label{0.1}
    \frac{\hat{l}}{l} \xrightarrow[]{a.s.} 1, \hspace{0.2cm} |<p_l,u_l>| \xrightarrow[]{a.s.} 1
\end{equation}
when finitely many of the eigenvalues grow with \(n,N,\) the rest tend to a fixed finite constant, and \(\mathbf{Z}_n \in \mathbb{R}^{N \times n}\) has i.i.d. entries whose distributions are independent of \(n\) with \(\mathbb{E}[z_{11}]=0, \mathbb{E}[z^2_{11}]=1, \mathbb{E}[z^4_{11}]<\infty.\) Their results require no assumptions on the relation between \(N\) and \(n,\) but rather on the growth of the considered eigenvalues \(l\) relative to \(N,n,\) or their ratio \(N/n,\) covering thus parts of the three regimes mentioned above. For instance, if \(l\) is \textit{separated} from the rest (i.e., \(|\frac{l}{l'}-1|>\epsilon_0>0\) for any eigenvalue \(l' \ne l\)) and \textit{divergent} (i.e., \(l=l(n,N) \to \infty\)) while \(N,n\) are comparable, then (\ref{0.1}) holds. Two further questions their work evoke would be how the sample eigenvalues fluctuate from the ground truth, and what the rate of convergence of the sample eigenvectors to the true ones is, as Shen et al. do not address these aspects.
\par
Two recent papers tackled these questions and obtained partial answers under separability and divergence conditions on the considered spikes \(l\) similar to the ones in Shen et al.~\cite{shen}. First, Wang and Fan~\cite{fan} considered the case in which the columns of \(\mathbf{Z}_n \in \mathbb{R}^{N \times n}\) are random vectors with independent subgaussian entries of uniformly bounded norms, mean zero, variance one, and found the asymptotic behavior of the eigenstructure of \(\mathbf{S}_n.\) Second, Cai et al.~\cite{cai} derived CLTs for the eigenvalues and eigenvector consistency when \(M\) can grow with \(n\) but slower than \(n^{1/6}\) within a more general model than the prevalent one: \(\Sigma=\Gamma \Gamma^T, \mathbf{X}_n=\Gamma \mathbf{Z}_n,\) where \(\Gamma \in \mathbb{R}^{N \times (N+l)}\) with \(l/N \to 0,\) and \(\mathbf{Z}_n \in \mathbb{R}^{(N+l) \times n}\) with i.i.d. entries of mean zero, variance one, and uniformly bounded fourth moment (the price paid for this level of generality is twofold: tracking an additional matrix, one containing right singular vectors of \(\Gamma,\) and assuming a certain averaging property of the columns of these matrices). Both frameworks make no rigid assumption on the relation between \(N\) and \(n,\) other than \(N>n,\) and \(N \geq cn\) for some \(c>0,\) respectively, and allow a fair amount of flexibility for the bulk eigenvalues, their conditions involving ratios of the type \(\frac{N}{nl}\) (for spikes \(l\)) being bounded or tending to zero.  Besides, while the CLT centering in theorem \(3.1\) in Wang and Fan~\cite{fan} is not explicit (it is a sum of two terms, one deterministic and one stochastically bounded), Cai et al.~\cite{cai} find a deterministic quantity for it, \(l+\frac{N-M}{n} \cdot \frac{l}{l-1},\) which differs from \(\phi(l),\) the almost sure limit of the empirical spiked eigenvalues for Johnstone's original model, solely in one respect, namely \(\frac{N-M}{n}\) substitutes \(\gamma.\)
\par
In this paper, we let \(M\) grow slightly slower than \(n\) within a relaxed version of the random matrix theory regime (i.e., \(\frac{N}{n} \in (\gamma^{-1},\gamma)\) for some \(\gamma \in (0,\infty)\)), and assume the i.i.d. entries of \(\mathbf{Z}_n \in \mathbb{R}^{N \times n}\) are subgaussian with uniformly bounded norms, mean zero, and variance one. Our main results are two CLTs for empirical eigenvalues \(\hat{l}\) whose deterministic counterparts \(l\) are separated divergent spikes, one statistical in nature (i.e., random centering), and one oracle (i.e., deterministic centering), together with eigenvector consistency rates under several growth regimes of \(l.\) Although our framework does not cover the high-dimensional low-sample size regime as Wang and Fan~\cite{fan} and Cai et al.~\cite{cai} implicitly do, it has the advantage of allowing the number of spikes to be much larger than what had been previously considered. Particularly, the statistical CLT holds for separated divergent eigenvalues as long as \(M=o(n^{2/3}),\) while the oracle CLT requires \(M\) to grow slightly slower than \(n.\) Furthermore, a quasi transition emerges: namely, if \(M=o(\sqrt{n}),\) then up to a large extent each such empirical eigenvalue \(\hat{l}\) has its CLT unaffected by the rest of the spikes, whereas this does not occur in the remaining parts of the two ranges.

\subsection{Model and Eigenstructure Representation}\label{1.1}

In this subsection, we introduce the employed model, and the key equations relating the eigenstructure of the sample covariance matrix to its deterministic counterpart. 
\par
For \(n,N=N(n),M=M(n) \in \mathbb{N},\) and fixed \(K,\gamma \in (0,\infty),\) consider
\vspace{0.2cm}
\(\newline A1. \hspace{0.1cm} \Sigma_n=\mathbf{U}_n diag(l^{(n)}_1,l^{(n)}_2, \hspace{0.05cm} ... \hspace{0.05cm},l^{(n)}_M,1,1,\hspace{0.05cm} ... \hspace{0.05cm},1)\mathbf{U}_n^T \in \mathbb{R}^{N \times N}\) deterministic, \(l^{(n)}_1 \geq l^{(n)}_2 \geq ... \geq l^{(n)}_M \geq 1,\) and \(\mathbf{U}_n= [u_1^{(n)} \hspace{0.1cm} u_2^{(n)} \hspace{0.1cm} ... \hspace{0.1cm} u_N^{(n)}] \in \mathbb{R}^{N \times N}\) orthogonal,
\(\newline A2. \hspace{0.2cm} \mathbf{X}_n=\mathbf{U}_ndiag(\sqrt{l^{(n)}_1},\sqrt{l^{(n)}_2},\hspace{0.05cm} ... \hspace{0.05cm},\sqrt{l^{(n)}_M},1,1, \hspace{0.05cm} ... \hspace{0.05cm}, 1) \mathbf{Z}_n \in \mathbb{R}^{N \times n},\) \(\mathbf{Z}_n=(z^{(n)}_{ij})_{1 \leq i \leq N, 1 \leq j \leq n}\) with entries i.i.d. subgaussian random variables, \(\mathbb{E}[z^{(n)}_{11}]=0, \hspace{0.1cm} \mathbb{E}[(z^{(n)}_{11})^2]=1, \hspace{0.1cm} ||z^{(n)}_{11}||_{\psi_2} \leq K,\) where \(||\cdot||_{\psi_2}\) is given by definition \(2.5.6\) in Vershynin~\cite{vershynin}:
\[||X||_{\psi_2}=\inf{\{t>0: \mathbb{E}[\exp(X^2/t^2)] \leq 2\}},\]
\(A3. \hspace{0.2cm} \frac{1}{\gamma} \leq \frac{N}{n} \leq \gamma,\)
\(\newline A4. \hspace{0.2cm} \lim_{n \to \infty}{\frac{M(n)}{n}}=0.\)

\vspace{0.2cm}
\par 
Since any symmetric matrix is diagonalizable, \(A1\) is merely adopting a spiked model with bulk eigenvalues identical to one, spikes at least one, and considering an orthonormal basis \(\{u_1^{(n)}, u_2^{(n)}, \hspace{0.1cm} ... \hspace{0.1cm}, u_N^{(n)}\}\) of \(\mathbb{R}^{N}\) formed of eigenvectors of \(\Sigma_n.\) As mentioned earlier, the independence condition in \(A2\) is prevalent when handling random matrices whose columns are i.i.d. samples, tantamount to imposing that \[\mathbf{Z}_n:=diag(\frac{1}{\sqrt{l^{(n)}_1}},\frac{1}{\sqrt{l^{(n)}_2}},\hspace{0.05cm} ... \hspace{0.05cm},\frac{1}{\sqrt{l^{(n)}_M}},1,1,\hspace{0.05cm} ... \hspace{0.05cm}, 1)\mathbf{U}_n^T\mathbf{X}_n \in \mathbb{R}^{N \times n}\] 
has i.i.d. entries, a strengthening of \(\frac{1}{n}\mathbb{E}[\mathbf{Z}_n\mathbf{Z}_n^T]=I_N\) (equivalent to \(\frac{1}{n}\mathbb{E}[\mathbf{X}_n\mathbf{X}_n^T]=\Sigma_n\)), while subgaussianity, technical in nature, provides a toolkit for controlling quantitatively a wide range of events involving high-dimensional random vectors and matrices (the literature on concentration inequalities for such random variables is extensive: see, for instance, Vershynin~\cite{vershynin}). 
\par
Condition \(A3\) reflects that the dimension of the observations and their number increase at the same rate, which can be interpreted as a relaxation of the random matrix theory regime. A looser control on \(\frac{N}{n}\) than in the classical regimes described in the introduction suffices primarily because the investigated spikes grow to infinity with \(n,\) and this growth annihilates the fluctuations originating in the bulk: namely, these latter contributions stay bounded whenever \(\frac{N}{n}\) is bounded from above, while sometimes dividing with \(\frac{N}{n}\) is needed, which explains the imposed lower bound. Lastly, \(A4,\) stating that the number of spikes increases slower than the number of samples (or equivalently, than the number of features), is the leading force pushing most of the errors to zero. It will be seen that this condition suffices for eigenstructure consistency (Proposition~\ref{prop1} and Proposition~\ref{prop2}), whereas our proof technique for refining the convergences behind it necessitates a larger gap between \(M=M(n)\) and \(n\) than \(\lim_{n \to \infty}{\frac{M(n)}{n}}=0.\)
\par
Moreover, the dimension of the samples, \(N,\) and the number of spikes, \(M,\) are taken to be functions of the number of samples, \(n,\) and for the sample covariance matrix \(\mathbf{S}_n=\frac{1}{n}\mathbf{X}_n\mathbf{X}_n^T,\)
\[\hat{l}^{(n)}_1 \geq \hat{l}^{(n)}_2 \geq ... \geq \hat{l}^{(n)}_M\]
denote its largest \(M\) eigenvalues with corresponding linearly independent unit eigenvectors 
\[p^{(n)}_1,p^{(n)}_2, \hspace{0.05cm} ...  \hspace{0.05cm}, p^{(n)}_M\]
(\(\mathbf{S}_n\) is symmetric and so diagonalizable, which justifies why these \(M\) unit eigenvectors can be chosen thus). Because our goal is comparing the eigenstructures of \(\mathbf{S}_n\) and \(\Sigma_n,\) we assume for simplicity in what follows that \(\mathbf{U}_n=I_N,\) state the results nevertheless in full generality (i.e., using \(u_1^{(n)}, u_2^{(n)},\hspace{0.05cm}...\hspace{0.05cm}, u_N^{(n)}\)) whenever necessary, and drop the superscripts marking the dependency on \(n\) to keep notation as light as possible. 
\par
This paper builds on Paul's approach in \cite{paul}, and consequently the starting point is the block decomposition of \(\mathbf{S}_n\) that provides two crucial equations, (\ref{4}) and (\ref{5}), upon which everything else relies. Keeping most of his notation, split \(\mathbf{S}_n\) into four rectangular matrices by dividing each column of \(\mathbf{Z}_n\) into two vectors with entry indices given by \(A=\{1,2,\hspace{0.05cm} ... \hspace{0.05cm}, M\}, B=\{M+1,M+2, \hspace{0.05cm} ... \hspace{0.05cm}, N\}:\) that is,
\[
\mathbf{Z}_n= \begin{pmatrix}
Z_A\\
Z_B
\end{pmatrix},
\hspace{0.1cm}
\mathbf{X}_n=\Sigma_n^{1/2}\mathbf{Z}_n= \begin{pmatrix}
\Lambda^{1/2} Z_A\\
Z_B
\end{pmatrix},
\hspace{0.1cm}
\mathbf{S}_n=
\begin{pmatrix}
S_{AA} & S_{AB}\\
S_{BA} & S_{BB}
\end{pmatrix},\] 
for \(Z_A \in \mathbb{R}^{M \times n},Z_B \in \mathbb{R}^{(N-M) \times n},\) and  \(\Lambda=diag(l_1,l_2,\hspace{0.05cm} ... \hspace{0.05cm},l_M).\) Denote by \(\mathcal{M} \in \mathbb{R}^{(N-M) \times (N-M)}\) the diagonal matrix containing the eigenvalues of \(S_{BB}=\frac{1}{n}Z_BZ_B^T,\)
\[\frac{1}{\sqrt{n}}Z_B=V\mathcal{M}^{1/2}H^T, \hspace{0.5cm} T=\frac{1}{\sqrt{n}}H^TZ_A^T=[t_1 \hspace{0.1cm} t_2 \hspace{0.1cm} ... \hspace{0.1cm} t_M]\]
for \(V \in \mathbb{R}^{(N-M) \times (N-M)}\) orthogonal, and \(H \in \mathbb{R}^{n \times (N-M)}\) with its first \(\min(N-M,n)\) columns forming an orthonormal set and the remaining zero. Such a decomposition can always be derived from an SVD, \(\frac{1}{\sqrt{n}}Z_B=V_0\mathcal{M}_0H_0^T\): let \(V=V_0, \mathcal{M}=\mathcal{M}_0\mathcal{M}_0^T\); if \(N-M \leq n,\) take \(H^T\) to be \(H_0^T\) with its last \(n-(N-M)\) rows removed (as \(\mathcal{M}_0H_0^T=\mathcal{M}^{1/2}H^T\)); if \(N-M>n,\) let \(H^T\) be given by concatenating \(H_0^T\) and \((N-M)-n\) zero rows (again, \(\mathcal{M}_0H_0^T=\mathcal{M}^{1/2}H^T\)). This decomposition has the advantage of \(\mathcal{M}^{1/2}\) being diagonal (implying that \(\frac{1}{\sqrt{n}}Z_B^T=H\mathcal{M}^{1/2}V^T,\) which will be used to deduce (\ref{4}) and (\ref{5})), while still giving an eigendecomposition for \(S_{BB}=\frac{1}{n}Z_BZ_B^T=V\mathcal{M} V^T\) as it can be easily seen that \(\mathcal{M}^{1/2} H^TH=\mathcal{M}^{1/2}.\)
\par
Denote by \(p_\nu^T=[p_{A,\nu}^T \hspace{0.1cm} p_{B,\nu}^T]\) where \(p_{A,\nu} \in \mathbb{R}^{M}, p_{B,\nu} \in \mathbb{R}^{N-M},\) and \(R_\nu=||p_{B,\nu}||\) for \(\nu \in \{1,2, \hspace{0.05cm} ... \hspace{0.05cm}, M\}.\) Because \(p_1, p_2, \hspace{0.05cm} ... \hspace{0.05cm}, p_M\) are linearly independent unit eigenvectors of \(\mathbf{S}_n,\)
\begin{equation}\label{1}
    S_{AA}p_{A,\nu}+S_{AB}p_{B,\nu}=\hat{l}_\nu p_{A,\nu}
\end{equation}
\begin{equation}\label{2}
    S_{BA}p_{A,\nu}+S_{BB}p_{B,\nu}=\hat{l}_\nu p_{B,\nu}
\end{equation}
\begin{equation}\label{3}
    p_{A,\nu}^Tp_{A,\nu'}+p_{B,\nu}^Tp_{B,\nu'}=\delta_{\nu, \nu'}
\end{equation}
for all \(\nu, \nu' \in \{1,2, \hspace{0.05cm} ...\hspace{0.05cm}, M\}.\)
\par
If \(\hat{l}_\nu I-S_{BB}\) is invertible (which will be the case for divergent eigenvalues \(l_\nu\)), then \(R_\nu<1\) (otherwise, \(p_{A,\nu}=0,||p_{B,\nu}||=1,\) from which \((\hat{l}_\nu I-S_{BB})p_{B,\nu}=\mathbf{0}\) using (\ref{2}), absurd), and for \(a_\nu=\frac{p_{A,\nu}}{||p_{A,\nu}||}=\frac{p_{A,\nu}}{\sqrt{1-R_\nu^2}}\) the following identities ensue:
\begin{equation}\label{4}
    (S_{AA}+\Lambda^{1/2}T^T\mathcal{M}(\hat{l}_\nu I- \mathcal{M})^{-1}T\Lambda^{1/2})a_\nu=\hat{l}_\nu a_\nu
\end{equation}
\begin{equation}\label{5}
    a^T_\nu\Lambda^{1/2}T^T\mathcal{M}(\hat{l}_\nu I- \mathcal{M})^{-2}T\Lambda^{1/2} a_\nu=\frac{R_\nu^2}{1-R_\nu^2}.
\end{equation}
To see why this is so, suppose \(\hat{l}_\nu I-S_{BB}\) is invertible. Then (\ref{2}) yields
\(p_{B,\nu}=(\hat{l}_\nu I-S_{BB})^{-1}S_{BA}p_{A,\nu},\) which together with (\ref{1}), upon dividing both sides with \(\sqrt{1-R_\nu^2}>0\) gives (\ref{4}) since \(S_{BA}^T=S_{AB}=\frac{1}{n}\Lambda^{1/2}Z_AZ_B^T,S_{BB}=\frac{1}{n}Z_BZ_B^T=\frac{1}{n}V\mathcal{M}V^T,\) \(\frac{1}{\sqrt{n}}Z_B=V\mathcal{M}^{1/2}H^T, T=\frac{1}{\sqrt{n}}H^TZ_A^T,\) as well as (\ref{5}) using (\ref{3}) for \((\nu,\nu),\) and \(||p_{A,\nu}||^2=1-R_\nu^2.\)  
\par
We conclude this section by stressing the significance of Paul's method in \cite{paul} (reproduced above) for our results: its simplicity may only be equaled by its far-reaching power which this paper seeks to exploit further. In light of this statement, a brief digression into the different methods adopted in the more recent literature is relevant. On the one hand, Wang and Fan~\cite{fan} study the eigenvalues of the sample covariance matrix \(\mathbf{S}_n\) with the aid of its dual
\[\mathbf{S}^{D}_{n}=\frac{1}{n}\mathbf{X}_n^T\mathbf{X}_n=\frac{1}{n}\sum_{1 \leq i \leq M}{l_iZ_i^TZ_i}+\frac{1}{n}\sum_{M+1 \leq i \leq N}{Z_i^TZ_i}:=\mathcal{A}_n+\mathcal{B}_n,\]
where \(Z_1, Z_2, \hspace{0.05cm} ... \hspace{0.05cm}, Z_N \in \mathbb{R}^{1 \times n}\) are the rows of \(\mathbf{Z}_n.\) The authors obtain first a CLT for the eigenvalues of \(\mathcal{A}_n\) using ideas from Anderson~\cite{anderson} and then account for the lower order contribution of \(\mathcal{B}_n,\) which yields the stochastically bounded component in the centering. On the other hand, Cai et al. in \cite{cai} start with an alternative definition of the eigenvalues of \(\mathbf{S}_n,\) namely they are the zeros of its characteristic polynomial. This furnishes an equation for the empirical eigenvalues in the form of the determinant of an \(M \times M\) matrix being zero, and a CLT arises via an entry-wise control of this random object, whose potentially large dimension in turn imposes the number of spikes to grow slower than \(n^{1/6}.\)

\section{Discussion of Results}\label{sec2}

We begin by considering the consistency between divergent eigenvalues and their empirical counterparts as well as for their corresponding eigenvectors. The former, a generalization of its analogue in Shen et al.~\cite{shen}, always holds for the model described above, while the latter, tightly related to \(l_2\)-convergence, requires separation.

\begin{prop}\label{prop1}
    Suppose \(A1-A4\) hold and
    \[\lim_{n \to \infty}{l^{(n)}_{\nu(n)}}=\infty\]
    for some \(\nu(n) \in \{1,2,\hspace{0.05cm} ... \hspace{0.05cm},M(n)\}\) for all \(n \in \mathbb{N}.\)
    Then as \(n \to \infty,\) jointly for all \(k \in \{1,2, \hspace{0.05cm}... \hspace{0.05cm}, \nu(n)\},\)
    \[\frac{\hat{l}^{(n)}_k}{l^{(n)}_k} \xrightarrow[]{a.s.} 1.\]
\end{prop}

\begin{prop}\label{prop2}
    Suppose \(A1-A4\) hold
    \[\liminf_{n \to \infty}{\frac{l_{\nu(n)}^{(n)}}{l_{\nu(n)+1}^{(n)}}} >1+\epsilon_0, \hspace{0.5cm} \liminf_{n \to \infty}{\frac{l_{\nu(n)-1}^{(n)}}{l_{\nu(n)}^{(n)}}}>1+\epsilon_0, \hspace{0.5cm}  \lim_{n \to \infty}{l^{(n)}_{\nu(n)}}=\infty\] 
    for some \(\nu(n) \in \{1,2, \hspace{0.05cm} ... \hspace{0.05cm}, M(n)\}\) for all \(n \in \mathbb{N},\) and \(\epsilon_0>0\) (where by convention \(l_0=l_{M(n)+1}=1\)). Then as \(n \to \infty,\)
    \[<p^{(n)}_{\nu(n)},u^{(n)}_{\nu(n)}>^2 \xrightarrow[]{a.s.} 1.\]
\end{prop}

\par
Our focus will be the study of eigenvalues satisfying the conditions in Proposition~\ref{prop2} together with their corresponding eigenvectors, and thus we address next the relation between such separated divergent spikes and their empirical counterparts. Assume without loss of generality that
\[<p_\nu,u_\nu>=\sqrt{1-R_\nu^2} \cdot <a_\nu,e_\nu> \geq 0.\]
Then Proposition~\ref{prop2} yields \(a_\nu \approx e_\nu,\) which in conjunction with (\ref{4}) provides a link between \(\hat{l}_\nu\) and \(l_\nu:\) 
\begin{equation}\label{9990}
    a_\nu^T(S_{AA}+\Lambda^{1/2}T^T\mathcal{M}(\hat{l}_\nu I- \mathcal{M})^{-1}T\Lambda^{1/2})a_\nu=\hat{l}_\nu
\end{equation}
because the left-hand side term is, roughly speaking, \(l_\nu \cdot (1+t_\nu^T\mathcal{M}(\hat{l}_\nu I- \mathcal{M})^{-1}t_\nu).\) This representation of \(\hat{l}_\nu\) is the pillar of our main CLT for \(\frac{\hat{l}_\nu}{l_\nu}.\) Once \(a_\nu \approx e_\nu\) is sharpened, we get a hold of the fluctuations of \(\frac{\hat{l}_\nu}{l_\nu}-1\) with the peculiarity that our analysis leads us to a more exotic centering than usual: a sum of one random term and one deterministic, which is the content of our first theorem.

\begin{theorem}\label{th2}
Suppose \(A1-A3\) hold,
\[\liminf_{n \to \infty}{\frac{l_{\nu(n)}^{(n)}}{l_{\nu(n)+1}^{(n)}}} >1+\epsilon_0, \hspace{0.5cm} \liminf_{n \to \infty}{\frac{l_{\nu(n)-1}^{(n)}}{l_{\nu(n)}^{(n)}}}>1+\epsilon_0, \hspace{0.5cm} \lim_{n \to \infty}{l^{(n)}_{\nu(n)}}=\infty, \hspace{0.5cm} \lim_{n \to \infty}{\frac{l^{(n)}_{\nu(n)}}{M(n)/\sqrt{n}}}=\infty,\]
\[\lim_{n \to \infty}{\frac{\sqrt{\log{n}}}{ \log{\frac{n}{M(n)}}}}=0, \hspace{0.5cm} \mathbb{E}[(z^{(n)}_{11})^4] \geq 1+\delta_0\]
for \(\nu(n) \in \{1,2, \hspace{0.05cm} ... \hspace{0.05cm}, M(n)\}\) (where by convention \(l_0=l_{M(n)+1}=1\)), and \(\epsilon_0,\delta_0>0.\) Then for \(\nu = \nu(n),\) as \(n \to \infty,\)
\begin{equation}\label{6}
    \sqrt{\frac{n}{\mathbb{E}[(z^{(n)}_{11})^4]-1}} \cdot (\frac{\hat{l}_\nu}{l_\nu}-1-\frac{1}{n}tr(\mathcal{M}(\hat{l}_\nu I-\mathcal{M})^{-1})-x_{n,\nu,(l_i)_{1 \leq i \leq M}}) \Rightarrow N(0,1),
\end{equation}
where \(x_{n,\nu,(l_i)_{1 \leq i \leq M}}\) is deterministic with
\(|x_{n,\nu,(l_i)_{1 \leq i \leq M}}| \leq c(\epsilon_0) \cdot \frac{M(n)}{n}.\) 
\end{theorem}

\par
A few observations are in order regarding the three additional constraints imposed (on the minimal growth of the considered eigenvalues, the number of spikes, and the fourth moment of the underlying i.i.d. variables). The first condition (meaningful only when \(\limsup_{n \to \infty}{\frac{M(n)}{\sqrt{n}}}=\infty\) because \(\lim_{n \to \infty}{l^{(n)}_{\nu(n)}}=\infty\)), \(\lim_{n \to \infty}{\frac{l^{(n)}_{\nu(n)}}{M(n)/\sqrt{n}}}=\infty,\) is somehow necessary: if the number of spikes is very large and the considered eigenvalues are relatively small, then new contributions, more difficult to control than the rest, arise: specifically, cross-terms such as \[\sum_{k \ne \nu, 1 \leq k \leq M}{\frac{1}{n}z_\nu^Tz_k \cdot t_k\mathcal{M}(\hat{l}_\nu I -\mathcal{M})^{-1}t_\nu}\]
have be dealt with as they are no longer negligible (see proof of Theorem~\ref{th2}). The second condition is likely suboptimal, yet necessary for the approach taken here. As mentioned earlier, (\ref{9990}) is the foundation of this CLT, and its left-hand side term becomes a sum of three terms involving \(a_\nu-e_\nu\) from \(a_\nu=e_\nu+(a_\nu-e_\nu).\) To control this difference, we decompose its entries into series containing entries of random matrices, a representation allowing us to handle it more easily. This in turn translates into expressing our object of interest \(\frac{\hat{l}_\nu}{l_\nu}-1\) as a series, and when we scrutinize its summands, a new difficulty emerges: each has small variance, but these variances can generate a non-summable series, which would render the decomposition tactic futile. To bypass this issue, we truncate the (first) series, and this clipping is beneficial if both the tail and the bounds for the finite part are negligible, leading to the growth restraint on the number of spikes. The third condition is needed primarily for concluding that \(\frac{1}{\sqrt{n}}\sum_{1 \leq i \leq n}{(z_{\nu i}^2-1)},\) appropriately normalized, converges to a standard normal: notice that in the extreme case in which \(\mathbb{E}[z^4_{11}]=1\) (equivalent to \(\mathbb{E}[(z^2_{11}-1)^2]=0,\) or \(\mathbb{P}(z_{11}=1)=\mathbb{P}(z_{11}=-1)=\frac{1}{2}\)), such normalization cannot exist as the sum is zero with probability one.
\par
Our second theorem offers an algebraic description of the deterministic component \(x_{n,\nu,(l_i)_{1 \leq i \leq M}}\) in the centering above: it is a zero of a polynomial whose degree depends on how close \(\frac{n}{M}\) and \(n\) are at a logarithmic scale. This result comes for free from the proof of Theorem~\ref{th2}, rendering the left-hand side in (\ref{6}) less mysterious than it might otherwise seem.

\begin{theorem}\label{th4}
Under the assumptions of Theorem~\ref{th2} and \(\frac{M(n)}{n} \leq c_0(\epsilon_0),\) \(x=x_{n,\nu,(l_i)_{1 \leq i \leq M}}\) satisfies the following polynomial equation
\[x=\overline{O}+\sum_{1 \leq j \leq 2s^2+2s}{\overline{O}_jx^j},\]
where \(s=\left \lfloor{\frac{8\log{n}}{\log{\frac{n}{M}}}}\right \rfloor, |\overline{O}| \leq c(\epsilon_0) \cdot \frac{M(n)}{n}, |\overline{O}_j| \leq c(\epsilon_0)^{j+1} \cdot \frac{M(n)}{n}.\) Moreover, these coefficients are given by 
\[\sum_{0 \leq j \leq s}{(\sum_{0 \leq i \leq 2s}{(2a_i+b_i+c_i)z^i})(\sum_{0 \leq i \leq 2s}{b_iz^i})^j}=\overline{O}+\sum_{1 \leq j \leq 2s^2+2s}{\overline{O}_{j}z^j},\]
for \(\Tilde{I} \in \mathbb{R}^{M \times M}\) with all entries \(1,\) \(\mathcal{R}_\nu=diag(\frac{1}{l_1-l_\nu},\hspace{0.1cm} ... \hspace{0.1cm} ,\frac{1}{l_{\nu-1}-l_\nu},0,\frac{1}{l_{\nu+1}-l_\nu},\hspace{0.1cm}  ... \hspace{0.1cm} ,\frac{1}{l_M-l_\nu}) \in \mathbb{R}^{M \times M},\)
\begin{equation}\label{129}
    \mathcal{M}_{\nu}(z)=\sum_{0 \leq j \leq s}{\frac{1}{n^{j+1}}\Lambda^{1/2}(-\mathcal{R}_\nu \Lambda^{1/2} \Tilde{I} \Lambda^{1/2}+nzl_\nu \mathcal{R}_\nu)^j\mathcal{R}_\nu \Lambda^{1/2}\Tilde{I}},
\end{equation}
\begin{equation}\label{130}
    -\sum_{k \ne \nu}{(\mathcal{M}_{\nu}(z))_{k \nu}}=\sum_{0 \leq i \leq s}{a_iz^i},
\end{equation}
\begin{equation}\label{131}
    -nl_\nu(\mathcal{M}_{\nu}(z)^T\Lambda^{-1}\mathcal{M}_{\nu}(z))_{\nu \nu}=\sum_{0 \leq i \leq 2s}{b_iz^i},
\end{equation}
\begin{equation}\label{132}
    n(\mathcal{M}_{\nu}(z)^T\mathcal{M}_{\nu}(z))_{\nu \nu}+(\mathcal{M}_{\nu}(z)^T\Tilde{I}\mathcal{M}_{\nu}(z))_{\nu \nu}=\sum_{0 \leq i \leq 2s}{c_iz^i},
\end{equation}
with the equations defining \(\overline{O},\overline{O}_j,a_i,b_i,c_i\) interpreted as equalities of polynomials in \(z\) (by convention, \(a_j=0\) for \(s< j \leq 2s\)).

\end{theorem}

\par
In spite of our result suggesting that computing \(x=x_{n,\nu,(l_i)_{1 \leq i \leq M}}\) is infeasible since it is a root of a polynomial with degree \(2s^2+2s \geq 2 \cdot 8^2+2 \cdot 8=144,\) note that, up to an extent, its value is irrelevant: Slutsky's lemma implies any \(\Tilde{x}=x+o(\frac{1}{\sqrt{n}})\) can substitute \(x\) in (\ref{6}). In particular, if \(\lim_{n \to \infty}{\frac{M(n)}{n^{1-1/(2k_0)}}}=0\) for some fixed \(k_0 \in \mathbb{N}\) (i.e., \(s\) can be chosen independently of \(n\)), then such an \(\Tilde{x}\) can be obtained as a sum of \(k_0-1\) terms, one summand at a time with the \(k^{th}\) being of order \((M/n)^k\) (when \(k_0=1,\) this approximation is zero). To see this, notice first that for \(n\) large enough, \(s \leq 16k_0.\) Take next \(x_1=x-\overline{O},\) for which the binomial theorem yields
\begin{equation}\label{5001}
    x_1=\overline{O}'+\sum_{1 \leq j \leq 2s^2+2s}{\overline{O}'_jx^j_1}
\end{equation}
with 
\[|\overline{O}'|=|\sum_{1 \leq j \leq 2s^2+2s}{\overline{O}^j \cdot \overline{O}_j}| \leq c_1(\epsilon_0) \cdot \frac{M^2}{n^2}, \hspace{0.5cm} |\overline{O}'_j|= |\sum_{j \leq l \leq 2s^2+2s}{\overline{O}_l \cdot \binom{l}{j} \cdot \overline{O}^{l-j}}| \leq c_1(\epsilon_0) \cdot \frac{M}{n},\]
where \(c_1(\epsilon_0)=c_1(\epsilon_0, k_0).\) These bounds and \(|x_1| \leq |x|+|\overline{O}| \leq c'(\epsilon_0) \cdot \frac{M}{n}\) imply that for \(n\) sufficiently large
\[|x_1| \leq |\overline{O}'|+\sum_{1 \leq j \leq 2 s^2+2s}{c_1(\epsilon_0) \cdot \frac{M}{n} \cdot |x_1|^j} \leq |\overline{O}'|+2c_1(\epsilon_0) \cdot \frac{M}{n} \cdot |x_1| \leq |\overline{O}'|+\frac{|x_1|}{2},\]
from which \(x_1=O(c_1(\epsilon_0) \cdot \frac{M^2}{n^2}).\)
Iterating this step \(k_0-1\) times renders \[x=\overline{O}^{(1)}+\overline{O}^{(2)}+...+\overline{O}^{(k_0-1)}+x_{k_0-1}\] 
with \(\overline{O}^{(i)}=O(c_{i}(\epsilon_0) \cdot \frac{M^i}{n^i}),\) and \(x_{k_0-1}\) satisfying a polynomial equation analogous to (\ref{5001}) with the constant term \(O(c_{k_0}(\epsilon_0) \cdot \frac{M^{k_0}}{n^{k_0}}),\) from which first
\[|x_{k_0-1}| \leq |x|+|\overline{O}^{(1)}|+|\overline{O}^{(2)}|+...+|\overline{O}^{(k_0-1)}| \leq c'_{k_0}(\epsilon_0) \cdot \frac{M}{n},\]
and then, reasoning as for \(x_1,\) 
\[x_{k_0-1}=O(c_{k_0}(\epsilon_0) \cdot \frac{M^{k_0}}{n^{k_0}})=o(\frac{1}{\sqrt{n}}).\] 
Hence we can take
\[\Tilde{x}=\overline{O}^{(1)}+\overline{O}^{(2)}+...+\overline{O}^{(k_0-1)}.\]

\subsection{Statistical and Oracle CLTs}\label{2.1}

With a complete understanding of the CLT in (\ref{6}) under our belt, we now force its centering to be either purely random or purely deterministic.
\par
Consider first finding a statistical version of (\ref{6}), which is tantamount to devising an estimator of \(x_{n,\nu,(l_i)_{1 \leq i \leq M}}.\) A corollary of the definitions of \(a_i, b_i, c_i,\) and the rationale in the last paragraph of the previous subsection is 
\[x_{n,\nu,(l_i)_{1 \leq i \leq M}}=\frac{1}{n}\sum_{k \ne \nu}{\frac{l_k}{l_k-l_\nu}}+O(c_1(\epsilon_0) \cdot \frac{M^2}{n^2}),\] 
which furnishes the desired estimator for a smaller range of \(M\) than that in Theorem~\ref{th2}, \(\lim_{n \to \infty}{\frac{M(n)}{n^{3/4}}}=0\) (\(\sqrt{n} \cdot \frac{M^2}{n^2} \to 0\)), after showing that its empirical variant, \(\frac{1}{n}\sum_{k \ne \nu}{\frac{\hat{l}_k}{\hat{l}_k-\hat{l}_\nu}},\) is close enough to its deterministic counterpart, a task not considerably difficult.

\begin{theorem}\label{th1}
Suppose \(A1-A3\) hold,
\[\liminf_{n \to \infty}{\frac{l_{\nu(n)}^{(n)}}{l_{\nu(n)+1}^{(n)}}} >1+\epsilon_0, \hspace{0.5cm} \liminf_{n \to \infty}{\frac{l_{\nu(n)-1}^{(n)}}{l_{\nu(n)}^{(n)}}}>1+\epsilon_0, \hspace{0.5cm} \lim_{n \to \infty}{l^{(n)}_{\nu(n)}}=\infty, \hspace{0.5cm} \lim_{n \to \infty}{\frac{l^{(n)}_{\nu(n)}}{M(n)/\sqrt{n}}}=\infty,\]
\[\lim_{n \to \infty}{\frac{M(n)}{n^{2/3}}}=0, \hspace{0.5cm} \mathbb{E}[(z^{(n)}_{11})^4] \geq 1+\delta_0\]
for \(\nu(n) \in \{1,2, \hspace{0.05cm} ... \hspace{0.05cm}, M(n)\}\) (where by convention \(l_0=l_{M(n)+1}=1\)), and some \(\epsilon_0,\delta_0>0.\) Then as \(n \to \infty,\) for \(\nu = \nu(n),\)
\begin{equation}\label{301}
    \sqrt{\frac{n}{\mathbb{E}[(z^{(n)}_{11})^4]-1}} \cdot (\frac{\hat{l}_\nu}{l_\nu}-1-\frac{1}{n}tr(\mathcal{M}(\hat{l}_\nu I-\mathcal{M})^{-1})-\frac{1}{n}\sum_{k \ne \nu}{\frac{\hat{l}_k}{\hat{l}_k-\hat{l}_\nu}}) \Rightarrow N(0,1),
\end{equation}
where almost surely \(\hat{l}_k \ne \hat{l}_\nu\) for all \(k \ne \nu.\)
\end{theorem}

\par
Let us expand on the statement and implications of this theorem, a statistical CLT that, to best of our knowledge, has no close relative in the existing literature. The reason for which the last condition mentioned above, \(\lim_{n \to \infty}{\frac{M(n)}{n^{3/4}}}=0,\) has to be strengthened to \(\lim_{n \to \infty}{\frac{M(n)}{n^{2/3}}}=0\) is making the substitution from deterministic to empirical possible: in other words,
\[\sqrt{n} \cdot (\frac{1}{n}\sum_{k \ne \nu}{\frac{\hat{l}_k}{\hat{l}_k-\hat{l}_\nu}}-\frac{1}{n}\sum_{k \ne \nu}{\frac{l_k}{l_k-l_\nu}}) \xrightarrow[]{p} 0\]
is needed, and this adds a new toll on the growth of \(M.\) Regarding the purely empirical centering in (\ref{301}), notice it can be computed once an estimator for \(M\) is chosen because the latter turns 
\[1+\frac{1}{n}tr(\mathcal{M}(\hat{l}_\nu I-\mathcal{M})^{-1})+\frac{1}{n}\sum_{k \ne \nu}{\frac{\hat{l}_k}{\hat{l}_k-\hat{l}_\nu}}\] 
explicit and produces an estimator of \(\mathbb{E}[z_{11}^4]-1,\) making the left-hand side term in (\ref{301}) fit for tasks such as obtaining confidence intervals for the true eigenvalue \(l_\nu\). However, in this paper we do not pursue further this question and leave it instead to future research. 
\par
Next, the oracle CLT is derived anew from (\ref{6}), by converting the centering into a deterministic quantity, which requires modifying \(\frac{1}{n}tr(\mathcal{M}(\hat{l}_\nu I -\mathcal{M})^{-1}).\) Note that (\ref{6}) virtually renders the true eigenvalue as a function of its empirical counterpart:
\[l_\nu \approx \frac{\hat{l}_\nu}{1+\frac{1}{n}tr(\mathcal{M}(\hat{l}_\nu I -\mathcal{M})^{-1})+x_{n,\nu,(l_i)_{1 \leq i \leq M}}},\]
while, at a high level, a deterministic centering asks for inverting this relationship (i.e., expressing \(\hat{l}_\nu\) in terms of the ground truth). In our situation, the key towards this inversion is the Marchenko-Pastur result in \cite{marchenkopastur} which allows replacing \(\frac{1}{n}tr(\mathcal{M}(\hat{l}_\nu I -\mathcal{M})^{-1}),\) a rational function with random coefficients, by a non-random function of \(\hat{l}_\nu,\) even though the empirical distributions behind these traces generally do not converge to any law (since we are not operating under the random matrix theory regime). It is also worth emphasizing that an ingenuous quasi inversion such as 
\(\hat{l}_\nu \approx l_\nu (1+\frac{1}{n}tr(\mathcal{M}(l_\nu I -\mathcal{M})^{-1})+x_{n,\nu,(l_i)_{1 \leq i \leq M}})\)
is generally not sharp enough because \(l_\nu\) can grow to infinity much slower than \(n\) does. Take, for instance, a simpler yet closely related function to the relation between \(l_\nu\) and \(\hat{l}_\nu,\) 
\[y=y(x)=\frac{x}{1+\frac{1}{x}}, \hspace{0.2cm} x>0;\]
for \(y\) large, \(x \approx y\) yields the naive proxy \(y \cdot (1+\frac{1}{y})=y+1,\) whereas the actual inverse is \(\frac{y+\sqrt{y^2+4y}}{2}=y+1+\frac{2}{\sqrt{y^2+4y}+y+2},\) about \(O(\frac{1}{y})\) away from \(y+1,\) an error that can explode once multiplied by, say, \(\frac{\sqrt{n}}{y}.\)

\begin{theorem}\label{th3}
Under the assumptions of Theorem~\ref{th2}, as \(n \to \infty,\)
\begin{equation}\label{302}
    \sqrt{\frac{n}{\mathbb{E}[(z^{(n)}_{11})^4]-1}} \cdot (\frac{\hat{l}_\nu}{l_\nu}-1-\frac{N-M}{n} \cdot \frac{1}{l_\nu-1}-x_{n,\nu,(l_i)_{1 \leq i \leq M}}) \Rightarrow N(0,1).
\end{equation}
\end{theorem}

\par
This result extends theorem \(2.2\) in Cai et al.~\cite{cai}, which assumes \(\lim_{n \to \infty}{\frac{M(n)}{n^{1/6}}}=0.\) To see why this is the case, note the conclusion of their theorem can be formulated as \[\sqrt{\frac{n}{\mathbb{E}[(z^{(n)}_{11})^4]-1}} \cdot (\frac{\hat{l}_\nu}{\theta_\nu}-1) \Rightarrow N(0,1)\] for \(\theta_\nu=l_\nu(1+\frac{N-M}{n} \cdot \frac{1}{l_\nu-1});\) because \(x_{n,\nu,(l_i)_{1 \leq i \leq M}}=O(c(\epsilon_0) \cdot  \frac{M}{n})=o(\frac{1}{\sqrt{n}})\)) is negligible, (\ref{302}) can be rewritten as \[\sqrt{\frac{n}{\mathbb{E}[(z^{(n)}_{11})^4]-1}} \cdot (\frac{\hat{l}_\nu}{l_\nu}-1-\frac{N-M}{n} \cdot \frac{1}{l_\nu-1}) \Rightarrow N(0,1),\] 
and lastly,
\[\frac{\hat{l}_\nu}{\theta_\nu}-1=\frac{l_\nu}{\theta_\nu} \cdot (\frac{\hat{l}_\nu}{l_\nu}-1-\frac{N-M}{n} \cdot \frac{1}{l_\nu-1}), \hspace{0.4cm} \frac{l_\nu}{\theta_\nu} \to 1.\]
\par
It is worth noticing that when \(\lim_{n \to \infty}{\frac{M}{\sqrt{n}}}=0,\) our CLTs for separated divergent eigenvalues are virtually unaffected by the rest of the spikes: that is,
\[\sqrt{\frac{n}{\mathbb{E}[(z^{(n)}_{11})^4]-1}} \cdot (\frac{\hat{l}_\nu}{l_\nu}-1-\frac{1}{n}tr(\mathcal{M}(\hat{l}_\nu I-\mathcal{M})^{-1})) \Rightarrow N(0,1)\]
\[\sqrt{\frac{n}{\mathbb{E}[(z^{(n)}_{11})^4]-1}} \cdot (\frac{\hat{l}_\nu}{l_\nu}-1-\frac{N-M}{n} \cdot \frac{1}{l_\nu-1}) \Rightarrow N(0,1)\]
since \(x_{n,\nu,(l_i)_{1 \leq i \leq M}}=o(\frac{1}{\sqrt{n}}),\) and the proof of Theorem~\ref{th1} yields \(\frac{\hat{l}_k}{|\hat{l}_k-\hat{l}_\nu|} \leq c(\epsilon_0)\) for all \(k \ne \nu\) with probability tending to one, from which 
\[\sqrt{n} \cdot |\frac{1}{n}\sum_{k \ne \nu}{\frac{\hat{l}_k}{\hat{l}_k-\hat{l}_\nu}}| \leq \frac{Mc(\epsilon_0)}{\sqrt{n}}=o(1).\] 
In contrast, when \(M\) grows at least as fast as \(\sqrt{n},\) both \(\frac{1}{n}\sum_{k \ne \nu}{\frac{\hat{l}_k}{\hat{l}_k-\hat{l}_\nu}}\) and \(x_{n,\nu,(l_i)_{1 \leq i \leq M}}\) may make non-negligible contributions to the CLTs they appear in, an effect unveiled by our results and specific to a (relatively) large number of spikes.

\subsection{Consistency Rates for Eigenvectors}\label{2.2}

Our last two theorems consider the consistency rates of the empirical eigenvectors to their true counterparts, the sixth being a statistical version of the fifth:

\begin{theorem}\label{th5}
Suppose \(A1-A4\) hold,
\[\liminf_{n \to \infty}{\frac{l_{\nu(n)}^{(n)}}{l_{\nu(n)+1}^{(n)}}} >1+\epsilon_0, \hspace{0.5cm} \liminf_{n \to \infty}{\frac{l_{\nu(n)-1}^{(n)}}{l_{\nu(n)}^{(n)}}}>1+\epsilon_0, \hspace{0.5cm} \lim_{n \to \infty}{l^{(n)}_{\nu(n)}}=\infty\]
for some \(\nu(n) \in \{1,2, \hspace{0.05cm} ... \hspace{0.05cm}, M(n)\}\) for all \(n \in \mathbb{N},\) and some \(\epsilon_0>0\) (where by convention \(l_0=l_{M(n)+1}=1\)).

\par
\vspace{0.3cm}
\((a)\) If \(\lim_{n \to \infty}{\frac{l^{(n)}_{\nu(n)}}{n/M(n)}}=0,\) 
then as \(n \to \infty,\) 
\[l^{(n)}_{\nu(n)}(1-<p^{(n)}_{\nu(n)},u^{(n)}_{\nu(n)}>^2)-\frac{N(n)}{n} \xrightarrow[]{p} 0.\]
\par
If additionally
\(\sum_{n \geq 1}{\exp(-\frac{an}{l^{(n)}_{\nu(n)}})}<\infty, \forall a>0,\) then as \(n \to \infty,\) 
\[l^{(n)}_{\nu(n)}(1-<p^{(n)}_{\nu(n)},u^{(n)}_{\nu(n)}>^2)-\frac{N(n)}{n} \xrightarrow[]{a.s.} 0.\]

\vspace{0.3cm}
\((b)(i)\) If \(\lim_{n \to \infty}{\frac{l^{(n)}_{\nu(n)}}{n/M(n)}}=c_\nu \in (0,\infty),\) \(M(n)=M, \nu(n)=\nu\) are fixed, and
\[\lim_{n \to \infty}{\frac{l^{(n)}_{k}l^{(n)}_{\nu}}{(l^{(n)}_{k}-l^{(n)}_{\nu})^2}}=c_{k \nu} \in [0,\infty), k \ne \nu,\]
then as \(n \to \infty,\)
\[l_{\nu(n)}^{(n)}(1-<p^{(n)}_{\nu(n)},u^{(n)}_{\nu(n)}>^2)-\frac{N(n)}{n} \Rightarrow \frac{c_\nu}{M} X_{(c_{k \nu})_{k \ne \nu}},\]
where
\[X_{(c_{k \nu})_{k \ne \nu}}=\sum_{k \ne \nu}{c_{k\nu}y^2_{k\nu}},\]
for \((y_{k \nu})_{k \ne \nu}\) standard normal variables, mutually independent.
In addition, \(y_{ij}=y_{ji}\) for \(i \ne j,\) and \((y_{ij})_{1 \leq i<j \leq M}\) are mutually independent.

\vspace{0.2cm}
\par
\((b)(ii)\) If \(\lim_{n \to \infty}{\frac{l^{(n)}_{\nu(n)}}{n/M(n)}}=c_\nu \in (0,\infty),\) \(\lim_{n \to \infty}{M(n)}=\infty,\) then as \(n \to \infty,\)
\[l_{\nu(n)}^{(n)}(1-<p^{(n)}_{\nu(n)},u^{(n)}_{\nu(n)}>^2)-\frac{l^{(n)}_{\nu(n)}}{n}\sum_{k \ne \nu}{\frac{l^{(n)}_kl^{(n)}_{\nu(n)}}{(l^{(n)}_k-l^{(n)}_{\nu(n)})^2}}-\frac{N(n)}{n} \xrightarrow[]{p} 0.\]
Moreover, if in addition \(\sum_{n \geq 1}{\exp(-a \cdot M(n))}<\infty, \forall a>0,\) then the convergence is almost surely.

\vspace{0.3cm}
\((c)(i)\) If \(\lim_{n \to \infty}{\frac{l^{(n)}_{\nu(n)}}{n/M(n)}}=\infty,\) \(M(n)=M, \nu(n)=\nu\) are fixed, and
\[\lim_{n \to \infty}{\frac{l^{(n)}_kl^{(n)}_\nu}{(l^{(n)}_k-l^{(n)}_\nu)^2}}=c_{k\nu} \in [0,\infty), k \in \{1,2, \hspace{0.05cm} ... \hspace{0.05cm}, M\}, k \ne \nu,\]
then as \(n \to \infty,\)
\[n(1-<p^{(n)}_\nu,u^{(n)}_\nu>^2) \Rightarrow X_{(c_{k \nu})_{k \ne \nu}},\]
where \(X_{(c_{k \nu})_{k \ne \nu}}\) is defined as in \((b)(i).\)

\vspace{0.2cm}
\par
\((c)(ii)\) If \(\lim_{n \to \infty}{\frac{l^{(n)}_{\nu(n)}}{n/\sqrt{M(n)}}}=c_\nu \in (0,\infty),\) \(\lim_{n \to \infty}{M(n)}=\infty, \lim_{n \to \infty}{\frac{M(n)}{\sqrt{n}}}=0,\) and
\[\lim_{n \to \infty}{\frac{1}{M(n)}\sum_{k \ne \nu(n)}{\frac{(l^{(n)}_kl^{(n)}_{\nu(n)})^2}{(l^{(n)}_k-l^{(n)}_{\nu(n)})^4}}}=\sigma_\nu>0,\]
then as \(n \to \infty,\)
\[l_{\nu(n)}^{(n)}(1-<p^{(n)}_{\nu(n)},u^{(n)}_{\nu(n)}>^2)-\frac{l^{(n)}_{\nu(n)}}{n}\sum_{k \ne \nu}{\frac{l^{(n)}_kl^{(n)}_{\nu(n)}}{(l^{(n)}_k-l^{(n)}_{\nu(n)})^2}}-\frac{N(n)}{n} \Rightarrow c_\nu N(0,2\sigma_\nu).\]
\end{theorem}

\vspace{0.2cm}
\textit{Remark:}  When \(M(n)=M\) is fixed, the conditions on \(\nu(n)\) can be relaxed in \((b)(i),\) and \((c)(i):\) instead of having \(\nu(n)\) independent of \(n,\) it suffices to have a permutation \(\tau_1,\tau_2, \hspace{0.05cm} ...\hspace{0.05cm}, \tau_{M-1}\) of \(\{k: k \in \{1,2, \hspace{0.05cm} ... \hspace{0.05cm}, M\}, k \ne {\nu(n)}\}\) for which \(\lim_{n \to \infty}{\frac{l^{(n)}_{\tau(i)}l^{(n)}_{\nu(n)}}{(l^{(n)}_{\tau(i)}-l^{(n)}_{\nu(n)})^2}}=c_{i\nu} \in [0,\infty), i \in \{1,2, \hspace{0.05cm} ... \hspace{0.05cm}, M-1\}.\) Additionally, to avoid dealing with sets potentially empty, \(\{k: k \ne \nu, 1 \leq k \leq M\},\) assume \(M \geq 2,\) which does not restrict generality.

\begin{theorem}\label{th6}
In Theorem~\ref{th5}, the true eigenvalues on the left-hand side terms of the convergences can be replaced by their empirical counterparts.
\end{theorem}

\par
Wang and Fan~\cite{fan} consider as well the consistency between the empirical and the true eigenvectors: \(M\) is fixed in their framework, and thus the authors formulate their results in terms of the differences \(\frac{p_{A,\nu}}{||p_{A,\nu}||}-u_{A,\nu},\) which, after being appropriately centered and normalized, are asymptotically multivariate normal. Since our setting allows \(M\) to depend on \(n,\) we investigate instead inner products and obtain their limiting behavior in terms of the underlying model. Nevertheless, the series representations obtained for each of the entries of these differences can be employed to recoup their asymptotic behavior when \(M\) is fixed: see beginning of section~\ref{sec4} for more details.
\par
We conclude this section with some remarks on Theorem~\ref{th5}: it will become self-evident from its proof that the ratio \(\frac{l_\nu}{n/M}\) heavily influences the behavior of the inner products under consideration. The three parts of the theorem, \((a), (b), (c),\) handle regions of three different regimes, the limit of this ratio as \(n\) tends to infinity being \(0,\) a positive constant, and infinity, respectively, revealing also interesting connections between the eigenvectors and the number of spikes: namely, when the eigenvalues under consideration are relatively large (parts \((b)\) and \((c)\)), whether the number of spikes stays fixed or grows to infinity with \(n\) plays a crucial role in the fluctuations of the corresponding inner products, situation which does not occur when the eigenvalues are relatively small (part \((a)\)). In the former case, more stringent conditions are imposed on the eigenvalues (pointwise limits of ratios), whereas in the latter case, there is more room for flexibility to the degree that only some averages of functions of their ratios have to converge. 
\par
These phenomena are mainly explained by the fact that what we are looking at is determined by the length of a projection (what is denoted by \(R_\nu\)), and a dot product of two (relatively) low-dimensional vectors (\(a_\nu = \frac{p_{A,\nu}}{||p_{A,\nu}||},u_{\nu,A} \in \mathbb{R}^{M},\) whereas \(p_\nu,u_\nu \in \mathbb{R}^{N}\)). The size of \(R_\nu^2\) can be derived from (\ref{5}), while the square of the dot product will be seen to be primarily given by a sum of about \(M\) random variables: 
\[\frac{l_k l_\nu}{(l_k-l_\nu)^2}(\frac{1}{\sqrt{n}}z_k^Tz_\nu)^2, \hspace{0.1cm} 1 \leq k \leq M, k \ne \nu,\] 
where \((z_k)_{1 \leq k \leq M} \in \mathbb{R}^n\) are the rows of \(Z_A.\) It will be shown that \(\frac{1}{\sqrt{n}}z_k^Tz_\nu, k \ne \nu\) are asymptotically mutually independent standard normal variables, and thus whether \(M\) is fixed or grows to infinity with \(n\) dictates the behavior of this sum and the employed normalizations, which in turn indicate the growth conditions required of the spike \(l_\nu.\) In the former case, all the individual ratios must converge to some limits as \(n\) tends to infinity, whereas in the latter, apart from imposing an average of functions of these fractions to be convergent, there are two normalization options, as it happens with a sum of \(q_n\) centered i.i.d. random variables with \(\lim_{n \to \infty}{q_n}=\infty:\) either by \(q_n,\) generating an SLLN-type result (part \((b)(ii)\)), or by \(\sqrt{q_n},\) producing a CLT-type result (part \((c)(ii)\)).

\section{Core Ideas Behind the Proofs}\label{sec3'}

In this section, heuristic and sketchy in nature, we offer a bird's-eye view of the central steps that led us to the results presented above in the hope it will offer the reader a sense of our approach and facilitate his traversal of the proofs. We begin with our propositions and continue with the theorems (the first two in \ref{3'.1}, and the rest in \ref{3'.2}).
\par
Proposition~\ref{prop1} requires an extension of the method in Shen et al.~\cite{shen} whose setup differs from ours in at least one fundamental way, a fixed number of spikes. We consider as well the dual of \(\mathbf{S}_n,\)
\[\mathbf{S}^{D}_n=\frac{1}{n}\mathbf{X}_n^T\mathbf{X}_n=\frac{1}{n}\sum_{1 \leq i \leq M}{l_iZ_i^TZ_i}+\frac{1}{n}\sum_{M+1 \leq i \leq N}{Z_i^TZ_i},\]
where \(Z_1, Z_2, \hspace{0.05cm} ... \hspace{0.05cm}, Z_N \in \mathbb{R}^{1 \times n}\) are the rows of \(\mathbf{Z}_n,\) matrix that has the advantage over \(\mathbf{S}_n\) of separating the eigenvalues of \(\Sigma_n\) in the following sense:
\[(\mathbf{S}^{D}_n)_{ij}=\frac{1}{n}(\mathbf{Z}_n^T\Sigma_n\mathbf{Z}_n)_{ij}=\frac{1}{n}\sum_{1 \leq k \leq M}{l_kz_{ki}z_{kj}}+\frac{1}{n}\sum_{M+1 \leq k \leq N}{z_{ki}z_{kj}},\]
whereas
\[(\mathbf{S}_n)_{ij}=\frac{1}{n}(\Sigma_n^{1/2}\mathbf{Z}_n\mathbf{Z}_n^T\Sigma_n^{1/2})_{ij}=\frac{\sqrt{l_il_j}}{n} 
\sum_{1 \leq k \leq n}{z_{ik}z_{jk}}.\] 
With the aid of \(\mathbf{S}^{D}_n\) and the following inequality (theorem \(4.6.1,\) Vershynin~\cite{vershynin}), our argument runs smoothly: for all \(t \geq 0\) and any random matrix \(A \in \mathbb{R}^{p \times q}\) whose entries are independent, of mean zero, subgaussian with \(\max_{1 \leq i \leq p,1 \leq j \leq q}{||A_{ij}||_{\psi_2}} \leq K,\)
\[\sqrt{p}-cK^2(\sqrt{q}+t) \leq s_q(A) \leq s_1(A) \leq \sqrt{p}+cK^2(\sqrt{q}+t)\]
with probability at least \(1-2\exp(-t^2),\) where \(s_i(\cdot)\) is the \(i^{th}\) largest singular value. These bounds will supply us throughout the proofs not only with the size of the eigenvalues of the random matrices we must deal with, but also with a quantitative control of their fluctuations.  
\par
Proposition~\ref{prop2} relies on the analysis of two separate components because
\[<p_\nu,u_\nu>^2=(1-R_\nu^2) \cdot <a_\nu,e_\nu>^2.\]
For the inner product, the key is the decomposition of the projection of \(a_\nu\) onto \(e_\nu\) from Paul~\cite{paul}:
\[\mathcal{P}_\nu^\perp a_\nu=-\mathcal{R}_\nu \mathcal{D}_\nu a_\nu+(\hat{l}_\nu-l_\nu)\mathcal{R}_\nu a_\nu\]
for 
\[\mathcal{D}_\nu=S_{AA}-\Lambda+\Lambda^{1/2}T^T\mathcal{M}(\hat{l}_\nu I- \mathcal{M})^{-1}T\Lambda^{1/2},\]
yielding \(<a_\nu,e_\nu>^2 \to 1,\) which further implies with equation (\ref{5}) that \(R_\nu^2 \to 0:\) roughly speaking,
\[\frac{R_\nu^2}{1-R_\nu^2}=a^T_\nu\Lambda^{1/2}T^T\mathcal{M}(\hat{l}_\nu I- \mathcal{M})^{-2}T\Lambda^{1/2} a_\nu \approx e^T_\nu\Lambda^{1/2}T^T\mathcal{M}(\hat{l}_\nu I- \mathcal{M})^{-2}T\Lambda^{1/2} e_\nu=\]
\begin{equation}\label{4010}
    =l_\nu t_\nu^T\mathcal{M}(\hat{l}_\nu I- \mathcal{M})^{-2}t_\nu \leq l_\nu \cdot \frac{c}{\hat{l}_\nu^2}=O(\frac{1}{l_\nu})=o(1).
\end{equation}

\subsection{Main Ingredients for an Eigenvalue CLT}\label{3'.1}

In this subsection, we discuss the observations that conduced us to the CLT described in (\ref{6}), and some of the steps turning our formal arguments into rigorous justifications.  
\par
As mentioned in section~\ref{sec2}, equation (\ref{9990}), derived from (\ref{4}), represents the cornerstone of Theorem~\ref{th2}. Since we expect \(a_\nu \approx e_\nu,\) rewrite (\ref{9990}) in the following form:
\begin{equation}\label{3008}
    \sqrt{n} \cdot \frac{\hat{l}_\nu}{l_\nu}=\frac{\sqrt{n}}{l_\nu} \cdot e_\nu^T\Lambda^{1/2}F_\nu\Lambda^{1/2} e_\nu+\frac{2\sqrt{n}}{l_\nu} \cdot  e_\nu^T\Lambda^{1/2}F_\nu\Lambda^{1/2}(a_\nu-e_\nu)+\frac{\sqrt{n}}{l_\nu} \cdot (a_\nu-e_\nu)^T\Lambda^{1/2}F_\nu \Lambda^{1/2} (a_\nu-e_\nu)
\end{equation}
where
\[F_\nu=T^T\mathcal{M}(\hat{l}_\nu I- \mathcal{M})^{-1}T+\frac{1}{n}Z_AZ_A^T.\] 
The first component on the right-hand side of (\ref{3008}) will generate the normal distribution in the CLT:
\[\frac{\sqrt{n}}{l_\nu} \cdot e_\nu^T\Lambda^{1/2}F_\nu \Lambda^{1/2} e_\nu=\frac{1}{\sqrt{n}}(z_\nu^T z_\nu+z_\nu^T H\mathcal{M} (\hat{l}_\nu I-\mathcal{M})^{-1}H^T z_\nu) \approx \sqrt{n}(1+\frac{1}{n}tr(\mathcal{M} (\hat{l}_\nu I-\mathcal{M})^{-1}))+\mathcal{D}_n\]
with 
\[\frac{1}{\sqrt{\mathbb{E}[z_{11}^4]-1}}\mathcal{D}_n:=\frac{1}{\sqrt{n(\mathbb{E}[z_{11}^4]-1})}(z_\nu^T z_\nu-n)=\frac{1}{\sqrt{n(\mathbb{E}[z_{11}^4]-1)}}\sum_{1 \leq i \leq n}{(z_{\nu i}^2-1)} \Rightarrow N(0,1),\] 
while an analysis of the remaining two calls for comprehending \(a_\nu-e_\nu.\) Before examining this difference, we must point out that if \(\lim_{n \to \infty}{\frac{M(n)}{\sqrt{n}}}=0,\) then it can be shown the last two terms in (\ref{3008}) are \(o_p(1),\) yielding immediately a CLT: put differently, what follows in this subsection can be circumvented for this range of \(M(n)\) by employing
\[\frac{\sqrt{n}}{l_\nu} \cdot  |e_\nu^T\Lambda^{1/2}F_\nu \Lambda^{1/2}(a_\nu-e_\nu)| \leq \sqrt{n} \cdot \frac{||a_\nu-e_\nu||^2}{2}+ \sqrt{\frac{n}{l_\nu}} \cdot  ||(F_\nu-I)e_\nu|| \cdot ||\Lambda^{1/2}(a_\nu-e_\nu)||,\]
\[0 \leq \frac{\sqrt{n}}{l_\nu} \cdot (a_\nu-e_\nu)^T\Lambda^{1/2}F_\nu \Lambda^{1/2} (a_\nu-e_\nu) \leq \frac{\sqrt{n}}{l_\nu} \cdot ||F_\nu|| \cdot ||\Lambda^{1/2}(a_\nu-e_\nu)||^2.\]
Inasmuch as even this scenario requires some care (see part \((b)\) of Lemma~\ref{lemma1}, (\ref{25}), and (\ref{33})), we choose not to dwell on it and proceed instead with our scrutiny of \(a_\nu-e_\nu.\)
\par
It can be easily seen that
\begin{equation}\label{4020}
    a_\nu-e_\nu=-\mathcal{R}_\nu \mathcal{D}_\nu e_\nu+r_\nu,
\end{equation}
\[r_\nu=(<a_\nu,e_\nu>-1)e_\nu-M_\nu (a_\nu-e_\nu),\]
where \(M_\nu=\mathcal{R}_\nu \mathcal{D}_\nu-(\hat{l}_\nu-l_\nu)\mathcal{R}_\nu,\) implying that for \(k \ne \nu,\)
\[(a_\nu-e_\nu)_k=(-\mathcal{R}_\nu \mathcal{D}_\nu e_\nu)_k+((-M_\nu) (a_\nu-e_\nu))_k,\]
which further renders a series expansion
\[(a_\nu-e_\nu)_k=(\frac{||a_\nu-e_\nu||^2}{2}-1)\sum_{j \geq 0}{((-M_\nu)^j\mathcal{R}_\nu \mathcal{D}_\nu e_\nu)_k}:=(\frac{||a_\nu-e_\nu||^2}{2}-1) \Sigma_{0,k}.\]
The final touch is recalling \(a_\nu\) has length one, leading to 
\[||a_\nu-e_\nu||^2-\frac{||a_\nu-e_\nu||^4}{4}=\frac{\sum_{k \ne \nu}{\Sigma_{0,k}^2}}{1+\sum_{k \ne \nu}{\Sigma_{0,k}^2}}.\]
As an aside, note that the direction of \(a_\nu\) is unique (\(\hat{l}_\nu\) has multiplicity one), which is reflected in the symmetry displayed by this difference:
\[||a_\nu-e_\nu||^2-\frac{||a_\nu-e_\nu||^4}{4}=\frac{1}{4} \cdot ||a_\nu-e_\nu||^2 \cdot (4-||a_\nu-e_\nu||^2)=\frac{1}{4} \cdot ||a_\nu-e_\nu||^2 \cdot ||-a_\nu-e_\nu||^2.\] 
\par
After rearranging its right-hand side terms, (\ref{3008}) can be rewritten using the series \(\Sigma_{0,k}, k \ne \nu\) instead of \(a_\nu-e_\nu.\) Formally, since \(M_\nu=\mathcal{R}_\nu \mathcal{D}_\nu-(\frac{\hat{l}_\nu}{l_\nu}-1)l_\nu \mathcal{R}_\nu,\) these sums look like \(\sum_{i \geq 0}{\alpha_i(\frac{\hat{l}_\nu}{l_\nu}-1)^i},\) which in turn suggests (\ref{3008}) can be reformulated as
\begin{equation}\label{3009}
    \sqrt{n} \cdot (\frac{\hat{l}_\nu}{l_\nu}-1)=\sqrt{n} \cdot c_{tr}+\mathcal{D}_n+\sqrt{n}\sum_{i \geq 0}{\alpha_i(\frac{\hat{l}_\nu}{l_\nu}-1)^i}+o_p(1)
\end{equation}
for \(c_{tr}=\frac{1}{n} tr(\mathcal{M} (\hat{l}_\nu I-\mathcal{M})^{-1}) \xrightarrow[]{a.s.} 0.\) Given that the size of a power series at zero is dictated by its first term when the argument is sufficiently small, we can speculate this is the case in (\ref{3009}) as well because \(\frac{\hat{l}_\nu}{l_\nu}-1 \xrightarrow[]{a.s.} 0.\) This leads us to the following heuristic: if \(\alpha_i, i \geq 0\) and \(c_{tr}\) are small, then (\ref{3009}) can be maneuvered into a CLT for \(\frac{\hat{l}_\nu}{l_\nu}-1\) in three steps:
\[1. \hspace{0.2cm} \sqrt{n} \cdot (\frac{\hat{l}_\nu}{l_\nu}-1-c_{tr})=\mathcal{D}_n+\sqrt{n}\sum_{i \geq 0}{\alpha_i(\frac{\hat{l}_\nu}{l_\nu}-1-c_{tr})^i}+o_p(1)\]
\[2. \hspace{0.2cm}\sqrt{n} \cdot (\frac{\hat{l}_\nu}{l_\nu}-1-c_{tr}-\overline{\alpha})=\mathcal{D}_n+\sqrt{n}\sum_{i \geq 1}{\alpha_i(\frac{\hat{l}_\nu}{l_\nu}-1-c_{tr}-\overline{\alpha})^i}+o_p(1)\]
for some small \(\overline{\alpha} \in \mathbb{R}\) with \(\overline{\alpha}=\sum_{i \geq 0}{\alpha_i\overline{\alpha}^i},\) giving \(\frac{\hat{l}_\nu}{l_\nu}-1-c_{tr}-\overline{\alpha} \to 0,\) and thus
\[3. \hspace{0.2cm} \sqrt{n} \cdot (\frac{\hat{l}_\nu}{l_\nu}-1-c_{tr}-\overline{\alpha})=\mathcal{D}_n+o_p(1) \cdot \sqrt{n} \cdot (\frac{\hat{l}_\nu}{l_\nu}-1-c_{tr}-\overline{\alpha})+o_p(1),\]
which in conjunction with \(\frac{1}{\sqrt{\mathbb{E}[z_{11}^4]-1}}\mathcal{D}_n \Rightarrow N(0,1)\) and Slutsky's lemma yields
\[\sqrt{\frac{n}{\mathbb{E}[z_{11}^4]-1}} \cdot (\frac{\hat{l}_\nu}{l_\nu}-1-c_{tr}-\overline{\alpha}) \Rightarrow N(0,1).\]
Making this daydream mathematically rigorous completes the proof of Theorem~\ref{th2}. Specifically, we truncate the series \(\Sigma_{0,k}, k \ne \nu,\) leading to \(\overline{\alpha}\) being the root of a polynomial (the content of Theorem~\ref{th4}), and show there exist deterministic coefficients \(\alpha_i,\) growing sufficiently slow and rendering an identity analogous to (\ref{3009}) (since \(M_\nu\) is random, it is not a priori clear that this is even possible).
\par
Consider first the coefficients \(\alpha_i.\) Notice that the sums \(\Sigma_{0,k}, k \ne \nu\) are up to a large extent determined by the powers of \(M_\nu.\) Inasmuch as 
\[M_\nu=\mathcal{R}_\nu \mathcal{D}_\nu-(\hat{l}_\nu-l_\nu)\mathcal{R}_\nu=\mathcal{R}_\nu \Lambda^{1/2}(\Tilde{D}+\Tilde{T})\Lambda^{1/2}-(\frac{\hat{l}_\nu}{l_\nu}-1) \cdot l_\nu \mathcal{R}_\nu\]
where
\[\Tilde{D}=\frac{1}{n}Z_AZ_A^T-I, \hspace{0.2cm} \Tilde{T}=T^T\mathcal{M}(\hat{l}_\nu I- \mathcal{M})^{-1}T,\]
the component \(\mathcal{R}_\nu \Lambda^{1/2}\Tilde{T}\Lambda^{1/2}\) seems the most delicate out of the three: it depends on \(\hat{l}_\nu,\) which in turn is dependent on everything else. This is a serious impediment because generally speaking sums of products of random vectors and matrices are successfully handled when the two are independent. Besides, note that the multinomial expansions for these powers generate two categories of terms, those containing at least a factor of \(\tilde{T}\) (call them \textit{first type}), and the rest (\textit{second type}). Given the potential intricacy of dealing with \(\hat{l}_\nu\) in the first type terms, we introduce the growth condition of \(l_\nu\) which turns them negligible (i.e., their size is \(o_p(1)\) in (\ref{3008}); moreover, this threshold is optimal in the sense that below it these summands contribute to the right-hand side term in (\ref{3008})) and allows us to focus on the second type terms. Once we unpack these components, we see their main constituents are products of the form
\begin{equation}\label{3010}
    \Tilde{D}_{\nu k_1}\Tilde{D}_{k_1 k_2}...\Tilde{D}_{k_{j-1} k_j} \Tilde{D}_{k_j \nu}
\end{equation}
for \(k_1, \hspace{0.05cm} ... \hspace{0.05cm}, k_j \in \{1,2, \hspace{0.05cm} ... \hspace{0.05cm}, M\} -\{\nu\},\) most of which have expectation \(\frac{1}{n^j}:\) namely, this is the case whenever \(k_1, \hspace{0.05cm} ... \hspace{0.05cm}, k_j\) are pairwise distinct because for \(\kappa_1 \ne \kappa_2,\)
\[\Tilde{D}_{\kappa_1 \kappa_2}=\frac{1}{n}z_{\kappa_1}^Tz_{\kappa_2}=\frac{1}{n}\sum_{1 \leq i \leq n}{z_{\kappa_1 i} z_{\kappa_2 i}}\] 
and thus the terms in (\ref{3010}) with non-zero expectation have the same index \(i\) in all its \(j+1\) factors, from which 
\[\mathbb{E}[\Tilde{D}_{\nu k_1}\Tilde{D}_{k_1 k_2}...\Tilde{D}_{k_{j-1} k_j} \Tilde{D}_{k_j \nu}]=\frac{1}{n^{j+1}} \cdot n=\frac{1}{n^j},\]
centering which will define the polynomial coefficients \(\alpha_i.\) 
\par
Consider next the errors, which naturally suggest us to analyze the following sums
\begin{equation}\label{4040}
    \sum_{(k_1, k_2, ... , k_j),(k'_1, k'_2, ..., k'_t), k_i,k'_i \ne \nu}{|\mathbb{E}[(\Tilde{D}_{k_1 \nu}\Tilde{D}_{k_1 k_2}...\Tilde{D}_{k_{j-1} k_j} \Tilde{D}_{k_j \nu}-\frac{1}{n^{j}}) \cdot (\Tilde{D}_{k'_1 \nu}\Tilde{D}_{k'_1 k'_2}...\Tilde{D}_{k'_{t-1} k'_t} \Tilde{D}_{k'_t \nu}-\frac{1}{n^{t}})]|},
\end{equation}
tightly connected with the squares of our objects of interest. Namely, we show the expectations of the latter are \(o(\frac{1}{n}),\) which would then render the desired conclusion using Chebyshev's inequality. Most of these terms are negligible: for \(k_1,k_2, \hspace{0.05cm} ... \hspace{0.05cm}, k_j, k'_1, k'_2, \hspace{0.05cm} ...\hspace{0.05cm}, k'_t\) pairwise distinct, 
\[\mathbb{E}[(\Tilde{D}_{k_1 \nu}\Tilde{D}_{k_1 k_2}...\Tilde{D}_{k_{j-1} k_j} \Tilde{D}_{k_j \nu}-\frac{1}{n^{j}}) \cdot (\Tilde{D}_{k'_1 \nu}\Tilde{D}_{k'_1 k'_2}...\Tilde{D}_{k'_{t-1} k'_t} \Tilde{D}_{k'_t \nu}-\frac{1}{n^{t}})]=\]
\[=\frac{1}{n^{j+t+2}}\sum_{1 \leq i, j \leq n}{\mathbb{E}[z_{\nu i}^2z_{\nu j}^2]}-\frac{1}{n^{j+t}}=\frac{1}{n^{j+t+1}} \cdot (\mathbb{E}[z_{11}^4]-1),\]
and since their number of such tuples is at most \(M^{j+t},\) their contribution is at most \(\frac{M^{j+t}}{n^{j+t+1}} \cdot (\mathbb{E}[z_{11}^4]-1)=o(\frac{1}{n}).\) Overall, it will be seen the sum in (\ref{4040}) is 
\[O(\frac{c(j,t)}{n} \cdot (\frac{M}{n})^{(j+t)/2}):\]
because \(c(j,t)\) increases with \(j+t,\) these bounds pile up when \(j\) or \(t\) is large. To avoid such a blowup, we truncate the series \(\Sigma_{0,k}, k \ne \nu\) and obtain bounds for (\ref{4040}) with relatively small \(j,t,\) finishing thus the proof of our first two theorems.

\subsection{Refining CLT Centerings, and Eigenvector Consistency}\label{3'.2}

With the foundations laid by the CLT in (\ref{6}) and the consistency in Proposition~\ref{prop2}, we proceed to polish the centering of the former and derive the rates of the latter. Before embarking on these modifications, a couple of remarks on altering the CLT centerings are in order: when \(\lim_{n \to \infty}{\frac{M(n)}{\sqrt{n}}}=0,\) Theorem~\ref{th2} already yields an empirical centering (\(x_{n,\nu,(l_i)_{1 \leq i \leq M}}\) can be replaced by zero), whereas no alike simplification occurs in Theorem~\ref{th3} for a relatively small number of spikes (the heart of the matter is the trace component which is not more manageable in this part of the range of \(M(n)\) than in the rest of it).
\par
Theorem~\ref{th1} aims for an empirical centering, and the substitution mentioned in subsection~\ref{2.1} 
\[\frac{1}{n}\sum_{k \ne \nu}{\frac{l_k}{l_k-l_\nu}} \hspace{0.2cm} \longrightarrow \hspace{0.2cm} \frac{1}{n}\sum_{k \ne \nu}{\frac{\hat{l}_k}{\hat{l}_k-\hat{l}_\nu}}\]
is implemented in two phases, replacing first \(l_\nu\) by \(\hat{l}_\nu,\) and second \(l_k, k \ne \nu\) by \(\hat{l}_k, k \ne \nu.\) The first swap occurs almost effortlessly due to the fluctuations of \(\frac{\hat{l}_\nu}{l_\nu}-1\) captured by the CLT in Theorem~\ref{th2}. The second, nonetheless, is more intricate because it deals with \(\hat{l}_k, k \ne \nu,\) random variables depending on each other in a complex way. Thus, to accomplish this change, we impose the fluctuations of each to be sufficiently small, leading to the stronger assumption on the growth of \(M.\)
\par
Theorem~\ref{th3} finds a deterministic replacement \(\xi_\nu\) for the trace component:
\begin{equation}\label{4060}
    \sqrt{n} \cdot (\frac{1}{n}tr(\mathcal{M}(\hat{l}_\nu I-\mathcal{M})^{-1})-\xi_\nu) \xrightarrow[]{p} 0.
\end{equation}
Since the first term is closely related to the Stieltjes transform of the empirical distribution of \(S_{BB},\) 
\[m_{F_n}(z)=\frac{1}{N-M}tr((zI-S_{BB})^{-1})=\frac{1}{N-M}tr((zI-\mathcal{M})^{-1})\]
for \(z \in \mathbb{C}-\{m_1, m_2,\hspace{0.05cm} ... \hspace{0.05cm}, m_{N-M}: \mathcal{M}=diag(m_1, m_2, \hspace{0.05cm} ... \hspace{0.05cm}, m_{N-M})\},\) rewrite (\ref{4060}) for \(\newline \gamma_n=\frac{N-M}{n}\) as
\begin{equation}\label{4050}
    \sqrt{n} \cdot (\gamma_n \cdot (\hat{l}_\nu m_{F_n}(\hat{l}_\nu)-1)-\xi_\nu) \xrightarrow[]{p} 0.
\end{equation}
A particular case of the Marchenko-Pastur result in \cite{marchenkopastur}, \(m_{F_n}(z) \xrightarrow[]{p} m_{\gamma_0}(z)\) when \(\gamma_n \to \gamma_0\) for \(z \in \mathbb{C}^{+},\)  where for \(\gamma_0>0,\)
\[m_{\gamma_0}(z)=\frac{z+\gamma_0-1-\sqrt{(z-\gamma_0+1)^2-4z}}{2z\gamma_0}, \hspace{0.2cm} z \ne 0,\]
suggests looking for a deterministic estimate \(\xi_0\) of \(\hat{l}_\nu\) and then take 
\[\xi_\nu = \gamma_n \cdot (\xi_0 \cdot m_{\gamma_n}(\xi_0)-1).\] Ignoring everything in our CLT that depends on the rest of the spikes and employing the approximation \(m_{F_n} \approx m_{\gamma_n},\) we have
\[l_\nu \approx \frac{\hat{l}_\nu}{1+\frac{1}{n}tr(\mathcal{M}(\hat{l}_\nu I-\mathcal{M})^{-1})},\]
\[l_\nu \approx \frac{\hat{l}_\nu}{1+\gamma_n \cdot (\hat{l}_\nu \cdot m_{\gamma_n}(\hat{l}_\nu)-1)}.\]
Therefore, what is needed is inverting the function 
\[f(x)=\frac{x}{1+\gamma_n \cdot (x \cdot \frac{x+\gamma_n-1-\sqrt{(x-\gamma_n+1)^2-4x}}{2\gamma_n x}-1)},\]
which is to say, expressing \(\hat{l}_\nu\) in terms of \(l_\nu.\) Notice that the square root is the main source of complexity in \(f:\) to dispense with it, take \(x-\gamma_n+1=a+b,x=ab\) with \(a \leq b,\) from which
\[x=b \cdot (1+\frac{\gamma_n}{b-1}), \hspace{0.3cm}  f(x)=\frac{x}{1+\gamma_n \cdot (x \cdot \frac{a+b+2\gamma_n-2-(b-a)}{2\gamma_n x}-1)}=\frac{x}{a}=b.\]
In other words, the following candidates present themselves
\[\hat{l}_\nu \approx \xi_0=l_\nu \cdot (1+\frac{\gamma_n}{l_\nu-1}), \hspace{0.2cm} \xi_\nu=\gamma_n \cdot (\xi_0 \cdot m_n(\xi_0)-1)=\frac{\gamma_n}{l_\nu-1}.\]
For this value of \(\xi_\nu,\) after manipulating (\ref{4050}), it remains to show that for deterministic \(\newline z=z_n \to \infty,\)
\[\sqrt{n} \cdot (m_{F_n}(z)-m_n(z)) \xrightarrow[]{p} 0,\]
which can be justified using one of the ideas behind the proof of the Marchenko-Pastur result in \cite{marchenkopastur} employed anew by Ledoit and Péché in \cite{ledoitpeche}.
\par
Theorem~\ref{th5} sharpens the convergences obtained in Proposition~\ref{prop2},
\[R_\nu^2 \xrightarrow[]{a.s.} 0, \hspace{0.2cm} <a_\nu, e_\nu>^2 \xrightarrow[]{a.s.} 1\] 
since the object of interest remains the same
\[<p_\nu,u_\nu>^2=(1-R_\nu^2) \cdot <a_\nu,e_\nu>^2.\]
Specifically, we find that \(l_\nu R_\nu^2-\frac{N}{n} \xrightarrow[]{a.s.} 0:\) given (\ref{4010}), we expect \(R_\nu^2=O(\frac{1}{l_\nu}),\) and the ratio \(\frac{N}{n}\) arises from 
\[t_\nu^T\mathcal{M}t_\nu-\frac{N}{n}=\frac{1}{n} z_\nu^T H\mathcal{M} H^T z_\nu-\frac{N}{n} \xrightarrow[]{a.s.} 0,\] 
a consequence of Hanson-Wright inequality (\ref{9}) and
\[\frac{1}{n}tr(H\mathcal{M} H^T)=\frac{1}{n}tr(\mathcal{M})=\frac{1}{n}tr(S_{BB})=\frac{N-M}{n}+\frac{1}{n^2}\sum_{M+1 \leq i \leq N, 1 \leq j \leq n}{(z_{ij}^2-1)}=\frac{N}{n}+o_{a.s.}(1).\] 
\par
As regards the second constituent, \(a_\nu \approx e_\nu\) and (\ref{4020}) yield \(||a_\nu-e_\nu|| \approx ||\mathcal{R}_\nu \mathcal{D}_\nu e_\nu||,\) 
\[<a_\nu,e_\nu>^2=\frac{1}{4}(||a_\nu-e_\nu||^2-2)^2=1+\frac{||a_\nu-e_\nu||^4}{4}-||a_\nu-e_\nu||^2 \approx 1-||\mathcal{R}_\nu \mathcal{D}_\nu e_\nu||^2+o(||\mathcal{R}_\nu \mathcal{D}_\nu e_\nu||^2).\]
Let us attend the dominating term,
\[||\mathcal{R}_\nu \mathcal{D}_\nu e_\nu||^2=\sum_{k \ne \nu}{\frac{l_kl_\nu}{(l_k-l_\nu)^2}(\frac{1}{n}z_k^Tz_\nu+t_k^T\mathcal{M}(\hat{l}_\nu I-\mathcal{M})^{-1}t_\nu)^2},\] 
for which we forecast that
\[\sum_{k \ne \nu}{\frac{l_kl_\nu}{(l_k-l_\nu)^2}(\frac{1}{n}z_k^Tz_\nu)^2}\]
generates the main contribution and notice that Linderberg's CLT yields for \(k \ne \nu,\) 
\[\frac{1}{\sqrt{n}}z_k^Tz_\nu=\frac{1}{\sqrt{n}}\sum_{1 \leq i \leq n}{z_{ki} z_{\nu i}} \Rightarrow N(0,1).\]
\par
If \(M\) is fixed, then computing the limits of the moments of
\begin{equation}\label{2021}
    \sum_{k \ne \nu}{\frac{l_kl_\nu}{(l_k-l_\nu)^2}(\frac{1}{\sqrt{n}}z_k^Tz_\nu)^2}
\end{equation}
together with Carleman's condition render that asymptotically (\ref{2021}) is a weighted sum of mutually independent chi-distributions with weights given by the limits of the ratios 
\[\frac{l_kl_\nu}{(l_k-l_\nu)^2}, \hspace{0.1cm} k \ne \nu,\] 
explaining as well as how the size of \(l_\nu\) relative to \(n\) comes into play in \((b)(i)\) and \((c)(i).\) Moreover, Crámer-Wold theorem implies  \(\frac{1}{\sqrt{n}}z_k^Tz_\nu, k \ne \nu, 1 \leq k \leq M\) are asymptotically mutually independent, which in conjunction with our analysis of \(a_\nu-e_\nu\) in Theorem~\ref{th2} can be used to infer the limiting behavior of this difference. 
\par
If \(M\) grows to infinity with \(n,\) then tackling (\ref{2021}) requires a different approach. In virtue of what is known about large sums of i.i.d. random variables, two normalizations of (\ref{2021}) are presumably fruitful
\begin{equation}\label{2022}
    \frac{1}{M}\sum_{k \ne \nu}{\frac{l_kl_\nu}{(l_k-l_\nu)^2} \cdot (\frac{1}{\sqrt{n}}z_k^Tz_\nu)^2}
\end{equation}
\begin{equation}\label{2023}
    \frac{1}{\sqrt{M}}\sum_{k \ne \nu}{\frac{l_kl_\nu}{(l_k-l_\nu)^2} \cdot ((\frac{1}{\sqrt{n}}z_k^Tz_\nu)^2-1)}
\end{equation}
yielding results in the spirit of an SLLN and a CLT, respectively. We treat (\ref{2022}) with the aid of moments and (\ref{2023}) with a martingale CLT. Since \(\{z_k, k \ne \nu, z_\nu\}\) are independent, the conditional expectation of (\ref{2023}) on any subset of it can be easily obtained: in particular, for any \(C \subset \{1,2, \hspace{0.05cm} ... \hspace{0.05cm}, M\}-\{\nu\},\)
\[\mathbb{E}[\frac{1}{\sqrt{M}}\sum_{k \ne \nu}{\frac{l_kl_\nu}{(l_k-l_\nu)^2} \cdot ((\frac{1}{\sqrt{n}}z_k^Tz_\nu)^2-1)} \hspace{0.05cm} | \hspace{0.05cm} z_\nu, z_c, c \in C]=\]
\[=\frac{1}{\sqrt{M}}\sum_{k \in C }{\frac{l_kl_\nu}{(l_k-l_\nu)^2} \cdot ((\frac{1}{\sqrt{n}}z_k^Tz_\nu)^2-1)}+\frac{1}{\sqrt{M}}\sum_{k \ne \nu, k \not \in C}{\frac{l_kl_\nu}{(l_k-l_\nu)^2} \cdot (\frac{1}{n}||z_\nu||^2-1)},\]
leading to a martingale representation of (\ref{2023}) amenable to computations. Finally, Theorem~\ref{th6} can be obtained along similar lines of reasoning to the ones giving Theorem~\ref{th1}.
\par 
The rest of the paper in organized as follows: section~\ref{sec3} contains the justification of the eigenstructure consistency for divergent spikes; the next three sections consist of the proofs of the theorems:
\begin{itemize}
    \item Section~\ref{sec4}: Theorems \ref{th2}, and \ref{th4},
    \item Section~\ref{sec5}: Theorems \ref{th1}, and \ref{th3},
    \item Section~\ref{sec6}: Theorems \ref{th5}, and \ref{th6};
\end{itemize}
section~\ref{sec7} presents three auxiliary lemmas; the appendix contains the concentration inequalities employed repeatedly throughout the proofs.
\par
\vspace{0.2cm}
\textbf{Acknowledgements:} I would like to thank professors George Papanicolaou and Lenya Ryzhik for their comments and suggestions, especially for the feedback concerning the expository aspects of this paper.   

\section{Eigenstructure Consistency}\label{sec3}

Subsections \ref{3.1} and \ref{3.2} contain the proofs of Proposition~\ref{prop1}, and Proposition~\ref{prop2}, respectively. Both rely on inequality (\ref{7}), while the latter uses also identities (\ref{4}) and (\ref{5}) which provide connections between \(R_\nu^2, \hat{l}_\nu,\) and multidimensional objects, matrices such as \(S_{AA},\Lambda,\) and the vector \(a_\nu.\)

\subsection{Proof of Proposition~\ref{prop1}}\label{3.1}

We will use ideas from Shen et al.~\cite{shen} to show this result (the authors assume, among other conditions, that \(M(n)\) is fixed). Consider the dual matrices of \(\mathbf{S}_n,S_{AA},S_{BB},\)
\[\mathbf{S}^D_n=\frac{1}{n}\mathbf{X}_n^T\mathbf{X}_n=\frac{1}{n}\sum_{1 \leq i \leq N}{l_iZ_i^TZ_i}, \hspace{0.1cm} \mathcal{A}=\frac{1}{n}\sum_{1 \leq i \leq M}{l_iZ_i^TZ_i}, \hspace{0.1cm}  \mathcal{B}=\frac{1}{n}\sum_{M+1 \leq i \leq N}{Z_i^TZ_i},\]
where \(Z_i \in \mathbb{R}^{1 \times n}, 1 \leq i \leq N\) are the rows of \(\mathbf{Z}_n.\)
\par
Denote by \(\lambda_k(J)\) the \(k^{th}\) largest eigenvalue of a symmetric matrix \(J.\) We prove first that for \(1 \leq k \leq \nu,\)
\begin{equation}\label{11}
    \frac{\lambda_k(\mathcal{A})}{l_k} \xrightarrow[]{a.s.} 1.
\end{equation}
\par
Inequality (\ref{7}) yields
\[\frac{\lambda_1(\mathcal{A})}{l_1} \leq \lambda_1(\frac{1}{n}\sum_{1 \leq i \leq M}{Z_i^TZ_i}) \leq \frac{1}{n}(\sqrt{n}+cK^2(\sqrt{M}+n^{1/4}))^2\]
with probability at least \(1-2\exp(-\sqrt{n}),\) which implies with Borel-Cantelli lemma that almost surely 
\[\limsup_{n \to \infty}{\frac{\lambda_1(\mathcal{A})}{l_1}} \leq 1;\]
since \(|d_{11}| \leq \lambda_1(D)\) for \(D \in \mathbb{R}^{p \times p},\) and 
\[\frac{1}{l_1}(S_{AA})_{11}=\frac{1}{n}\sum_{1 \leq i \leq n}{z_{1i}^2} \xrightarrow[]{a.s.} 1\] 
from Borel-Cantelli lemma and
\begin{equation}\label{726}
    \mathbb{E}[(\frac{1}{n}\sum_{1 \leq i \leq n}{(z_{1i}^2}-1))^4]=\frac{1}{n^3}\mathbb{E}[(z_{11}^2-1)^4]+\frac{6}{n^4} \cdot \binom{n}{2} (\mathbb{E}[(z_{11}^2-1)^2])^2 \leq \frac{c(K)}{n^2},
\end{equation}
we obtain (\ref{11}) for \(k=1:\)
\begin{equation}\label{2024}
    \frac{\lambda_1(\mathcal{A})}{l_1}=\frac{\lambda_1(S_{AA})}{l_1} \xrightarrow[]{a.s.} 1.
\end{equation}
\par 
Consider now \(2 \leq k \leq \nu\) for which
\begin{equation}\label{781}
    \lambda_k(\frac{1}{n}\sum_{1 \leq i \leq k}{Z_i^TZ_i}) \leq \frac{\lambda_k(\mathcal{A})}{l_k} \leq \lambda_1(\frac{1}{n}\sum_{k \leq i \leq M}{Z_i^TZ_i}).
\end{equation}
To justify this chain of inequalities, notice the lower bound is clear from \(l_i \geq l_k\) for \(i \leq k,\) while for the upper bound
\begin{equation}\label{782}
    \lambda_k(\mathcal{A}) \leq \lambda_k(\frac{1}{n}\sum_{1 \leq i \leq k-1}{l_iZ^T_iZ_i})+\lambda_1(\frac{1}{n}\sum_{k \leq i \leq M}{l_iZ^T_iZ_i}) \leq l_k \cdot \lambda_1(\frac{1}{n}\sum_{k \leq i \leq M}{Z_i^TZ_i}),
\end{equation}
where the first claim ensues from Courant-Fischer min-max formula for symmetric matrices \(A \in \mathbb{R}^{p \times p}\)
\[\lambda_k(A)=\min_{w_1,w_2, \hspace{0.05cm} ... \hspace{0.05cm}, w_{k-1}}{\hspace{0.05cm}\max_{||w||=1, w \perp w_1, w_2, \hspace{0.05cm} ... \hspace{0.05cm}, w_{k-1}}{w^TAw}}\]
for \(1 \leq k \leq p, w_1, w_2, \hspace{0.05cm} ... \hspace{0.05cm}, w_k, w \in \mathbb{R}^p\) (see for instance theorem \(A.7\) in Bai and Silverstein~[\ref{baisilvbook}]), and the second follows from
\[\lambda_k(\frac{1}{n}\sum_{1 \leq i \leq k-1}{l_iZ^T_iZ_i})=0\] 
(the matrix \(\frac{1}{n}\sum_{1 \leq i \leq k-1}{l_iZ^T_iZ_i}\) has rank at most \(k-1<k\)), and \(l_i \leq l_k\) for \(i \geq k.\)
\par
For any \(t>0, n \geq n(t), 2 \leq k \leq M,\) (\ref{7}) entails each of the following two events occurs with probability at least \(1-2\exp(-\sqrt{n})\)
\[\lambda_k(\frac{1}{n}\sum_{1 \leq i \leq k}{Z^T_iZ_i}) \geq \frac{1}{n}(\sqrt{n}-cK^2(\sqrt{k}+n^{1/4}))^2 \geq \frac{1}{n}(\sqrt{n}-cK^2(\sqrt{M}+n^{1/4}))^2 \geq 1-t,\]
\[\lambda_1(\frac{1}{n}\sum_{k \leq i \leq M}{Z_i^TZ_i}) \leq \frac{1}{n}(\sqrt{n}+cK^2(\sqrt{M-k+1}+n^{1/4}))^2 \leq \frac{1}{n}(\sqrt{n}+cK^2(\sqrt{M}+n^{1/4}))^2 \leq 1+t,\]
which together with (\ref{781}) renders
\[\mathbb{P}(|\frac{\lambda_k(\mathcal{A})}{l_k}-1| \geq t) \leq 4\exp(-\sqrt{n}).\]
Since \(4M\exp(-\sqrt{n}) \leq \exp(-n^{1/4})\) for \(n\) large enough, (\ref{11}) ensues from this last inequality, Borel-Cantelli lemma, and (\ref{2024}).
\par
Finally, notice that with probability one, \(||\mathcal{B}||\) is bounded (again using (\ref{7})), and thus the conclusion of the proposition follows from (\ref{11}), \(\lim_{n \to \infty}{l^{(n)}_k}=\infty,\) and
\begin{equation}\label{12}
    \frac{\lambda_k(\mathcal{A})}{l_k} \leq \frac{\hat{l}_k}{l_k}=\frac{\lambda_k(\mathcal{A}+\mathcal{B})}{l_k} \leq \frac{\lambda_k(\mathcal{A})+\lambda_1(\mathcal{B})}{l_k}.
\end{equation}

\subsection{Proof of Proposition~\ref{prop2}}\label{3.2}

For \(e_1,e_2, \hspace{0.05cm} ... \hspace{0.05cm}, e_M \in \mathbb{R}^M\) the standard basis,
\begin{equation}\label{49}
    <p_\nu,u_\nu>^2=(1-R_\nu^2)<a_\nu,e_\nu>^2,
\end{equation}
and so the conclusion is equivalent to
\[R_\nu^2 \xrightarrow[]{a.s.} 0, \hspace{0.5cm} <a_\nu, e_\nu>^2 \xrightarrow[]{a.s.} 1.\]
In light of Proposition~\ref{prop1} and 
\[||S_{BB}||=||\mathcal{M}|| \leq c_{K,\gamma}:=2cK^2(1+\sqrt{\gamma})^2\] 
almost surely (from (\ref{7})), it follows that \(\hat{l}_\nu I-S_{BB}\) is invertible, and thus (\ref{4}) and (\ref{5}) hold with probability one. Having justified the validity of these equations, we begin with the second convergence, \(<a_\nu, e_\nu>^2,\) and continue with the first, \( R_\nu^2.\)

\vspace{0.5cm}
\par
\textit{Second term, \(<a_\nu,e_\nu>^2:\)} As in Paul~\cite{paul},
\begin{equation}\label{791}
    \mathcal{P}_\nu^\perp a_\nu=-\mathcal{R}_\nu \mathcal{D}_\nu a_\nu+(\hat{l}_\nu-l_\nu)\mathcal{R}_\nu a_\nu,
\end{equation}
for \(\mathcal{P}^\perp_\nu=I_M-e_\nu e_\nu^T,\) \(\mathcal{R}_\nu \in \mathbb{R}^{M \times M}\) the diagonal matrix whose \(k^{th}\) diagonal entry \(\frac{1}{l_k-l_\nu}\) for \(k \ne \nu,\) and \(0\) for \(k=\nu,\) and 
\[\mathcal{D}_\nu=S_{AA}-\Lambda+\Lambda^{1/2}T^T\mathcal{M}(\hat{l}_\nu I- \mathcal{M})^{-1}T\Lambda^{1/2}\]
since 
\[\mathcal{P}^\perp_\nu a_\nu=(I_M-e_\nu e_\nu ^T)a_\nu=\mathcal{R}_\nu (\Lambda-l_\nu I)a_\nu=\mathcal{R}_\nu (-\mathcal{D}_\nu+(\mathcal{D}_\nu+\Lambda)-l_\nu I)a_\nu\]
from equation (\ref{4}), \((\mathcal{D}_\nu+\Lambda)a_\nu=\hat{l}_\nu a_\nu.\)
\par
It suffices to show that \(\alpha_\nu := ||\mathcal{R}_\nu\mathcal{D}_\nu||+|\hat{l}_\nu-l_\nu|\cdot ||\mathcal{R}_\nu|| \xrightarrow[]{a.s.} 1\) because this and (\ref{791}) will render
\begin{equation}\label{2025}
    <a_\nu,e_\nu>^2=1-(\mathcal{P}_\nu ^\perp a_\nu)^2 \xrightarrow[]{a.s.} 1.
\end{equation}
Notice that
\[|\hat{l}_\nu-l_\nu|\cdot ||\mathcal{R}_\nu||=|\hat{l}_\nu-l_\nu|\cdot \frac{1}{\min_{k \ne \nu, k \leq M}{|l_k-l_\nu|}}=\frac{|\frac{\hat{l}_\nu}{l_\nu}-1|}{\min_{k \ne \nu, k \leq M}{|\frac{l_k}{l_\nu}-1|}} \xrightarrow[]{a.s.} 0,\]
\[||\mathcal{R}_\nu\mathcal{D}_\nu||=||\Lambda^{1/2} \mathcal{R}_\nu \Lambda^{1/2}(\frac{1}{n}Z_AZ_A^T-I_M+T^T\mathcal{M}(\hat{l}_\nu I- \mathcal{M})^{-1}T)|| \leq\]
\[\leq ||\Lambda^{1/2} \mathcal{R}_\nu \Lambda^{1/2}|| \cdot ||\frac{1}{n}Z_AZ_A^T-I_M+T^T\mathcal{M}(\hat{l}_\nu I- \mathcal{M})^{-1}T|| \xrightarrow[]{a.s.} 0,\] 
because 
\[||\Lambda^{1/2} \mathcal{R}_\nu \Lambda^{1/2}||=\max_{k \ne \nu, k \leq M}{\frac{l_k}{|l_k-l_\nu|}}=\max{(\frac{1}{\frac{l_\nu}{l_{\nu+1}}-1}, \frac{1}{1-\frac{l_\nu}{l_{\nu-1}}})} \leq c(\epsilon_0),\]
while inequality (\ref{7}) and Proposition~\ref{prop1} yield that almost surely \(||\frac{1}{n}Z_AZ_A^T|| \to 1, ||\mathcal{M}|| \leq c_{K,\gamma}\) for \(n\) large enough, and \(\hat{l}_\nu \to \infty,\) from which
\[||T^T\mathcal{M}(\hat{l}_\nu I- \mathcal{M})^{-1}T|| \leq ||TT^T|| \cdot ||\mathcal{M}(\hat{l}_\nu I- \mathcal{M})^{-1}|| \leq \frac{2c_{K,\gamma}}{\hat{l}_v-c_{K,\gamma}} \xrightarrow[]{a.s.} 0,\]
as \(||TT^T|| \leq ||\frac{1}{n}Z_AZ_A^T||.\)

\vspace{0.5cm}
\par
\textit{First term, \(R_\nu^2\):} Equations (\ref{4}) and (\ref{5}) give that with probability one
\begin{equation}\label{794}
    a^T_\nu \Lambda^{1/2}(\frac{1}{n}Z_AZ_A^T+T^T\mathcal{M}(\hat{l}_\nu I- \mathcal{M})^{-1}T)\Lambda^{1/2} a_\nu=\hat{l}_\nu,
\end{equation}
\begin{equation}\label{793}
    a^T_\nu \Lambda^{1/2}T^T\mathcal{M}(\hat{l}_\nu I- \mathcal{M})^{-2}T\Lambda^{1/2} a_\nu=\frac{R_\nu^2}{1-R_\nu^2}.
\end{equation}
\par
Because \(\mathcal{M}\) is diagonal with non-negative entries, and almost surely for \(n\) large \(||\mathcal{M}|| \leq c_{K,\gamma} \leq \frac{\hat{l}_\nu}{2},\) it follows from (\ref{793}) that
\[\frac{R_\nu^2}{1-R_\nu^2} \leq \frac{2}{\hat{l}_\nu} \cdot a^T_\nu \Lambda^{1/2}T^T\mathcal{M}(\hat{l}_\nu I- \mathcal{M})^{-1}T\Lambda^{1/2} a_\nu,\]
as \(\frac{x}{(l-x)^2} \leq \frac{2x}{l(l-x)}\) for \(0 \leq x \leq \frac{l}{2}, l>0.\) Furthermore, (\ref{794}) implies
\[\frac{1}{\hat{l}_\nu}a^T_\nu \Lambda^{1/2}T^T\mathcal{M}(\hat{l}_\nu I- \mathcal{M})^{-1}T\Lambda^{1/2} a_\nu=1-\frac{1}{\hat{l}_\nu}a^T_\nu \Lambda^{1/2}(\frac{1}{n}Z_AZ_A^T)\Lambda^{1/2} a_\nu \leq 1-\lambda_{\min}(\frac{1}{n}Z_AZ_A^T) \cdot \frac{l_\nu}{\hat{l}_\nu}<a_\nu,e_\nu>^2.\]
\par
Hence
\[0 \leq R_\nu^2 \leq \frac{R_\nu^2}{1-R_\nu^2} \leq 2 \cdot (1-\lambda_{\min}(\frac{1}{n}Z_AZ_A^T) \cdot \frac{l_\nu}{\hat{l}_\nu}<a_\nu,e_\nu>^2) \xrightarrow[]{a.s.} 0,\]
where we have used \(\lambda_{\min}(\frac{1}{n}Z_AZ_A^T) \xrightarrow[]{a.s.} 1\) (from (\ref{7})), Proposition~\ref{prop1}, and (\ref{2025}).

\section{Primary CLT for Eigenvalues}\label{sec4}

Subsection \ref{4.1} presents the crux of the proof of Theorem~\ref{th2}, and the full justification of this result will be covered in the coming eight subsections. An outline of the method employed is as follows: equation~(\ref{4}) offers a decomposition of \(\frac{\hat{l}_\nu}{l_\nu}-1\) in three terms depending on \(a_\nu-e_\nu.\) Once the entries of this difference are decomposed into series, three identities need to be shown, (\ref{21}), (\ref{22}), and (\ref{23}). Subsection \ref{4.2} presents the proofs of (\ref{21}) and (\ref{22}). The last missing piece, (\ref{23}), requires most of the work: we reduce it to finding polynomial representations in \(\frac{\hat{l}_\nu}{l_\nu}-1\) of three sums (\(\Tilde{\Sigma}_0, \Tilde{\Sigma}_1, \Tilde{\Sigma}_2\) defined below), up to \(o_p(\frac{1}{\sqrt{n}})\) errors, in subsection \ref{4.3}; in other words, (\ref{23}) ensues from three new identities, (\ref{27}), (\ref{28}), (\ref{29}). Next, for each \(\Tilde{\Sigma}_i, 0 \leq i \leq 2,\) we obtain these polynomial decompositions in three stages: 
\begin{itemize}
    \item Stage I: truncate the underlying series coming from the entries of \(a_\nu-e_\nu,\) being left with computing finitely many powers of a sum of three random matrices,
    
    \item Stage II: show the terms containing the first matrix (called \textit{first type}) in these multinomial expansions are negligible (i.e., \(o_p(\frac{1}{\sqrt{n}})\)),
    
    \item Stage III: prove the contribution of the rest of the summands (called \textit{second type}) can be replaced, up to an \(o_p(\frac{1}{\sqrt{n}})\) error, by deterministic quantities.
\end{itemize}
Subsections \ref{4.4}, \ref{4.5}, and \ref{4.6} present the justification of (\ref{27}), split in these three phases; subsection \ref{4.7} consists of a combinatorial result that allows us to control the errors from the third stage. As the rationales for (\ref{28}) and (\ref{29}) are similar to the one employed for (\ref{27}), the proofs of the former are succinctly covered in subsections \ref{4.8} and \ref{4.9}. Finally, subsection \ref{4.10} presents the proof of Theorem~\ref{th4}, a by-product of the identities employed for Theorem~\ref{th2}. Last but not least, the result of Wang and Fan~\cite{fan} regarding the asymptotic behavior \(a_\nu-e_\nu\) for \(M\) is fixed can be recovered from the expansions of \((a_\nu-e_\nu)_k, k \ne \nu\) given by (\ref{17}) and of \(\Tilde{\Sigma}_3\) ((\ref{67}) and (\ref{28})) (only the first terms in \((a_\nu-e_\nu)_k, k \ne \nu\) will contribute, and \((a_\nu-e_\nu)_\nu=-\frac{||a_\nu-e_\nu||^2}{2}\) can be recovered). 

\subsection{Proof of Theorem~\ref{th2}: A Three-Component Decomposition}\label{4.1}

Denote by 
\[F=T^T\mathcal{M}(\hat{l}_\nu I- \mathcal{M})^{-1}T+\frac{1}{n}Z_AZ_A^T.\] 
Arguing in the same vein as for Proposition~\ref{prop2}, both (\ref{4}) and (\ref{5}) hold almost surely, the former yielding
\[\hat{l}_\nu=a_\nu^T\Lambda^{1/2}F\Lambda^{1/2} a_\nu.\]
Because \(F\) is symmetric, this identity can be rewritten as
\begin{equation}\label{16}
    \sqrt{n} \cdot \frac{\hat{l}_\nu}{l_\nu}=\frac{\sqrt{n}}{l_\nu} \cdot e_\nu^T\Lambda^{1/2}F\Lambda^{1/2} e_\nu+\frac{2\sqrt{n}}{l_\nu} \cdot  e_\nu^T\Lambda^{1/2}F\Lambda^{1/2}(a_\nu-e_\nu)+\frac{\sqrt{n}}{l_\nu} \cdot (a_\nu-e_\nu)^T\Lambda^{1/2}F\Lambda^{1/2} (a_\nu-e_\nu).
\end{equation}
\par
We introduce next the key decompositions that furnish us with the means of handling the last two terms in the right-hand side above: for \(k \ne \nu,\)
\begin{equation}\label{17}
    (a_\nu-e_\nu)_k=(\frac{||a_\nu-e_\nu||^2}{2}-1)\sum_{j \geq 0}{((-M_\nu)^j\mathcal{R}_\nu \mathcal{D}_\nu e_\nu)_k}:=(\frac{||a_\nu-e_\nu||^2}{2}-1) \Sigma_{0,k}
\end{equation}
for \(M_\nu=\mathcal{R}_\nu \mathcal{D}_\nu-(\hat{l}_\nu-l_\nu)\mathcal{R}_\nu.\) To establish these equalities, note that, as in Paul~\cite{paul},
\begin{equation}\label{13}
    a_\nu-e_\nu=-\mathcal{R}_\nu \mathcal{D}_\nu e_\nu+r_\nu
\end{equation}
\begin{equation}\label{14}
    r_\nu=(<a_\nu,e_\nu>-1)e_\nu-\mathcal{R}_\nu \mathcal{D}_\nu (a_\nu-e_\nu)+(\hat{l}_\nu-l_\nu)\mathcal{R}_\nu(a_\nu-e_\nu)
\end{equation}
since
\[-\mathcal{R}_\nu \mathcal{D}_\nu e_\nu+r_\nu=(<a_\nu,e_\nu>-1)e_\nu+(-\mathcal{R}_\nu \mathcal{D}_\nu+(\hat{l}_\nu-l_\nu)\mathcal{R}_\nu) a_\nu\]
\[a_\nu-<a_\nu,e_\nu>e_\nu=\mathcal{P}_\nu^{\perp}a_\nu=\mathcal{R}_\nu(\Lambda-l_\nu I)a_\nu,\]
\[(\Lambda-l_\nu I) a_\nu = (-\mathcal{D}_\nu+(\hat{l}_\nu-l_\nu)I)a_\nu,\]
the last equation being a rearrangement of (\ref{4}). In other words,
\[a_\nu-e_\nu=-\mathcal{R}_\nu \mathcal{D}_\nu e_\nu+r_\nu, \hspace{0.2cm} r_\nu=(<a_\nu,e_\nu>-1)e_\nu-M_\nu(a_\nu-e_\nu),\]
from which for \(k \ne \nu\) and any \(m \in \mathbb{N},\)
\begin{equation}\label{61}
    (a_\nu-e_\nu)_k=(-\mathcal{R}_\nu \mathcal{D}_\nu e_\nu)_k+((-M_\nu)(a_\nu-e_\nu))_k,
\end{equation}
\begin{equation}\label{19}
    (-M_\nu)^m(a_\nu-e_\nu)=-(-M_\nu)^m\mathcal{R}_\nu \mathcal{D}_\nu e_\nu+\frac{||a_\nu-e_\nu||^2}{2}(-M_\nu)^{m-1}\mathcal{R}_\nu \mathcal{D}_\nu e_\nu+(-M_\nu)^{m+1}(a_\nu-e_\nu),
\end{equation}
where we have used 
\[<a_\nu,e_\nu>-1=-\frac{||a_\nu-e_\nu||^2}{2}, \hspace{0.2cm} M_\nu e_\nu = \mathcal{R}_\nu \mathcal{D}_\nu e_\nu;\]
lastly, (\ref{61}) and (\ref{19}) in conjunction with
\[|(M_\nu^{m}(a_\nu-e_\nu))_k| \leq ||M_\nu||^m \cdot ||a_\nu-e_\nu||, \hspace{0.2cm} ||M_\nu|| \xrightarrow[]{a.s.} 0,\] render (\ref{17}).
\par
Next, separate the terms containing \((a_\nu-e_\nu)_\nu=-\frac{||a_\nu-e_\nu||^2}{2}\) and \((a_\nu-e_\nu)_k\) for \(k \ne \nu\) in the last two right-hand side terms in (\ref{16}) while employing the freshly established representations in (\ref{17}):
\[\frac{2\sqrt{n}}{l_\nu} \cdot e_\nu^T\Lambda^{1/2}F\Lambda^{1/2}(a_\nu-e_\nu)=-\sqrt{n} \cdot ||a_\nu-e_\nu||^2F_{\nu \nu}+2(\frac{||a_\nu-e_\nu||^2}{2}-1)\sqrt{\frac{n}{l_\nu}}\sum_{k \ne \nu}{F_{k \nu} \sqrt{l}_k\Sigma_{0,k}},\]
\[\frac{\sqrt{n}}{l_\nu}(a_\nu-e_\nu)^T\Lambda^{1/2}F\Lambda^{1/2} (a_\nu-e_\nu)=\sqrt{n} \cdot \frac{||a_\nu-e_\nu||^4}{4} F_{\nu \nu}-||a_\nu-e_\nu||^2(\frac{||a_\nu-e_\nu||^2}{2}-1)\sqrt{\frac{n}{l_\nu}}\sum_{k \ne \nu}{F_{k \nu} \sqrt{l}_k\Sigma_{0,k}}+\]
\[+\frac{\sqrt{n}}{l_\nu}(\frac{||a_\nu-e_\nu||^2}{2}-1)^2 \sum_{k_1 \ne \nu, k_2 \ne \nu}{F_{k_1k_2}\sqrt{l_{k_1}}\Sigma_{0,k_1}\sqrt{l_{k_2}}\Sigma_{0,k_2}}.\]
Hence (\ref{16}) is equivalent to
\[\sqrt{n} \cdot \frac{\hat{l}_\nu}{l_\nu}=\sqrt{n} \cdot F_{\nu \nu}+\sqrt{n} \cdot (\frac{||a_\nu-e_\nu||^4}{4}-||a_\nu-e_\nu||^2)F_{\nu \nu}
-2\sqrt{n}\cdot(\frac{||a_\nu-e_\nu||^2}{2}-1)^2\Tilde{\Sigma}_1+\sqrt{n}\cdot(\frac{||a_\nu-e_\nu||^2}{2}-1)^2\Tilde{\Sigma}_2\]
where
\[\Tilde{\Sigma}_1=\frac{1}{\sqrt{l_\nu}}\sum_{k \ne \nu}{F_{k \nu} \sqrt{l}_k\Sigma_{0,k}}, \hspace{0.2cm} \Tilde{\Sigma}_2=\frac{1}{l_\nu}\sum_{k_1 \ne \nu, k_2 \ne \nu}{F_{k_1k_2}\sqrt{l_{k_1}}\Sigma_{0,k_1}\sqrt{l_{k_2}}\Sigma_{0,k_2}}.\]
Denote by 
\[\Tilde{\Sigma}_3=||a_\nu-e_\nu||^2-\frac{||a_\nu-e_\nu||^4}{4},\] 
from which (\ref{16}) becomes
\begin{equation}\label{20}
    \sqrt{n} \cdot \frac{\hat{l}_\nu}{l_\nu}=\sqrt{n} \cdot F_{\nu \nu}-\sqrt{n} \cdot \Tilde{\Sigma}_3 \cdot F_{\nu \nu}
    -2\sqrt{n} \cdot (1-\Tilde{\Sigma}_3)\Tilde{\Sigma}_1+\sqrt{n} \cdot (1-\Tilde{\Sigma}_3)\Tilde{\Sigma}_2.
\end{equation}
\par
The three identities completing the proof of the theorem are
\begin{equation}\label{21}
    \sqrt{n} \cdot F_{\nu \nu}=\sqrt{n} \cdot (1+\frac{1}{n}tr(\mathcal{M}(\hat{l}_\nu I-\mathcal{M})^{-1})+\sqrt{n} \cdot (\frac{1}{n}z_\nu^Tz_\nu-1)+o_p(1),
\end{equation}
\begin{equation}\label{22}
    \sqrt{n} \cdot \Tilde{\Sigma}_3 \cdot F_{\nu \nu}=\sqrt{n} \cdot \Tilde{\Sigma}_3+o_p(1),
\end{equation}
\begin{equation}\label{23}
    -\Tilde{\Sigma}_3-2(1-\Tilde{\Sigma}_3)\Tilde{\Sigma}_1+(1-\Tilde{\Sigma}_3)\Tilde{\Sigma}_2=\overline{O}+\sum_{1 \leq j \leq s}{(\frac{\hat{l}_\nu}{l_\nu}-1)^j\overline{O}_{j}}+o_p(\frac{1}{\sqrt{n}}),
\end{equation}
for some 
\[s=s(n,M(n)), \hspace{0.2cm} \lim_{n \to \infty}{\frac{s(n,M(n))}{n/M(n)}}=0,\]
and deterministic constants depending only on \(n,\nu,l_1,l_2, \hspace{0.05cm} ... \hspace{0.05cm}, l_M,\)  
\[\overline{O}, \overline{O}_j \in \mathbb{R}, \hspace{0.2cm} |\overline{O}| \leq c(\epsilon_0) \cdot \frac{M}{n}, \hspace{0.2cm} |\overline{O}_j| \leq c(\epsilon_0)^{j+1} \cdot \frac{M}{n}.\] 
In the remainder of this subsection, we explain how these equations lead to the claimed conclusion. 
\par
Let \(x=x_{n,\nu,(l_i)_{1 \leq i \leq M}} \in \mathbb{R}, |x| \leq 2c(\epsilon_0) \cdot \frac{M}{n}\) solve
\begin{equation}\label{24}
    x=\overline{O}+\sum_{1 \leq j \leq s}{x^j\overline{O}_{j}}:
\end{equation}
such a solution always exists for \(n\) sufficiently large because \(f(y):=y-\sum_{1 \leq j \leq s}{y^j\overline{O}_{j}}\) has 
\[|f(y)-y| \leq \sum_{j \geq 1}{|y|^j \cdot c(\epsilon_0)^{j+1} \cdot \frac{M}{n}}=c(\epsilon_0) \cdot \frac{M}{n} \cdot \frac{|y|c(\epsilon_0)}{1-|y|c(\epsilon_0)} \leq 2c^2(\epsilon_0) \cdot \frac{M}{n} \cdot |y| \leq \frac{|y|}{2}\]
for \(|y| \leq 2c(\epsilon_0) \cdot \frac{M}{n} \leq \frac{1}{2c(\epsilon_0)},\) which implies with the intermediate value theorem that the image of \(f\) on \(|y| \leq 2c(\epsilon_0) \cdot \frac{M}{n}\) contains the interval \([-c(\epsilon_0) \cdot \frac{M}{n},c(\epsilon_0) \cdot \frac{M}{n}]\) in which \(\overline{O}\) lies. Next, employing (\ref{21}), (\ref{22}), and (\ref{23}), rewrite (\ref{20}) as
\begin{equation}\label{610}
    \sqrt{n} \cdot (\frac{\hat{l}_\nu}{l_\nu}-1)=\sqrt{n} \cdot c_{tr}+\mathcal{D}_n+ \sqrt{n} \cdot \overline{O}+\sqrt{n}\sum_{1 \leq j \leq s}{(\frac{\hat{l}_\nu}{l_\nu}-1)^j\overline{O}_{j}}+o_p(1),
\end{equation}
where 
\[c_{tr}=\frac{1}{n}tr(\mathcal{M}(\hat{l}_\nu I-\mathcal{M})^{-1}), \hspace{0.2cm} \mathcal{D}_n=\sqrt{n} \cdot (\frac{1}{n}z_\nu^Tz_\nu-1).\] 
\par
Lastly, change (\ref{610}) thrice to obtain the desired CLT:

\vspace{0.2cm}
\begin{equation}\label{62}
    1. \hspace{0.2cm} \sqrt{n} \cdot (\frac{\hat{l}_\nu}{l_\nu}-1-c_{tr})=\sqrt{n} \cdot \overline{O}+\mathcal{D}_n+\sqrt{n}\sum_{1 \leq j \leq s}{(\frac{\hat{l}_\nu}{l_\nu}-1-c_{tr})^j\overline{O}_{j}}+o_p(1).
\end{equation}
To justify this step, in light of (\ref{610}), it suffices to show
\begin{equation}\label{5007}
    \sqrt{n}\sum_{1 \leq j \leq s}{((\frac{\hat{l}_\nu}{l_\nu}-1)^j-(\frac{\hat{l}_\nu}{l_\nu}-1-c_{tr})^j)\overline{O}_{j}}=o_p(1):
\end{equation}
note that \(|x^j-(x-y)^j| \leq |y| \cdot j(|x|+|y|)^{j-1}\) because \(|x|,|x-y| \leq |x|+|y|,\) from which
\[|\sqrt{n}\sum_{1 \leq j \leq s}{((\frac{\hat{l}_\nu}{l_\nu}-1)^j-(\frac{\hat{l}_\nu}{l_\nu}-1-c_{tr})^j)\overline{O}_{j}}| \leq \sqrt{n} \cdot c^2(\epsilon_0) \cdot \frac{M}{n} \cdot |c_{tr}| \sum_{1 \leq j \leq s}{j \cdot (c(\epsilon_0) \cdot (|\frac{\hat{l}_\nu}{l_\nu}-1|+|c_{tr}|))^{j-1}}.\]
Since \(\frac{\hat{l}_\nu}{l_\nu}-1,c_{tr} \xrightarrow[]{a.s.} 0,\) it follows that almost surely the sum is bounded (\(\sum_{j \geq 1}{jx^{j-1}}=\frac{1}{(1-x)^2}\) for \(|x|<1\)) which together with \(\frac{M}{\sqrt{n}} \cdot c_{tr} \xrightarrow[]{a.s.} 0\)
(with probability one, \(0 \leq c_{tr}=tr(\mathcal{M}(\hat{l}_\nu I-\mathcal{M})^{-1}) \leq \frac{c(K,\gamma)}{\hat{l}_\nu} \leq \frac{2c(K,\gamma)}{l_\nu}\)) completes the justification of (\ref{5007}) and consequently of (\ref{62}).

\vspace{0.5cm}
\begin{equation}\label{63}
    2. \hspace{0.2cm} \sqrt{n} \cdot (\frac{\hat{l}_\nu}{l_\nu}-1-c_{tr}-x)=\mathcal{D}_n+\sqrt{n}\sum_{1 \leq j \leq s}{(\frac{\hat{l}_\nu}{l_\nu}-1-c_{tr}-x)^j\overline{O}_{j1}}+o_p(1)
\end{equation}
for some deterministic \(\overline{O}_{j1}\) with \(|\overline{O}_{j1}| \leq c(\epsilon_0)^{j+1} \cdot \frac{M}{n}:\) the binomial theorem, (\ref{62}), and (\ref{24}) give
\[\sqrt{n}\sum_{1 \leq j \leq s}{((\frac{\hat{l}_\nu}{l_\nu}-1-c_{tr})^j-(\frac{\hat{l}_\nu}{l_\nu}-1-c_{tr}-x)^j-x^j)\overline{O}_{j}}=\sqrt{n}\sum_{1 \leq j \leq s-1}{(\frac{\hat{l}_\nu}{l_\nu}-1-c_{tr}-x)^j\overline{O}_{j1}}\]
with
\[|\overline{O}_{j1}|=|\sum_{j+1 \leq l \leq s}{x^{l-j}\binom{l}{j}\overline{O}_l}| \leq c(\epsilon_0)^{j+1} \cdot \frac{M}{n}\sum_{1 \leq l \leq s}{|x|^{l} \cdot (c(\epsilon_0)s)^l} \leq c(\epsilon_0)^{j+1} \cdot \frac{M}{n}\]
because 
\[\binom{l}{j} \leq l^{l-j} \leq s^{l-j}, \hspace{0.4cm} |x| \cdot c(\epsilon_0)s \leq 2c^2(\epsilon_0) \cdot s \cdot \frac{M}{n} \to 0,\]
and in particular, \(|x| \cdot c(\epsilon_0)s \leq \frac{1}{2}\) for \(n\) large enough.

\vspace{0.5cm}
\begin{equation}\label{64}
    3. \hspace{0.2cm} \sqrt{n} \cdot (\frac{\hat{l}_\nu}{l_\nu}-1-c_{tr}-x)=\mathcal{D}_n+\sqrt{n} \cdot (\frac{\hat{l}_\nu}{l_\nu}-1-c_{tr}-x)o_1+o_p(1)
\end{equation}
where \(o_1=o_p(1):\) using (\ref{63}) and \(|\frac{\hat{l}_\nu}{l_\nu}-1-c_{tr}-x| \leq |\frac{\hat{l}_\nu}{l_\nu}-1|+|c_{tr}|+|x| \xrightarrow[]{a.s.} 0,\) almost surely
\[|\sum_{1 \leq j \leq s}{(\frac{\hat{l}_\nu}{l_\nu}-1-c_{tr}-x)^{j-1}\overline{O}_{j1}}| \leq c^2(\epsilon_0) \cdot \frac{M}{n}\sum_{1 \leq j \leq s}{|c(\epsilon_0)(\frac{\hat{l}_\nu}{l_\nu}-1-c_{tr}-x)|^{j-1}} \leq  2c^2(\epsilon_0) \cdot \frac{M}{n}.\]

\vspace{0.3cm}
\par
Finally, (\ref{64}), Lindeberg's CLT, and two applications of Slutsky's lemma yield Theorem~\ref{th2}: first, 
\[\sqrt{\frac{n}{\mathbb{E}[z_{11}^4]-1}} \cdot (\frac{\hat{l}_\nu}{l_\nu}-1-c_{tr}-x)(1-o_1) \Rightarrow N(0,1),\]
as for \(T_{ni}=\frac{1}{\sqrt{n}} \cdot \frac{z_{\nu i}^2-1}{\sqrt{\mathbb{E}[z_{11}^4]-1}}\) for \(1 \leq i \leq n,\) \(\sum_{1 \leq i \leq n}{\mathbb{E}[T^2_{ni}]}=1,\) and for any \(\epsilon>0,\) 
\[\sum_{1 \leq i \leq n}{\mathbb{E}[T^2_{ni}\chi_{|T_{ni}| \geq \epsilon}]}=\frac{1}{\mathbb{E}[z_{11}^4]-1} \mathbb{E}[(z_{11}^2-1)^2\chi_{|z_{11}^2-1| \geq \epsilon \sqrt{n(\mathbb{E}[z_{11}^4]-1)}}] \leq \frac{1}{\delta_0} \cdot \frac{\mathbb{E}[(z_{11}^2-1)^4]}{n\delta_0\epsilon^2} \leq \frac{8(1+c_8(K))}{n(\delta_0 \epsilon)^2} \to 0,\]
where \(c_8(K) = \sup_{\mathbb{E}[x]=0,||x||_{\psi_2} \leq K}{\mathbb{E}[x^8]}<\infty\) from Vershynin~\cite{vershynin}, proposition \(2.5.2,\) and second, 
\[\sqrt{\frac{n}{\mathbb{E}[z_{11}^4]-1}} \cdot (\frac{\hat{l}_\nu}{l_\nu}-1-c_{tr}-x)=\frac{1}{1-o(1)} \cdot \sqrt{\frac{n}{\mathbb{E}[z_{11}^4]-1}} \cdot (\frac{\hat{l}_\nu}{l_\nu}-1-c_{tr}-x)(1-o_1) \Rightarrow N(0,1).\]

\subsection{First and Second Components}\label{4.2}

In this subsection, we justify (\ref{21}) and (\ref{22}).

\par
\vspace{0.2cm}
\textit{Proof of (\ref{21}):}
Denote by \(\hat{\mathcal{M}}_\nu=H\mathcal{M}(\hat{l}_\nu I- \mathcal{M})^{-1}H^T, \mathcal{M}_\nu=H\mathcal{M}(l_\nu I- \mathcal{M})^{-1}H^T.\) We show next that 
\[\frac{1}{\sqrt{n}}(z_\nu ^T  \hat{\mathcal{M}}_\nu z_\nu-tr(\hat{\mathcal{M}}_\nu)) \xrightarrow[]{p} 0,\]
which yields (\ref{21}) since \(F_{\nu \nu}=\frac{1}{n} z_\nu^T \hat{\mathcal{M}}_\nu z_\nu+\frac{1}{n}z_\nu^T z_\nu,\) and \(tr(\hat{\mathcal{M}}_\nu)=tr(\mathcal{M}(\hat{l}_\nu I- \mathcal{M})^{-1}).\)
\par
If \(|\frac{\hat{l}_\nu}{l_\nu}-1| \leq \frac{1}{8}, ||\mathcal{M}|| \leq c_{K,\gamma} \leq \frac{l_\nu}{2},\) then
\[z_\nu ^T  \hat{\mathcal{M}}_\nu z_\nu-tr(\hat{\mathcal{M}}_\nu)=\sum_{k \geq 0}{(l_\nu-\hat{l}_\nu)^k(z_\nu ^T H\mathcal{M}(l_\nu I-\mathcal{M})^{-(k+1)}H^T z_\nu-tr(\mathcal{M}(l_\nu I-\mathcal{M})^{-(k+1)})}\]
because for \(0 \leq x \leq \frac{l_\nu}{2},\) \(|l_\nu-\hat{l}_\nu|<\frac{l_\nu}{2} \leq l_\nu-x,\)
\[\frac{1}{\hat{l}_\nu-x}=\frac{1}{l_\nu-x} \cdot \frac{1}{1-\frac{l_\nu-\hat{l}_\nu}{l_\nu-x}}=\sum_{k \geq 0}{\frac{(l_\nu-\hat{l}_\nu)^k}{(l_\nu-x)^{k+1}}}.\]
\par
Hence for \(t>0,\) and \(n\) large enough so that \(l_\nu \geq 2c_{K,\gamma},\)
\[\mathbb{P}(\frac{1}{\sqrt{n}}|z_\nu ^T  \hat{\mathcal{M}}_\nu z_\nu-tr(\hat{\mathcal{M}}_\nu)| \geq 2t) \leq \mathbb{P}(|\frac{\hat{l}_\nu}{l_\nu}-1|>\frac{1}{8})+\mathbb{P}(||\mathcal{M}||>c_{K,\gamma})+\sum_{k \geq 0}{\mathbb{P}(\frac{|\sigma_k|}{\sqrt{n}} \geq \frac{t}{2^k}| \hspace{0.05cm} ||\mathcal{M}|| \leq c_{K,\gamma})},\]
for 
\[\sigma_k=(\frac{l_\nu}{8})^k \cdot (z_\nu ^T H\mathcal{M}(l_\nu I-\mathcal{M})^{-(k+1)}H^T z_\nu-tr(\mathcal{M}(l_\nu I-\mathcal{M})^{-(k+1)}).\] 
Hanson-Wright inequality (\ref{9}) gives
\[\mathbb{P}(\frac{1}{\sqrt{n}}|z_\nu ^T H\mathcal{M}(l_\nu I-\mathcal{M})^{-(k+1)}H^T z_\nu-tr(\mathcal{M}(l_\nu I-\mathcal{M})^{-(k+1)}| \geq \frac{4^kt}{l_\nu^k}|Z_B) \leq\]
\[\leq 2\exp(-c\min{(\frac{4^{2k}t^2}{K^4l_\nu^{2k} \cdot ||\mathcal{M}(l_\nu I-\mathcal{M})^{-(k+1)}||^2},\frac{4^{k}t\sqrt{n}}{K^2l_\nu^{k} \cdot ||\mathcal{M}(l_\nu I-\mathcal{M})^{-(k+1)}||})}),\]
and if \(||\mathcal{M}|| \leq c_{K,\gamma} \leq \frac{l_\nu}{2},\) then \(||\mathcal{M}(l_\nu I-\mathcal{M})^{-(k+1)}||\leq \frac{c_{K,\gamma}}{(\frac{l_\nu}{2})^{k+1}},\) from which
\[\mathbb{P}(\frac{|\sigma_k|}{\sqrt{n}} \geq \frac{t}{2^k}| \hspace{0.05cm} ||\mathcal{M}|| \leq c_{K,\gamma}) \leq 2\exp(-c\min{(\frac{4^{k}l_\nu^2t^2}{4K^4c_{K,\gamma}^2},\frac{4^{k}t\sqrt{n}l_\nu}{2K^2c_{K,\gamma}})}),\]
implying that
\[\mathbb{P}(\frac{1}{\sqrt{n}}|z_\nu ^T  \hat{\mathcal{M}}_\nu z_\nu-tr(\hat{\mathcal{M}}_\nu)| \geq 2t) \leq \mathbb{P}(|\frac{\hat{l}_\nu}{l_\nu}-1|>\frac{1}{8})+\mathbb{P}(||\mathcal{M}||>c_{K,\gamma})+\]
\[+2\sum_{k \geq 0}{\exp(-\frac{c4^{k}l_\nu^2t^2}{K^4c_{K,\gamma}^2})}+2\sum_{k \geq 0}{\exp(-\frac{c4^{k}t\sqrt{n}l_\nu}{K^2c_{K,\gamma}})} \to 0,\]
because \(||\mathcal{M}|| \leq c_{K,\gamma}, \frac{\hat{l}_\nu}{l_\nu}-1 \to 0\) almost surely from (\ref{7}) and Proposition~\ref{prop1}, \(l_\nu \to \infty,\) and for \(x \geq 1,\)
\[\sum_{k \geq 0}{\exp(-4^kx)} \leq \sum_{k \geq 0}{\exp(-2^kx)} \leq \sum_{k \geq 1}{\exp(-kx)}=\frac{\exp(-x)}{1-\exp(-x)} \leq 2\exp(-x).\]

\par
\vspace{0.3cm}
\textit{Proof of (\ref{22}):}
It suffices to show that
\begin{equation}\label{66}
    \sqrt{n} \cdot |F_{\nu \nu}-1| \cdot \beta_\nu^2 \xrightarrow[]{p} 0
\end{equation}
because Lemma~\ref{lemma2} yields \(||a_\nu-e_\nu|| \leq 2\beta_\nu\) almost surely, from which
\[|\sqrt{n} \cdot (F_{\nu \nu}-1) \cdot (||a_\nu-e_\nu||^2-\frac{||a_\nu-e_\nu||^4}{4})| \leq \sqrt{n} \cdot |F_{\nu \nu}-1| \cdot 4\beta_\nu^2.\]
\par
We prove next that for any sequences \(\delta_n, \delta^*_n\) with \(\delta_n \to \infty, \hspace{0.2cm} \delta^*_n \to 0,  \hspace{0.2cm} \frac{\sqrt{n}\delta^*_n}{l_\nu} \to 0,\)
\begin{equation}\label{25}
    \frac{n}{M\delta_n} \cdot \beta_\nu^2 \xrightarrow[]{p} 0, \hspace{0.4cm} \delta^*_n \cdot \sqrt{n} \cdot (F_{\nu \nu}-1) \xrightarrow[]{p} 0,
\end{equation}
from which (\ref{66}) ensues because then
\[\sqrt{n} \cdot |F_{\nu \nu}-1| \cdot \beta_\nu^2=(\frac{M\delta_n}{n} \cdot \sqrt{n} \cdot |F_{\nu \nu}-1|) \cdot \frac{n}{M\delta_n}\beta_\nu^2 \xrightarrow[]{p} 0\]
for \(\delta_n \to \infty\) such that \(\frac{M\delta_n}{n} \to 0, \frac{M\delta_n}{\sqrt{n}} \cdot \frac{1}{l_\nu} \to 0:\) 
for instance, \(\delta_n=\min{(\sqrt{\frac{n}{M}},\sqrt{\frac{l_\nu}{M/\sqrt{n}}})} \to \infty.\) 
\par
\textit{First convergence in (\ref{25}):}
\[\beta_\nu^2=\sum_{k \ne \nu}{\frac{l_kl_\nu}{(l_k-l_\nu)^2}(\frac{1}{n}z_k^Tz_\nu+t_k^T\mathcal{M}(\hat{l}_\nu I-\mathcal{M})^{-1}t_\nu)^2} \leq c(\epsilon_0)\sum_{k \ne \nu}{(\frac{1}{n} z_k^T z_\nu)^2}+c(\epsilon_0)\sum_{k \ne \nu}{(t_k^T\mathcal{M}(\hat{l}_\nu I- \mathcal{M})^{-1}t_\nu)^2};\]
as \(\mathbb{E}[\sum_{k \ne \nu}{(\frac{1}{n} z_k^T z_\nu)^2}]=\frac{M-1}{n},\) it follows that the first term under multiplication with \(\frac{n}{M\delta_n}\) tends to zero in probability, while for the second term, part \((b)\) of Lemma~\ref{lemma1} yields the result: \(\newline \frac{n}{M\delta_n}\sum_{k \ne \nu}{(t_k^T\mathcal{M}(\hat{l}_\nu I- \mathcal{M})^{-1}t_\nu)^2} \geq t\) implies for \(n\) large enough \(\newline \frac{n}{M}\sum_{k \ne \nu}{(t_k^T\mathcal{M}(\hat{l}_\nu I- \mathcal{M})^{-1}t_\nu)^2} \geq t\delta_n \geq t.\)

\vspace{0.2cm}
\textit{Second convergence in (\ref{25}):} (\ref{21}) gives
\[\delta^*_n \cdot\sqrt{n} \cdot (F_{\nu \nu}-1)=\delta^*_n \cdot\sqrt{n} \cdot (\frac{1}{n}z_\nu^Tz_\nu-1)+\delta^*_n \cdot\sqrt{n} \cdot c_{tr}+\delta^*_n \cdot o_p(1)=o_p(1)\]
where we have used Chebyshev's inequality, \(\delta^*_n \to 0\) in conjunction with
\[\mathbb{E}[n \cdot (\frac{1}{n}z_\nu^Tz_\nu-1)^2]=\mathbb{E}[z_{11}^4]-1 \leq c(K),\]
and almost surely
\[\delta^*_n \cdot \sqrt{n} c_{tr} \leq c(K,\gamma) \cdot \frac{\sqrt{n}\delta^*_n}{l_\nu}.\]

\subsection{Third Component: Reduction to Three Sums}\label{4.3}

In this subsection, we reduce the last identity, (\ref{23}), needed to finalize the proof of our first theorem to obtaining polynomial representations for three sums.
\par
Rewrite \(\Tilde{\Sigma}_3\) in terms of \(\Sigma_{0,k}, k \ne \nu:\) since
\[1=||a_\nu||^2=(\frac{||a_\nu-e_\nu||^2}{2}-1)^2+\sum_{k \ne \nu}{((a_\nu-e_\nu)_k)^2},\]
(\ref{17}) yields
\[(\frac{||a_\nu-e_\nu||^2}{2}-1)^2+(\frac{||a_\nu-e_\nu||^2}{2}-1)^2\sum_{k \ne \nu}{\Sigma_{0,k}^2}=1,\]
from which
\begin{equation}\label{67}
    \Tilde{\Sigma}_3=||a_\nu-e_\nu||^2-\frac{||a_\nu-e_\nu||^4}{4}=1-(\frac{||a_\nu-e_\nu||^2}{2}-1)^2=\frac{\sum_{k \ne \nu}{\Sigma_{0,k}^2}}{1+\sum_{k \ne \nu}{\Sigma_{0,k}^2}}:=\frac{\Tilde{\Sigma}_0}{1+\Tilde{\Sigma}_0}.
\end{equation}
\par
To obtain (\ref{23}), it suffices to find expansions of the type 
\[O+\sum_{1 \leq j \leq s}{(\frac{\hat{l}_\nu}{l_\nu}-1)^jO_{j}}+o_p(\frac{1}{\sqrt{n}})\] 
for \(\Tilde{\Sigma}_0,\Tilde{\Sigma}_1,\Tilde{\Sigma}_2\) with \(|O| \leq c(\epsilon_0) \cdot \frac{M}{n},|O_{j}| \leq c(\epsilon_0)^{j+1} \cdot \frac{M}{n}\) for some \(s=s(n,M) \in \mathbb{N}\) with \(\sqrt{n} \cdot (4c(\epsilon_0)\frac{M}{n})^s \to 0:\) then 
\[\Tilde{\Sigma}_2-2\Tilde{\Sigma}_1-\Tilde{\Sigma}_0=O+\sum_{1 \leq j \leq s}{(\frac{\hat{l}_\nu}{l_\nu}-1)^jO_{j}}+o_p(\frac{1}{\sqrt{n}}), \hspace{0.3cm} -\Tilde{\Sigma}_0=O'+\sum_{1 \leq j \leq s}{(\frac{\hat{l}_\nu}{l_\nu}-1)^jO'_{j}}+o_p(\frac{1}{\sqrt{n}}),\]
render
\[-\Tilde{\Sigma}_3-2\Tilde{\Sigma}_1(1-\Tilde{\Sigma}_3)+(1-\Tilde{\Sigma}_3)\Tilde{\Sigma}_2=\frac{\Tilde{\Sigma}_2-2\Tilde{\Sigma}_1-\Tilde{\Sigma}_0}{1+\Tilde{\Sigma}_0}=\]
\begin{equation}\label{26}
    =\sum_{0 \leq t \leq s}{(O+\sum_{1 \leq j \leq s}{(\frac{\hat{l}_\nu}{l_\nu}-1)^jO_{j}}) \cdot (O'+\sum_{1 \leq j \leq s}{(\frac{\hat{l}_\nu}{l_\nu}-1)^jO'_{j}})^t}+o_p(\frac{1}{\sqrt{n}})=\overline{O}+\sum_{1 \leq j \leq s^2+s}{(\frac{\hat{l}_\nu}{l_\nu}-1)^j\overline{O}_{j}}+o_p(\frac{1}{\sqrt{n}}),
\end{equation}
since
\[\frac{O+\sum_{1 \leq j \leq s}{(\frac{\hat{l}_\nu}{l_\nu}-1)^jO_{j}}+o_p(\frac{1}{\sqrt{n}})}{1-O'-\sum_{1 \leq j \leq s}{(\frac{\hat{l}_\nu}{l_\nu}-1)^jO'_{j}}-o_p(\frac{1}{\sqrt{n}})}=\frac{O+\sum_{1 \leq j \leq s}{(\frac{\hat{l}_\nu}{l_\nu}-1)^jO_{j}}}{1-O'-\sum_{1 \leq j \leq s}{(\frac{\hat{l}_\nu}{l_\nu}-1)^jO'_{j}}}+o_p(\frac{1}{\sqrt{n}})\]
from
\[\frac{o_p(\frac{1}{\sqrt{n}})(1-O'-\sum_{1 \leq j \leq s}{(\frac{\hat{l}_\nu}{l_\nu}-1)^jO'_{j}})+o_p(\frac{1}{\sqrt{n}})(O+\sum_{1 \leq j \leq s}{(\frac{\hat{l}_\nu}{l_\nu}-1)^jO_{j}})}{(1-O'-\sum_{1 \leq j \leq s}{(\frac{\hat{l}_\nu}{l_\nu}-1)^jO'_{j}}) \cdot (1-O'-\sum_{1 \leq j \leq s}{(\frac{\hat{l}_\nu}{l_\nu}-1)^jO'_{j}})+o_p(\frac{1}{\sqrt{n}}))}=o_p(\frac{1}{\sqrt{n}})\]
as \(O+\sum_{1 \leq j \leq s}{(\frac{\hat{l}_\nu}{l_\nu}-1)^jO_{j}} \xrightarrow[]{p} 0,O'+\sum_{1 \leq j \leq s}{(\frac{\hat{l}_\nu}{l_\nu}-1)^jO'_{j}} \xrightarrow[]{p} 0,\) and
\[\frac{O+\sum_{1 \leq j \leq s}{(\frac{\hat{l}_\nu}{l_\nu}-1)^jO_{j}}}{1-O'-\sum_{1 \leq j \leq s}{(\frac{\hat{l}_\nu}{l_\nu}-1)^jO'_{j}}}=o_p(\frac{1}{\sqrt{n}})+ \sum_{0 \leq t \leq s}{(O+\sum_{1 \leq j \leq s}{(\frac{\hat{l}_\nu}{l_\nu}-1)^jO_{j}}) \cdot (O'+\sum_{1 \leq j \leq s}{(\frac{\hat{l}_\nu}{l_\nu}-1)^jO'_{j}})^t},\]
because almost surely for \(n\) large enough,
\[|O'+\sum_{1 \leq j \leq s}{(\frac{\hat{l}_\nu}{l_\nu}-1)^jO'_{j}}| \leq c(\epsilon_0) \cdot \frac{M}{n} \sum_{0 \leq j \leq s}{(c(\epsilon_0) \cdot |\frac{\hat{l}_\nu}{l_\nu}-1|)^{j}} \leq 2c(\epsilon_0) \cdot \frac{M}{n},\]
entailing \(|O'+\sum_{1 \leq j \leq s}{(\frac{\hat{l}_\nu}{l_\nu}-1)^jO'_{j}}|^{s+1}=o_p(\frac{1}{\sqrt{n}}).\)
\par
Regarding the sizes of \(\overline{O}, \overline{O}_j, 1 \leq j \leq 2s^2+s,\) for \(0 \leq u \leq s^2+s,\) the absolute value of the coefficient of \((\frac{\hat{l}_\nu}{l_\nu}-1)^u\) in (\ref{26}) is at most
\begin{equation}\label{68}
    \sum_{0 \leq t \leq s, 0 \leq v \leq u}{|O_v| \cdot (\frac{M}{n})^t \cdot c(\epsilon_0)^{t+u} \cdot \binom{t+u-v-1}{t-1}}
\end{equation}
where \(O_0:=O,\) because the coefficient of \((\frac{\hat{l}_\nu}{l_\nu}-1)^{u-v}\) in \((O'+\sum_{1 \leq j \leq s}{(\frac{\hat{l}_\nu}{l_\nu}-1)^jO'_{j}})^t\) (seen as a polynomial in \(\frac{\hat{l}_\nu}{l_\nu}-1\)) is a sum, whose terms are each bounded in absolute value by 
\[|O'|^{t-(u-v)} \max{\prod_{\sum{j}=u-v}{|O'_j|}} \leq (c(\epsilon_0) \cdot \frac{M}{n})^{t-(u-v)} \cdot (c(\epsilon_0))^{2(u-v)} \cdot  (\frac{M}{n})^{u-v} \leq  (\frac{M}{n})^t \cdot c(\epsilon_0)^{t+u},\]
and their number is at most \(\Tilde{a}_{t,u-v}=\binom{t+u-v-1}{t-1}\) for 
\begin{equation}\label{4001}
    \Tilde{a}_{p,n}:=|\{(x_1, \hspace{0.05cm} ... \hspace{0.05cm}, x_p):x_1+...+x_p=n, x_i \in \mathbb{Z}_{\geq 0}\}|=\binom{p+n-1}{p-1}
\end{equation}
(by induction: \(\Tilde{a}_{1,n}=1,\Tilde{a}_{p,n}=\Tilde{a}_{p-1,n}+\Tilde{a}_{p,n-1}\) for \(p>1\)); the sum in (\ref{68}) is upper bounded by
\[c(\epsilon_0)^{u}\sum_{0 \leq t \leq s, 0 \leq v \leq u}{|O_v| \cdot (\frac{M}{n})^t \cdot c(\epsilon_0)^{t} \cdot 2^{t+u}}= 
(2c(\epsilon_0))^{u} (\sum_{0 \leq v \leq u}{|O_v|})\sum_{0 \leq t \leq s}{(2c(\epsilon_0) \cdot \frac{M}{n})^t} \leq \]
\[\leq 2 \cdot (2c(\epsilon_0))^{u} \cdot \sum_{0 \leq v \leq u}{|O_v|} \leq 2 \cdot (2c(\epsilon_0))^{u} \cdot \frac{M}{n}\cdot c(\epsilon_0) (c(\epsilon_0)+1)^{u+1}=\frac{M}{n} \cdot (2c(\epsilon_0)(c(\epsilon_0)+1))^{u+1},\]
using the elementary inequality \(\sum_{0 \leq i \leq u}{x^u} \leq (x+1)^{u+1}\) for \(x \geq 0, u \in \mathbb{Z}_{\geq 0}.\)  
(We abuse notation for the sake of simplicity and denote by \(s\) the index in (\ref{23}) when in reality it is \(s^2+s\).).
\vspace{0.2cm}
\par
Take \(s=s(n,M)=\left \lfloor{\frac{8\log{n}}{\log{\frac{n}{M}}}}\right \rfloor \in \mathbb{N},\) for which
\[s\log{\frac{n}{M}}-4\log{n} \to \infty, \hspace{0.2cm} \frac{s}{\log{\frac{n}{M}}} \to 0:\]
note that \(\sqrt{n} \cdot (4c(\epsilon_0)\frac{M}{n})^{2s} \to 0,\) and, in light of the previous paragraph, it suffices to show that the following expansions hold:
\begin{equation}\label{27}
    \Tilde{\Sigma}_1=\frac{1}{\sqrt{l_\nu}}\sum_{k \ne \nu}{F_{k \nu}\sqrt{l_k}\Sigma_{0,k}}=-a_0-\sum_{1 \leq m \leq s}{a_m(\frac{\hat{l}_\nu}{l_\nu}-1)^m}+o_p(\frac{1}{\sqrt{n}}),
\end{equation}
\begin{equation}\label{28}
    \Tilde{\Sigma}_0=\sum_{k \ne \nu}{\Sigma^2_{0,k}}=-b_0-\sum_{1 \leq m \leq 2s}{b_m(\frac{\hat{l}_\nu}{l_\nu}-1)^m}+o_p(\frac{1}{\sqrt{n}}),
\end{equation}
\begin{equation}\label{29}
    \Tilde{\Sigma}_2= \sum_{k \ne \nu}{\frac{l_k}{l_\nu}\Sigma^2_{0,k}}+\frac{1}{l_\nu}\sum_{k_1\ne \nu, k_2 \ne \nu}{F_{k_1 k_2} \sqrt{l_{k_1}}\Sigma_{0,k_1} \cdot \sqrt{l_{k_2}}\Sigma_{0,k_2}} = c_0+\sum_{1 \leq m \leq 2s}{c_m(\frac{\hat{l}_\nu}{l_\nu}-1)^m}+o_p(\frac{1}{\sqrt{n}}),
\end{equation}
together with \(|a_i|,|b_i|,|c_i| \leq c(\epsilon_0)^{i+1} \cdot \frac{M}{n}\) for \(0 \leq i \leq 2s.\)
\par
Recall that 
\[F=\Tilde{D}+\Tilde{T}+I, \hspace{0.2cm} -M_\nu=\mathcal{R}_\nu \Lambda^{1/2} (\Tilde{T}+\Tilde{D})\Lambda^{1/2}-(\frac{\hat{l}_\nu}{l_\nu}-1)l_\nu \mathcal{R}_\nu,\]
for 
\[\Tilde{T}=T^T\mathcal{M}(\hat{l}_\nu I- \mathcal{M})^{-1}T, \hspace{0.2cm} \Tilde{D}=\frac{1}{n}Z_AZ_A^T-I.\] 
Since each term in the sums underlying \(\Tilde{\Sigma}_i, 0 \leq i \leq 2\) arises from some 
\[\Sigma_{0,k}=\sum_{j \geq 0}{((-M_\nu)^j\mathcal{R}_\nu \mathcal{D}_\nu e_\nu)_k},\] 
after expanding \((-M_\nu)^j,\) any either contains some factor of \(\Tilde{T}\) (\textit{first type} terms) or does not (\textit{second type} terms). The forthcoming subsections consist of proving the former are negligible (i.e., \(o_p(\frac{1}{\sqrt{n}})\)) while the latter contribute solely through some expectations to the sum (with an error of order \(o_p(\frac{1}{\sqrt{n}})\)): these are the second and third stages mentioned at the beginning of this section while the first phase is the truncation of the series \(\Sigma_{0,k}, k \ne \nu.\)

\subsection{First Sum, Stage I: Truncation}\label{4.4}

In this subsection, we truncate the series in (\ref{27}) by establishing the following:
\begin{equation}\label{30}
    \frac{1}{\sqrt{l_\nu}}\sum_{k \ne \nu}{F_{k \nu}\sqrt{l_k}\Sigma_{0,k}}=\frac{1}{\sqrt{l_\nu}}\sum_{k \ne \nu}{\Tilde{D}_{k \nu}\sqrt{l_k}\Sigma_{0,k}}+o_p(\frac{1}{\sqrt{n}}),
\end{equation}
\begin{equation}\label{31}
    \frac{1}{\sqrt{l_\nu}}\sum_{k \ne \nu}{\Tilde{D}_{k \nu}\sqrt{l_k}\Sigma_{0,k}}=\frac{1}{\sqrt{l_\nu}}\sum_{0 \leq j \leq s}{\sum_{k \ne \nu}{\Tilde{D}_{k \nu}\sqrt{l_k}}((-M_\nu)^j\mathcal{R}_\nu \mathcal{D}_\nu e_\nu)_k}+o_p(\frac{1}{\sqrt{n}}).
\end{equation}
\par
\textit{Proof of (\ref{30}):} Recall that \(F=\Tilde{D}+\Tilde{T}+I,\) from which \(F_{k\nu}=\Tilde{D}_{k \nu}+\Tilde{T}_{k\nu}\) for \(k \ne \nu.\) Cauchy-Schwarz inequality and (\ref{17}) yield
\begin{equation}\label{931}
    |\sqrt{\frac{n}{l_\nu}}\sum_{k \ne \nu}{\Tilde{T}_{k \nu}\sqrt{l_k}\Sigma_{0,k}}| \leq \sqrt{\frac{n}{l_\nu}} \cdot \frac{||\Lambda^{1/2}(a_\nu-e_\nu)||}{|1-\frac{||a_\nu-e_\nu||^2}{2}|} \cdot \sqrt{\sum_{k \ne \nu}{\Tilde{T}^2_{k \nu}}}.
\end{equation}
Using (\ref{13}), and (\ref{14}),
\[\Lambda^{1/2}(a_\nu-e_\nu)=-\Lambda^{1/2}\mathcal{R}_\nu \mathcal{D}_\nu e_\nu -\frac{||a_\nu-e_\nu||^2}{2}\sqrt{l_\nu}e_\nu -\Lambda^{1/2}M_\nu(a_\nu-e_\nu),\]
from which
\[\frac{||\Lambda^{1/2}(a_\nu-e_\nu)||}{\sqrt{l_\nu}} \leq \frac{||\Lambda^{1/2}\mathcal{R}_\nu \mathcal{D}_\nu e_\nu||}{\sqrt{l_\nu}}+\frac{||a_\nu-e_\nu||^2}{2}+\frac{||\Lambda^{1/2}M_\nu(a_\nu-e_\nu)||}{\sqrt{l_\nu}};\]
since \(||M_\nu|| \leq \frac{1}{2}\) almost surely, \(||\Lambda^{1/2}M_\nu(a_\nu-e_\nu)|| \leq \frac{1}{2}||\Lambda^{1/2}(a_\nu-e_\nu)||,\) yielding
\begin{equation}\label{1010}
    \frac{||\Lambda^{1/2}(a_\nu-e_\nu)||}{\sqrt{l_\nu}} \leq 2 \cdot \frac{||\Lambda^{1/2}\mathcal{R}_\nu \mathcal{D}_\nu e_\nu||}{\sqrt{l_\nu}}+||a_\nu-e_\nu||^2 \leq 2 \cdot \frac{||\Lambda^{1/2}\mathcal{R}_\nu \mathcal{D}_\nu e_\nu||}{\sqrt{l_\nu}}+4\beta_\nu^2,
\end{equation}
employing Lemma~\ref{lemma2}. This gives that for any \(\delta_n \to \infty,\)
\begin{equation}\label{33}
    \sqrt{\frac{n}{M\delta_n}} \cdot \frac{||\Lambda^{1/2}(a_\nu-e_\nu)||}{\sqrt{l_\nu}} \xrightarrow[]{p} 0,
\end{equation}
arguing as for (\ref{25}). Because almost surely \(|1-\frac{||a_\nu-e_\nu||^2}{2}| \geq \frac{1}{2},\) it follows from (\ref{931}) that 
\[\sqrt{\frac{n}{l_\nu}}\sum_{k \ne \nu}{\Tilde{T}_{k \nu}\sqrt{l_k}\Sigma_{0,k}} \xrightarrow[]{p} 0\] 
because for \(\delta_n=\frac{l_\nu}{M/\sqrt{n}} \to \infty,\)
\begin{equation}\label{1001}
    M\delta_n\sum_{k \ne \nu}{\Tilde{T}^2_{k\nu}} \xrightarrow[]{p} 0
\end{equation}
from (\ref{909}) as \(M\delta_n\sum_{k \ne \nu}{\Tilde{T}^2_{k\nu}} \geq t\) is equivalent to \(\frac{n}{M}\sum_{k \ne \nu}{\Tilde{T}^2_{k\nu}} \geq \frac{tn}{M^2\delta_n},\) and \(\frac{tn}{M^2\delta_n} \cdot l_\nu^2 \to \infty\) for \(t>0.\) Hence,
\[\frac{1}{\sqrt{l_\nu}}\sum_{k \ne \nu}{F_{k \nu}\sqrt{l_k}\Sigma_{0,k}}=\frac{1}{\sqrt{l_\nu}}\sum_{k \ne \nu}{\Tilde{D}_{k \nu}\sqrt{l_k}\Sigma_{0,k}}+o_p(\frac{1}{\sqrt{n}}).\]

\vspace{0.2cm}
\par
\textit{Proof of (\ref{31}):} Note that for any \(j \geq 0,\)
\[|\frac{1}{\sqrt{l_\nu}}\sum_{k \ne \nu}{\Tilde{D}_{k \nu}\sqrt{l_k}((-M_\nu)^j\mathcal{R}_\nu \mathcal{D}_\nu e_\nu)_k}| \leq \frac{1}{\sqrt{l_\nu}} \sqrt{\sum_{k \ne \nu}{\Tilde{D}^2_{k \nu}}} \cdot ||\Lambda^{1/2}(-M_\nu)^j\mathcal{R}_\nu \mathcal{D}_\nu e_\nu|| \leq\]
\[\leq \sqrt{\sum_{k \ne \nu}{\Tilde{D}^2_{k \nu}}} \cdot ||M_\nu||^{j} \cdot \frac{||\Lambda^{1/2}\mathcal{R}_\nu \mathcal{D}_\nu e_\nu||}{\sqrt{l_\nu}},\]
from which almost surely
\begin{equation}\label{107}
    |\frac{\sqrt{n}}{\sqrt{l_\nu}}\sum_{k \ne \nu}{\Tilde{D}_{k \nu}\sqrt{l_k}\sum_{j>s}{((-M_\nu)^j\mathcal{R}_\nu \mathcal{D}_\nu e_\nu)_k}}| \leq 2\sqrt{n} \cdot ||M_\nu||^{s+1}\sqrt{\sum_{k \ne \nu}{\Tilde{D}^2_{k \nu}}} \cdot \frac{||\Lambda^{1/2}\mathcal{R}_\nu \mathcal{D}_\nu e_\nu||}{\sqrt{l_\nu}}
\end{equation}
because \(||M_\nu|| \xrightarrow[]{a.s.} 0.\) As for (\ref{25}), it can be shown that for any \(\delta_n \to \infty,\) 
\begin{equation}\label{1002}
    \frac{n}{M\delta_n}\sum_{k \ne \nu}{\Tilde{D}^2_{k \nu}} \xrightarrow[]{p} 0, \hspace{0.4cm} \frac{n}{M\delta_n} \cdot \frac{||\Lambda^{1/2}\mathcal{R}_\nu \mathcal{D}_\nu e_\nu||^2}{l_\nu} \xrightarrow[]{p} 0,
\end{equation}
providing
\begin{equation}\label{108}
    \frac{n}{M\delta_n}\sqrt{\sum_{k \ne \nu}{\Tilde{D}^2_{k \nu}}}\cdot \frac{||\Lambda^{1/2}\mathcal{R}_\nu \mathcal{D}_\nu e_\nu||}{\sqrt{l_\nu}} \xrightarrow[]{p} 0,
\end{equation}
and to obtain (\ref{31}), in light of (\ref{107}) and (\ref{108}), it suffices to prove that for \(\delta_n=\min{(\frac{n}{M},\frac{l_\nu}{M/\sqrt{n}}}) \to \infty,\)
\begin{equation}\label{110}
    \frac{M\delta_n}{\sqrt{n}} \cdot ||M_\nu||^{s+1} \xrightarrow[]{p} 0:
\end{equation}
using (\ref{12}) and (\ref{781}), with probability one,
\[||M_\nu|| \leq ||\mathcal{R}_\nu \mathcal{D}_\nu||+c(\epsilon_0) \cdot |\frac{\hat{l}_\nu}{l_\nu}-1| \leq ||\mathcal{R}_\nu \Lambda^{1/2} \Tilde{T}\Lambda^{1/2}||+||\mathcal{R}_\nu \Lambda^{1/2} \Tilde{D}\Lambda^{1/2}||+c(\epsilon_0) \cdot |\frac{\hat{l}_\nu}{l_\nu}-1| \leq \]
\[\leq c(K,\gamma,\epsilon_0) \cdot \frac{1}{l_\nu}+c(\epsilon_0) \cdot ||\frac{1}{n}Z_AZ_Z^T-I||+c(\epsilon_0) \cdot (\frac{\lambda_1(\mathcal{B})}{l_\nu}+|\lambda_\nu(\sum_{1 \leq j \leq \nu}{Z_j^TZ_j})-1|+|\lambda_1(\sum_{\nu \leq j \leq M}{Z_j^TZ_j})-1|),\]
and thus (\ref{7}) yields that in probability, for any \(\Tilde{\delta}_n \to \infty,\)
\begin{equation}\label{109}
    ||M_\nu|| \leq c(K,\gamma,\epsilon_0)(\frac{1}{l_\nu}+\sqrt{\frac{M\Tilde{\delta}_n}{n}}).
\end{equation}
Note that (\ref{109}) for \(\Tilde{\delta}_n=(\frac{n}{M})^{1/3} \to \infty\) gives (\ref{110}): for \(n\) large enough, it implies \(\log{\frac{1}{||M_\nu||}} \geq \frac{1}{4}\min{(\log{l_\nu},\log{\frac{n}{M}})},\) and
\[(s+1)\min{(\log{l_\nu},\log{\frac{n}{M}})}-4\log{\frac{M\delta_n}{\sqrt{n}}} \to \infty\] 
because \(\delta_n=\min{(\frac{n}{M},\frac{l_\nu}{M/\sqrt{n}}}) \to \infty,\)
\[(s+1)\log{l_\nu}-4\log{\frac{M\delta_n}{\sqrt{n}}} \geq (s+1)\log{l_\nu}-4\log{l_\nu} \geq 5\log{l_\nu} \to \infty,\]
\[(s+1)\log{\frac{n}{M}}-4\log{\frac{M\delta_n}{\sqrt{n}}} \geq (s+1)\log{\frac{n}{M}}-2\log{n} \to \infty.\]

\subsection{First Sum, Stage II: First Type Terms}\label{4.5}

Having established (\ref{30}) and (\ref{31}), turn to the remaining term in (\ref{27}):
\begin{equation}\label{111}
    \sum_{0 \leq j \leq s}{\frac{1}{\sqrt{l_\nu}}\sum_{k \ne \nu}{\Tilde{D}_{k \nu}\sqrt{l_k}((-M_\nu)^j\mathcal{R}_\nu \mathcal{D}_\nu e_\nu)_k}}.
\end{equation}
Denote by \(\Tilde{R}=l_\nu\mathcal{R}_\nu,\) and recall that
\[-M_\nu=-\mathcal{R}_\nu \mathcal{D}_\nu+(\hat{l}_\nu-l_\nu)\mathcal{R}_\nu=-\mathcal{R}_\nu \Lambda^{1/2} (\Tilde{T}+\Tilde{D})\Lambda^{1/2}+(\frac{\hat{l}_\nu}{l_\nu}-1)\Tilde{R};\]
take
\[-M_\nu=A_\nu+B_\nu, \hspace{0.2cm} A_\nu=-\mathcal{R}_\nu \Lambda^{1/2} \Tilde{T}\Lambda^{1/2}, \hspace{0.2cm} B_\nu=-\mathcal{R}_\nu \Lambda^{1/2} \Tilde{D}\Lambda^{1/2}+(\frac{\hat{l}_\nu}{l_\nu}-1)\Tilde{R}.\]
In this subsection, we show the terms in (\ref{111}) containing some factor of \(A_\nu\) (i.e., first type) are negligible.
\par
Since 
\[(A+B)^{j+1}-B^{j+1}=\sum_{0 \leq m \leq j}{(A+B)^{j-m}AB^m}\] 
for all \(j \geq 0, A,B \in \mathbb{R}^{p \times p}\) (by induction, and \((A+B)^{j+2}-B^{j+2}=((A+B)^{j+1}-B^{j+1})B+(A+B)^{j+1}A\)),
\[\sum_{1 \leq j \leq s}{\frac{1}{\sqrt{l_\nu}}\sum_{k \ne \nu}{\Tilde{D}_{k \nu}\sqrt{l_k}(((A_\nu+B_\nu)^j-B_\nu^j)\mathcal{R}_\nu \mathcal{D}_\nu e_\nu)_k}}=\]
\[=\sum_{1 \leq j \leq s, 0 \leq m \leq j-1}{\frac{1}{\sqrt{l_\nu}}\sum_{k \ne \nu}{\Tilde{D}_{k \nu}\sqrt{l_k}((A_\nu+B_\nu)^{j-m-1}A_\nu B_\nu^m\mathcal{R}_\nu \mathcal{D}_\nu e_\nu)_k}},\]
from which almost surely,
\[\sqrt{n}|\sum_{1 \leq j \leq s, 0 \leq m \leq j-1}{\frac{1}{\sqrt{l_\nu}}\sum_{k \ne \nu}{\Tilde{D}_{k \nu}\sqrt{l_k}(((A_\nu+B_\nu)^j-B_\nu^j)\mathcal{R}_\nu \mathcal{D}_\nu e_\nu)_k}}| \leq\]
\[\leq \sqrt{n} \cdot \sqrt{\sum_{k \ne \nu}{\Tilde{D}_{k\nu}^2}} \cdot \frac{||\Lambda^{1/2}\mathcal{R}_\nu \mathcal{D}_\nu e_\nu||}{\sqrt{l_\nu}} \cdot ||A_\nu||\sum_{1 \leq j \leq s, 0 \leq m \leq j-1}{||A_\nu+B_\nu||^{j-m-1} \cdot ||B_\nu||^m} \leq\]
\[\leq \sqrt{n} \cdot \sqrt{\sum_{k \ne \nu}{\Tilde{D}_{k\nu}^2}} \cdot \frac{||\Lambda^{1/2}\mathcal{R}_\nu \mathcal{D}_\nu e_\nu||}{\sqrt{l_\nu}}\cdot \frac{c(K,\gamma,\epsilon_0)}{l_\nu} \sum_{1 \leq j \leq s}{\frac{j}{2^{j-1}}} \leq \frac{4c(K,\gamma,\epsilon_0)}{l_\nu} \cdot \sqrt{n} \cdot \sqrt{\sum_{k \ne \nu}{\Tilde{D}_{k\nu}^2}} \cdot \frac{||\Lambda^{1/2}\mathcal{R}_\nu \mathcal{D}_\nu e_\nu||}{\sqrt{l_\nu}}\]
using \(||A_\nu|| \leq \frac{c(K,\gamma,\epsilon_0)}{l_\nu},||B_\nu|| \xrightarrow[]{a.s.} 0.\) This last inequality in conjunction with (\ref{108}) for \(\delta_n=\frac{l_\nu}{M/\sqrt{n}} \to \infty\) further yields
\begin{equation}\label{5080}
    \sum_{1 \leq j \leq s, 0 \leq m \leq j-1}{\frac{1}{\sqrt{l_\nu}}\sum_{k \ne \nu}{\Tilde{D}_{k \nu}\sqrt{l_k}(((A_\nu+B_\nu)^j-B_\nu^j)\mathcal{R}_\nu \mathcal{D}_\nu e_\nu)_k}}=o_p(\frac{1}{\sqrt{n}}),
\end{equation}
completing the analysis of the first type terms in (\ref{27}).

\subsection{First Sum, Stage III: Second Type Terms}\label{4.6}

This subsection concludes the proof of (\ref{27}): (\ref{30}), (\ref{31}), and (\ref{5080}) imply 
\[\tilde{\Sigma}_1=\sum_{0 \leq j \leq s}{\frac{1}{\sqrt{l_\nu}}\sum_{k \ne \nu}{\Tilde{D}_{k \nu}\sqrt{l_k}((-\mathcal{R}_\nu \Lambda^{1/2} \Tilde{D}\Lambda^{1/2}+(\frac{\hat{l}_\nu}{l_\nu}-1)\Tilde{R})^j\mathcal{R}_\nu \mathcal{D}_\nu e_\nu)_k}}+o_p(\frac{1}{\sqrt{n}}),\]
and we show next that
\begin{equation}\label{112}
    \sum_{0 \leq j \leq s}{\frac{1}{\sqrt{l_\nu}}\sum_{k \ne \nu}{\Tilde{D}_{k \nu}\sqrt{l_k}((-\mathcal{R}_\nu \Lambda^{1/2} \Tilde{D}\Lambda^{1/2}+(\frac{\hat{l}_\nu}{l_\nu}-1)\Tilde{R})^j\mathcal{R}_\nu \mathcal{D}_\nu e_\nu)_k}}
\end{equation}
is the polynomial in \(\frac{\hat{l}_\nu}{l_\nu}-1\) on the right-hand side of (\ref{27}) up to an \(o_p(\frac{1}{\sqrt{n}})\) error.
\par
Once again \(\mathcal{R}_\nu \mathcal{D}_\nu\) can be replaced by \(\mathcal{R}_\nu \Lambda^{1/2} \Tilde{D}\Lambda^{1/2}\) since this substitution generates a negligible error:
\[|\sum_{0 \leq j \leq s}{\frac{1}{\sqrt{l_\nu}}\sum_{k \ne \nu}{\Tilde{D}_{k \nu}\sqrt{l_k}(B_\nu^j\mathcal{R}_\nu \Lambda^{1/2}\Tilde{T} \Lambda^{1/2} e_\nu)_k}}|=|\sum_{0 \leq j \leq s}{\sum_{k \ne \nu}{\Tilde{D}_{k \nu}\sqrt{l_k}(B_\nu^j\mathcal{R}_\nu \Lambda^{1/2}\Tilde{T} e_\nu)_k}}| \leq\]
\[\leq \sum_{0 \leq j \leq s}{\sqrt{\sum_{k \ne \nu}{\Tilde{D}^2_{k \nu}}} \cdot ||B_\nu||^j \cdot ||\Lambda \mathcal{R}_\nu\Tilde{T} e_\nu||} \leq 2c(K,\gamma, \epsilon_0) \sqrt{\sum_{k \ne \nu}{\Tilde{D}^2_{k \nu}}} \cdot \sqrt{\sum_{k \ne \nu}{\Tilde{T}^2_{k \nu}}}=o_p(\frac{1}{\sqrt{n}})\]
employing (\ref{1001}) and (\ref{1002}) 
\[\frac{n}{M\delta_n}\sum_{k \ne \nu}{\Tilde{D}^2_{k \nu}} \xrightarrow[]{p} 0, \hspace{0.4cm} M\delta_n\sum_{k \ne \nu}{\Tilde{T}^2_{k \nu}} \xrightarrow[]{p} 0\] 
for \(\delta_n=\frac{l_\nu}{M/\sqrt{n}} \to \infty.\)
\par
Hence, the remaining contribution in (\ref{112}) is
\begin{equation}\label{113}
    \sum_{0 \leq j \leq s}{\sum_{k \ne \nu}{\Tilde{D}_{k \nu}\sqrt{l_k}((-\mathcal{R}_\nu \Lambda^{1/2} \Tilde{D}\Lambda^{1/2}+(\frac{\hat{l}_\nu}{l_\nu}-1)\Tilde{R})^j\mathcal{R}_\nu \Lambda^{1/2} \Tilde{D})_{k\nu}}}.
\end{equation}
Think of this sum as a polynomial in \(\frac{\hat{l}_\nu}{l_\nu}-1:\) then the sources of randomness behind its coefficients are products such as
\[\Tilde{D}_{k_1 \nu}\Tilde{D}_{k_1 k_2}...\Tilde{D}_{k_m \nu}\]
for \(k_1, k_2, \hspace{0.05cm} ...\hspace{0.05cm}, k_m \in \{1,2,\hspace{0.05cm} ... \hspace{0.05cm}, M\}-\{\nu\}, m \leq s,\) most of which have expectation \(\frac{1}{n^m}\) (whenever the \(m\)-tuple consists of pairwise distinct elements). The deterministic coefficients \(-a_0,\hspace{0.05cm} ... \hspace{0.05cm}, -a_s\) thus arise after replacing these products by this shared first moment, substitution that produces an \(o_p(\frac{1}{\sqrt{n}})\) error. In this subsection, we prove these two claims together with the bounds \(|a_m| \leq c(\epsilon_0)^{m+1} \cdot \frac{M}{n}\) separately for \(m=0\) and \(m>0\) inasmuch as the former case, slightly simpler than the latter, sheds light into the core of the argument.
\vspace{0.2cm}
\par
\textit{\(1. \hspace{0.2cm} m=0:\)} Consider the terms in (\ref{113}) with no factors of \(\frac{\hat{l}_\nu}{l_\nu}-1,\) that is:
\begin{equation}\label{141}
    \sum_{0 \leq j \leq s}{\sum_{k \ne \nu}{\Tilde{D}_{k \nu}\sqrt{l_k}((-\mathcal{R}_\nu \Lambda^{1/2} \Tilde{D}\Lambda^{1/2})^j \mathcal{R}_\nu \Lambda^{1/2} \Tilde{D} )_{k\nu}}}.
\end{equation}
Note that the sum corresponding to \(j, 0 \leq j \leq s\) is
\[-\sum_{(k_0,k_1,k_2,...,k_j), k_i \ne \nu}{ \frac{l_{k_0}l_{k_1}...l_{k_j}}{(l_\nu-l_{k_0})(l_\nu-l_{k_1})...(l_\nu-l_{k_j})}\Tilde{D}_{k_0 \nu}\Tilde{D}_{k_0 k_1}...\Tilde{D}_{k_{j-1} k_j} \Tilde{D}_{k_j \nu}},\]
because
\[(-\mathcal{R}_\nu \Lambda^{1/2} \Tilde{D}\Lambda^{1/2})_{\nu i}=0, \hspace{0.1cm} 1 \leq i \leq M, \hspace{0.4cm} (-\mathcal{R}_\nu \Lambda^{1/2} \Tilde{D}\Lambda^{1/2})_{ki}=\frac{\sqrt{l_kl_i}}{l_\nu-l_k}\Tilde{D}_{ki}, \hspace{0.1cm} k \ne \nu.\]
\par
Let 
\[a_{j,t}=\sum_{(k_1,k_2,...,k_j),(k'_1,k'_2,...,k'_t), k_i,k'_i \ne \nu}{|\mathbb{E}[(\Tilde{D}_{k_1 \nu}\Tilde{D}_{k_1 k_2}...\Tilde{D}_{k_{j-1} k_j} \Tilde{D}_{k_j \nu}-\frac{1}{n^{j}}) \cdot (\Tilde{D}_{k'_1 \nu}\Tilde{D}_{k'_1 k'_2}...\Tilde{D}_{k'_{t-1} k'_t} \Tilde{D}_{k'_t \nu}-\frac{1}{n^{t}})]|},\]
where \(1 \leq j,t \leq 2s+1,\) which will be shown to satisfy
\begin{equation}\label{34}
    a_{j,t} \leq \frac{(j+t)!}{n} \cdot (\frac{M}{n})^{\frac{j+t}{2}} \cdot (2^{k_0}(j+t))^{3(j+t)+4}
\end{equation}
for \[2^{k_0}=64+64c_1(K)+64(c_1(K))^2+64(c_1(K))^4, \hspace{0.2cm} c_1(K)=64+64c_0(K)+64(c_0(K))^2+64(c_0(K))^4,\] 
where \(c_0(K)=\max_{p \geq 1, \mathbb{E}[z]=0, ||z||_{\psi_2} \leq K}{\frac{||z||_{L^p}}{\sqrt{p}}}<\infty\) from proposition \(2.5.2\) in Vershynin~\cite{vershynin}. 
\par
Suppose (\ref{34}) holds, and let us see how it yields that (\ref{141}) is
\begin{equation}\label{117}
    -\sum_{(k_0,k_1,k_2,...,k_j), k_i \ne \nu, 0 \leq j \leq s}{\frac{1}{n^{j+1}} \cdot \frac{l_{k_0}l_{k_1}...l_{k_j}}{(l_\nu-l_{k_0})(l_\nu-l_{k_1})...(l_\nu-l_{k_j})}}+o_p(\frac{1}{\sqrt{n}}):
\end{equation}
notice that (\ref{117}) completes the analysis of the case \(m=0\) since the first term in it is \(-a_0\) using (\ref{129}) and (\ref{130}) while for \(n\) large enough
\[|\sum_{(k_0,k_1,k_2,...,k_j), k_i \ne \nu, 0 \leq j \leq s}{\frac{1}{n^{j+1}} \cdot  \frac{l_{k_0}l_{k_1}...l_{k_j}}{(l_\nu-l_{k_0})(l_\nu-l_{k_1})...(l_\nu-l_{k_j})}}|=|\sum_{0 \leq j \leq s}{(\frac{1}{n}\sum_{k \ne \nu}{\frac{l_k}{l_\nu-l_k}})^{j+1}}| \leq \]
\[\leq 2 \cdot \frac{1}{n}\sum_{k \ne \nu}{\frac{l_k}{|l_\nu-l_k|}} \leq c(\epsilon_0) \cdot \frac{M}{n}.\]
\par
Return now to (\ref{117}): for any \(1 \leq j \leq s+1,\) (\ref{34}) and \((2j)! \leq (2j)^{2j}=(4j^2)^j, 6j+4 \leq 10j\) yield
\[a_{j,j} \leq \frac{(2j)!}{n} \cdot (\frac{M}{n})^{j} \cdot (2^{k_0+1}j)^{6j+4} \leq \frac{1}{n} \cdot (\frac{M j^{10} \cdot 2^{10k_0+12}}{n})^{j},\]
\[n\mathbb{E}[(\sum_{(k_1,k_2,...,k_j), k_i \ne \nu}{ \frac{l_{k_1}...l_{k_j}}{(l_\nu-l_{k_1})...(l_\nu-l_{k_j})}(\Tilde{D}_{k_1 \nu}\Tilde{D}_{k_1 k_2}...\Tilde{D}_{k_{j-1} k_j} \Tilde{D}_{k_j \nu}-\frac{1}{n^{j}})})^2] \leq (\frac{Mj^{10}c(K,\epsilon_0)}{n})^{j},\]
where we have employed \(\frac{l_k}{|l_\nu-l_k|} \leq c(\epsilon_0), k \ne \nu,\) from which (\ref{117}) ensues because for \(t>0,\) 
\[S(j)=\sum_{(k_0,k_2,...,k_{j}), k_i \ne \nu}{ \frac{l_{k_0}...l_{k_j}}{(l_\nu-l_{k_0})...(l_\nu-l_{k_j})}(\Tilde{D}_{k_0 \nu}\Tilde{D}_{k_0 k_1}...\Tilde{D}_{k_{j-1} k_{j}} \Tilde{D}_{k_{j} \nu}-\frac{1}{n^{j+1}})},\]
\[\mathbb{P}(\sqrt{n}|\sum_{0 \leq j \leq s}{S(j)}| \geq t) \leq \sum_{0 \leq j \leq s}{\mathbb{P}(\sqrt{n}|S(j)| \geq \frac{t}{2^{j+1}})} \leq  \sum_{0 \leq j \leq s}{\frac{2^{2(j+1)}\mathbb{E}[nS_j^2]}{t^2}} \leq \]
\begin{equation}\label{999}
    \leq \frac{1}{t^2}\sum_{0 \leq j \leq s}{(\frac{4M(j+1)^{10}c(K, \epsilon_0)}{n})^{j+1}} \leq \frac{2}{t^2} \cdot \frac{4M(s+1)^{10}c(K, \epsilon_0)}{n} \to 0
\end{equation}
as
\begin{equation}\label{5019}
    \frac{s}{(n/M)^{1/10}} \leq \frac{16\log{n}}{(n/M)^{1/10}\log{(n/M)}} \leq \frac{160\log{n}}{(\log{(n/M)})^2} \to 0.
\end{equation}
\par
We proceed with the case \(m>0\) and present the proof of (\ref{34}) in the coming subsection.

\vspace{0.2cm}
\par
\textit{\(2. \hspace{0.2cm} m>0:\)}  Consider \(m\) with \(1 \leq m \leq s,\) and the terms having exactly a factor of \((\frac{\hat{l}_\nu}{l_\nu}-1)^m\) for some fixed \(j\) with \(m \leq j \leq s\) in
\[\sum_{0 \leq j \leq s}{\sum_{k \ne \nu}{\Tilde{D}_{k \nu}\sqrt{l_k}((-\mathcal{R}_\nu \Lambda^{1/2} \Tilde{D}\Lambda^{1/2}+(\frac{\hat{l}_\nu}{l_\nu}-1)\Tilde{R})^j\mathcal{R}_\nu \Lambda^{1/2} \Tilde{D})_{k\nu}}},\]
which are
\begin{equation}\label{941}
    \sum{\Tilde{D}_{k \nu}\sqrt{l_k}((-\mathcal{R}_\nu \Lambda^{1/2} \Tilde{D}\Lambda^{1/2})^{j_1}\Tilde{R}^{j_2}(-\mathcal{R}_\nu \Lambda^{1/2} \Tilde{D}\Lambda^{1/2})^{j_3}\Tilde{R}^{j_4}...\Tilde{R}^{j_{2l}}\mathcal{R}_\nu \Lambda^{1/2} \Tilde{D})_{k \nu}},
\end{equation}
where the summation is over \(k \ne \nu\) and integer tuples \((j_1,j_2,\hspace{0.05cm} ... \hspace{0.05cm}, j_{2l})\) with 
\[j_2+j_4+...+j_{2l}=m, j_2, \hspace{0.05cm} ... \hspace{0.05cm}, j_{2l-2}>0, j_{2l} \geq 0, j_1+j_3+...+j_{2l-1}=j-m, j_1 \geq 0, j_3,\hspace{0.05cm} ...\hspace{0.05cm}, j_{2l-1}>0.\] 
Consequently (\ref{941}) is a sum over such \((j_1,j_2,\hspace{0.05cm} ... \hspace{0.05cm}, j_{2l}),\) whose contribution is
\[-\sum{\frac{l_kl_{k_1}...l_{k_{j-m}}}{(l_\nu-l_k)(l_\nu-l_{k_1})...(l_\nu-l_{k_{j-m}})} \cdot (\frac{l_\nu}{l_{k_{j_1}-l_\nu}})^{j_2}...(\frac{l_\nu}{l_{k_{j_{2l-1}}-l_\nu}})^{j_{2l}}\Tilde{D}_{k\nu}\Tilde{D}_{kk_1}...\Tilde{D}_{k_{j-m}\nu}},\]
where the summation is over \((k, k_1,k_2, \hspace{0.05cm} ... \hspace{0.05cm}, k_{j-m}),k, k_i \ne \nu.\) Similarly to the case \(m=0,\) we prove next that by replacing these random products by \(\frac{1}{n^{j-m+1}},\) the expectation of most such terms, the error is negligible. That is, by putting together all \(j\) with \(m \leq j \leq s,\) and letting \(m\) vary in \(\{1,2, \hspace{0.05cm} ... \hspace{0.05cm}, s\}\) we show that these sums add up to 
\[-\sum_{1 \leq m \leq s}{\sum_{m \leq j \leq s}{\sum{\frac{1}{n^{j-m+1}} \cdot \frac{l_kl_{k_1}...l_{k_{j-m}}}{(l_\nu-l_k)(l_\nu-l_{k_1})...(l_\nu-l_{k_{j-m}})} \cdot (\frac{l_\nu}{l_{k_{j_1}-l_\nu}})^{j_2}...(\frac{l_\nu}{l_{k_{j_{2l-1}}-l_\nu}})^{j_{2l}}}}}+o_p(\frac{1}{\sqrt{n}}),\]
where the \(m^{th}\) constituent is \(-a_m\) (an easy consequence of (\ref{129}) and (\ref{130})) with \(|a_m| \leq c(\epsilon_0)^{m+1} \cdot \frac{M}{n},\) completing the analysis of the case \(m>0.\) 
\par
Consider first the deterministic coefficients \(a_m, 1 \leq m \leq s:\)
\[a_m=(\frac{1}{n}\sum_{k \ne \nu}{\frac{l_k}{l_\nu-l_k}})^{j-m-(l-1)}\cdot \frac{1}{n^{l}} \cdot  \frac{l_{k_{j_1}}}{l_\nu-l_{k_{j_1}}} (\frac{l_\nu}{l_{k_{j_1}-l_\nu}})^{j_2} \cdot ... \cdot \frac{l_{k_{j_{2l-1}}}}{l_\nu-l_{k_{j_{2l-1}}}}(\frac{l_\nu}{l_{k_{j_{2l-1}}-l_\nu}})^{j_{2l}}\]
where the summation is over \((k_{j_1},k_{j_3}, \hspace{0.05cm} ... \hspace{0.05cm}, k_{j_{2l-1}},j_{2}, \hspace{0.05cm} ... \hspace{0.05cm}, j_{2l})\) with \(k_i \ne \nu,\) and \[j_2+j_4+...+j_{2l}=m, \hspace{0.2cm} j_2, \hspace{0.05cm} ... \hspace{0.05cm}, j_{2l-2}>0,\hspace{0.05cm} j_{2l} \geq 0.\]
Since \(\frac{l_\nu}{|l_k-l_\nu|},\frac{l_k}{|l_\nu-l_k|} \leq c(\epsilon_0),\) for \(\sigma=\frac{1}{n}|\sum_{k \ne \nu}{\frac{l_k}{l_\nu-l_k}}|,\)
\[|a_m| \leq \sum_{m \leq j \leq s, 1 \leq l \leq \min{(m+1,j-m+1)}}{\sigma^{j-m-l+1}\cdot \frac{1}{n^l} \cdot c(\epsilon_0)^{m+l} \cdot M^l \cdot \binom{m}{l-1}}\]
because for each fixed \(l,\) there are at most \(M^{l}\) tuples \((k_{j_1},k_{j_3},\hspace{0.05cm} ...\hspace{0.05cm}, k_{j_{2l-1}}),\) and the number of tuples \((j_2, \hspace{0.05cm} ... \hspace{0.05cm}, j_{2l})\) is at most \(\Tilde{a}_{l,m-(l-1)}=\binom{(m-(l-1))+l-1}{l-1}=\binom{m}{l-1}\) from (\ref{4001}), and the constraints on \(l\) hold since \(m>0\) implies \(l>0\) (a sum of \(l\) non-negative numbers is \(m\)), \(j_2+j_4+...+j_{2l}=m \geq l-1,\) and \(l \leq j-m+1\) as a vestige of \(j_1+j_3+...+j_{2l-1}=j-m \geq l-1.\) Therefore, for each fixed \(l \geq 1,\) the conditions on \(j\) are \(l+m-1 \leq j \leq s,\) entailing with \(\sigma \leq \frac{Mc(\epsilon_0)}{n} \leq \frac{1}{2}\) for \(n\) sufficiently large
\[|a_m| \leq c(\epsilon_0)^m \sum_{1 \leq l \leq m+1}{\frac{\sigma^{l+m-1}-\sigma^{s+1}}{1-\sigma}  \cdot \sigma^{-m-l+1}\cdot (\frac{Mc(\epsilon_0)}{n})^l \cdot \binom{m}{l-1}} \leq\]
\[\leq \frac{c(\epsilon_0)^m}{1-\sigma}\sum_{1 \leq l \leq m+1}{(\frac{Mc(\epsilon_0)}{n})^l \cdot \binom{m}{l-1}}=\frac{Mc(\epsilon_0)}{n} \cdot \frac{c(\epsilon_0)^m}{1-\sigma}(1+\frac{Mc(\epsilon_0)}{n})^{m} \leq (2c(\epsilon_0)(1+c(\epsilon_0)))^{m+1} \cdot \frac{M}{n},\]
yielding the claim \(|a_m| \leq c(\epsilon_0)^{m+1} \cdot \frac{M}{n}.\)
\par
Consider now the overall error given by
\[\sum_{1 \leq m \leq s}{S'(m)},\]
\[S'(m):=\sum{\frac{l_kl_{k_1}...l_{k_{j-m}}}{(l_\nu-l_k)(l_\nu-l_{k_1})...(l_\nu-l_{k_{j-m}})} \cdot (\frac{l_\nu}{l_{k_{j_1}-l_\nu}})^{j_2}...(\frac{l_\nu}{l_{k_{j_{2l-1}}-l_\nu}})^{j_{2l}}(\Tilde{D}_{k\nu}\Tilde{D}_{kk_1}...\Tilde{D}_{k_{j-m}\nu}-\frac{1}{n^{j-m+1}})}:\]
to conclude it is \(o_p(\frac{1}{\sqrt{n}}),\) it suffices to show
\[\sum_{1 \leq m \leq s}{2^{2m}\mathbb{E}[n(S'(m))^2]}=o(1),\]
and use the analogue of (\ref{999}). 
\par
Fix \(m\) and \(k,k_1, \hspace{0.05cm} ...\hspace{0.05cm}, k_{j-m};\) then
\[|\frac{l_kl_{k_1}...l_{k_{j-m}}}{(l_\nu-l_k)(l_\nu-l_{k_1})...(l_\nu-l_{k_{j-m}})} \cdot (\frac{l_\nu}{l_{k_{j_1}-l_\nu}})^{j_2}...(\frac{l_\nu}{l_{k_{j_{2l-1}}-l_\nu}})^{j_{2l}}| \leq c(\epsilon_0)^{j+1},\]
there are at most \(\binom{j-m+1}{l}+\binom{j-m+1}{l-1}=\binom{j-m+2}{l}\) ways of choosing \(k_{j_1},\hspace{0.05cm} ...\hspace{0.05cm}, k_{j_{2l-1}}\) (either \(l\) or \(l-1\) unordered positions have to be chosen depending on \(j_{1}\) being nonzero or not), and at most \(\binom{m}{l-1}\) tuples \((j_2, \hspace{0.05cm} ... \hspace{0.05cm}, j_{2l})\) for each fixed \(l.\) Thus, the number of apparitions of each such product is upper bounded by
\[\sum_{l \geq 0}{\binom{j-m+2}{l} \binom{m}{m-l+1}}=\binom{j+2}{m+1},\]
from which
\[\mathbb{E}[(S'(m))^2] \leq \sum_{m \leq j_1,j_2 \leq s}{c(\epsilon_0)^{j_1+1} \cdot c(\epsilon_0)^{j_2+1} \cdot \binom{j_1+2}{m+1} \binom{j_2+2}{m+1}a_{j_1-m+1,j_2-m+1}}.\]
Note that
\[\sum_{j_1+j_2=j}{\binom{j_1+2}{m+1} \binom{j_2+2}{m+1}} \leq \sum_{j_1+j_2=j}{2^{j_1+2} \cdot 2^{j_2+2}}=2^{j+4}(j+1) \leq 2^{2j+4},\]
and for \(j=j_1+j_2-2m+2\)
\[\mathbb{E}[(S'(m))^2] \leq \sum_{2 \leq j \leq 2s}{c(\epsilon_0)^{j+2m}2^{2j+4m}\max_{l_1+l_2=j}{a_{l_1,l_2}}} \leq (4c(\epsilon_0))^{2m} \sum_{2 \leq j \leq 2s}{c(\epsilon_0)^{j}2^{2j} \cdot \frac{j!}{n} \cdot (\frac{M}{n})^{j/2} \cdot (2^{k_0}j)^{7j}} \leq\]
\[\leq \frac{(4c(\epsilon_0))^{2m}}{n} \cdot \sum_{2 \leq j \leq 2s}{(c(K,\epsilon_0)\sqrt{\frac{M}{n}}j^8)^{j}} \leq \frac{(4c(\epsilon_0))^{2m}}{n} \sum_{2 \leq j \leq 2s}{(c(K,\epsilon_0)\sqrt{\frac{M}{n}}s^8)^{j}}  \leq \frac{(4c(\epsilon_0))^{2m}}{n} \cdot \frac{2M}{n}s^{16}c(K,\epsilon_0),\]
where the second inequality and the last ensue from (\ref{34}) and an analogue of (\ref{5019}), respectively.
\par
Finally, this last bound implies 
\[\sum_{1 \leq m \leq s}{2^{2m}\mathbb{E}[n(S'(m))^2]} \leq 2c(K,\epsilon_0) \cdot \frac{Ms^{16}}{n} \sum_{1 \leq m \leq s}{(8c(\epsilon_0))^{2m}} \leq 2c(K,\epsilon_0) \cdot \frac{Ms^{16}}{n} \cdot (8c(\epsilon_0)+1)^{2s}=o(1)\]
from 
\[2s \cdot 8c(\epsilon_0)+16\log{s}-\log{\frac{n}{M}} \leq 16c(\epsilon_0) \cdot \frac{9\log{n}}{\log{\frac{n}{M}}} +16 \cdot \frac{9\log{n}}{\log{\frac{n}{M}}} -\log{\frac{n}{M}} \leq -\frac{1}{2} \log{\frac{n}{M}} \to -\infty\]
for \(n\) sufficiently large.

\subsection{Controlling Errors: A Combinatorial Lemma}\label{4.7}

In this subsection, we justify (\ref{34}) for which the key is the following inequality:
\begin{equation}\label{35}
    a_{j,t} \leq \frac{j}{n}a_{j-1,t}+\frac{t}{n}a_{j,t-1}+\frac{M^{\frac{t+j}{2}} \cdot (\frac{t+j}{2}+1)^{t+j}}{n^{\frac{t+j}{2}+1}} \cdot (\alpha(j)\alpha(t)+\alpha(j)+\alpha(t)+1)
\end{equation}
for \(2 \leq j,t \leq 2s+1,\) where 
\[\alpha(k)=(c_1(K)(k+1))^{\frac{k+1}{2}}+\frac{1}{\sqrt{n}}(c_1(K)(k+1))^{2k+2}.\]
\par
To see why (\ref{35}) holds, notice there are two types of pairs of tuples \((k_1, k_2,\hspace{0.05cm} ...\hspace{0.05cm}, k_j),(k'_1, k'_2, \hspace{0.05cm}... \hspace{0.05cm}, k'_t):\)
\par
\(1.\) There exists an index \(k_i\) or \(k'_{i'}\) such that its multiplicity in the multiset \(k_1, k_2, \hspace{0.05cm} ...\hspace{0.05cm}, k_j, k'_1, k'_2, \hspace{0.05cm} ...\hspace{0.05cm}, k'_t\) is one (choose \(i\) or \(i'\) to be minimal). Suppose \(k_i\) has the property just stated: then
\[\mathbb{E}[\Tilde{D}_{k_1 \nu}\Tilde{D}_{k_1 k_2}...\Tilde{D}_{k_{j-1} k_j} \Tilde{D}_{k_j \nu}]=\frac{1}{n}\mathbb{E}[\Tilde{D}_{k_1 \nu}\Tilde{D}_{k_1 k_2}...\Tilde{D}_{k_{i-2} k_{i-1}}\Tilde{D}_{k_{i-1} k_{i+1}}\Tilde{D}_{k_{i+1} k_{i+2}} ...\Tilde{D}_{k_{j-1} k_j} \Tilde{D}_{k_j \nu}]\]
since each term in the product will have two factors of the form \(z_{k_i l}, 1 \leq l \leq n,\) and they have to be identical to make a non-zero contribution to the expectation (that is, \(l_1=l_2\)), and thus after conditioning on \(z_{k_i},\) the product above ensues; a similar identity holds for the cross expectation, and so these pairs contribute at most \(\frac{j}{n}a_{j-1,t}+\frac{t}{n}a_{j,t-1}.\)
\par
\(2.\) Each value in the multiset \(k_1, k_2, \hspace{0.05cm} ...\hspace{0.05cm}, k_j, k'_1, k'_2,\hspace{0.05cm} ...\hspace{0.05cm}, k'_t\) has multiplicity at least two. Note that
\begin{equation}\label{72}
    \mathbb{E}[\Tilde{D}^2_{k_1 \nu}\Tilde{D}^2_{k_1 k_2}...\Tilde{D}^2_{k_{j-1} k_j} \Tilde{D}^2_{k_j \nu}] \leq \frac{\alpha^2(j)}{n^{j+1}}
\end{equation}
because
\[\mathbb{E}[\Tilde{D}^2_{k_1 \nu}\Tilde{D}^2_{k_1 k_2}...\Tilde{D}^2_{k_{j-1} k_j} \Tilde{D}^2_{k_j \nu}] \leq \frac{1}{j+1}(\mathbb{E}[\Tilde{D}^{2j+2}_{k_1 \nu}]+\mathbb{E}[\Tilde{D}^{2j+2}_{k_1 k_2}]+...+\mathbb{E}[\Tilde{D}^{2j+2}_{k_j \nu}]) \leq\]
\begin{equation}\label{2029}
    \leq \frac{1}{n^{j+1}}\max{(\mathbb{E}[(\frac{1}{\sqrt{n}} z_{1}^Tz_{2})^{2j+2}],\mathbb{E}[(\frac{1}{\sqrt{n}}||z_{1}||^2-\sqrt{n})^{2j+2}])} \leq \frac{1}{n^{j+1}}((c_1(K)(j+1))^{j+1}+\frac{1}{n} (c_1(K)(j+1))^{4j+4}):
\end{equation}
the products in
\[\mathbb{E}[(\frac{1}{\sqrt{n}} z_{1}^Tz_{2})^{2j+2}]\]
are of the form \((\frac{1}{\sqrt{n}}z_{1i_1}z_{2i_1})^{m_1} \cdot ... \cdot (\frac{1}{\sqrt{n}}z_{1i_p}z_{2i_p})^{m_p}\) for some 
\[i_1<...<i_p, \hspace{0.1cm} m_1+m_2+...+m_p=2j+2, \hspace{0.1cm} m_1,\hspace{0.05cm} ...\hspace{0.05cm}, m_p>0;\]
if \(m_i=1\) for some \(i,\) then the expectation of this product is zero; if some \(m_i>2,\) then \(p \leq j,\) and so the overall contribution of such terms is at most
\[\frac{1}{n^{j+1}}\sum_{1 \leq p \leq j}{\binom{n}{p} \sum_{m_1+m_2+...+m_p=2j+2}{\binom{2j+2}{m_1,m_2,\hspace{0.05cm} ... \hspace{0.05cm}, m_p}} \mathbb{E}[z_{11}^{4j+4}]} \leq \frac{1}{n^{j+1}}\sum_{1 \leq p \leq j}{n^{j}p^{2j+2} \mathbb{E}[z_{11}^{4j+4}]} \leq \]
\begin{equation}\label{69}
    \leq \frac{1}{n}\mathbb{E}[z_{11}^{4j+4}] \cdot \frac{(j+1)^{2j+3}}{2j+3} \leq \frac{1}{n}(c_0(K)\sqrt{4(j+1)})^{4(j+1)} \cdot (j+1)^{2j+2} \leq \frac{1}{n}(c_1(K)(j+1))^{4j+4}
\end{equation}
where we have used there are \(\binom{n}{p}\) ways of choosing \(i_1< i_2< ... <i_p,\) and the weighted arithmetic-geometric mean inequality
\begin{equation}\label{5084}\tag{AM-GM}
    x_1^{n_1}...x_p^{n_p} \leq \frac{1}{n_1+...+n_p}(n_1x_1^{n_1+...+n_p}+...+n_px_p^{n_1+...+n_p})
\end{equation}
for \(x_i \geq 0, n_i>0;\) lastly, when \(m_1=m_2=...=m_p=2, p=j+1,\) the expectation of each such term is one and their overall contribution will be 
\begin{equation}\label{70}
    \frac{1}{n^{j+1}}\binom{n}{j+1} \binom{2j+2}{2} \binom{2j}{2} ... \binom{2}{2}=\frac{(n-1)(n-2)...(n-j)}{n^j} (2j+2)!! \leq (2j+2)^{j+1};
\end{equation}
a similar analysis holds for \(\mathbb{E}[(\frac{1}{\sqrt{n}}||z_{1}||^2-\sqrt{n})^{2j+2}]\) with the distinction that in the analogue (\ref{69}), \(\mathbb{E}[z_{11}^{4j+4}]\) is replaced by 
\[\mathbb{E}[(z_{11}^2-1)^{2j+2}] \leq 2^{2j+2}\mathbb{E}[z_{11}^{4j+4}] \leq 2^{2j+2} (c_0(K)\sqrt{4j+4})^{4j+4} \leq (c_1(K)(j+1))^{2j+2},\] 
employing \((x+y)^{2j+2} \leq 2^{2j+2}\max{(x^{2j+2},y^{2j+2})},\)
and in the surrogate of (\ref{70}), \(\mathbb{E}[z_{11}^2]^{j+1}=1\) is substituted by 
\[\mathbb{E}[(z_{11}^2-1)^{2}]^{j+1} \leq \mathbb{E}[z^4_{11}]^{j+1} \leq (\sqrt{4}c_0(K))^{4(j+1)} \leq (\frac{c_1(K)}{2})^{j+1}.\]
Cauchy-Schwarz inequality yields 
\[|\mathbb{E}[(X-a)(Y-b)]| \leq |ab|+\sqrt{\mathbb{E}[X^2] \cdot \mathbb{E}[Y^2]}+|b|\sqrt{\mathbb{E}[X^2]}+|a|\sqrt{\mathbb{E}[Y^2]},\] 
and together with (\ref{72}) it implies the contribution of each such pair is at most
\[\frac{1}{n^{t+j}}+\frac{\alpha(j)}{n^{t+\frac{j+1}{2}}}+\frac{\alpha(t)}{n^{j+\frac{t+1}{2}}}+\frac{\alpha(j) \cdot \alpha(t)}{n^{\frac{t+j}{2}+1}} \leq \frac{\alpha(j)\alpha(t)+\alpha(j)+\alpha(t)+1}{n^{\frac{t+j}{2}+1}},\]
while their number is at most \(M^{\frac{t+j}{2}} \cdot (\frac{t+j}{2}+1)^{t+j}\) because for each \((t+j)\)-tuple of this form, let \(m_1,m_2, \hspace{0.05cm} ... \hspace{0.05cm}, m_l\) be the multiplicities of the distinct values in the tuple (increasingly ordered): then \(l \leq \frac{t+j}{2}\) from \(\newline m_i \geq 2, m_1+...+m_l=t+j,\) and thus their number is at most
\[\sum_{1 \leq l \leq \frac{t+j}{2},m_1+...+m_l=t+j}{\binom{t+j}{m_1, m_2, \hspace{0.05cm}...\hspace{0.05cm}, m_l}\binom{M}{l}} \leq \sum_{1 \leq l \leq \frac{t+j}{2}}{l^{t+j}\binom{M}{l}} \leq \]
\begin{equation}\label{2027}
    \leq \sum_{1 \leq l \leq \frac{t+j}{2}}{l^{t+j}M^l} \leq \sum_{1 \leq l \leq \frac{t+j}{2}}{l^{t+j}M^{(t+j)/2}} \leq M^{(t+j)/2} \cdot \frac{(\frac{t+j}{2}+1)^{t+j+1}}{t+j+1} \leq M^{(t+j)/2} \cdot (\frac{t+j}{2}+1)^{t+j},
\end{equation}
completing the justification of (\ref{35}).
\par
To finish the proof of (\ref{34}), we get upper bounds for \(a_{j,t}\) given by some \(\Tilde{a}_{j+t}.\) Notice that any sequence for which the following hold
\begin{equation}\label{37}
    \Tilde{a}_{l} \geq \frac{l}{n}\Tilde{a}_{l-1}+\frac{1}{n} \cdot (\frac{M}{n})^{l/2} \cdot (c_1(K)l)^{3l+4}, \hspace{0.5cm} \Tilde{a}_{l} \geq a_{1,l-1}
\end{equation}
for \(l \geq 2\) works using (\ref{35}), \((\frac{l}{2}+1)^{l} \leq l^l, \alpha(j)\alpha(t) \leq  c_1(K)^{2j+2t+4} \cdot (j+t)^{2j+2t+4}\) from
\begin{equation}\label{2028}
    \alpha(k)=(c_1(K)(k+1))^{\frac{k+1}{2}}+\frac{1}{\sqrt{n}}(c_1(K)(k+1))^{2k+2} \leq (c_1(K)(k+1))^{2k+2}
\end{equation}
and the upper bound of \(\alpha(j) \alpha(t)\) reused for \(\alpha(j)+\alpha(t)+1\) inasmuch as \(1+x+y \leq xy\) for \(x,y \geq 3.\) 
\par
Denote by
\[\Tilde{a}_l=\frac{l!}{n} \cdot (\frac{M}{n})^{l/2} \cdot \overline{a}_l,\] 
and we show inductively that \(\overline{a}_{l}=(2^{k_0}l)^{3l+4}\) satisfies (\ref{37}). For the base case, take \(\tilde{a}_2=a_{1,1}\) with
\[a_{1,1}=\sum_{k_1 \ne \nu,k_2 \ne \nu}{|\mathbb{E}[(\Tilde{D}^2_{\nu k_1}-\frac{1}{n})(\Tilde{D}^2_{\nu k_2}-\frac{1}{n})]|} \leq \frac{M^2(2c_0(K))^4}{n^3}+\frac{M\alpha^2(1)}{n^2} \leq \frac{2}{n} \cdot \frac{M}{n} \cdot 2^{k_0}\]
since for \(k_1 \ne k_2,\) 
\[\mathbb{E}[\Tilde{D}^2_{\nu k_1}\Tilde{D}^2_{\nu k_2}]=\frac{1}{n^4} \sum_{1 \leq i,j \leq n}{\mathbb{E}[z_{\nu i}^2 z_{\nu j}^2]}= \frac{1}{n^2}+\frac{\mathbb{E}[z_{11}^4]-1}{n^3} \leq \frac{1}{n^2}+\frac{(2c_0(K))^4}{n^3}\]
and for \(k_1=k_2,\) \(\mathbb{E}[\Tilde{D}^4_{\nu k_1}] \leq \frac{\alpha^2(1)}{n^2}\) from (\ref{72}). We proceed with the induction step, which requires the two inequalities in (\ref{37}): the first is equivalent to
\[\overline{a}_{l} \geq \frac{\overline{a}_{l-1}}{\sqrt{nM}}
+\frac{(c_1(K)l)^{3l+4}}{l!},\]
which holds true as \(k_0 \geq 3,\)
\[\frac{\overline{a}_{l-1}}{\sqrt{nM}}=\frac{(2^{k_0}(l-1))^{3l+1}}{\sqrt{nM}} \leq \frac{(2^{k_0}l)^{3l+4}}{2},\]
\[\frac{(c_1(K)l)^{3l+4}}{l!} \leq \frac{(c_1(K)l)^{3l+4}}{2} \leq \frac{(2^{k_0}l)^{3l+4}}{2}.\]
For the second inequality, notice that
\[a_{1,l}=\sum_{k,k_1, \hspace{0.05cm} ... \hspace{0.05cm}, k_l \ne \nu}{|\mathbb{E}[(\Tilde{D}_{\nu k}^2-\frac{1}{n}) \cdot (\Tilde{D}_{\nu k_1}...\Tilde{D}_{k_l \nu}-\frac{1}{n^l})]|} \leq \]
\begin{equation}\label{2031}
    \leq \frac{c_1(K)}{n} \cdot (\frac{M}{n})^{l+1}+\frac{l}{n}a_{1,l-1}+ \alpha(l+2) \cdot (\frac{l}{2}+1)^{l} \cdot \frac{M^{\frac{l}{2}+1}}{n^{\frac{l}{2}+2}}+2l\alpha(l+3)\cdot (\frac{l}{2}+1)^{l} \cdot \frac{M^{\frac{l}{2}}}{n^{\frac{l+3}{2}}}
\end{equation}
since each \((k,k_1, \hspace{0.05cm} ... \hspace{0.05cm}, k_l)\) falls in one of the following four categories:
\vspace{0.2cm}
\(\newline 1. \hspace{0.2cm} k,k_1, \hspace{0.05cm} ... \hspace{0.05cm}, k_l\) are pairwise distinct: by conditioning on \(z_k,z_{k_1},\hspace{0.05cm} ... \hspace{0.05cm}, z_{k_l},\)
\[\mathbb{E}[(\Tilde{D}_{\nu k}^2-\frac{1}{n}) \cdot (\Tilde{D}_{\nu k_1}...\Tilde{D}_{k_l \nu}-\frac{1}{n^l})]=\mathbb{E}[\Tilde{D}_{\nu k}^2\Tilde{D}_{\nu k_1}...\Tilde{D}_{k_l \nu}]-\frac{1}{n}\mathbb{E}[\Tilde{D}_{\nu k_1}...\Tilde{D}_{k_l \nu}]=\]
\[=\frac{1}{n^{l+3}}\sum_{1 \leq i,j \leq n}{\mathbb{E}[z^2_{\nu i}z^2_{\nu j}]} -\frac{1}{n^{l+1}}=\frac{1}{n^{l+2}}(\mathbb{E}[z_{11}^4]-1) \leq \frac{c_1(K)}{n^{l+2}}.\] 
\vspace{0.2cm}
\(\newline 2.\) There exists \(k_i\) such that its multiplicity in \(k,k_1, \hspace{0.05cm} ... \hspace{0.05cm}, k_l\) is one: by conditioning on \(z_{k_i}\) (with \(k_i\) minimal), their overall contribution is at most \(\frac{l}{n}a_{1,l-1}.\) 
\vspace{0.2cm}
\(\newline 3. \hspace{0.2cm} k\) has multiplicity one in \(k,k_1, \hspace{0.05cm} ... \hspace{0.05cm}, k_l\) whereas \(k_1, \hspace{0.05cm} ... \hspace{0.05cm}, k_l\) have multiplicity at least two:
\[\mathbb{E}[\Tilde{D}_{\nu k}^2\Tilde{D}_{\nu k_1}...\Tilde{D}_{k_l \nu}]-\frac{1}{n}\mathbb{E}[\Tilde{D}_{\nu k_1}...\Tilde{D}_{k_l \nu}]=\frac{1}{n}\mathbb{E}[\Tilde{D}_{\nu \nu}\Tilde{D}_{\nu k_1}...\Tilde{D}_{k_l \nu}] \leq \frac{\alpha(l+2)}{n^{\frac{l}{2}+2}}\] 
by conditioning on \(z_k,\)
\[\mathbb{E}[\Tilde{D}_{\nu k}^2\Tilde{D}_{\nu k_1}...\Tilde{D}_{k_l \nu}]=\mathbb{E}[\frac{1}{n^2}||z_\nu||^2\Tilde{D}_{\nu k_1}...\Tilde{D}_{k_l \nu}]\]
in conjunction with Cauchy-Schwarz inequality and (\ref{2029}). There are at most \(M \cdot M^{l/2} \cdot (\frac{l}{2}+1)^l\) such tuples: \(k\) can take at most \(M\) values, and recall the bound (\ref{2027}) for the number of \((k_1, \hspace{0.05cm} ... \hspace{0.05cm}, k_l).\)
\vspace{0.2cm}
\(\newline 4. \hspace{0.2cm} k=k_i\) for some \(i,\) and \(k_1, \hspace{0.05cm} ... \hspace{0.05cm}, k_l\) have multiplicity at least two:
\[|\mathbb{E}[\Tilde{D}_{\nu k}^2\Tilde{D}_{\nu k_1}...\Tilde{D}_{k_l \nu}]-\frac{1}{n}\mathbb{E}[\Tilde{D}_{\nu k_1}...\Tilde{D}_{k_l \nu}]| \leq |\mathbb{E}[\Tilde{D}_{\nu k}^2\Tilde{D}_{\nu k_1}...\Tilde{D}_{k_l \nu}]|+\frac{1}{n}|\mathbb{E}[\Tilde{D}_{\nu k_1}...\Tilde{D}_{k_l \nu}]| \leq \frac{2\alpha(l+3)}{n^{\frac{l+3}{2}}}\]
using (\ref{2029}) and \(\alpha(l+3) \geq \alpha(l+1).\) The number of such tuples is at most \(l \cdot M^{l/2} \cdot (\frac{l}{2}+1)^l.\) 
\vspace{0.2cm}
\par
Then by induction, 
\[a_{1,l} \leq \frac{(l+1)!}{n} \cdot (\frac{M}{n})^{\frac{l+1}{2}} \cdot (2^{k_0}(l+1))^{3l+4}:\]
the base case \(l=1\) was already covered above, and the induction step follows from (\ref{2031}) and 
\[1. \hspace{0.2cm} \frac{c_1(K)}{n} \cdot (\frac{M}{n})^{l+1} \leq \frac{1}{4} \cdot \frac{(l+1)!}{n} \cdot (\frac{M}{n})^{\frac{l+1}{2}} \cdot (2^{k_0}(l+1))^{3l+4},\]
\[2. \hspace{0.2cm} \frac{l}{n}a_{1,l-1} \leq \frac{(l+1)!}{n^2} \cdot (\frac{M}{n})^{\frac{l}{2}} \cdot (2^{k_0}l)^{3l+1} \leq \frac{1}{4} \cdot \frac{(l+1)!}{n} \cdot (\frac{M}{n})^{\frac{l+1}{2}} \cdot (2^{k_0}(l+1))^{3l+4},\]
\[3. \hspace{0.2cm} \alpha(l+2) \cdot (\frac{l}{2}+1)^{l} \cdot \frac{M^{\frac{l}{2}+1}}{n^{\frac{l}{2}+2}} \leq  \alpha(l+2) \cdot (\frac{l}{2}+1)^{l} \cdot \frac{1}{n} \cdot (\frac{M}{n})^{\frac{l+1}{2}} \leq \frac{1}{4} \cdot \frac{(l+1)!}{n} \cdot (\frac{M}{n})^{\frac{l+1}{2}} \cdot (2^{k_0}(l+1))^{3l+4}\]
from (\ref{2028}) and \(l+3 \leq 2(l+1)\) which give
\[\alpha(l+2) \cdot (\frac{l}{2}+1)^{l} \leq (2c_1(K)(l+1))^{2l+6} \cdot  (\frac{l}{2}+1)^{l} \leq \frac{(l+1)!}{4} \cdot (2^{k_0}(l+1))^{3l+4}\] 
because
\[\frac{(l+1)!}{4} \cdot (2^{k_0}(l+1))^{3l+4} \geq \frac{(l+1)^2}{8} \cdot (4c_1(K)(l+1))^{3l+4}=\frac{(l+1)^2}{8} \cdot (4c_1(K)(l+1))^{2l+6} \cdot (4c_1(K)(l+1))^{l-2} \geq\]
\[\geq \frac{(l+1)^2}{8} \cdot (4c_1(K)(l+1))^{2l+6} \cdot (l+1)^{l-2} \geq (2c_1(K)(l+1))^{2l+6} \cdot (l+1)^l,\]
\[4. \hspace{0.2cm} 2l\alpha(l+3)\cdot (\frac{l}{2}+1)^{l} \cdot \frac{M^{\frac{l}{2}}}{n^{\frac{l+3}{2}}} \leq 2l\alpha(l+3)\cdot (\frac{l}{2}+1)^{l} \cdot \frac{1}{n} \cdot (\frac{M}{n})^{\frac{l+1}{2}} \leq \frac{1}{4} \cdot \frac{(l+1)!}{n} \cdot (\frac{M}{n})^{\frac{l+1}{2}} \cdot (2^{k_0}(l+1))^{3l+4},\]
from (\ref{2028}) and \(l+4 \leq 2(l+1)\) yielding
\[2l\alpha(l+3)\cdot (\frac{l}{2}+1)^{l} \leq  (2c_1(K)(l+1))^{2l+8} \cdot (l+2)^{l} \leq (2c_1(K)(l+1))^{3l+8} \leq \frac{(l+1)!}{4} \cdot (2^{k_0}(l+1))^{3l+4}\]
inasmuch as
\[\frac{(l+1)!}{4} \cdot (2^{k_0}(l+1))^{3l+4} \geq \frac{(l+1)^2}{8} \cdot (2c_1(K)(l+1))^{3l+4} \cdot (2c_1(K))^{3l+4} \geq\]
\[\geq (2c_1(K)(l+1))^{3l+4} \cdot (l+1)^5 \cdot (c_1(K))^{3l+4} \geq (2c_1(K)(l+1))^{3l+8}.\]

\subsection{Second Sum: Stages I, II, and III}\label{4.8}

In this subsection, we justify (\ref{28}): the argument is similar to the one for (\ref{27}) presented in \ref{4.4}-\ref{4.7}, and thus we skip some of the details.
\par
First, we truncate first the series underlying \(\Sigma_{0,k}\) for \(k \ne \nu:\)
\[\Tilde{\Sigma}_0=\sum_{k \ne \nu}{(\sum_{0 \leq j \leq s}{((-M_\nu)^j\mathcal{R}_\nu \mathcal{D}_\nu e_\nu)_k}+\sum_{j>s }{((-M_\nu)^j\mathcal{R}_\nu \mathcal{D}_\nu e_\nu)_k})^2}.\]
Recall that \(||M_\nu|| \xrightarrow[]{a.s.} 0,\) and note that
\[\sqrt{\sum_{k \ne \nu}{(\sum_{j>s}{((-M_\nu)^j\mathcal{R}_\nu \mathcal{D}_\nu e_\nu)_k})^2}}=||\sum_{j>s}{(-M_\nu)^j\mathcal{R}_\nu \mathcal{D}_\nu e_\nu}|| \leq \sum_{j>s}{||M_\nu||^j \cdot \beta_\nu} \leq 2||M_\nu||^{s+1} \cdot \beta_\nu,\]
which together with (\ref{17}) and Lemma~\ref{lemma2} gives almost surely for \(k \ne \nu,\)
\[\sqrt{\sum_{k \ne \nu}{(\sum_{0 \leq j \leq s}{((-M_\nu)^j\mathcal{R}_\nu \mathcal{D}_\nu e_\nu)_k})^2}} \leq \frac{||a_\nu-e_\nu||}{|1-\frac{||a_\nu-e_\nu||^2}{2}|}+\sqrt{\sum_{k \ne \nu}{(\sum_{j>s}{((-M_\nu)^j\mathcal{R}_\nu \mathcal{D}_\nu e_\nu)_k})^2}} \leq \]
\begin{equation}\label{1012}
    \leq 2||a_\nu-e_\nu||+2\beta_\nu \leq 6\beta_\nu.
\end{equation}
Cauchy-Schwarz inequality then yields
\[|\Tilde{\Sigma}_0-\sum_{k \ne \nu}{(\sum_{0 \leq j \leq s}{((-M_\nu)^j\mathcal{R}_\nu \mathcal{D}_\nu e_\nu)_k})^2}| \leq 
4||M_\nu||^{2s+2} \cdot \beta_\nu^2+24\cdot ||M_\nu||^{s+1} \cdot \beta_\nu^2=o_p(\frac{1}{\sqrt{n}})\]
because \(\sqrt{n} \cdot ||M_\nu||^{s+1} \cdot \beta_\nu^2 \xrightarrow[]{p} 0\) from (\ref{25}) and (\ref{110}).
\par
Second, we show the first type terms in
\[\sum_{k \ne \nu}{(\sum_{0 \leq j \leq s}{((-M_\nu)^j\mathcal{R}_\nu \mathcal{D}_\nu e_\nu)_k})^2}\]
are negligible. Recall that
\[-M_\nu=A_\nu+B_\nu, \hspace{0.2cm} A_\nu=-\mathcal{R}_\nu \Lambda^{1/2} \Tilde{T}\Lambda^{1/2}, \hspace{0.2cm} B_\nu=-\mathcal{R}_\nu \Lambda^{1/2} \Tilde{D}\Lambda^{1/2}+(\frac{\hat{l}_\nu}{l_\nu}-1)\Tilde{R},\]
\[\sum_{0 \leq j \leq s}{(A_\nu+B_\nu)^j\mathcal{R}_\nu \mathcal{D}_\nu e_\nu}-\sum_{0 \leq j \leq s}{B_\nu^j\mathcal{R}_\nu \mathcal{D}_\nu e_\nu}=\sum_{1 \leq j \leq s, 0 \leq m \leq j-1}{(A_\nu+B_\nu)^{j-1-m}A_\nu B_\nu^m\mathcal{R}_\nu \mathcal{D}_\nu e_\nu},\]
\[||\sum_{0 \leq j \leq s}{B_\nu^j\mathcal{R}_\nu \mathcal{D}_\nu e_\nu}|| \leq \sum_{0 \leq j \leq s}{||B_\nu||^j \cdot \beta_\nu} \leq 2\beta_\nu,\]
\[||\sum_{1 \leq j \leq s, 0 \leq m \leq j-1}{(A_\nu+B_\nu)^{j-1-m}A_\nu B_\nu^m\mathcal{R}_\nu \mathcal{D}_\nu e_\nu}|| \leq\]
\[\leq \sum_{1 \leq j \leq s, 0 \leq m \leq j-1}{(||A_\nu||+||B_\nu||)^{j-1-m} \cdot ||A_\nu|| \cdot || B_\nu||^m \cdot \beta_\nu} \leq ||A_\nu|| \cdot \beta_\nu \sum_{j \geq 1}{\frac{j}{2^{j-1}}}=4||A_\nu|| \cdot \beta_\nu,\]
as \(||B_\nu||, ||A_\nu|| \xrightarrow[]{a.s.} 0,\) from which the desired claim ensues:
\[|\sum_{k \ne \nu}{(\sum_{0 \leq j \leq s}{(A_\nu+B_\nu)^j\mathcal{R}_\nu \mathcal{D}_\nu e_\nu}})_k^2-\sum_{k \ne \nu}{(\sum_{0 \leq j \leq s}{B_\nu^j\mathcal{R}_\nu \mathcal{D}_\nu e_\nu})_k^2}| \leq 16\beta_\nu^2||A_\nu||^2+16\beta_\nu^2||A_\nu||=o_p(\frac{1}{\sqrt{n}})\]
since \(\sqrt{n} \cdot \beta_\nu^2||A_\nu|| \xrightarrow[]{p} 0\) from (\ref{25}) for \(\delta_n=\frac{l_\nu}{M/\sqrt{n}} \to \infty,\) and \(||A_\nu|| \leq \frac{c(K,\gamma,\epsilon_0)}{l_\nu}\) almost surely. 
\par
Third, we turn to the remaining contribution in (\ref{28})
\begin{equation}\label{1004}
    \sum_{0 \leq j_1,j_2 \leq s, k \ne \nu}{(M(j_1))_k \cdot (M(j_2))_k},
\end{equation}
where 
\[M(j)=(-\mathcal{R}_\nu \Lambda^{1/2} \Tilde{D}\Lambda^{1/2}+(\frac{\hat{l}_\nu}{l_\nu}-1)\Tilde{R})^{j_1}\mathcal{R}_\nu \Lambda^{1/2}\Tilde{D}\Lambda^{1/2} e_\nu,\]
because
\[|\sum_{k \ne \nu}{(\sum_{0 \leq j \leq s}{B_\nu^j\mathcal{R}_\nu \mathcal{D}_\nu e_\nu})_k^2}-\sum_{k \ne \nu}{(\sum_{0 \leq j \leq s}{B_\nu^j\mathcal{R}_\nu \Lambda^{1/2}\Tilde{D}\Lambda^{1/2} e_\nu})_k^2}| \leq ||\sum_{0 \leq j \leq s}{B_\nu^j\mathcal{R}_\nu \Lambda^{1/2}\Tilde{T}\Lambda^{1/2} e_\nu}||^2+\]
\[+2 \cdot ||\sum_{0 \leq j \leq s}{B_\nu^j\mathcal{R}_\nu \Lambda^{1/2}\Tilde{D}\Lambda^{1/2} e_\nu}|| \cdot ||\sum_{0 \leq j \leq s}{B_\nu^j\mathcal{R}_\nu \Lambda^{1/2}\Tilde{T}\Lambda^{1/2} e_\nu}|| \leq c(\epsilon_0) \cdot (||\Tilde{T}e_\nu||^2+2||\Tilde{T}e_\nu|| \cdot ||\Tilde{D}e_\nu||)=o_p(\frac{1}{\sqrt{n}})\]
since (\ref{909}) gives \(\sqrt{n} \cdot ||\Tilde{T}e_\nu||^2 \xrightarrow[]{p} 0\) as for \(t>0,\) \(\frac{t\sqrt{n}}{M} \cdot l_\nu^2 \to \infty,\) while (\ref{1001}) and (\ref{1002}) imply 
\[\sqrt{n} \cdot ||\Tilde{T}e_\nu|| \cdot ||\Tilde{D}e_\nu|| \xrightarrow[]{p} 0.\]
As before, we begin with the terms in (\ref{1004}) that have no factor of \(\frac{\hat{l}_\nu}{l_\nu}-1:\)
\[\sum_{0 \leq j_1,j_2 \leq s, k \ne \nu}{((-\mathcal{R}_\nu \Lambda^{1/2} \Tilde{D}\Lambda^{1/2})^{j_1+1})_{k\nu} \cdot ((-\mathcal{R}_\nu \Lambda^{1/2} \Tilde{D}\Lambda^{1/2})^{j_2+1})_{k\nu}}=\]
\[=\sum_{0 \leq j_1,j_2 \leq s}{\frac{l_kl_\nu}{(l_\nu-l_k)^2} \cdot \frac{l_{k_1}...l_{k_{j_1+j_2}}}{(l_\nu-l_{k_1})...(l_\nu-l_{k_{j_1+j_2}})}\Tilde{D}_{kk_1}...\Tilde{D}_{k_{j_1} \nu} \Tilde{D}_{kk_{j_1+1}}...\Tilde{D}_{k_{j_1+j_2}\nu}}=\]
\begin{equation}\label{1003}
    =\sum_{1 \leq j \leq 2s+1}{(\sum_{1 \leq m \leq \min{(s+1,j)}}{\frac{l_\nu}{l_\nu-l_{k_m}}}) \cdot \frac{l_{k_1}...l_{k_{j}}}{(l_\nu-l_{k_1})...(l_\nu-l_{k_{j}})}\Tilde{D}_{k_1\nu}\Tilde{D}_{k_1 k_2}...\Tilde{D}_{k_{j}\nu}},
\end{equation}
where the last equality comes from a change of summation: instead of \((j_1,j_2,k,k_1,\hspace{0.05cm} ... \hspace{0.05cm}, k_{j_1+j_2})\) with \(k, k_i \ne \nu, 0 \leq j_1,j_2 \leq s\) use \((j_1+j_2+1,k_{j_1}, \hspace{0.05cm} ... \hspace{0.05cm}, k_1, k, k_{j_1+1}, \hspace{0.05cm} ... \hspace{0.05cm}, k_{j_1+j_2});\) because
\[|(\sum_{1 \leq m \leq \min{(s+1,j)}}{\frac{l_\nu}{l_\nu-l_{k_m}}}) \cdot \frac{l_{k_1}...l_{k_{j}}}{(l_\nu-l_{k_1})...(l_\nu-l_{k_{j}})}| \leq (s+1) \cdot c(\epsilon_0)^{j+1},\]
an analogous argument to the one in subsection \(6.6\) yields that the sum in (\ref{1003}) is
\[\sum_{1 \leq j \leq 2s+1}{\frac{1}{n^{j}}(\sum_{1 \leq m \leq \min{(s+1,j)}}{\frac{l_\nu}{l_\nu-l_{k_m}}}) \cdot \frac{l_{k_1}...l_{k_{j}}}{(l_\nu-l_{k_1})...(l_\nu-l_{k_{j}})}}+o_p(\frac{1}{\sqrt{n}}),\]
with the first term being \(b_0\) (employing (\ref{1003}) again) and
\[|\sum_{1 \leq j \leq 2s+1}{\frac{1}{n^{j}}(\sum_{1 \leq m \leq \min{(s+1,j)}}{\frac{l_\nu}{l_\nu-l_{k_m}}}) \cdot \frac{l_{k_1}...l_{k_{j}}}{(l_\nu-l_{k_1})...(l_\nu-l_{k_{j}})}}| \leq \]
\[\leq \sum_{1 \leq j \leq 2s+1}{\frac{1}{n^{j}}(\sum_{1 \leq m \leq j}{\frac{l_\nu}{|l_\nu-l_{k_m}|}}) \cdot \frac{l_{k_1}...l_{k_{j}}}{|l_\nu-l_{k_1}| \cdot ... \cdot |l_\nu-l_{k_{j}}|}} \leq \frac{1}{n}\sum_{k \ne \nu}{\frac{l_kl_\nu}{(l_k-l_\nu)^2}} \sum_{j \geq 0}{(\frac{1}{n}\sum_{k \ne \nu}{\frac{l_k}{|l_\nu-l_k|}})^j} \leq\]
\[\leq \frac{2}{n}\sum_{k \ne \nu}{\frac{l_kl_\nu}{(l_k-l_\nu)^2}} \leq c(\epsilon_0) \cdot \frac{M}{n}.\]
\par
We consider next the terms in (\ref{1004}) with a factor of exactly \((\frac{\hat{l}_\nu}{l_\nu}-1)^m\) for \(1 \leq m \leq 2s:\)
\[\sum_{0 \leq j_1,j_2 \leq s}{\frac{l_kl_\nu}{(l_k-l_\nu)^2} \cdot \frac{l_{k_1}...l_{k_{j_1+j_2-m}}}{(l_\nu-l_{k_1})...(l_\nu-l_{k_{j_1+j_2-m_1-m_2}})} \cdot (\frac{l_\nu}{l_{k_{j_1}}-l_\nu})^{j_2}...(\frac{l_\nu}{l_{k_{j_{2l_1+2l_2-1}}}-l_\nu})^{j_{2l_1+2l_2}}} \cdot\]
\[\cdot \Tilde{D}_{kk_1}...\Tilde{D}_{k_{j_1-m_1}\nu}\Tilde{D}_{kk_{j_1-m_1+1}}...\Tilde{D}_{k_{j_1+j_2-m}\nu},\]
with the constraints
\[m_1+m_2=m, m_1 \leq j_1, m_2 \leq j_2, 1 \leq l_1 \leq \min{(m_1+1,j_1-m_1+1)},1 \leq l_2 \leq \min{(m_2+1,j_2-m_2+1)},\]
\[j_{2}+...+j_{2l_1}=m_1,j_{2l_1+2}+...+j_{2l_1+2l_2}=m_2, j_{2},\hspace{0.05cm} ... \hspace{0.05cm}, j_{2l_1-2},j_{2l_1+2},\hspace{0.05cm} ...\hspace{0.05cm}, j_{2l_1+2l_2-2}>0, j_{2l_1},j_{2l_1+2l_2} \geq 0.\]
The analysis of these terms is completed arguing in the same vein as in subsection \(6.6\) for the case \(m>0:\)  using a change of summation the \(m^{th}\) contribution is
\[\sum_{0 \leq j_1,j_2 \leq s}{\frac{1}{n^{j_1+j_2-m+1}} \cdot \frac{l_kl_\nu}{(l_k-l_\nu)^2} \cdot \frac{l_{k_1}...l_{k_{j_1+j_2-m}}}{(l_\nu-l_{k_1})...(l_\nu-l_{k_{j_1+j_2-m}})} \cdot (\frac{l_\nu}{l_{k_{j_1}}-l_\nu})^{j_2}...(\frac{l_\nu}{l_{k_{j_{2l_1+2l_2-1}}}-l_\nu})^{j_{2l_1+2l_2}}}+S''(m),\]
the first term being \(b_m\) from (\ref{1003}) with \(|b_m| \leq c(\epsilon_0)^{m+1} \cdot \frac{M}{n},\) and 
\[\sum_{1 \leq m \leq 2s}{\mathbb{E}[n(S''(m))^2]}=o(1).\]

\subsection{Third Sum: Stages I, II, and III}\label{4.9}

In this subsection, we justify the last identity left to complete the proof of Theorem~\ref{th2}, (\ref{29}): because the argument is in the same vein as the one employed for (\ref{27}) in \ref{4.4}-\ref{4.7}, we omit some of the details:
\[\Tilde{\Sigma}_2=\sum_{k \ne \nu}{\frac{l_k}{l_\nu}\Sigma_{0,k}^2}+
\frac{1}{l_\nu}\sum_{k_1,k_2 \ne \nu}{(\Tilde{T}_{k_1k_2}+\Tilde{D}_{k_1k_2})\sqrt{l_{k_1}}\Sigma_{0,k_1}\sqrt{l_{k_2}}\Sigma_{0,k_2}}.\]
\par
For the first component, the analysis for \(\Tilde{\Sigma}_0\) yields a decomposition for it as well, giving the first sum in (\ref{132}), and hence we consider the second component. Before employing the three-step proof previously used, note that
\[|\frac{1}{l_\nu}\sum_{k_1,k_2 \ne \nu}{\Tilde{T}_{k_1k_2}\sqrt{l_{k_1}}\Sigma_{0,k_1}\sqrt{l_{k_2}}\Sigma_{0,k_2}}| \leq \frac{1}{l_\nu} \cdot ||\Tilde{T}|| \cdot \sum_{k \ne \nu}{l_k\Sigma_{0,k}^2}= \frac{||\Lambda^{1/2}(a_\nu-e_\nu)||^2}{l_\nu(\frac{||a_\nu-e_\nu||^2}{2}-1)^2} \cdot ||\Tilde{T}||=o_p(\frac{1}{\sqrt{n}}),\]
using (\ref{33}) and \(||\Tilde{T}|| \leq c(\epsilon_0) \cdot \frac{1}{l_\nu}\) almost surely, yielding
\[\frac{1}{l_\nu}\sum_{k_1,k_2 \ne \nu}{(\Tilde{T}_{k_1k_2}+\Tilde{D}_{k_1k_2})\sqrt{l_{k_1}}\Sigma_{0,k_1}\sqrt{l_{k_2}}\Sigma_{0,k_2}}=\frac{1}{l_\nu}\sum_{k_1,k_2 \ne \nu}{\Tilde{D}_{k_1k_2}\sqrt{l_{k_1}}\Sigma_{0,k_1}\sqrt{l_{k_2}}\Sigma_{0,k_2}}+o_p(\frac{1}{\sqrt{n}}).\]
\par
Next, we proceed with the first stage, truncation of the series at \(s:\)
\[|\frac{1}{l_\nu}\sum_{k_1,k_2 \ne \nu}{\Tilde{D}_{k_1k_2}\sqrt{l_{k_1}}\Sigma_{0,k_1}\sqrt{l_{k_2}}\Sigma_{0,k_2}}-\frac{1}{l_\nu}\sum_{k_1,k_2 \ne \nu}{\Tilde{D}_{k_1k_2}\sqrt{l_{k_1}}\Sigma_{s,k_1}\sqrt{l_{k_2}}\Sigma_{s,k_2}}| \leq \]
\begin{equation}\label{2099}
    \leq 2||\Tilde{D}|| \cdot \frac{1}{l_\nu} \cdot \sqrt{\sum_{k \ne \nu}{l_k\Sigma_{s,k}^2}} \cdot \sqrt{\sum_{k \ne \nu}{l_k(\Sigma_{0,k}-\Sigma_{s,k})^2}}+||\Tilde{D}|| \cdot \frac{1}{l_\nu} \cdot \sum_{k \ne \nu}{l_k(\Sigma_{0,k}-\Sigma_{s,k})^2}=o_p(\frac{1}{\sqrt{n}}),
\end{equation}
where \(\Sigma_{s,k}=\sum_{0 \leq j \leq s}{((-M_\nu)^j\mathcal{R}_\nu \mathcal{D}_\nu e_\nu)_k},\) using that almost surely
\[\frac{1}{\sqrt{l_\nu}} \cdot \sqrt{\sum_{k \ne \nu}{l_k\Sigma^2_{s,k}}} \leq 2 \cdot \frac{||\Lambda^{1/2}(a_\nu-e_\nu)||}{\sqrt{l_\nu}}+2 \cdot \frac{||\Lambda^{1/2}\mathcal{R}_\nu \mathcal{D}_\nu e_\nu||}{\sqrt{l_\nu}} \leq 4 \cdot \frac{||\Lambda^{1/2}\mathcal{R}_\nu \mathcal{D}_\nu e_\nu||}{\sqrt{l_\nu}}+8\beta_\nu^2,\]
\[\sqrt{\sum_{k \ne \nu}{(\Sigma_{0,k}-\Sigma_{s,k})^2}} \leq 2||M_\nu||^{s+1} \cdot \frac{||\Lambda^{1/2}(a_\nu-e_\nu)||}{\sqrt{l_\nu}} \leq ||M_\nu||^{s+1}(4 \cdot \frac{||\Lambda^{1/2}\mathcal{R}_\nu \mathcal{D}_\nu e_\nu||}{\sqrt{l_\nu}}+8\beta_\nu^2)\]
where the first inequalities in each chain are obtained arguing in the same vein as for (\ref{1012}), the second arise from (\ref{1010}), entailing with (\ref{25}), (\ref{1002}), and (\ref{110}) the bound in (\ref{2099}).  
\par
We continue with the second step, removing the first type terms in the truncated sums:
\[|\frac{1}{l_\nu}\sum_{k_1,k_2 \ne \nu}{\Tilde{D}_{k_1k_2}\sqrt{l_{k_1}}\Sigma_{s,k_1}\sqrt{l_{k_2}}\Sigma_{s,k_2}}-\frac{1}{l_\nu}\sum_{k_1,k_2 \ne \nu}{\Tilde{D}_{k_1k_2}\sqrt{l_{k_1}} \cdot \overline{\Sigma}_{s,k_1}\sqrt{l_{k_2}} \cdot \overline{\Sigma}_{s,k_2}}| \leq \]
\[\leq 2||\Tilde{D}|| \cdot \frac{\sqrt{\sum_{k \ne \nu}{l_k\overline{\Sigma}_{s,k}^2}} \cdot \sqrt{\sum_{k \ne \nu}{l_k(\Sigma_{s,k}-\overline{\Sigma}_{s,k})^2}}}{l_\nu}+||\Tilde{D}|| \cdot \frac{\sum_{k \ne \nu}{l_k(\Sigma_{s,k}-\overline{\Sigma}_{s,k})^2}}{l_\nu}=o_p(\frac{1}{\sqrt{n}})\]
where \(\overline{\Sigma}_{s,k}=\sum_{0 \leq j \leq s}{(B_\nu^j\mathcal{R}_\nu \mathcal{D}_\nu e_\nu)_k},\) since
\[\frac{1}{\sqrt{l_\nu}} \cdot \sqrt{\sum_{k \ne \nu}{l_k\overline{\Sigma}_{s,k}^2}} \leq 2 \cdot \frac{||\Lambda^{1/2}\mathcal{R}_\nu \mathcal{D}_\nu e_\nu||}{\sqrt{l_\nu}},\]
\[\frac{1}{\sqrt{l_\nu}} \cdot \sqrt{\sum_{k \ne \nu}{l_k(\Sigma_{s,k}-\overline{\Sigma}_{s,k})^2}} \leq 4 \cdot ||A_\nu||\cdot \frac{||\Lambda^{1/2}\mathcal{R}_\nu \mathcal{D}_\nu e_\nu||}{\sqrt{l_\nu}}\]
from
\[||\sum_{0 \leq j \leq s}{\Lambda^{1/2}B_\nu^j\mathcal{R}_\nu \mathcal{D}_\nu e_\nu}|| \leq \sum_{0 \leq j \leq s}{||B_\nu||^j \cdot ||\Lambda^{1/2}\mathcal{R}_\nu \mathcal{D}_\nu e_\nu||},\]
\[\sqrt{\sum_{k \ne \nu}{l_k(\Sigma_{s,k}-\overline{\Sigma}_{s,k})^2}}= ||\sum_{1 \leq j \leq s, 0 \leq m \leq j-1}{\Lambda^{1/2}(A_\nu+B_\nu)^{j-1-m}A_\nu B_\nu^m\mathcal{R}_\nu \mathcal{D}_\nu e_\nu}|| \leq\]
\[\leq \sum_{1 \leq j \leq s, 0 \leq m \leq j-1}{(||A_\nu||+||B_\nu||)^{j-1-m}||A_\nu|| \cdot || B_\nu||^m \cdot ||\Lambda^{1/2}\mathcal{R}_\nu \mathcal{D}_\nu e_\nu||} \leq\]
\[\leq ||A_\nu|| \cdot ||\Lambda^{1/2}\mathcal{R}_\nu \mathcal{D}_\nu e_\nu|| \sum_{1 \leq j}{\frac{j}{2^{j-1}}}=4||A_\nu||\cdot ||\Lambda^{1/2}\mathcal{R}_\nu \mathcal{D}_\nu e_\nu||,\]
and (\ref{1002}), \(||A_\nu|| \leq \frac{c(K,\gamma,\epsilon_0)}{l_\nu}\) almost surely. Similarly, we can replace \(\mathcal{R}_\nu \mathcal{D}_\nu\) by \(\mathcal{R}_\nu \Lambda^{1/2}\Tilde{D}\Lambda^{1/2}\) in \(\overline{\Sigma}_{s,k},\) and we are left with
\begin{equation}\label{1202}
    \frac{1}{l_\nu}\sum_{k_1,k_2 \ne \nu}{\Tilde{D}_{k_1k_2}\sqrt{l_{k_1}}\sum_{0 \leq j \leq s}{((C_\nu+(\frac{\hat{l}_\nu}{l_\nu}-1)\Tilde{R})^jC_\nu)_{k_1\nu}}\sqrt{l_{k_2}}\sum_{0 \leq j \leq s}{(C_\nu+(\frac{\hat{l}_\nu}{l_\nu}-1)\Tilde{R})^jC_\nu)_{k_2\nu}}}
\end{equation}
for \(C_\nu=-\mathcal{R}_\nu \Lambda^{1/2}\Tilde{D}\Lambda^{1/2}.\)
\par
Consider first the terms with no factor of \(\frac{\hat{l}_\nu}{l_\nu}-1:\)
\[\frac{1}{l_\nu}\sum_{k_1,k_2 \ne \nu}{\Tilde{D}_{k_1k_2}\sqrt{l_{k_1}}\sum_{0 \leq j \leq s}{(C_\nu^{j+1})_{k_1\nu}}\sqrt{l_{k_2}}\sum_{0 \leq j \leq s}{(C_\nu^{j+1})_{k_2\nu}}}=\]
\[=\sum_{0 \leq j_1,j_2 \leq s}{\frac{l_{k_1}l_{k_2}l_{\kappa_{1}}...l_{\kappa_{j_1+j_2}}}{(l_\nu-l_{k_1})(l_\nu-l_{k_2})(l_\nu-l_{\kappa_1})...(l_\nu-l_{\kappa_{j_1+j_2}})}\Tilde{D}_{k_1\kappa_1}...\Tilde{D}_{\kappa_{j_1}\nu} \Tilde{D}_{k_2\kappa_{j_1+1}}...\Tilde{D}_{\kappa_{j_1+j_2}\nu}}=\]
\[=\sum_{0 \leq j_1,j_2 \leq s}{\frac{1}{n^{j_1+j_2+2}} \cdot \frac{l_{k_1}l_{k_2}l_{\kappa_{1}}...l_{\kappa_{j_1+j_2}}}{(l_\nu-l_{k_1})(l_\nu-l_{k_2})(l_\nu-l_{\kappa_1})...(l_\nu-l_{\kappa_{j_1+j_2}})}}+o_p(\frac{1}{\sqrt{n}}),\]
where again a change of summation is used, from \((k_1,k_2,\kappa_1, \hspace{0.05cm} ... \hspace{0.05cm}, \kappa_{j_1+j_2})\) to \((\kappa_{j_1}, \hspace{0.05cm} ...\hspace{0.05cm}, \kappa_1, k_1, k_2, \kappa_{j_1+1}, \hspace{0.05cm}... \hspace{0.05cm}, \kappa_{j_1+j_2}),\) and the analysis of the error for \(\Tilde{\Sigma}_0,\) with the first term coming from the second component in (\ref{132}) and
\[\sum_{0 \leq j_1,j_2 \leq s}{\frac{1}{n^{j_1+j_2+2}} \cdot \frac{l_{k_1}l_{k_2}l_{\kappa_{1}}...l_{\kappa_{j_1+j_2}}}{(l_\nu-l_{k_1})(l_\nu-l_{k_2})(l_\nu-l_{\kappa_1})...(l_\nu-l_{\kappa_{j_1+j_2}})}}=\]
\[=\sum_{2 \leq j \leq 2s+2}{\min{(j-1,s+1) \cdot (\frac{1}{n}\sum_{k \ne \nu}{\frac{l_k}{l_k-l_\nu}})^j}}=O(c(\epsilon_0) \cdot \frac{M}{n}),\]
as each \(j=(j_1+j_2+2)\)-tuple appears \(\min{(j-1,s+1)}\) times (there is a bijection between these appearances and the solutions \((j_1,j_2)\) of \(j=j_1+j_2+2, 0 \leq j_1, j_2 \leq s\)). 
\par
Lastly, consider the terms containing a factor of \((\frac{\hat{l}_\nu}{l_\nu}-1)^m\) for \(1 \leq m \leq 2s:\) 
\[\sum{\frac{l_{k_1}l_{k_2}l_{\kappa_{1}}...l_{\kappa_{j_1+j_2-m_1-m_2}}}{(l_\nu-l_{k_1})(l_\nu-l_{k_2})(l_\nu-l_{\kappa_1})...(l_\nu-l_{\kappa_{j_1+j_2-m_1-m_2}})} \cdot (\frac{l_\nu}{l_{\kappa_{\overline{j}_1}}})^{\overline{j}_2} ...(\frac{l_\nu}{l_{\kappa_{2\overline{j}_{2l_1+2l_2-1}}}})^{\overline{j}_{2l_1+2l_2}} \cdot }\]
\begin{equation}\label{1200}
    \cdot \Tilde{D}_{k_1\kappa_1}...\Tilde{D}_{\kappa_{j_1-m_1}\nu} \Tilde{D}_{k_2\kappa_{j_1-m_1+1}}...\Tilde{D}_{\kappa_{j_1-m_1+j_2-m_2}\nu},
\end{equation}
where the summation is over \((j_1,j_2,m_1,m_2, k_1, k_2, \kappa_{1}, \hspace{0.05cm} ...\hspace{0.05cm}, \kappa_{j_1+j_2-m_1-m_2},\overline{j}_1, \hspace{0.05cm} ...\hspace{0.05cm}, \overline{j}_{2l_1+2l_2})\) with 
\[k_1, k_2 \ne \nu, \kappa_i \ne \nu, m=m_1+m_2, \overline{j}_1+...+\overline{j}_{2l_1-1}=j_1-m_1,\overline{j}_2+...+\overline{j}_{2l_1}=m_1,\]
\[\overline{j}_{2l_1+1}+...+\overline{j}_{2l_1+2l_2-1}=j_2-m_2, \overline{j}_{2l_1+2}+...+\overline{j}_{2l_1+2l_2}=m_2,\]
\[\overline{j}_1 \geq 0, \overline{j}_2,\hspace{0.05cm} ...\hspace{0.05cm}, \overline{j}_{2l_1-1}>0, \overline{j}_{2l_1} \geq 0, \overline{j}_{2l_1+1} \geq 0, \overline{j}_{2l_1+2}, \hspace{0.05cm} ...\hspace{0.05cm}, \overline{j}_{2l_1+2l_2-1}>0, \overline{j}_{2l_1+2l_2} \geq 0.\]
Proceeding as in the case \(m>0\) in subsection \ref{4.6}, it follows from (\ref{1202}) that (\ref{1200}) is
\[\sum{\frac{1}{n^{j_1+j_2-m+2}} \cdot \frac{l_{k_1}l_{k_2}l_{\kappa_{1}}...l_{\kappa_{j_1+j_2-m_1-m_2}}}{(l_\nu-l_{k_1})(l_\nu-l_{k_2})(l_\nu-l_{\kappa_1})...(l_\nu-l_{\kappa_{j_1+j_2-m_1-m_2}})} \cdot (\frac{l_\nu}{l_{\kappa_{\overline{j}_1}}})^{\overline{j}_2} ...(\frac{l_\nu}{l_{\kappa_{2\overline{j}_{2l_1+2l_2-1}}}})^{\overline{j}_{2l_1+2l_2}}}+\Tilde{S}(m),\]
with the first term \(c'_m\) coming from the second sum in (\ref{132}) and
\[|c'_m| \leq c(\epsilon_0)^{m+1} \cdot \frac{M}{n}, \hspace{0.2cm} \sum_{1 \leq m \leq 2s}{\mathbb{E}[n(\Tilde{S}(m))^2]}=o(1).\]

\subsection{Proof of Theorem~\ref{th4}}\label{4.10}

Recall (\ref{24}) and (\ref{23}): \(x=x_{n,\nu,(l_i)_{1 \leq i \leq M}} \in \mathbb{R}, |x| \leq 2c(\epsilon_0) \cdot \frac{M}{n}\) satisfies
\[x=\overline{O}+\sum_{1 \leq j \leq s}{x^j\overline{O}_{j}},\]
where
\[-\Tilde{\Sigma}_3-2(1-\Tilde{\Sigma}_3)\Tilde{\Sigma}_1+(1-\Tilde{\Sigma}_3)\Tilde{\Sigma}_2=\overline{O}+\sum_{1 \leq j \leq s}{(\frac{\hat{l}_\nu}{l_\nu}-1)^j\overline{O}_{j}}+o_p(\frac{1}{\sqrt{n}}).\]
(\ref{26}) with (\ref{27}), (\ref{28}), (\ref{29}) then yield
\[\sum_{0 \leq j \leq s}{(\sum_{0 \leq i \leq 2s}{(2a_i+b_i+c_i)z^i})(\sum_{0 \leq i \leq 2s}{b_iz^i})^j}=\overline{O}+\sum_{1 \leq j \leq 2s^2+2s}{z^j\overline{O}_{j}}.\]

\section{Switching between Empirical and Deterministic Centering}\label{sec5}

Subsection~\ref{5.1} presents the proof of Theorem~\ref{th1}, and subsection~\ref{5.2} contains the justification of Theorem~\ref{th3}. Both results build on the CLT given by Theorem~\ref{th2}, whose description is complemented by Theorem~\ref{th4}, and as such their proofs are concerned only with altering its centering.

\subsection{Proof of Theorem~\ref{th1}}\label{5.1}

Theorem~\ref{th4} gives 
\[x=2a_0+b_0+c_0+o(\frac{1}{\sqrt{n}})\]
since \(O(\frac{M^2}{n^2})=o(\frac{1}{\sqrt{n}}),\) and thus
\[a_0=\frac{1}{n}\sum_{k \ne \nu}{\frac{l_k}{l_\nu-l_k}}+o(\frac{1}{\sqrt{n}}), \hspace{0.1cm} b_0=-\frac{1}{n}\sum_{k \ne \nu}{\frac{l_kl_\nu}{(l_\nu-l_k)^2}}+o(\frac{1}{\sqrt{n}}), \hspace{0.1cm} c_0=\frac{1}{n}\sum_{k \ne \nu}{\frac{l_k}{l_\nu-l_k}}+\frac{1}{n}\sum_{k \ne \nu}{\frac{l_k^2}{(l_\nu-l_k)^2}}+o(\frac{1}{\sqrt{n}}),\]
from subsections \(\ref{4.7}, \ref{4.8}, \ref{4.9},\) from which
\[x=2a_0+b_0+c_0+o(\frac{1}{\sqrt{n}})=\frac{1}{n}\sum_{k \ne \nu}{\frac{l_k}{l_\nu-l_k}}+o(\frac{1}{\sqrt{n}}).\]
\par
Then by employing Theorem~\ref{th2}, it suffices to prove that 
\begin{equation}\label{103}
    \frac{1}{\sqrt{n}} \cdot (\sum_{k \ne \nu}{\frac{l_k}{l_\nu-l_k}}-\sum_{k \ne \nu}{\frac{\hat{l}_k}{\hat{l}_\nu-\hat{l}_k}}) \xrightarrow[]{p} 0.
\end{equation}
\par
We show first that almost surely, \(l_k \ne \hat{l}_\nu\) and \(\hat{l}_k \ne \hat{l}_\nu\) for all \(k \ne \nu.\) As \(\frac{\hat{l}_\nu}{l_\nu} \xrightarrow[]{a.s.} 1,\) the first claim follows from the separability assumption. For \(k<\nu,\) the second claim is clear because \(\frac{\hat{l}_k}{l_k} \xrightarrow[]{a.s.} 1\) from Proposition~\ref{prop1}. For \(k>\nu,\) if \(\hat{l}_\nu=\hat{l}_k,\) then 
\begin{equation}\label{5081}
    |\frac{\hat{l}_k}{l_k}-1| \geq \frac{\epsilon_0}{(1+\epsilon_0)(2+\epsilon_0)} \cdot \frac{l_\nu}{l_k}
\end{equation}
or
\begin{equation}\label{5082}
    |\frac{\hat{l}_\nu}{l_\nu}-1| \geq \frac{\epsilon_0}{2+\epsilon_0}:
\end{equation}
from
\[\frac{\hat{l}_k}{l_k}-1=\frac{\hat{l}_\nu
}{l_\nu} \cdot (\frac{l_\nu}{l_k}-1)+(\frac{\hat{l}_\nu
}{l_\nu}-1),\] 
if \(|\frac{\hat{l}_\nu}{l_\nu}-1|<\frac{\epsilon_0}{2+\epsilon_0},\) then 
\[\frac{\hat{l}_k}{l_k}-1 \geq \frac{2}{2+\epsilon_0 }\cdot (\frac{l_\nu}{l_k}-1)-\frac{\epsilon_0}{2+\epsilon_0} \geq \frac{2}{2+\epsilon_0 }\cdot \frac{\epsilon_0}{1+\epsilon_0} \cdot \frac{l_\nu}{l_k}-\frac{\epsilon_0}{2+\epsilon_0} \geq \frac{\epsilon_0}{(1+\epsilon_0)(2+\epsilon_0)} \cdot \frac{l_\nu}{l_k}.\] 
Proposition~\ref{prop1} yields that (\ref{5082}) occurs with probability zero with the same conclusion ensuing for (\ref{5081}) from a union bound, Borel-Cantelli lemma, and the following inequality holding for all \(k>\nu, \newline c_1(\epsilon_0)=\frac{\epsilon_0}{(1+\epsilon_0)(2+\epsilon_0)},\) and \(n \geq n(\epsilon_0),\)
\[\mathbb{P}(|\frac{\hat{l}_k}{l_k}-1| \geq c_1(\epsilon_0) \cdot \frac{l_\nu}{l_k}) \leq 6\exp(-\sqrt{n}):\]
(\ref{12}), (\ref{781}) give for \(t>0\)
\[\mathbb{P}(|\frac{\hat{l}_k}{l_k}-1| \geq t) \leq \mathbb{P}(\lambda_1(\mathcal{B}) \geq \frac{tl_k}{2})+\mathbb{P}(|\frac{\lambda_k(\mathcal{A})}{l_k}-1| \geq \frac{t}{2}),\]
\begin{equation}\label{105}
    \mathbb{P}(|\frac{\hat{l}_k}{l_k}-1| \geq t) \leq \mathbb{P}(\lambda_1(\mathcal{B}) \geq \frac{tl_k}{2})+ \mathbb{P}(|\lambda_k(\frac{1}{n}\sum_{1 \leq i \leq k}{Z_i^TZ_i})-1| \geq \frac{t}{2})+\mathbb{P}(|\lambda_1(\frac{1}{n}\sum_{k \leq i \leq M}{Z_i^TZ_i})-1| \geq \frac{t}{2}),
\end{equation}
which in conjunction with \(\frac{l_\nu}{l_k} \geq 1\) implies
\[\mathbb{P}(|\frac{\hat{l}_k}{l_k}-1| \geq c_1(\epsilon_0) \cdot \frac{l_\nu}{l_k}) \leq \mathbb{P}(\lambda_1(\mathcal{B}) \geq \frac{c_1(\epsilon_0)}{2} \cdot l_\nu)+\mathbb{P}(|\lambda_1(\frac{1}{n}\sum_{1 \leq i \leq k}{Z_i^TZ_i})-1| \geq \frac{c_1(\epsilon_0)}{2})+\]
\[+\mathbb{P}(|\lambda_1(\frac{1}{n}\sum_{k \leq i \leq M}{Z_i^TZ_i})-1| \geq \frac{c_1(\epsilon_0)}{2}) \leq \mathbb{P}(\lambda_1(\mathcal{B}) \geq 2c_{K,\gamma})+4\exp(-\sqrt{n}) \leq 6\exp(-\sqrt{n})\]
because from inequality (\ref{7}) with probability at least \(1-2\exp(-\sqrt{n}),\)
\[(1-cK^2(\sqrt{\frac{k}{n}}+n^{-1/4}))^2 \leq  \lambda_1(\frac{1}{n}\sum_{1 \leq i \leq k}{Z_i^TZ_i}) \leq (1+cK^2(\sqrt{\frac{k}{n}}+n^{-1/4}))^2,\]
from which for \(n \geq n(\epsilon_0),\)
\begin{equation}\label{202}
|\lambda_1(\frac{1}{n}\sum_{1 \leq i \leq k}{Z_i^TZ_i})-1| \leq 3cK^2(\sqrt{\frac{k}{n}}+n^{-1/4}) \leq 3cK^2(\sqrt{\frac{M}{n}}+n^{-1/4}) \leq  \frac{c_1(\epsilon_0)}{4},
\end{equation}
and similarly 
\[\mathbb{P}(|\lambda_1(\frac{1}{n}\sum_{k \leq i \leq M}{Z_i^TZ_i})-1| \geq \frac{c_1(\epsilon_0)}{2}) \leq 2\exp(-\sqrt{n}).\]
\par
Next, to obtain (\ref{103}) (whose left-hand side term is well-defined almost surely from the previous paragraph), we show
\begin{equation}\label{39}
    \frac{1}{\sqrt{n}} \cdot (\sum_{k \ne \nu}{\frac{l_k}{l_k-l_\nu}}-\sum_{k \ne \nu}{\frac{l_k}{l_k-\hat{l}_\nu}}) \xrightarrow[]{p} 0,    
\end{equation}
\begin{equation}\label{40}
    \frac{1}{\sqrt{n}} \cdot (\sum_{k \ne \nu}{\frac{l_k}{l_k-\hat{l}_\nu}}-\sum_{k \ne \nu}{\frac{\hat{l}_k}{\hat{l}_k-\hat{l}_\nu}}) \xrightarrow[]{p} 0.   
\end{equation}
\par
We begin with (\ref{39}):
\[\frac{l_k}{l_k-l_\nu}-\frac{l_k}{l_k-\hat{l}_\nu}=-\frac{l_k(\hat{l}_\nu-l_\nu)}{(l_k-l_\nu)(l_k-\hat{l}_\nu)}=\frac{\frac{l_k}{l_\nu} \cdot (\frac{\hat{l}_\nu}{l_\nu}-1)}{(\frac{l_k}{l_\nu}-1)(\frac{\hat{l}_\nu}{l_\nu}-\frac{l_k}{l_\nu})},\]
from which Proposition~\ref{prop1} and the separability assumption give that almost surely
\[|\frac{l_k}{l_k-l_\nu}-\frac{l_k}{l_k-\hat{l}_\nu}| \leq c(\epsilon_0) \cdot |\frac{\hat{l}_\nu}{l_\nu}-1|.\]
Since
\[\frac{M}{\sqrt{n}} \cdot |\frac{\hat{l}_\nu}{l_\nu}-1| \leq \frac{M}{\sqrt{n}} \cdot \frac{1}{n}tr(\mathcal{M}(\hat{l}_\nu I -\mathcal{M})^{-1})+\frac{M}{\sqrt{n}} \cdot \frac{1}{n}|\sum_{k \ne \nu}{\frac{l_k}{l_k-l_\nu}}|+\]
\[+\frac{M}{n} \cdot \sqrt{n} \cdot |\frac{\hat{l}_\nu}{l_\nu}-1-\frac{1}{n}tr(\mathcal{M}(\hat{l}_\nu I -\mathcal{M})^{-1})-\frac{1}{n}\sum_{k \ne \nu}{\frac{l_k}{l_\nu-l_k}}| \xrightarrow[]{p} 0\]
from
\[\frac{M}{\sqrt{n}} \cdot \frac{1}{n}|tr(\mathcal{M}(\hat{l}_\nu I -\mathcal{M})^{-1})| \leq c(K,\gamma) \cdot \frac{M/\sqrt{n}}{l_\nu} \to 0, \hspace{0.2cm} \frac{M}{\sqrt{n}} \cdot \frac{M}{n}=\frac{M^2}{n^{3/2}} \to 0, \hspace{0.2cm} \frac{M}{n} \to 0,\]
and Theorem~\ref{th2}, (\ref{39}) follows.
\par
We continue with (\ref{40}):
\begin{equation}\label{797}
    \frac{l_k}{l_k-\hat{l}_\nu}-\frac{\hat{l}_k}{\hat{l}_k-\hat{l}_\nu}=\frac{\hat{l}_\nu(\hat{l}_k-l_k)}{(l_k-\hat{l}_\nu)(\hat{l}_k-\hat{l}_\nu)}=(\frac{\hat{l}_k}{l_k}-1) \cdot \frac{\hat{l}_\nu l_k}{(l_k-\hat{l}_\nu)(\hat{l}_k-\hat{l}_\nu)},
\end{equation}
and prove that with probability one
\begin{equation}\label{104}
    |\frac{\hat{l}_\nu l_k}{(l_k-\hat{l}_\nu)(\hat{l}_k-\hat{l}_\nu)}|= |\frac{\frac{\hat{l}_\nu}{l_\nu} \cdot \frac{l_k}{l_\nu} }{(\frac{l_k}{l_\nu}-\frac{\hat{l}_\nu}{l_\nu}) \cdot (\frac{l_k}{l_\nu} \cdot \frac{\hat{l}_k}{l_k}-\frac{\hat{l}_\nu}{l_\nu})}| \leq c(\epsilon_0,K,\gamma) \cdot \min{(\frac{l_k}{l_\nu},\frac{l_\nu}{l_k})}.
\end{equation}
\par
If \(l_k \geq \frac{l_\nu}{4(c_{K,\gamma}+1)},\) Proposition~\ref{prop1} yields \(\frac{\hat{l}_k}{l_k} \xrightarrow[]{a.s.} 1, \frac{\hat{l}_\nu}{l_\nu} \xrightarrow[]{a.s.} 1\) jointly for all such \(k \ne \nu,\) and so (\ref{104}) follows from \[|\frac{\frac{\hat{l}_\nu}{l_\nu}}{(\frac{l_k}{l_\nu}-\frac{\hat{l}_\nu}{l_\nu}) \cdot (\frac{l_k}{l_\nu} \cdot \frac{\hat{l}_k}{l_k}-\frac{\hat{l}_\nu}{l_\nu})}| \leq c(\epsilon_0,K,\gamma), \hspace{0.2cm} |\frac{\frac{\hat{l}_\nu}{l_\nu} \cdot (\frac{l_k}{l_\nu})^2 }{(\frac{l_k}{l_\nu}-\frac{\hat{l}_\nu}{l_\nu}) \cdot (\frac{l_k}{l_\nu} \cdot \frac{\hat{l}_k}{l_k}-\frac{\hat{l}_\nu}{l_\nu})}| \leq c(\epsilon_0,K,\gamma).\] If \(l_k \leq \frac{l_\nu}{4(c_{K,\gamma}+1)},\) then almost surely for all such \(k,\) \(\frac{\hat{l}_k}{l_k} \leq 2c_{K,\gamma}+2\) from (\ref{105}) (\(\frac{(2c_{K,\gamma}+1)l_k}{2} \geq \frac{2c_{K,\gamma}+1}{2} >c_{K,\gamma}\)). Hence with probability one
\[|\frac{\frac{\hat{l}_\nu}{l_\nu} }{(\frac{l_k}{l_\nu}-\frac{\hat{l}_\nu}{l_\nu}) \cdot (\frac{l_k}{l_\nu} \cdot \frac{\hat{l}_k}{l_k}-\frac{\hat{l}_\nu}{l_\nu})}| 
\leq \frac{2}{(\frac{1}{1+\epsilon_0/2}-\frac{1}{1+\epsilon_0}) \cdot \frac{1}{6}},\]
from 
\[\frac{l_k}{l_\nu} \cdot \frac{\hat{l}_k}{l_k}-\frac{\hat{l}_\nu}{l_\nu} \leq \frac{1}{2}-\frac{2}{3}=-\frac{1}{6}, \hspace{0.4cm} \frac{l_k}{l_\nu}-\frac{\hat{l}_\nu}{l_\nu} \leq \frac{1}{1+\epsilon_0}-\frac{1}{1+\epsilon_0/2}\] 
concluding the proof of (\ref{104}).
\par
Now we return to (\ref{40}), and show that
\begin{equation}\label{106}
    \mathbb{P}(|\frac{\hat{l}_k}{l_k}-1| \cdot \min{(\frac{l_k}{l_\nu},\frac{l_\nu}{l_k})} \geq \frac{t\sqrt{n}}{M}) \leq 6\exp(-\sqrt{n}),
\end{equation}
for all \(t>0, n \geq n(t):\) if \(l_k \geq l_\nu,\) then (\ref{105}) gives
\[\mathbb{P}(|\frac{\hat{l}_k}{l_k}-1| \geq \frac{tl_k\sqrt{n}}{Ml_\nu}) \leq \mathbb{P}(|\frac{\hat{l}_k}{l_k}-1| \geq \frac{t\sqrt{n}}{M}) \leq \mathbb{P}(\lambda_1(\mathcal{B}) \geq \frac{tl_k\sqrt{n}}{2M})+\mathbb{P}(|\lambda_1(\frac{1}{n}\sum_{1 \leq i \leq k}{Z_i^TZ_i})-1| \geq \frac{t\sqrt{n}}{2M})+\]
\[+\mathbb{P}(|\lambda_1(\frac{1}{n}\sum_{k \leq i \leq M}{Z_i^TZ_i})-1| \geq \frac{t\sqrt{n}}{2M}) \leq \mathbb{P}(\lambda_1(\mathcal{B}) \geq 2c_{K,\gamma})+4\exp(-\sqrt{n}) \leq 6\exp(-\sqrt{n}),\]
using (\ref{202}) and 
\[(\frac{1}{n^{1/4}}+\sqrt{\frac{M}{n}}) \cdot \frac{M}{\sqrt{n}}=\frac{M}{n^{3/4}}+\frac{M^{3/2}}{n} \to 0;\]
if \(l_k<l_\nu,\) then
\[\mathbb{P}(\frac{l_k}{l_\nu} \cdot |\frac{\hat{l}_k}{l_k}-1| \geq \frac{t\sqrt{n}}{M}) \leq \mathbb{P}(\lambda_1(\mathcal{B}) \geq \frac{tl_\nu \sqrt{n}}{2M})+\mathbb{P}(|\lambda_1(\frac{1}{n}\sum_{1 \leq i \leq k}{Z_i^TZ_i})-1| \geq \frac{t\sqrt{n}}{M})+\]
\[+\mathbb{P}(|\lambda_1(\frac{1}{n}\sum_{k \leq i \leq M}{Z_i^TZ_i})-1| \geq \frac{t\sqrt{n}}{M}) \leq \mathbb{P}(\lambda_1(\mathcal{B}) \geq 2c_{K,\gamma})+4\exp(-\sqrt{n}) \leq 6\exp(-\sqrt{n}).\]
\par
Finally, from (\ref{797}), (\ref{104}), (\ref{106}), with probability one
\[\frac{1}{\sqrt{n}} \cdot |\sum_{k \ne \nu}{\frac{l_k}{l_k-\hat{l}_\nu}}-\sum_{k \ne \nu}{\frac{\hat{l}_k}{\hat{l}_k-\hat{l}_\nu}}| \leq \frac{c(\epsilon_0,K,\gamma)}{\sqrt{n}} \sum_{k \ne \nu}{ |\frac{\hat{l}_k}{l_k}-1| \cdot \min{(\frac{l_k}{l_\nu},\frac{l_\nu}{l_k})}},\]
and for any \(t>0, n \geq n(t)\)
\[ \sum_{k \ne \nu}{\mathbb{P}(|\frac{\hat{l}_k}{l_k}-1| \cdot \min{(\frac{l_k}{l_\nu},\frac{l_\nu}{l_k})} \geq \frac{t\sqrt{n}}{c(\epsilon_0,K,\gamma)M})} \leq 6M\exp(-\sqrt{n}),\]
yielding (\ref{40}).

\subsection{Proof of Theorem~\ref{th3}}\label{5.2}

Denote by
\begin{equation}\label{50}
    \overline{l}_\nu=l_\nu(1+\frac{N-M}{n} \cdot \frac{1}{l_\nu-1}).
\end{equation}
The conclusion can be derived from the following convergence
\begin{equation}\label{38}
    \sqrt{n} \cdot (\frac{1}{l_\nu-1}-\frac{1}{N-M}tr(\mathcal{M}(\overline{l}_\nu I -\mathcal{M})^{-1})) \xrightarrow[]{p} 0:
\end{equation}
let us see first how it leads to the desired result and second why it is correct. 
\par
First, suppose (\ref{38}) holds. Theorem~\ref{th2} gives
\[\sqrt{\frac{n}{\mathbb{E}[(z^{(n)}_{11})^4]-1}} \cdot (\frac{\hat{l}_\nu}{l_\nu}-1-\frac{1}{n}tr(\mathcal{M}(\hat{l}_\nu I -\mathcal{M})^{-1})-x_{n,\nu,(l_i)_{1 \leq i \leq M}}) \Rightarrow N(0,1),\]
and therefore, from Slutsky's lemma, it suffices to prove that
\[\sqrt{n} \cdot (\frac{1}{n}tr(\mathcal{M}(\hat{l}_\nu I -\mathcal{M})^{-1})-\frac{N-M}{n} \cdot \frac{1}{l_\nu-1}) \xrightarrow[]{p} 0,\]
or equivalently
\[\sqrt{n} \cdot (\frac{1}{N-M}tr(\mathcal{M}(\hat{l}_\nu I -\mathcal{M})^{-1})-\frac{1}{l_\nu-1}) \xrightarrow[]{p} 0,\]
which in turn can be rewritten, employing (\ref{38}), as
\begin{equation}\label{727}
    \sqrt{n} \cdot (\frac{1}{N-M}tr(\mathcal{M}(\hat{l}_\nu I -\mathcal{M})^{-1})-\frac{1}{N-M}tr(\mathcal{M}(\overline{l}_\nu I -\mathcal{M})^{-1})) \xrightarrow[]{p} 0.
\end{equation}
Take \(x_n \in \mathbb{R}\) such that
\begin{equation}\label{80}
    l_\nu=\frac{\hat{l}_\nu}{1+\frac{1}{n}tr(\mathcal{M}(\hat{l}_\nu I -\mathcal{M})^{-1})+x_{n,\nu,(l_i)_{1 \leq i \leq M}}+\frac{x_n}{\sqrt{n}}},
\end{equation}
and from Theorem~\ref{th1},
\[\frac{x_n}{\sqrt{\mathbb{E}[z_{11}^4]-1}} \Rightarrow N(0,1).\]
(\ref{38}) yields
\begin{equation}\label{81}
    \frac{\sqrt{n}}{l_\nu} \cdot (\frac{\overline{l}_\nu}{1+\frac{1}{n}tr(\mathcal{M}(\overline{l}_\nu I -\mathcal{M})^{-1})}-l_\nu)= \frac{N-M}{n} \cdot \frac{\sqrt{n}\cdot (\frac{1}{l_\nu-1}-\frac{1}{N-M}tr(\mathcal{M}(\overline{l}_\nu I -\mathcal{M})^{-1}))}{1+\frac{1}{n}tr(\mathcal{M}(\overline{l}_\nu I -\mathcal{M})^{-1})} \xrightarrow[]{p} 0,
\end{equation}
since \(\frac{1}{n}tr(\mathcal{M}(\overline{l}_\nu I -\mathcal{M})^{-1}) \xrightarrow[]{a.s.} 0;\) rewrite (\ref{81}) using (\ref{80}) as
\[\frac{\sqrt{n}}{l_\nu} \cdot (\frac{\overline{l}_\nu}{1+\frac{1}{n}tr(\mathcal{M}(\overline{l}_\nu I -\mathcal{M})^{-1})}-\frac{\hat{l}_\nu}{1+\frac{1}{n}tr(\mathcal{M}(\hat{l}_\nu I -\mathcal{M})^{-1})+x_{n,\nu,(l_i)_{1 \leq i \leq M}}+\frac{x_n}{\sqrt{n}}}) \xrightarrow[]{p} 0,\]
from which 
\[y_n:=\frac{\sqrt{n}}{l_\nu} \cdot (\frac{\overline{l}_\nu}{1+\frac{1}{n}tr(\mathcal{M}(\overline{l}_\nu I -\mathcal{M})^{-1})}-\frac{\hat{l}_\nu}{1+\frac{1}{n}tr(\mathcal{M}(\hat{l}_\nu I -\mathcal{M})^{-1})})+\]
\[+\frac{\hat{l}_\nu}{l_\nu} \cdot \frac{\sqrt{n}x_{n,\nu,(l_i)_{1 \leq i \leq M}}+x_n}{(1+\frac{1}{n}tr(\mathcal{M}(\hat{l}_\nu I -\mathcal{M})^{-1})) \cdot (1+\frac{1}{n}tr(\mathcal{M}(\hat{l}_\nu I -\mathcal{M})^{-1})+x_{n,\nu,(l_i)_{1 \leq i \leq M}}+\frac{x_n}{\sqrt{n}})} \xrightarrow[]{p} 0.\]
Note that for 
\[f(x)=\frac{x}{1+\frac{1}{n}tr(\mathcal{M}(xI -\mathcal{M})^{-1})},\] 
and \(x>c(\gamma,K) \cdot ||\mathcal{M}||,\) \[f'(x)=\frac{1+\frac{1}{n}tr(\mathcal{M}(xI -\mathcal{M})^{-1})+x \cdot \frac{1}{n}tr(\mathcal{M}(xI -\mathcal{M})^{-2})}{(1+\frac{1}{n}tr(\mathcal{M}(xI -\mathcal{M})^{-1}))^2}=\frac{1+\frac{1}{n}tr((2xI-\mathcal{M})\mathcal{M}(xI -\mathcal{M})^{-2})}{(1+\frac{1}{n}tr(\mathcal{M}(xI -\mathcal{M})^{-1}))^2} \geq \frac{1}{2},\]
from which, almost surely, using Proposition~\ref{prop1} and the mean-value theorem,
\begin{equation}\label{51}
    \frac{\sqrt{n} \cdot |\overline{l}_\nu-\hat{l}_\nu|}{l_\nu^2} \leq \frac{2|y_n|}{l_\nu}+\frac{4}{l_\nu} \cdot \frac{\sqrt{n} \cdot |x_{n,\nu,(l_i)_{1 \leq i \leq M}}|+|x_n|}{(1+\frac{1}{n}tr(\mathcal{M}(\hat{l}_\nu I -\mathcal{M})^{-1})) \cdot (1+\frac{1}{n}tr(\mathcal{M}(\hat{l}_\nu I -\mathcal{M})^{-1})-\frac{|x_n|}{\sqrt{n}}-c(\epsilon_0) \cdot \frac{M}{n})} \xrightarrow[]{p} 0
\end{equation}
because \(\frac{1}{n}tr(\mathcal{M}(\hat{l}_\nu I -\mathcal{M})^{-1}) \xrightarrow[]{a.s.} 0,\)
\[y_n \xrightarrow[]{p} 0, \hspace{0.4cm} \frac{\sqrt{n} \cdot |x_{n,\nu,(l_i)_{1 \leq i \leq M}}|}{l_\nu} \leq c(\epsilon_0) \cdot \frac{M/\sqrt{n}}{l_\nu} \to 0, \hspace{0.4cm}  \frac{x_n}{\sqrt{\mathbb{E}[z_{11}^4]-1}} \Rightarrow N(0,1), \hspace{0.4cm} \mathbb{E}[z_{11}^4]-1 \leq c(K), \hspace{0.4cm} l_\nu \to \infty,\]
Finally, note (\ref{727}) ensues from (\ref{51}) since almost surely for \(n\) large enough, \(\newline ||\mathcal{M}|| \leq c_{K,\gamma} \leq \frac{l_\nu}{2}, \overline{l}_\nu, \hat{l}_\nu \geq \frac{l_\nu}{2},\) and for \(0 \leq x \leq c_{K,\gamma} \leq \frac{l_\nu}{2},\)
\[|\frac{x}{\hat{l}_\nu-x}-\frac{x}{\overline{l}_\nu-x}| \leq \frac{4c_{K,\gamma} \cdot |\overline{l}_\nu-\hat{l}_\nu|}{l_\nu^2}.\] 
\par
Second, we proceed with the proof of (\ref{38}), which is equivalent to
\begin{equation}\label{52}
    \sqrt{n} \cdot (\frac{l_\nu}{l_\nu-1}-\overline{l}_\nu \cdot \frac{1}{N-M}tr((\overline{l}_\nu I-\mathcal{M})^{-1})) \xrightarrow[]{p} 0.
\end{equation}
Denote by 
\[m_n(z)=\frac{z+\gamma_n-1-\sqrt{(z-\gamma_n+1)^2-4z}}{2\gamma_n z}\] 
for \(z \ne 0,\gamma_n=\frac{N-M}{n}\) (the Stieltjes transform of a Marchenko-Pastur law), for which 
\[m_n(\overline{l}_\nu)=\frac{l_\nu}{l_\nu-1} \cdot \frac{1}{\overline{l}_\nu}\] 
in virtue of
\[(\overline{l}_\nu-\gamma_n+1)^2-4\overline{l}_\nu=(l_\nu+1+\frac{\gamma_n}{l_\nu-1})^2-4l_\nu(1+\frac{\gamma_n}{l_\nu-1})=(l_\nu-1-\frac{\gamma_n}{l_\nu-1})^2,\] 
and rewrite (\ref{52}) as
\[\sqrt{n} \cdot \overline{l}_\nu(m_n(\overline{l}_\nu)-\frac{1}{N-M}tr((\overline{l}_\nu I-\mathcal{M})^{-1})) \xrightarrow[]{p} 0.\]
\par
Lastly, we prove that for deterministic \(z=z_n \to \infty,\) 
\[m_{F_n}(z)=\frac{1}{N-M}tr((z I-\mathcal{M})^{-1})=\frac{1}{N-M}tr((zI-S_{BB})^{-1}),\]
\begin{equation}\label{53}
    \sqrt{n} \cdot z(m_n(z)-m_{F_n}(z)) \xrightarrow[]{p} 0,
\end{equation}
which will complete the proof of (\ref{52}), by adopting the argument in lemma 2 from Ledoit and Péché~\cite{ledoitpeche}.
\par
Denote by \(Z_1,Z_2, \hspace{0.05cm} ... \hspace{0.05cm}, Z_{n} \in \mathbb{R}^{N-M}\) the columns of \(Z_B,\) and recall that \(S_{BB}=\frac{1}{n}Z_BZ_B^T,\) from which almost surely for \(n\) large enough such that \(z>c_{K,\gamma},\) \(zI-S_{BB}\) is invertible and so
\begin{equation}\label{740}
    -1+zm_{F_n}(z)=\frac{1}{n(N-M)}\sum_{1 \leq i \leq n}{Z_i^T(zI-S_{BB})^{-1}Z_i}
\end{equation}
from multiplying on the right by \((zI-S_{BB})^{-1}\) the identity
\[-(zI-S_{BB})+zI=\frac{1}{n}\sum_{1 \leq i \leq n}{Z_iZ_i^T},\]
taking the trace of each side, and dividing by \(N-M.\)
Take \(R_i=(zI-S_{BB}+\frac{1}{n}Z_iZ_i^T)^{-1}\) for \(1 \leq i \leq n,\) and recall the resolvent identity \[R_i-(zI-S_{BB})^{-1}=R_i \frac{1}{n}Z_iZ_i^T (zI-S_{BB})^{-1},\]
which by multiplication on the left by \(Z_i^T\) and by \(Z_i\) on the right yields 
\[Z_i^T(zI-S_{BB})^{-1}Z_i=\frac{Z_i^TR_iZ_i}{1-\frac{1}{n}Z_i^TR_iZ_i}=-n+\frac{n}{1-\frac{1}{n}Z_i^TR_iZ_i}\]
because 
\begin{equation}\label{5083}
    |\frac{1}{n}Z_i^TR_iZ_i| \leq \frac{1}{n}||Z_i||^2 \cdot ||R_i|| \leq 2 \cdot \frac{2}{z}<1
\end{equation}
almost surely (inequality (\ref{7}) yields \(||S_{BB}-\frac{1}{n}Z_iZ_i^T|| \leq c(K,\gamma)<\frac{z}{2}\) with probability one), from which (\ref{740}) becomes
\begin{equation}\label{54}
    -1+zm_{F_n}(z)=-\frac{1}{\gamma_n}
    +\frac{1}{N-M}\sum_{1 \leq i \leq n}{\frac{1}{1-\frac{1}{n}Z_i^TR_iZ_i}}=-\frac{1}{\gamma_n}
    +\frac{1}{\gamma_n} \cdot \frac{1}{1-\gamma_n m_{F_n}(z)}+\Delta_1+\Delta_2
\end{equation}
where
\[\Delta_1=\frac{1}{N-M}\sum_{1 \leq i \leq n}{(\frac{1}{1-\frac{1}{n}Z_i^TR_iZ_i}-\frac{1}{1-\frac{1}{n}tr(R_i)})},\]
\[\Delta_2=\frac{1}{N-M}\sum_{1 \leq i \leq n}{(\frac{1}{1-\frac{1}{n}tr(R_i)}-\frac{1}{1-\frac{1}{n}tr((zI-S_{BB})^{-1})})},\]
both being small in probability as we show next.

\vspace{0.2cm}
\(\Delta_1=o_p(\frac{1}{\sqrt{n}}):\) Almost surely
\[|\Delta_1| \leq \frac{1}{N-M}\sum_{1 \leq i \leq n}{4|\frac{1}{n}Z_i^TR_iZ_i-\frac{1}{n}tr(R_i)|}\]
because from (\ref{5083})
\[\frac{1}{n}Z_i^TR_iZ_i \leq \frac{4}{z} \leq \frac{1}{2}, \hspace{0.2cm} \frac{1}{n}tr(R_i) \leq \frac{N-M}{n} \cdot ||R_i|| \leq 2\gamma \cdot \frac{2}{z} \leq \frac{1}{2}.\] 
\par
It suffices to prove that
\[\mathbb{E}[(\frac{1}{n}Z_1^TR_1Z_1-\frac{1}{n}tr(R_1))^2 \chi_{ ||S_{BB}-\frac{1}{n}Z_1Z_1^T|| \leq \frac{z}{2}}] \leq (\frac{N-M}{n})^2 \cdot \frac{4}{z^2n}\]
since this yields for the events \(A_{i,n}=\{||S_{BB}-\frac{1}{n}Z_iZ_i^T|| \leq \frac{z}{2}\}, 1 \leq i \leq n,\)
\[n\mathbb{E}[(\frac{1}{N-M}\sum_{1 \leq i \leq n}{|\frac{1}{n}Z_i^TR_iZ_i-\frac{1}{n}tr(R_i)|})^2\chi_{A_{1,n}\cap A_{2,n} \cap ...\cap A_{n,n}}] \leq\] 
\[\leq n\mathbb{E}[\frac{n}{(N-M)^2}\sum_{1 \leq i \leq n}{(\frac{1}{n}Z_i^TR_iZ_i-\frac{1}{n}tr(R_i))^2}\chi_{A_{1,n}\cap A_{2,n} \cap ...\cap A_{n,n}}] \leq\]
\[\leq \frac{n^2}{(N-M)^2} \cdot \sum_{1 \leq i \leq n} {\mathbb{E}[(\frac{1}{n}Z_i^TR_iZ_i-\frac{1}{n}tr(R_i))^2\chi_{A_{i,n}}]} \leq \frac{4}{z^2} \to 0\]
with which the following simple fact can be used:

\vspace{0.2cm}
If \(|X_n| \leq Y_n\) almost surely, \(\lim_{n \to \infty}{\mathbb{E}[Y_n^2\chi_{A_n}]}=0, \lim_{n \to \infty}{\mathbb{P}(A_n)}=1,\) then \(X_n \xrightarrow[]{p} 0.\)

\vspace{0.2cm}
\textit{Proof:} For any \(t>0,\) 
\[\mathbb{P}(|X_n| \geq t) \leq \mathbb{P}(A_n^c)+\mathbb{P}(|X_n| \geq t, A_n) \leq \mathbb{P}(A_n^c)+\mathbb{P}(|X_n|>Y_n, A_n)+\mathbb{P}(Y_n \geq t, A_n) \leq \]
\[\leq \mathbb{P}(A_n^c)+\mathbb{P}(|X_n|>Y_n)+\frac{\mathbb{E}[Y_n^2\chi_{A_n}]}{t^2} \xrightarrow[]{n \to \infty} 0.\]

\vspace{0.2cm}
\par
Because \(Z_1=[z_{(M+1)1} \hspace{0.1cm} z_{(M+2)1} \hspace{0.1cm} ... \hspace{0.1cm} z_{N1}]^T\) and \(S_{BB}-\frac{1}{n}Z_1Z_1^T, R_1:=(r_{ij})_{1 \leq i,j \leq N-M}\) are independent,
\[\mathbb{E}[(\frac{1}{N-M}Z_1^TR_1Z_1-\frac{1}{N-M}tr(R_1))^2 \hspace{0.05cm} |S_{BB}-\frac{1}{n}Z_1Z_1^T]=\frac{1}{(N-M)^2}\mathbb{E}[\sum_{M+1 \leq i,j \leq N}{((z_{i1}z_{j1}-\delta_{i,j})r_{ij}})^2|R_1]=\]
\[=\frac{1}{(N-M)^2}(\mathbb{E}[(z_{11}^2-1)^2]\sum_{1 \leq i \leq N-M}{r^2_{ii}}+4\sum_{i < j}{r_{ij}}^2) \leq \frac{c(K)}{(N-M)^2} \cdot ||R_1||_F^2 \leq \frac{c(K)}{N-M} \cdot ||R_1||^2\]
since the cross term corresponding to \((i,j),(i',j')\) is zero in expectation unless \(\{i,j\}=\{i',j'\}\) as multisets and \(R_1\) is symmetric, from which
\[\mathbb{E}[(\frac{1}{N-M}Z_1^TR_1Z_1-\frac{1}{N-M}tr(R_1))^2 \chi_{||S_{BB}-\frac{1}{n}Z_1Z_1^T|| \leq \frac{z}{2}}] \leq \frac{c(K)}{N-M} \cdot \frac{4}{z^2}\]
because \(||S_{BB}-\frac{1}{n}Z_1Z_1^T|| \leq \frac{z}{2}\) implies 
\[||R_1||=\frac{1}{z-||S_{BB}-\frac{1}{n}Z_1Z_1^T||} \leq \frac{2}{z}.\]

\vspace{0.2cm}
\(\Delta_2=o_p(\frac{1}{\sqrt{n}}):\) For any \(t>0,\)
\[\mathbb{P}(\sqrt{n}|\Delta_2| \geq t) \leq n\mathbb{P}(|\frac{1}{1-\frac{1}{n}tr(R_1)}-\frac{1}{1-\frac{1}{n}tr((zI-S_{BB})^{-1})}| \geq \frac{t \cdot \frac{N-M}{n}}{\sqrt{n}}) \leq\]
\[\leq n \mathbb{P}(||R_1|| \geq \frac{n}{2(N-M)})+n\mathbb{P}(||zI-S_{BB}|| \leq \frac{N-M}{2n})+n\mathbb{P}(|\frac{1}{n}tr(R_1)-\frac{1}{n}tr((zI-S_{BB})^{-1})| \geq \frac{4t \cdot \frac{N-M}{n}}{\sqrt{n}})\]
which tends to zero as \(n \to \infty,\) using (\ref{7}) for the first two terms while for the third, for \(B \in \mathbb{R}^{p \times p}, u \in \mathbb{R}^{p}\) with \(B, B+uu^T\) invertible,
\[u^TB^{-1}(B+uu^T)=(1+u^TB^{-1}u)u^T \Longleftrightarrow u^TB^{-1}=(1+u^TB^{-1}u)u^T(B+uu^T)^{-1},\]
\[B^{-1}-(B+uu^T)^{-1}=(B+uu^T)^{-1}uu^TB^{-1}=(1+u^TB^{-1}u) \cdot (B+uu^T)^{-1}uu^T(B+uu^T)^{-1},\]
\[tr(B^{-1})-tr((B+uu^T)^{-1})=(1+u^TB^{-1}u) \cdot u^T(B+uu^T)^{-2}u,\]
from which for \(n\) large enough,
\[|tr((zI-S_{BB})^{-1})-tr((zI-R_1)^{-1})| \leq (1+\frac{1}{n}||Z_1||^2 \cdot ||zI-S_{BB}||^{-1}) \cdot \frac{1}{n}||Z_1||^2 \cdot ||zI-S_{BB}+\frac{1}{n}Z_1Z_1^T||^{-2} \xrightarrow[]{a.s.} 0,\]
yielding
\[n\mathbb{P}(|\frac{1}{n}tr(R_1)-\frac{1}{n}tr((zI-S_{BB})^{-1})| \geq \frac{4t \cdot \frac{N-M}{n}}{\sqrt{n}}) \leq n\mathbb{P}(\frac{1}{n}||Z_1||^2 \geq 2)+n\mathbb{P}(||S_{BB}||>c_{K,\gamma}) \to 0.\]
\par
Finally, 
\[m_n(z)=\frac{z+\gamma_n-1-\sqrt{(z-\gamma_n+1)^2-4z}}{2\gamma_n z}\] 
leads to 
\begin{equation}\label{55}
    -1+zm_n(z)=-\frac{1}{\gamma_n}+\frac{1}{\gamma_n} \cdot \frac{1}{1-\gamma_n m_n(z)}
\end{equation}
from
\[(z-\gamma_n+1)^2-4z=(z+\gamma_n-1)^2-4\gamma_nz,\]
and
\[\gamma_n zm^2_n(z)-m_n(z)(z+\gamma_n-1)+1=0\] 
or equivalently 
\[(\gamma_n m_n(z)-1) \cdot (zm_n(z)-1)=-m_n(z),\] 
while (\ref{54}) and (\ref{55}) give
\[z(m_{F_n}(z)-m_n(z))=\frac{1}{\gamma_n} \cdot (\frac{1}{1-\gamma_n m_{F_n}(z)}-\frac{1}{1-\gamma_n m_n(z)})+\Delta_1+\Delta_2,\]
from which
\[z(m_{F_n}(z)-m_n(z)) \cdot (1-\frac{1}{z(1-\gamma_n m_{F_n}(z)) \cdot (1-\gamma_n m_{n}(z))})=\Delta_1+\Delta_2.\]
Because \(\sqrt{n}\Delta_1,\sqrt{n}\Delta_2 \xrightarrow[]{p} 0,\) and 
\[\frac{1}{z(1-\gamma_n m_{F_n}(z)) \cdot (1-\gamma_n m_{n}(z))} \xrightarrow[]{a.s.} 0,\] 
(\ref{53}) ensues.

\section{Consistency of Eigenvectors}\label{sec6}

The proof of Theorem~\ref{th5} is covered in the first four subsections: \((a)\) in \ref{6.1}, \((b)(i), (c)(i)\) in \ref{6.2},  \((b)(ii)\) in \ref{6.3}, and \((c)(ii)\) in \ref{6.4}. The justification of Theorem~\ref{th6}, in the same spirit as the proof of Theorem~\ref{th1}, is presented in \ref{6.5}.

\subsection{Proof of Theorem~\ref{th5}, part (a)}\label{6.1}

Notice that
\begin{equation}\label{60}
    l_\nu(1-<p_\nu,u_\nu>^2)=l_\nu(1-(1-R_\nu^2)<a_\nu,e_\nu>^2)=\frac{l_{\nu}}{\hat{l}_{\nu}} \cdot (\hat{l}_\nu R_\nu^2-\frac{N}{n})+l_\nu (1-R_\nu^2) (1-<a_\nu, e_\nu>^2)+\frac{N}{n}\cdot \frac{l_{\nu}}{\hat{l}_{\nu}}.
\end{equation}
\par
Since Proposition~\ref{prop1}, Proposition~\ref{prop2}, and Lemma~\ref{lemma3} yield
\begin{equation}\label{720}
    \frac{l_{\nu}}{\hat{l}_{\nu}} \xrightarrow[]{a.s.} 1, \hspace{0.5cm} R_\nu^2 \xrightarrow[]{a.s.} 0 \hspace{0.5cm} \hat{l}_\nu R_\nu^2-\frac{N}{n} \xrightarrow[]{a.s.} 0,
\end{equation}
respectively, it suffices to show that
\begin{equation}\label{2097}
    l_\nu (1-<a_\nu, e_\nu>^2) \xrightarrow[]{p} 0,
\end{equation}
if \(\lim_{n \to \infty}{\frac{l^{(n)}_\nu}{n/M(n)}}=0,\) with the convergence holding almost surely if \( \sum_{n \in \mathbb{N}}{\exp(-\frac{an}{l^{(n)}_\nu})}<\infty, \forall a>0\) is also satisfied. As
\[l_\nu(1-<a_\nu,e_\nu>^2) \leq  2l_\nu(1-<a_\nu,e_\nu>)=l_\nu||a_\nu-e_\nu||^2,\]
for \(t>0,\)
\[\mathbb{P}(l_\nu(1-<a_\nu,e_\nu>^2) \geq t) \leq \mathbb{P}(||a_\nu-e_\nu|| \geq 2\beta_\nu)+\mathbb{P}(l_\nu \beta_\nu^2 \geq \frac{t}{4}).\]
Moreover, Lemma~\ref{lemma2} entails 
\[\sum_{n \to \infty}{\mathbb{P}(||a_\nu-e_\nu|| \geq 2||\beta_\nu||)}<\infty,\] 
and we show next \(\lim_{n \to \infty}{\mathbb{P}(l_\nu \beta_\nu^2 \geq \frac{t}{4})}=0,\) with \(\sum_{n \in \mathbb{N}}{\mathbb{P}(l_\nu \beta_\nu^2 \geq \frac{t}{4})}<\infty\) if the second condition holds as well, from which (\ref{2097}) and its almost sure variant will ensue. 
\par
With the aid of
\[\beta_\nu^2 \leq c(\epsilon_0)\sum_{k \ne \nu}{(\frac{1}{n}z_k^Tz_\nu)^2}+c(\epsilon_0)\sum_{k \ne \nu}{(t_k^T\mathcal{M}(\hat{l}_\nu I -\mathcal{M})^{-1}t_\nu)^2},\]
it follows that for \(n \geq n(\epsilon_0,t)\)
\begin{equation}\label{708}
    \mathbb{P}(l_\nu \beta_\nu^2 \geq \frac{t}{4}) \leq
    \mathbb{P}(l_\nu \sum_{k \ne \nu}{(\frac{1}{n}z_k^Tz_\nu)^2} \geq \frac{t}{8c(\epsilon_0)})+\mathbb{P}(\frac{n}{M}\sum_{k \ne \nu}{(t_k^T\mathcal{M}(\hat{l}_\nu I -\mathcal{M})^{-1}t_\nu)^2} \geq \frac{tn}{8c(\epsilon_0)Ml_\nu}).
\end{equation}
\par
Inequality (\ref{909}) from part \((b)\) of Lemma~\ref{lemma1} yields the second term in (\ref{708}) as a function of \(n \in \mathbb{N}\) is summable as:
\[\frac{tn}{Ml_\nu} \cdot l_\nu^2 \geq 2, \hspace{0.5cm}  Ml_\nu^4 \cdot (\frac{tn}{Ml_\nu})^2 \geq t^2n, \hspace{0.5cm} Ml_\nu^2 \cdot \frac{tn}{Ml_\nu} \geq tn.\]
\par
Regarding the first term in (\ref{708}), for \(t>0,\)
\[\mathbb{P}(l_\nu \sum_{k \ne \nu}{(\frac{1}{n}z_k^Tz_\nu)^2} \geq 2t|z_\nu) \leq \mathbb{P}(l_\nu \cdot \frac{M}{n^2}||z_\nu||^2 \geq t|z_\nu)+\mathbb{P}(l_\nu |\sum_{k \ne \nu}{((\frac{1}{n}z_k^Tz_\nu)^2-\frac{1}{n^2}||z_\nu||^2)}| \geq t|z_\nu).\]
Inequality (\ref{707}) gives
\[\mathbb{P}(l_\nu |\sum_{k \ne \nu}{((\frac{1}{n}z_k^Tz_\nu)^2-\frac{1}{n^2}||z_\nu||^2)}| \geq t|z_\nu) \leq 2\exp(-c\min{(\frac{t^2}{Ml_\nu^2 \cdot (\frac{1}{n}||z_\nu||)^4},\frac{t}{l_\nu \cdot (\frac{1}{n}||z_\nu||)^2})}).\]
Therefore, for \(n\) large enough such that \(l_\nu \leq \frac{n}{M} \cdot \frac{t}{2},\)
\[\mathbb{P}(l_\nu \sum_{k \ne \nu}{(\frac{1}{n}z_k^Tz_\nu)^2} \geq t) \leq \mathbb{P}(\frac{1}{n}||z_\nu||^2 \geq 2)+2\exp(-c\min{(\frac{t^2n^2}{4Ml_\nu^2},\frac{tn}{2l_\nu})}) \leq \]
\[\leq \mathbb{P}(\frac{1}{n}||z_\nu||^2 \geq 2)+2\exp(-\frac{ctn}{2l_\nu}) \leq \mathbb{P}(\frac{1}{n}||z_\nu||^2 \geq 2)+2\exp(-\frac{ct}{\frac{2l_\nu}{n/M}}),\]
from which we conclude together with \(\lim_{n \to \infty}{\frac{l^{(n)}_\nu}{n/M}}=0\) that the bounds in (\ref{708}) tend to zero, and that they are summable, if additionally \(\sum_{n \in \mathbb{N}}{\exp(-\frac{an}{l^{(n)}_\nu})}<\infty, \forall a>0.\)

\subsection{Proof of Theorem~\ref{th5}, parts (b)(i), (c)(i)}\label{6.2}

Notice that
\[n(1-<p_\nu,u_\nu>^2)=nR_\nu^2+n(1-R_\nu^2) (1-<a_\nu, e_\nu>^2).\]
\par
From Lemma~\ref{lemma3}, Proposition~\ref{prop1}, Proposition~\ref{prop2},
\[nR_\nu^2-\frac{N}{l_\nu}=\frac{n}{l_\nu} \cdot \frac{l_\nu}{\hat{l}_\nu} \cdot (\hat{l}_\nu R_\nu^2-\frac{N}{n})+\frac{N}{l_\nu}(1-\frac{l_\nu}{\hat{l}_\nu}) \xrightarrow[]{a.s.} 0, \hspace{0.5cm} R_\nu^2 \xrightarrow[]{a.s.} 0.\]
It suffices for both  \((b)(i),\) and \((c)(i)\) to show
\[n(1-<a_\nu, e_\nu>^2) \Rightarrow X_{(c_{k\nu})_{k \ne \nu}}\]
as the claims ensue then from Slutsky's lemma:
\[n(1-<p_\nu,u_\nu>^2)-\frac{N}{l_\nu} \Rightarrow X_{(c_{k\nu})_{k \ne \nu}},\]
which can be rewritten as
\[l_\nu(1-<p_\nu,u_\nu>^2)-\frac{N}{n} \Rightarrow \frac{c_\nu}{M}X_{(c_{k\nu})_{k \ne \nu}}.\]
\par
Note that
\[n(1-<a_\nu, e_\nu>^2)=2n(1-<a_\nu,e_\nu>)-n(1-<a_\nu,e_\nu>)^2=n||a_\nu-e_\nu||^2-\frac{n}{4}||a_\nu-e_\nu||^4.\]
Lemma~\ref{lemma2} gives that almost surely \(n||a_\nu-e_\nu||^4 \leq 16n\beta_\nu^4,\) and from (\ref{25}), \(\sqrt{n} \cdot \beta_\nu^2 \xrightarrow[]{p} 0.\) Therefore, by employing Slutsky's lemma once again and Lemma~\ref{lemma2}, it suffices to show that 
\[n\beta_\nu^2 \Rightarrow X_{(c_{k \nu})_{{k \ne \nu}}}.\]
\par
Recall that
\[n\beta_\nu^2=\sum_{k \ne \nu}{\frac{1}{(l_k-l_\nu)^2}\mathcal{D}_{\nu k}^2}=\sum_{k \ne \nu}{\frac{l_kl_\nu}{(l_k-l_\nu)^2}(\frac{1}{\sqrt{n}}z_k^Tz_\nu)^2}+\sum_{k \ne \nu}{\frac{l_kl_\nu}{(l_k-l_\nu)^2}(\sqrt{n} t_k^T\mathcal{M}(\hat{l}_\nu I -\mathcal{M})^{-1} t_\nu)^2}+\]
\begin{equation}\label{41}
    +2\sum_{k \ne \nu}{\frac{l_kl_\nu}{(l_k-l_\nu)^2} \cdot  \frac{1}{\sqrt{n}}z_k^T z_\nu} \cdot \sqrt{n} t_k^T\mathcal{M}(\hat{l}_\nu I -\mathcal{M})^{-1} t_\nu.
\end{equation}
Part \((b)\) of Lemma~\ref{lemma1} yields the second term in (\ref{41}) tends to zero in probability as
\begin{equation}\label{58}
    \sum_{k \ne \nu}{\frac{l_kl_\nu}{(l_k-l_\nu)^2}(\sqrt{n} t_k^T\mathcal{M}(\hat{l}_\nu I -\mathcal{M})^{-1} t_\nu)^2} \leq Mc(\epsilon_0)\cdot \frac{n}{M} \sum_{k \ne \nu}{( t_k^T\mathcal{M}(\hat{l}_\nu I -\mathcal{M})^{-1} t_\nu)^2} \xrightarrow[]{p} 0.
\end{equation}
Because
\[\mathbb{E}[\sum_{k \ne \nu}{\frac{l_kl_\nu}{(l_k-l_\nu)^2}(\frac{1}{\sqrt{n}}z_k^Tz_\nu)^2}]=\sum_{k \ne \nu}{\frac{l_kl_\nu}{(l_k-l_\nu)^2}} \leq Mc(\epsilon_0),\]
it follows the third term in (\ref{41}) also tends to zero in probability from Cauchy-Schwarz inequality, (\ref{58}), and the following simple fact:

\vspace{0.5cm}
If \(X_n,Y_n \geq 0\) are random variables defined on the same space and \(X_n \xrightarrow[]{p} 0,\mathbb{E}[Y_n] \leq c,\) then \begin{equation}\label{42}
    X_nY_n \xrightarrow[]{p} 0.
\end{equation}
\par
\textit{Proof:} For any \(t,s>0,\)
\[\mathbb{P}(X_nY_n \geq t) \leq \mathbb{P}(Y_n \geq s)+\mathbb{P}(X_n \geq \frac{t}{s}) \leq \frac{c}{s}+\mathbb{P}(X_n \geq \frac{t}{s}),\]
from which \(\limsup_{n \to \infty}{\mathbb{P}(X_nY_n \geq t)} \leq \frac{c}{s},\) and taking \(s \to \infty\) concludes the justification of the claim.
\vspace{0.5cm}

\par
We turn now to the main contribution in (\ref{41}):
\[\sum_{k \ne \nu}{\frac{l_kl_\nu}{(l_k-l_\nu)^2}(\frac{1}{\sqrt{n}}z_k^Tz_\nu)^2}.\]
\par
Because \(\lim_{n \to \infty}{\frac{l_kl_\nu}{(l_k-l_\nu)^2}}=c_{k\nu},\) it suffices to prove that 
\begin{equation}\label{43}
    T_n:=\sum_{k \ne \nu}{c_{k\nu}(\frac{1}{\sqrt{n}}z_k^Tz_\nu)^2} \Rightarrow X_{(c_{k\nu})_{k \ne \nu}},
\end{equation}
where 
\[X_{(c_{k\nu})_{k \ne \nu}}=\sum_{k \ne \nu}{c_{k \nu}y_{k \nu}^2},\] 
with \((y_{k \nu})_{k \ne \nu}\) mutually independent standard normal distributions, because
\[\sum_{k \ne \nu}{(\frac{l_kl_\nu}{(l_k-l_\nu)^2}-c_{k\nu})(\frac{1}{\sqrt{n}}z_k^Tz_\nu)^2} \xrightarrow[]{p} 0\]
as for any \(t>0, \epsilon>0,\)
\[\mathbb{P}(\sum_{k \ne \nu}{|\frac{l_kl_\nu}{(l_k-l_\nu)^2}-c_{k\nu}| \cdot (\frac{1}{\sqrt{n}}z_k^Tz_\nu)^2} \geq t) \leq \mathbb{P}(\sum_{k \ne \nu}{(\frac{1}{\sqrt{n}}z_k^Tz_\nu)^2} \geq \frac{t}{\epsilon}) \leq \frac{M\epsilon}{t},\]
using \(\mathbb{E}[\sum_{k \ne \nu}{(\frac{1}{\sqrt{n}}z_k^Tz_\nu)^2}]=M-1.\)
\par
If \(c_{k\nu}=0\) for all \(k \ne \nu,\) (\ref{43}) is clear. Suppose next that \(\max_{k \ne \nu}{c_{k \nu}}>0.\) We use the moment convergence theorem and Carleman's condition (lemmas \(B1, B3,\) respectively in Bai and Silverstein~\cite{baisilvbook}), stated below, to show (\ref{43}) holds:

\vspace{0.3cm}
A sequence of distributions \(F_n\) converges weakly to a limit if each \(F_n\) has moments of all orders; for all \(k \in \mathbb{N},\) their \(k^{th}\) moments converge to some finite limit \(\beta_k;\) and for any right-continuous nondecreasing functions \(F,G\) with the same moment sequence, \(F=G+c,\) for some \(c \in \mathbb{R}.\) 

\vspace{0.3cm}
If \((\beta_k)\) is the moment sequence of a distribution function \(F\) and \(\sum_{k \geq 1}{\beta_{2k}^{-1/(2k)}}=\infty,\) then \(F\) is uniquely determined by its moment sequence.

\vspace{0.5cm}
Namely, we show that
\begin{equation}\label{709}
    \lim_{n \to \infty}{\mathbb{E}[T_n^m]}=\beta_m, m \in \mathbb{N},\hspace{0.5cm} \sum_{m \geq 1}{\beta_{2m}^{-1/(2m)}}=\infty,
\end{equation}
for
\[\beta_m=\mathbb{E}[X_{(c_{k\nu})_{k \ne \nu}}^m]=\sum{\binom{m}{(q_k)_{k \ne \nu}}{\prod_{k \ne \nu}{\frac{(2q_k)!}{2^{q_k}q_k!}}c_{k\nu}^{q_k}}},\]
where the sum is over non-negative integer sequences \((q_k)_{k \ne \nu}\) with \(\sum_{k \ne \nu}{q_k}=m.\)
\par
Fix \(m \in \mathbb{N}.\) Then
\[\mathbb{E}[T_n^m]=\mathbb{E}[(\sum_{k \ne \nu}{c_{k\nu}(\frac{1}{\sqrt{n}}z_k^Tz_\nu)^2})^m].\]
Consider the term in this last sum corresponding to some sequence \((q_k)_{k \ne \nu}\) with \(\sum_{k \ne \nu}{q_k}=m:\)
\[\binom{m}{(q_k)_{k \ne \nu}}\prod_{k \ne \nu}{c_{k\nu}^{q_k}} \cdot \frac{1}{n^m} \mathbb{E}[\prod_{k \ne \nu}{(\sum_{1 \leq j \leq n}{z_{kj}z_{\nu j})^{2q_k}}}].\]
\par
For the first claim in (\ref{709}), it is enough to show that as \(n \to \infty,\)
\begin{equation}\label{710}
    \frac{1}{n^m} \mathbb{E}[\prod_{k \ne \nu}{(\sum_{1 \leq j \leq n}{z_{kj}z_{\nu j})^{2q_k}}}] \to \prod_{k \ne \nu}{\frac{(2q_k)!}{2^{q_k}q_k!}}.
\end{equation}
In this product, the factor corresponding to \(k \ne \nu\) will contribute to the terms that form up the expectation products of the type
\[c(m_1,m_2, \hspace{0.05cm} ... \hspace{0.05cm}, m_{t})(z_{kj_1}z_{\nu j_1})^{m_1}...(z_{kj_t}z_{\nu j_t})^{m_t}\]
for \(m_1+...+m_t=2q_k,m_i>0,\) and \(j_1, \hspace{0.05cm} ... \hspace{0.05cm}, j_t\) pairwise distinct. If some \(m_i=1,\) then the expectation of any term in which this product appears will be zero by independence (\(z_{k j_i}\) shows up only once in each such term, is independent of the rest, and has mean zero). Hence we are left with the products containing \(m_1,\hspace{0.05cm} ... \hspace{0.05cm}, m_t \geq 2.\) For any \(1 \leq t \leq m\) and fixed \(k \ne \nu,\) the number of such products having some \(m_i \geq 3\) is at most 
\[c(m)n^{t+\sum_{k' \ne k}{q_{k'}}} \leq c(m)n^{q_k-1+\sum_{k' \ne k}{q_k}} \leq c(m)n^{m-1},\] 
with an overall contribution in the left-hand side term in (\ref{710}) upper bounded by
\[\frac{1}{n^m} \cdot m \cdot M \cdot c(m)n^{m-1}\mathbb{E}[z_{11}^{4m}] \to 0\] 
from the arithmetic-geometric mean inequality (\ref{5084}) and \(\mathbb{E}[|z_{11}|^l] \leq c(K,l)<\infty,\) for all \(l \in \mathbb{N}.\)  
\par
Therefore, the only such products surviving in the limit are the ones for which \(m_1=...=m_t=2.\) Moreover, the number of products for which some \(z_{\nu j}\) appears in at least two factors corresponding to some \(k \ne k'\) is at most \(M^2 \cdot c(m)n^{m-1}:\) for each tuple \((q_k)_{k \ne \nu}\) with \(\sum_{k \ne \nu}{q_k}=m,\) denote the subsets of \(\{1,2, \hspace{0.05cm} ... \hspace{0.05cm}, n\}\) corresponding to \((q_k)_{k \ne \nu}\) by \((S_{q_k})_{k \ne \nu}\) (i.e., \(S_{q_k}=\{j_1,j_2, \hspace{0.05cm} ... \hspace{0.05cm}, j_t\}\)); each fixed \((q_k)_{k \ne \nu}\) and pair \((k_1,k_2), k_1 \ne k_2\) for which \(|S_{q_{k_1}} \cap S_{q_{k_2}}|>0\) generate at most \(m \cdot n^{m-1}\) products in the left-hand side term of (\ref{710}) because the collection \(\mathcal{S}\) of ordered tuples of subsets of \(\{1,2,\hspace{0.05cm} ... \hspace{0.05cm}, n\}\) with sizes \((q_k)_{k \ne \nu, k \ne k_1, k \ne k_2}, q_{k_1}-1, q_{k_2},\) respectively has 
\[|\mathcal{S}| \leq n^{-1+\sum_{k \ne \nu}{q_k}}=n^{m-1},\] 
and for the function taking 
\[\{(S_{q_k})_{k \ne \nu}: |S_{q_{k_1}} \cap S_{q_{k_2}}|>0, |S_{q_k}|=q_k, S_{q_k} \subset \{1,2,\hspace{0.05cm} ... \hspace{0.05cm}, n\}\}\] 
into \(\mathcal{S}\) by removing from \(S_{q_{k_1}}\) the smallest element in \(S_{q_{k_1}} \cap S_{q_{k_2}}\), the preimage of any tuple has size at most \(m\) (one element from \(S_{q_{k_2}}\) needs to be added to \(S_{q_{k_1}},\) and \(|S_{q_{k_2}}|= q_{k_2} \leq m\)). Consequently, the terms that will be nonzero in the limit will have pairwise distinct random variables from \(z_{\nu j}, 1 \leq j \leq n,\) and so their expectations will disappear, each being one, from which we get that as \(n \to \infty,\)
\[\frac{1}{n^m} \mathbb{E}[\prod_{k \ne \nu}{(\sum_{1 \leq j \leq n}{z_{kj}z_{\nu j})^{2q_k}}}] \to \frac{1}{n^m}\prod_{k \ne \nu}{\frac{(2q_k)!}{2^{q_k}} \binom{n-q_1-...-q_{k-1}}{q_k}} \to \prod_{k \ne \nu}{\frac{(2q_k)!}{2^{q_k}q_k!}},\]
where \(q_0=q_\nu=0.\) This concludes the proof of 
\[\lim_{n \to \infty}{\mathbb{E}[T_n^m]}=\beta_m, m \in \mathbb{N},\hspace{0.5cm} \beta_m=\mathbb{E}[X_{(c_{k\nu})_{k \ne \nu}}^m]=\sum{\binom{m}{(q_k)_{k \ne \nu}}{\prod_{k \ne \nu}{\frac{(2q_k)!}{2^{q_k}q_k!}}c_{k\nu}^{q_k}}}.\]
\par
The second result in (\ref{709}) follows from
\[\sum_{(q_k),\sum{q_k}=m}{\prod_{k \ne \nu}{\frac{(2q_k)!}{2^{q_k}q_k!}c_{k \nu}^{q_k}}} \leq
\sum_{(q_k),\sum{q_k}=m}{\prod_{k \ne \nu}{2^{q_k}q_k!c_{k \nu}^{q_k}}} \leq (2\max_{k \ne \nu}{c_{k \nu}})^m \sum_{(q_k),\sum{q_k}=m}{\prod_{k \ne \nu}{q_k!}} \leq (2\max_{k \ne \nu}{c_{k \nu}})^m \cdot m^M m!\]
because each product is at most \(m!\) and their number is at most \(m^M,\) yielding
\[0<\beta_m^{1/m} \leq 2\max_{k \ne \nu}{c_{k \nu}} \cdot  (m^M m!)^{1/m} \leq 2\max_{k \ne \nu}{c_{k \nu}} \cdot 2^{M}m.\] 
Finally, the claims about \((y_{ij})_{i \ne j}\) ensue from \(\frac{1}{\sqrt{n}}z_i^Tz_j \Rightarrow y_{ij}.\)

\subsection{Proof of Theorem~\ref{th5}, part (b)(ii)}\label{6.3}

Using (\ref{60}) and (\ref{720}), it suffices to prove
\[l_\nu(1-<a_\nu,e_\nu>^2)-\frac{l_\nu}{n} \sum_{k \ne \nu}{\frac{l_kl_\nu}{(l_k-l_\nu)^2}} \xrightarrow[]{p} 0,\]
with the convergence being almost sure if additionally \(\sum_{n \in \mathbb{N}}{\exp(-a \cdot M(n))}<\infty, \forall a>0.\) 
Because
\[l_\nu(1-<a_\nu,e_\nu>^2)=2l_\nu(1-<a_\nu,e_\nu>)-l_\nu(1-<a_\nu,e_\nu>)^2=l_\nu||a_\nu-e_\nu||^2-\frac{l_\nu}{4}||a_\nu-e_\nu||^4,\]
and from Lemma~\ref{lemma2} \(\frac{||a_\nu-e_\nu||}{\beta_\nu} \xrightarrow[]{a.s.} 1,\) we show that
\begin{equation}\label{721}
    l_\nu \beta_\nu^2 -\frac{l_\nu}{n} \sum_{k \ne \nu}{\frac{l_kl_\nu}{(l_k-l_\nu)^2}} \xrightarrow[]{p} 0,
\end{equation}
which will complete the proof of \((b)(ii)\) since this will yield \(l_\nu \beta_\nu^2\) bounded in probability, from which 
\[l_\nu||a_\nu-e_\nu||^2-l_\nu \beta_\nu^2=l_\nu \beta_\nu^2 \cdot (\frac{||a_\nu-e_\nu||^2}{\beta_\nu^2}-1) \xrightarrow[]{p} 0,\] 
while (\ref{721}) with Lemma~\ref{lemma2} imply
\[\beta_\nu^2=O(c(\epsilon_0) \cdot \frac{M}{n})+ o_p(\frac{1}{l_\nu}),\]
and so
\[\frac{l_\nu}{4}||a_\nu-e_\nu||^4 \xrightarrow[]{p} 0\] 
(these convergences hold almost surely if \(\sum_{n \in \mathbb{N}}{\exp(-a \cdot M(n))}<\infty, \forall a>0\) holds as well). 
\par
Recall that
\[l_\nu \beta_\nu^2=l_\nu\sum_{k \ne \nu}{\frac{l_kl_\nu}{(l_k-l_\nu)^2}(\frac{1}{n}z_k^Tz_\nu)^2}+l_\nu\sum_{k \ne \nu}{\frac{l_kl_\nu}{(l_k-l_\nu)^2}(t_k^T \mathcal{M}(\hat{l}_\nu I -\mathcal{M})^{-1} t_\nu)^2}+\]
\begin{equation}\label{723}
    +2l_\nu\sum_{k \ne \nu}{\frac{l_kl_\nu}{(l_k-l_\nu)^2} \cdot \frac{1}{n}z_k^Tz_\nu \cdot t_k^T \mathcal{M}(\hat{l}_\nu I -\mathcal{M})^{-1} t_\nu}.
\end{equation}
Inequality (\ref{909}) and \(\frac{l_kl_\nu}{(l_k-l_\nu)^2} \leq c(\epsilon_0), k \ne \nu, l_\nu \leq 2c_\nu \cdot \frac{n}{M}\) give the second term on the right-hand side of (\ref{723}) tends to zero almost surely as the upper bounds will be summable as functions of \(n \in \mathbb{N}: Ml_\nu^2 \geq Ml_\nu \geq \frac{c_\nu}{2} \cdot n.\) We show next that
\begin{equation}\label{722}
    l_\nu\sum_{k \ne \nu}{\frac{l_kl_\nu}{(l_k-l_\nu)^2}(\frac{1}{n}z_k^Tz_\nu)^2}-\frac{l_\nu}{n} \sum_{k \ne \nu}{\frac{l_kl_\nu}{(l_k-l_\nu)^2}} \xrightarrow[]{p} 0,
\end{equation}
which together with Cauchy-Schwarz inequality and 
\[\frac{l_\nu}{n} \sum_{k \ne \nu}{\frac{l_kl_\nu}{(l_k-l_\nu)^2}} \leq c_\nu c(\epsilon_0)\] 
will yield that the third term on the right-hand side of (\ref{723}) converges in probability to zero (almost surely if \(\sum_{n \in \mathbb{N}}{\exp(-a \cdot M(n))}<\infty, \forall a>0\) holds as well), concluding the justification of (\ref{721}).
\par
(\ref{722}) can be rewritten as
\[\frac{1}{M}\sum_{k \ne \nu}{\frac{l_kl_\nu}{(l_k-l_\nu)^2} \cdot ((\frac{1}{\sqrt{n}}z_k^Tz_\nu)^2-1)} \xrightarrow[]{p} 0,\]
with the convergence holding almost surely if additionally \(\sum_{n \in \mathbb{N}}{\exp(-a \cdot M(n))}<\infty, \forall a>0.\) Note that
\[\mathbb{P}(|\frac{1}{M}\sum_{k \ne \nu}{\frac{l_kl_\nu}{(l_k-l_\nu)^2} \cdot ((\frac{1}{\sqrt{n}}z_k^Tz_\nu)^2-1)}| \geq t) \leq \mathbb{P}(\frac{1}{M}\sum_{k \ne \nu}{|(\frac{1}{\sqrt{n}}z_k^Tz_\nu)^2-1|} \geq \frac{t}{c(\epsilon_0)}) \leq\]
\[\leq \mathbb{P}(\frac{1}{M}\sum_{k \ne \nu}{|(\frac{1}{\sqrt{n}}z_k^Tz_\nu)^2-\frac{1}{n}||z_\nu||^2|} \geq \frac{t}{2c(\epsilon_0)})+\mathbb{P}(|\frac{1}{n}||z_\nu||^2-1| \geq \frac{t}{2c(\epsilon_0)}),\]
and
\[\mathbb{P}(\frac{1}{M}\sum_{k \ne \nu}{|(\frac{1}{\sqrt{n}}z_k^Tz_\nu)^2-\frac{1}{n}||z_\nu||^2|} \geq t) \leq \mathbb{P}(\frac{1}{n}||z_\nu||^2 \geq 2)+2\exp(-c\min{(\frac{Mt^2}{4},\frac{tM}{2})})\]
from inequality (\ref{707})
\[\mathbb{P}(\frac{1}{M}\sum_{k \ne \nu}{|(\frac{1}{\sqrt{n}}z_k^Tz_\nu)^2-\frac{1}{n}||z_\nu||^2|} \geq t|z_\nu) \leq 2\exp(-c\min{(\frac{Mt^2}{(\frac{1}{n}||z_1||^2)^2},\frac{tM}{\frac{1}{n}||z_1||^2})}).\]
Both conclusions then follow by employing Hanson-Wright inequality (\ref{9}) and Borel-Cantelli lemma for the almost sure convergence.

\subsection{Proof of Theorem~\ref{th5}, part (c)(ii)}\label{6.4}

From (\ref{60}), (\ref{720}), and Proposition~\ref{prop1}, it suffices to show
\[l_\nu(1-<a_\nu,e_\nu>^2)-\frac{l_\nu}{n}\sum_{k \ne \nu}{\frac{l_kl_\nu}{(l_k-l_\nu)^2}} \Rightarrow c_\nu N(0,2\sigma_\nu),\]
which is a consequence of Slutsky's lemma and the following identities proven below:

\begin{equation}\label{44}
    l_\nu||a_\nu-e_\nu||^2-\frac{l_\nu}{4}||a_\nu-e_\nu||^4=l_\nu\beta_\nu^2+o_p(1),
\end{equation}
\begin{equation}\label{45}
    l_\nu\beta_\nu^2-\frac{l_\nu}{n}\sum_{k \ne \nu}{\frac{l_kl_\nu}{(l_k-l_\nu)^2}}=\frac{l_\nu}{n}\sum_{k \ne \nu}{\frac{l_kl_\nu}{(l_k-l_\nu)^2}[(\frac{1}{\sqrt{n}}z_k^Tz_\nu)^2-1]}+o_p(1),
\end{equation}
\begin{equation}\label{46}
    T_n:=\frac{1}{\sqrt{M}}\sum_{k \ne \nu}{\frac{l_kl_\nu}{(l_k-l_\nu)^2}[(\frac{1}{\sqrt{n}}z_k^Tz_\nu)^2-1]} \Rightarrow N(0,2\sigma_\nu).
\end{equation}

\vspace{0.3cm} \textit{Proof of (\ref{44}):} Recall (\ref{13}) and (\ref{14}),
\[a_\nu-e_\nu=-\mathcal{R}_\nu \mathcal{D}_\nu e_\nu+r_\nu,\]
\[r_\nu=(<a_\nu,e_\nu>-1)e_\nu-\mathcal{R}_\nu \mathcal{D}_\nu (a_\nu-e_\nu)+(\hat{l}_\nu-l_\nu)\mathcal{R}_\nu(a_\nu-e_\nu),\]
and denote by
\[\alpha_\nu=||\mathcal{R}_\nu \mathcal{D}_\nu||+|\hat{l}_\nu-l_\nu| \cdot ||\mathcal{R}_\nu||,
\beta_\nu=||\mathcal{R}_\nu \mathcal{D}_\nu e_\nu||.\]
Then
\[||a_\nu-e_\nu||^2=\beta_\nu^2+||r_\nu||^2-2r_\nu \cdot \mathcal{R}_\nu \mathcal{D}_\nu e_\nu,\]
\[||r_\nu||^2=\frac{||a_\nu-e_\nu||^4}{4}+||\mathcal{R}_\nu \mathcal{D}_\nu (a_\nu-e_\nu)-(\hat{l}_\nu-l_\nu)\mathcal{R}_\nu(a_\nu-e_\nu)||^2,\]
as \(e_\nu,\mathcal{R}_\nu v\) are orthogonal for all \(v \in \mathbb{R}^{M},\) from which
\[||a_\nu-e_\nu||^2-\frac{||a_\nu-e_\nu||^4}{4}=\beta_\nu^2-2r_\nu \cdot \mathcal{R}_\nu \mathcal{D}_\nu e_\nu+||\mathcal{R}_\nu \mathcal{D}_\nu (a_\nu-e_\nu)+(\hat{l}_\nu-l_\nu)\mathcal{R}_\nu(a_\nu-e_\nu)||^2.\]
Notice that
\[|2r_\nu \cdot \mathcal{R}_\nu \mathcal{D}_\nu e_\nu|=|2(\mathcal{R}_\nu \mathcal{D}_\nu (a_\nu-e_\nu)-(\hat{l}_\nu-l_\nu)\mathcal{R}_\nu(a_\nu-e_\nu)) \cdot \mathcal{R}_\nu \mathcal{D}_\nu e_\nu| \leq 2\beta_\nu \alpha_\nu ||a_\nu-e_\nu||,\]
\[||\mathcal{R}_\nu \mathcal{D}_\nu (a_\nu-e_\nu)+(\hat{l}_\nu-l_\nu)\mathcal{R}_\nu(a_\nu-e_\nu)||^2 \leq \alpha_\nu^2||a_\nu-e_\nu||^2.\]
As Lemma~\ref{lemma2} gives \(\frac{||a_\nu-e_\nu||}{\beta_\nu} \xrightarrow[]{a.s.} 1\) and \(\alpha_\nu \xrightarrow[]{a.s.} 0\) from Proposition~\ref{prop2}, (\ref{44}) will ensue from
\[l_\nu \beta_\nu^2 \alpha_\nu \xrightarrow[]{p} 0.\]
\par
As it will be shown below, \(\beta_\nu^2=\tau_\nu+\frac{1}{l_\nu}w_n\) for \(w_n \Rightarrow c_\nu N(0,2\sigma_\nu), \tau_\nu=\frac{1}{n}\sum_{k \ne \nu}{\frac{l_kl_\nu}{(l_k-l_\nu)^2}}.\) Then
\[l_\nu \beta^2_\nu \alpha_\nu=l_\nu \tau_\nu \alpha_\nu+\alpha _\nu w_n.\]
The second term converges to zero in probability because \(\alpha_\nu \xrightarrow[]{a.s.} 0, w_n \Rightarrow c_\nu N(0,2\sigma_\nu),\) while for the first
\[0 \leq l_\nu \tau_\nu \alpha_\nu \leq \frac{2c_\nu n}{\sqrt{M}} \cdot c(\epsilon_0)\frac{M}{n} \cdot \alpha_\nu=2c_\nu \cdot c(\epsilon_0) \cdot \sqrt{M}\alpha_\nu.\]
To conclude, we prove \(\sqrt{M}\alpha_\nu \xrightarrow[]{p} 0:\) note that from the analysis in Proposition~\ref{prop2},
\begin{equation}\label{82}
    \alpha_\nu \leq c(\epsilon_0) \cdot (||\frac{1}{n}Z_AZ_A^T-I||+||T\mathcal{M}(\hat{l}_\nu I- \mathcal{M})^{-1}T^T||+|\frac{\hat{l}_\nu}{l_\nu}-1|),
\end{equation}
and under multiplication by \(\sqrt{M}\) each of the three terms on the right tend to zero in probability:
\vspace{0.2cm}
\par
\(1.\) First term: inequality (\ref{7}) yields
\[\frac{1}{n}(\sqrt{n}-CK^2(\sqrt{M}+(n/M)^{1/4}))^2 \leq \lambda_{\min}(\frac{1}{n}Z_AZ_A^T) \leq \lambda_{\max}(\frac{1}{n}Z_AZ_A^T) \leq \frac{1}{n}(\sqrt{n}+CK^2(\sqrt{M}+(n/M)^{1/4}))^2\]
with probability at least \(1-4\exp(-\sqrt{n/M}).\) Since
\[|(1 \pm CK^2(\sqrt{M/n}+(nM)^{-1/4}))^2-1| \leq 3CK^2(\sqrt{M/n}+(nM)^{-1/4}) \leq \frac{t}{\sqrt{M}},\]
for all \(t>0,\) and \(n \geq n(t),\) it follows that \(\sqrt{M} \cdot ||\frac{1}{n}Z_AZ_A^T-I|| \xrightarrow[]{p} 0.\)
\vspace{0.2cm}
\par
\(2.\) Second term: for \(n\) sufficiently large, almost surely
\[\sqrt{M} \cdot ||T\mathcal{M}(\hat{l}_\nu I- \mathcal{M})^{-1}T^T|| \leq \sqrt{M} \cdot ||\frac{1}{n}Z_AZ_A^T|| \cdot ||\mathcal{M}(\hat{l}_\nu I- \mathcal{M})^{-1}|| \leq \sqrt{M} \cdot \frac{4c_{K,\gamma}}{l_\nu} \leq \frac{2c_{K,\gamma}}{c_\nu} \cdot \frac{M}{n},\]
and thus \(\sqrt{M} \cdot ||T\mathcal{M}(\hat{l}_\nu I- \mathcal{M})^{-1}T^T|| \xrightarrow[]{p} 0.\)
\vspace{0.2cm}
\par
\(3.\) Third term: from (\ref{12}),
\[|\frac{\hat{l}_\nu}{l_\nu}-1| \leq |\frac{\lambda_\nu(\mathcal{A})}{l_\nu}-1|+\frac{\lambda_1(\mathcal{B})}{l_\nu}.\]
As \(||\mathcal{B}|| \leq c_{K,\gamma}\) almost surely, \(\sqrt{M} \cdot \frac{\lambda_1(\mathcal{B})}{l_\nu} \leq \frac{2c_{K,\gamma}}{c_\nu} \cdot \frac{M}{n},\)
while (\ref{781}) gives
\[|\frac{\lambda_\nu(\mathcal{A})}{l_\nu}-1| \leq |\lambda_\nu(\frac{1}{n}\sum_{1 \leq i \leq \nu}{Z_i^TZ_i})-1|+|\lambda_1(\frac{1}{n}\sum_{\nu \leq i \leq M}{Z_i^TZ_i})-1|;\]
reasoning as above for the first term, it follows that upon multiplication with \(\sqrt{M}\) these upper bounds tend to zero in probability, and so the third term does as well.

\vspace{0.5cm} \textit{Proof of (\ref{45}):} Recall that
\[l_\nu\beta_\nu^2-\frac{l_\nu}{n}\sum_{k \ne \nu}{\frac{l_kl_\nu}{(l_k-l_\nu)^2}}=\frac{l_\nu}{n}\sum_{k \ne \nu}{\frac{l_kl_\nu}{(l_k-l_\nu)^2}[(\frac{1}{\sqrt{n}}z_k^Tz_\nu)^2-1]}+l_\nu\sum_{k \ne \nu}{\frac{l_kl_\nu}{(l_k-l_\nu)^2}(t_k^T \mathcal{M}(\hat{l}_\nu I-\mathcal{M})^{-1} t_\nu)^2}+\]
\[+2l_\nu\sum_{k \ne \nu}{\frac{l_kl_\nu}{(l_k-l_\nu)^2} \cdot \frac{1}{n}z_k^Tz_\nu \cdot t_k^T \mathcal{M}(\hat{l}_\nu I-\mathcal{M})^{-1} t_\nu}.\]
\par
\textit{Second term:} Inequality (\ref{909}) yields this sum tends to zero in probability as for \(n\) sufficiently large,
\[l_\nu\sum_{k \ne \nu}{\frac{l_kl_\nu}{(l_k-l_\nu)^2}(t_k^T \mathcal{M}(\hat{l}_\nu I-\mathcal{M})^{-1} t_\nu)^2} \leq 2c_\nu \cdot  c(\epsilon_0) \cdot \frac{n}{\sqrt{M}}\sum_{k \ne \nu}{(t_k^T \mathcal{M}(\hat{l}_\nu I-\mathcal{M})^{-1} t_\nu)^2} \xrightarrow[]{p} 0,\]
and
\[\frac{t}{\sqrt{M}} \cdot l_\nu^2 \geq \frac{tc_\nu^2n^2}{4M^{3/2}} \geq 2, \hspace{0.5cm} Ml_\nu^4 \cdot (\frac{t}{\sqrt{M}})^2=t^2l_\nu^4, \hspace{0.5cm} Ml_\nu^2 \cdot \frac{t}{\sqrt{M}} \geq tl_\nu^2.\]
\par
\textit{Third term:} Cauchy-Schwarz inequality implies the square of this sum is upper bounded by
\[4l^2_\nu \sum_{k \ne \nu}{\frac{l_kl_\nu}{(l_k-l_\nu)^2}(t_k^T \mathcal{M}(\hat{l}_\nu I-\mathcal{M})^{-1} t_\nu)^2} \cdot \sum_{k \ne \nu}{\frac{l_kl_\nu}{(l_k-l_\nu)^2} \cdot (\frac{1}{n}z_k^Tz_\nu)^2}.\]
Again using (\ref{909}), it can be shown that
\[\sqrt{M} l_\nu\sum_{k \ne \nu}{\frac{l_kl_\nu}{(l_k-l_\nu)^2}(t_k^T \mathcal{M}(\hat{l}_\nu I-\mathcal{M})^{-1} t_\nu)^2} \xrightarrow[]{p} 0,\]
because for \(n\) sufficiently large,
\[Ml_\nu^4 \cdot (\frac{t}{M})^2=t^2 \cdot \frac{l_\nu^4}{M} \geq \frac{t^2c_\nu^2n^4}{4M^3} \geq t^2n, \hspace{0.5cm} Ml_\nu^2 \cdot \frac{t}{M} \geq tl_\nu^2, \hspace{0.5cm} \frac{t}{M} \cdot l_\nu^2 \geq \frac{tc_\nu^2n^2}{4M^{2}} \geq 2,\]
and since
\[\mathbb{E}[\frac{l_\nu}{\sqrt{M}}\sum_{k \ne \nu}{\frac{l_kl_\nu}{(l_k-l_\nu)^2} \cdot (\frac{1}{n}z_k^Tz_\nu)^2}]=\frac{l_\nu}{n\sqrt{M}}\sum_{k \ne \nu}{\frac{l_kl_\nu}{(l_k-l_\nu)^2} \leq c(\epsilon_0) \cdot \frac{\sqrt{M}l_\nu}{n}} \leq 2c_\nu \cdot c(\epsilon_0),\]
the third term in the decomposition above tends to zero in probability employing the simple fact (\ref{42}).

\vspace{0.5cm} \textit{Proof of (\ref{46}):}
We use a CLT for martingales from Billingsley~\cite{billingsley} (theorem \(35.12\)) stated below:

\vspace{0.5cm}
Let for each \(n \in \mathbb{N},\) \((Y_{ni})_{1 \leq i \leq k_n}\) be a family of martingale differences for the filtration \(\mathcal{F}_{n1} \subset \mathcal{F}_{n2} \subset ... \subset \mathcal{F}_{nk_n}\) that have finite second moments. Suppose that for some \(\sigma>0,\) as \(n \to \infty,\)
\[\sum_{1 \leq i \leq k_n}{\mathbb{E}[Y_{ni}^2|\mathcal{F}_{n(i-1)}]} \xrightarrow[]{p} \sigma^2,\]
and for any \(\epsilon>0,\)
\[\sum_{1 \leq i \leq k_n}{\mathbb{E}[Y_{ni}^2 \cdot \chi_{|Y_{ni}| \geq \epsilon}]} \to 0.\]
Then 
\[\sum_{1 \leq i \leq k_n}{Y_{ni}} \Rightarrow N(0,\sigma^2).\]
\vspace{0.5cm}

\par
For simplicity, denote by 
\[c_{k}=\frac{l_kl_\nu}{(l_k-l_\nu)^2}, k \ne \nu;\] 
consider the nested sigma algebras
\[\sigma(z_\nu),\sigma(z_\nu,z_l, l \ne \nu, l \leq k), k \ne \nu, 1 \leq k \leq M,\]
denote them by
\[\mathcal{F}_{n1} \subset \mathcal{F}_{n2} \subset ... \subset \mathcal{F}_{nM},\]
and take the martingale differences
\[Y_{n1}=\mathbb{E}[T_n|\mathcal{F}_{n1}], Y_{ni}=\mathbb{E}[T_n|\mathcal{F}_{ni}]-\mathbb{E}[T_n|\mathcal{F}_{n(i-1)}], 2 \leq i \leq M.\]
\par
Then
\[\sum_{1 \leq i \leq M}{\mathbb{E}[Y^2_{ni}]}=\mathbb{E}[(\sum_{1 \leq i \leq M}{Y_{ni}})^2]=\mathbb{E}[T_n^2] \to 2\sigma_\nu>0\]
since
\[\mathbb{E}[T_n^2]=\frac{1}{M}\sum_{k \ne \nu}{c_k^2(\mathbb{E}[(\frac{1}{\sqrt{n}}z_k^Tz_\nu)^2-1)^2]-2)}+\frac{2}{M}\sum_{k \ne \nu}{c_k^2}+\frac{1}{M}\sum_{k \ne l, k \ne \nu,l \ne \nu}{c_kc_l\cdot\frac{1}{n}(\mathbb{E}[z_{11}^4]-1)},\]
\[|\mathbb{E}[T_n^2]-\frac{2}{M}\sum_{k \ne \nu}{c_k^2}| \leq \frac{c(K)}{Mn}\sum_{k \ne \nu}{c_k^2}+\frac{1}{M} \cdot M^2c(\epsilon_0) \cdot \frac{c(K)}{n} \leq 2c(K)c(\epsilon_0) \cdot \frac{M}{n} \to 0,\]
as by conditioning on \(z_k,z_l\) for \(k \ne l\)
\[\mathbb{E}[((\frac{1}{\sqrt{n}}z_{k}^Tz_\nu)^2-1) \cdot ((\frac{1}{\sqrt{n}}z_{l}^Tz_\nu)^2-1)]=\mathbb{E}[(\frac{1}{n}||z_\nu||^2-1)^2]=\frac{1}{n^2}\mathbb{E}[||z_\nu||^4]-1=\frac{1}{n}(\mathbb{E}[z_{11}^4]-1),\]
and
\[\mathbb{E}[((\frac{1}{\sqrt{n}}z_k^Tz_\nu)^2-1)^2]=\mathbb{E}[\frac{1}{n^2}(z_k^Tz_\nu)^4]-1=\frac{1}{n}\mathbb{E}([z_{11}^4])^2+\frac{1}{n^2} \cdot 6\binom{n}{2}-1 \leq 2+\frac{c(K)}{n}.\]
\par
Thus, it suffices to show that
\begin{equation}\label{5087}
    \sum_{1 \leq i \leq M}{(\mathbb{E}[Y_{ni}^2|\mathcal{F}_{n(i-1)}]-\mathbb{E}[Y_{ni}^2])} \xrightarrow[]{p} 0,
\end{equation}
\begin{equation}\label{5088}
    \sum_{1 \leq i \leq M}{\mathbb{E}[Y_{ni}^2 \cdot \chi_{|Y_{ni}| \geq \epsilon}]} \to 0.
\end{equation}
\par
Notice that
\[\mathbb{E}[T_n|\mathcal{F}_{n1}]=\frac{1}{\sqrt{M}}\sum_{k \ne \nu}{c_k} \cdot (\frac{1}{n}||z_\nu||^2-1),\]
\[\mathbb{E}[T_n|\mathcal{F}_{ni}]=\frac{1}{\sqrt{M}}\sum_{k \in S_i}{c_k((\frac{1}{\sqrt{n}}z_\nu^Tz_k)^2-1)}+\frac{1}{\sqrt{M}}\sum_{k \ne \nu, k \not \in S_i}{c_i(\frac{1}{n}||z_\nu||^2-1)}, \hspace{0.1cm} 2 \leq i \leq M,\]  
for \(S_i\) the set of the smallest \(i-1\) elements of \(\{k, k \ne \nu\},\) from which 
\begin{equation}\label{725}
    Y_{n1}=\frac{1}{\sqrt{M}}\sum_{k \ne \nu}{c_k} (\frac{1}{n}||z_\nu||^2-1), Y_{ni}=\frac{c_{\alpha(i)}}{\sqrt{M}}[((\frac{1}{\sqrt{n}}z_\nu^Tz_{\alpha(i)})^2-1)-(\frac{1}{n}||z_\nu||^2-1)], \hspace{0.1cm} 2 \leq i \leq M,
\end{equation}
for \(\alpha(i)=i-1, \hspace{0.1cm} 2 \leq i \leq \nu,\) and \(\alpha(i)=i,  \hspace{0.1cm}\nu < i \leq M.\)
\par
We proceed with (\ref{5087}). Using that \(z_l, l \ne \nu\) are i.i.d. and independent of \(z_\nu,\) it follows that for \(2 \leq k \leq M,\)
\[\mathbb{E}[Y_{nk}^2|\mathcal{F}_{n(k-1)}]=\frac{c_{\alpha(k)}^2}{M} \cdot (\mathbb{E}[((\frac{1}{\sqrt{n}}z_\nu^Tz_{\nu-1})^2-1)^2|z_\nu]-(\frac{1}{n}||z_\nu||^2-1)^2)=\]
\[=\frac{c_{\alpha(k)}^2}{M} \cdot (\mathbb{E}[((\frac{1}{\sqrt{n}}z_\nu^Tz_{\nu-1})^2)^2|z_\nu]+1-2 \cdot \frac{1}{n}||z_\nu||^2-(\frac{1}{n}||z_\nu||^2-1)^2)=\]
\[=\frac{c_{\alpha(k)}^2}{M}(\frac{\mathbb{E}[z_{11}^4]}{n^2}\sum_{1 \leq i \leq n}{z_{\nu i}^4}+\frac{6}{n^2}\sum_{1 \leq i<j \leq n}{z_{\nu i}^2z_{\nu j}^2}-(\frac{1}{n}||z_\nu||^2)^2),\]
which gives 
\[|\sum_{2 \leq k \leq M}{(\mathbb{E}[Y_{nk}^2|\mathcal{F}_{n(k-1)}]-\mathbb{E}[Y_{nk}^2])}| \leq c(\epsilon_0) \cdot \frac{\mathbb{E}[z_{11}^4]}{n^2}|\sum_{1 \leq i \leq n}{(z_{\nu i}^4-\mathbb{E}[z_{11}^4])}|+\]
\[+c(\epsilon_0) \cdot \frac{6}{n^2}|\sum_{1 \leq i<j \leq n}{(z_{\nu i}^2z_{\nu j}^2-1)}|+ c(\epsilon_0) \cdot |\frac{1}{n^2}||z_\nu||^4-\mathbb{E}[\frac{1}{n^2}||z_\nu||^4]|.\]
Consider each term in this sum:
\vspace{0.2cm}
\[1. \hspace{0.2cm} \frac{\mathbb{E}[z_{11}^4]}{n^2}\sum_{1 \leq i \leq n}{(z_{\nu i}^4-\mathbb{E}[z_{11}^4])} \xrightarrow[]{a.s.} 0,\]
as its second moment is bounded from above by \(\frac{c(K)}{n^3}\) since \(\mathbb{E}[z_{11}^4],\mathbb{E}[z_{11}^8] \leq c(K),\) and the independence of \(z_{\nu i}, 1 \leq i \leq n;\)

\vspace{0.2cm}
\[2. \hspace{0.2cm} \frac{6}{n^2}\sum_{1 \leq i<j \leq n}{(z_{\nu i}^2z_{\nu j}^2-1)}= 3(\frac{1}{n}\sum_{1 \leq i \leq n}{z_{\nu i}^2})^2-3+\frac{3}{n}-\frac{3}{n^2}\sum_{1 \leq i \leq n}{z_{\nu i}^4} \xrightarrow[]{p} 0,\]
using that the expectation of the last term is bounded from above by \(\frac{c(K)}{n}\) and
\[3(\frac{1}{n}\sum_{1 \leq i \leq n}{z_{\nu i}^2})^2-3 \xrightarrow[]{a.s.} 0\]
from Hanson-Wright inequality (\ref{9});

\vspace{0.2cm}
\[3. \hspace{0.2cm} \frac{1}{n^2}||z_\nu||^4-\mathbb{E}[\frac{1}{n^2}||z_\nu||^4]=(\frac{1}{n}||z_\nu||^2)^2-1-\frac{1}{n}(\mathbb{E}[z_{11}^4]-1) \xrightarrow[]{a.s.} 0.\] 

\vspace{0.2cm}
\par
In what regards (\ref{5088}), notice that for \(2 \leq k \leq M,\)
\[Y_{nk}^2 \leq \frac{c(\epsilon_0)}{M} \cdot (((\frac{1}{\sqrt{n}}z_\nu^Tz_{\alpha(k)})^2-1)-(\frac{1}{n}||z_\nu||^2-1))^2,\]
from which
\begin{equation}\label{5089}
    \sum_{2 \leq k \leq M}{\mathbb{E}[Y_{nk}^2\chi_{|Y_{nk}| \geq \epsilon}]} \leq c(\epsilon_0)\mathbb{E}[H_{n}^2\chi_{|H_{n}| \geq \epsilon \sqrt{\frac{M}{c(\epsilon_0)}}}],
\end{equation}
for 
\[H_n=(\frac{1}{\sqrt{n}}z_1^Tz_2)^2-\frac{1}{n}||z_1||^2.\]
Using \((a+b)^4 \leq 8(a^4+b^4),\)
\[H_n^4 \leq 8(\frac{1}{\sqrt{n}}z_1^Tz_2)^8+8(\frac{1}{n}||z_1||^2)^4,\]
which together with \(\frac{1}{\sqrt{n}}z_1^Tz_2 \Rightarrow N(0,1)\) (employing Lindeberg's CLT as in \ref{4.1}) and (\ref{726}) yields for \(n\) large enough
\[\mathbb{E}[H_n^4] \leq 8 \cdot 2 \cdot 8!!+8 \cdot \mathbb{E}[(\frac{1}{n}||z_1||^2)^4] \leq 8 \cdot 2 \cdot 8!!+8 \cdot (8+8\mathbb{E}[(\frac{1}{n}||z_1||^2-1)^4]) \leq 20 \cdot 8!!,\]
implying
\begin{equation}\label{5091}
    \mathbb{E}[H_{n}^2\chi_{|H_{n}| \geq \epsilon \sqrt{\frac{M}{c(\epsilon_0)}}}] \leq \frac{\mathbb{E}[H_n^4]}{\frac{M\epsilon^2}{c(\epsilon_0)}} \leq \frac{c(\epsilon_0)}{\epsilon^2} \cdot \frac{20  \cdot 8!!}{M} \to 0.
\end{equation}
Lastly, \(|Y_{n1}| \leq c(\epsilon_0) \cdot \sqrt{M} \cdot |\frac{1}{n}||z_1||^2-1|,\) (\ref{725}) and (\ref{726}) give
\begin{equation}\label{5090}
    \mathbb{E}[Y_{n1}^2\chi_{|Y_{n1}|\geq \epsilon}] \leq \frac{\mathbb{E}[Y_{n1}^4]}{\epsilon^2} \leq \frac{ (c(\epsilon_0))^4 \cdot M^2}{\epsilon^2} \cdot \frac{c(K)}{n^2} \to 0.
\end{equation}
Hence (\ref{5088}) ensues from (\ref{5089}), (\ref{5091}), and (\ref{5090}). 

\subsection{Proof of Theorem~\ref{th6}}\label{6.5}

Proposition~\ref{prop1} yields the claims of the theorem for \((a),(b)(i)\) because if \(a_n,b_n,c_n\) are random variables defined on the same space with 
\[a_n \xrightarrow[]{a.s} 1, \hspace{0.2cm} b_n \Rightarrow b, \hspace{0.2cm} c_n \xrightarrow[]{a.s.} c \in \mathbb{R},\] 
then \(a_nb_n \Rightarrow b\) (from Slutsky's lemma), and \(a_nc_n \xrightarrow[]{a.s.} c.\) Hence it suffices to consider \((b)(ii)\) and \((c)(ii).\) As the former is simpler than the latter, we prove them in reverse order.

\par
\vspace{0.2cm}
\((c)(ii)\) Theorem~\ref{th5}, Proposition~\ref{prop1}, and Slutsky's lemma give
\[\hat{l}_\nu(1-<p_\nu,u_\nu>^2)-\frac{\hat{l}_\nu}{n}\sum_{k \ne \nu}{\frac{l_kl_\nu}{(l_k-l_\nu)^2}} \Rightarrow \gamma+c_\nu N(0,2\sigma_\nu),\]
and by applying Slutsky's lemma anew, it suffices to show
\[\frac{\hat{l}_\nu}{n}\sum_{k \ne \nu}{\frac{l_kl_\nu}{(l_k-l_\nu)^2}}-\frac{\hat{l}_\nu}{n}\sum_{k \ne \nu}{\frac{\hat{l}_k\hat{l}_\nu}{(\hat{l}_k-\hat{l}_\nu)^2}} \xrightarrow[]{p} 0\]
(the second sum is defined almost surely reasoning as in the proof of Theorem~\ref{th1}), which in turn is equivalent, due to Proposition~\ref{prop1}, to
\begin{equation}\label{700}
    \frac{1}{\sqrt{M}} \cdot (\sum_{k \ne \nu}{\frac{l_kl_\nu}{(l_k-l_\nu)^2}}-\sum_{k \ne \nu}{\frac{\hat{l}_k\hat{l}_\nu}{(\hat{l}_k-\hat{l}_\nu)^2}}) \xrightarrow[]{p} 0.
\end{equation}
Similarly to the argument for Theorem~\ref{th1}, we use two intermediary steps to obtain (\ref{700}): 
\begin{equation}\label{47}
    \frac{1}{\sqrt{M}} \cdot (\sum_{k \ne \nu}{\frac{l_kl_\nu}{(l_k-l_\nu)^2}}-\sum_{k \ne \nu}{\frac{l_k\hat{l}_\nu}{(l_k-\hat{l}_\nu)^2}}) \xrightarrow[]{p} 0,    
\end{equation}
\begin{equation}\label{48}
    \frac{1}{\sqrt{M}} \cdot (\sum_{k \ne \nu}{\frac{l_k\hat{l}_\nu}{(l_k-\hat{l}_\nu)^2}}-\sum_{k \ne \nu}{\frac{\hat{l}_k\hat{l}_\nu}{(\hat{l}_k-\hat{l}_\nu)^2}}) \xrightarrow[]{p} 0,   
\end{equation}
with these new sums almost surely well-defined.
\par
We begin with the proof of (\ref{47}). Since 
\begin{equation}\label{5093}
    \frac{ax}{(a-x)^2}-\frac{ay}{(a-y)^2}=\frac{a(x-y)(a^2-xy)}{(a-x)^2 \cdot (a-y)^2},
\end{equation}
it follows that for \(k \ne \nu,\)
\[\frac{l_kl_\nu}{(l_k-l_\nu)^2}-\frac{l_k\hat{l}_\nu}{(l_k-\hat{l}_\nu)^2}=\frac{l_k(l_\nu-\hat{l}_\nu)(l_k^2-l_\nu \hat{l}_\nu)}{(l_k-l_\nu)^2 \cdot (l_k-\hat{l}_\nu)^2}=(1-\frac{\hat{l}_\nu}{l_\nu}) \cdot \frac{\frac{l_k}{l_\nu} \cdot ((\frac{l_k}{l_\nu})^2-\frac{\hat{l}_\nu}{l_\nu})}{(\frac{l_k}{l_\nu}-1)^2 \cdot (\frac{l_k}{l_\nu}-\frac{\hat{l}_\nu}{l_\nu})^2},\]
from which Proposition~\ref{prop1} and \(\frac{l_k}{l_\nu} \in [0, \frac{1}{1+\epsilon_0}] \cup [1+\epsilon_0,\infty)\) yield with probability one
\begin{equation}\label{230}
    |\frac{l_kl_\nu}{(l_k-l_\nu)^2}-\frac{l_k\hat{l}_\nu}{(l_k-\hat{l}_\nu)^2}| \leq c(\epsilon_0) \cdot |\frac{\hat{l}_\nu}{l_\nu}-1|,
\end{equation}
\[|\frac{1}{\sqrt{M}} \cdot (\sum_{k \ne \nu}{\frac{l_kl_\nu}{(l_k-l_\nu)^2}}-\sum_{k \ne \nu}{\frac{l_k\hat{l}_\nu}{(l_k-\hat{l}_\nu)^2}})| \leq c(\epsilon_0) \cdot \sqrt{M} \cdot |\frac{\hat{l}_\nu}{l_\nu}-1| \xrightarrow[]{p} 0,\]
where we have used the proof of (\ref{82}) for the last convergence.
\par
We continue with (\ref{48}): (\ref{5093}) gives for \(k \ne \nu,\)
\begin{equation}\label{701}
    \frac{l_k\hat{l}_\nu}{(l_k-\hat{l}_\nu)^2}-\frac{\hat{l}_k\hat{l}_\nu}{(\hat{l}_k-\hat{l}_\nu)^2}=\frac{\hat{l}_\nu (l_k-\hat{l}_k)(\hat{l}_\nu^2-l_k\hat{l}_k)}{(l_k-\hat{l}_\nu)^2 \cdot (\hat{l}_k-\hat{l}_\nu)^2}=(1-\frac{\hat{l}_k}{l_k}) \cdot
\frac{\hat{l}_\nu l_k (\hat{l}_\nu^2-l_k\hat{l}_k)}{(l_k-\hat{l}_\nu)^2 \cdot (\hat{l}_k-\hat{l}_\nu)^2}.
\end{equation}
Arguing as for (\ref{104}), with probability one,
\begin{equation}\label{211}
    |\frac{\hat{l}_\nu l_k (\hat{l}_\nu^2-l_k\hat{l}_k)}{(l_k-\hat{l}_\nu)^2 \cdot (\hat{l}_k-\hat{l}_\nu)^2}|=|\frac{\frac{\hat{l}_\nu}{l_\nu} \cdot \frac{l_k}{l_\nu} \cdot ((\frac{\hat{l}_\nu}{l_\nu})^2-\frac{\hat{l}_k}{l_k} \cdot (\frac{l_k}{l_\nu})^2)}{(\frac{l_k}{l_\nu}-\frac{\hat{l}_\nu}{l_\nu})^2 \cdot (\frac{l_k}{l_\nu} \cdot \frac{\hat{l}_k}{l_k}-\frac{\hat{l}_\nu}{l_\nu})^2}| \leq c(\epsilon_0,K,\gamma) \cdot \min{(\frac{l_k}{l_\nu},\frac{l_\nu}{l_k})},
\end{equation}
while reasoning as for (\ref{106}) yields that for any \(t>0, n \geq n(t),\)
\begin{equation}\label{221}
    \mathbb{P}(|\frac{\hat{l}_k}{l_k}-1| \cdot \min{(\frac{l_k}{l_\nu},\frac{l_\nu}{l_k})} \geq \frac{t}{\sqrt{M}}) \leq 6\exp(-\sqrt{n}),
\end{equation}
(using 
\[\frac{l_\nu}{\sqrt{M}} \geq \frac{c_\nu}{2} \cdot \frac{n}{M} \to \infty, \hspace{0.2cm} \sqrt{M} \cdot (\frac{1}{n^{1/4}}+\sqrt{\frac{M}{n}})=\frac{M^{1/2}}{n^{1/4}}+\frac{M}{n^{1/2}} \to 0.)\]
Finally, (\ref{701}),
(\ref{211}), and (\ref{221}) give for any \(t>0, n \geq n(t),\)
\[\mathbb{P}(\frac{1}{\sqrt{M}} \cdot |\sum_{k \ne \nu}{\frac{l_k\hat{l}_\nu}{(l_k-\hat{l}_\nu)^2}}-\sum_{k \ne \nu}{\frac{\hat{l}_k\hat{l}_\nu}{(\hat{l}_k-\hat{l}_\nu)^2}}| \geq t) \leq \sum_{k \ne \nu}{\mathbb{P}(|\frac{\hat{l}_k}{l_k}-1| \cdot \min{(\frac{l_k}{l_\nu},\frac{l_\nu}{l_k})} \geq \frac{t}{c(\epsilon_0,K,\gamma)\sqrt{M}})},\]
upper bounded by \(6M\exp(-\sqrt{n}),\) from which (\ref{48}) follows.

\par
\vspace{0.2cm}
\((b)(ii)\) The arguments are analogous to the ones employed in \((c)(ii):\) it suffices to show
\[\frac{\hat{l}_\nu}{n}\sum_{k \ne \nu}{\frac{l_kl_\nu}{(l_k-l_\nu)^2}}-\frac{\hat{l}_\nu}{n}\sum_{k \ne \nu}{\frac{\hat{l}_k\hat{l}_\nu}{(\hat{l}_k-\hat{l}_\nu)^2}} \xrightarrow[]{p} 0,\]
which we split into
\begin{equation}\label{227}
    \frac{1}{M} \cdot (\sum_{k \ne \nu}{\frac{l_kl_\nu}{(l_k-l_\nu)^2}}-\sum_{k \ne \nu}{\frac{l_k\hat{l}_\nu}{(l_k-\hat{l}_\nu)^2}}) \xrightarrow[]{p} 0,    
\end{equation}
\begin{equation}\label{228}
    \frac{1}{M} \cdot (\sum_{k \ne \nu}{\frac{l_k\hat{l}_\nu}{(l_k-\hat{l}_\nu)^2}}-\sum_{k \ne \nu}{\frac{\hat{l}_k\hat{l}_\nu}{(\hat{l}_k-\hat{l}_\nu)^2}}) \xrightarrow[]{p} 0. 
\end{equation}
(\ref{227}) is immediate from (\ref{230}) and Proposition~\ref{prop1}, while (\ref{228}) is an easy consequence of the analogue of (\ref{221}):
\begin{equation}\label{231}
    \mathbb{P}(|\frac{\hat{l}_k}{l_k}-1| \cdot \min{(\frac{l_k}{l_\nu},\frac{l_\nu}{l_k})} \geq t) \leq 6\exp(-\sqrt{n}),
\end{equation}
for any \(t>0,n \geq n(t).\)

\section{Auxiliary Lemmas}\label{sec7}

In virtue of the approach taken in our proofs, the random quantities
\[R_\nu^2, \hspace{0.1cm} ||\mathcal{R}_\nu \mathcal{D}_\nu e_\nu||,||a_\nu-e_\nu||, \hspace{0.1cm} \sum_{k \ne \nu}{(\frac{1}{n}z_k^Tz_\nu)^2}, \hspace{0.1cm} \sum_{k \ne \nu}{(t_k^T\mathcal{M}(\hat{l}_\nu I-\mathcal{M})^{-1}t_\nu)^2},\]
where \(z_1, \hspace{0.05cm} ... \hspace{0.05cm}, z_M \in \mathbb{R}^n\) denote the rows of \(Z_A,\) and \(t_1, \hspace{0.05cm}... \hspace{0.05cm}, t_M \in \mathbb{R}^{N-M}\) the columns of \(T,\) govern both the fluctuations of the eigenvalues and the consistency rates of eigenvectors. Consequently, we collect some results regarding their behavior that play a central role in our analysis.
\par
Lemma~\ref{lemma1} is primarily concerned with \(\Tilde{T}=T^T\mathcal{M}(\hat{l}_\nu I-\mathcal{M})^{-1}T,\) an essential matrix for Theorem~\ref{th2} (see subsection (\ref{3'.1})) as well as for Theorem~\ref{th5}: its justification relies on switching from \((\hat{l}_\nu I-\mathcal{M})^{-1}\) to \((l_\nu I-\mathcal{M})^{-1}\) since this allows us to employ concentration inequalities for sums depending on the entries of
\[T^T\mathcal{M}(l_\nu I-\mathcal{M})^{-1}T=\frac{1}{n} Z_AH\mathcal{M}(l_\nu I-\mathcal{M})^{-1}H^TZ_A^T,\]
the key observation being that \(Z_A\) and \(H\mathcal{M}(l_\nu I-\mathcal{M})^{-1}H^T\) are independent. Lemma~\ref{lemma2} relates \(||a_\nu-e_\nu||\) to \(||\mathcal{R}_\nu \mathcal{D}_\nu e_\nu||,\) an object more amenable to scrutiny (anew \(\tilde{T}\) comes into play as it directly influences \(||\mathcal{R}_\nu \mathcal{D}_\nu e_\nu||\)): its proof follows from a decomposition of \(a_\nu-e_\nu\) from Paul~\cite{paul}. Lemma~\ref{lemma3}, used repeatedly in Theorem~\ref{th5}, renders the rate of \(R_\nu^2 \to 0\) explicit: its justification is a careful analysis of (\ref{5}). 
\par
We state next these results and present their proofs in subsections \ref{7.1}, \ref{7.2}, and \ref{7.3}, respectively.

\begin{lemma}\label{lemma1}
Under the assumptions of Proposition~\ref{prop1}, as \(n \to \infty,\)

\vspace{0.2cm}
\[(a) \hspace{0.2cm} l_\nu^4\sum_{k \ne \nu}{(t_k^T\mathcal{M}(\hat{l}_\nu I-\mathcal{M})^{-2}t_\nu)^2} \xrightarrow[]{a.s.} 0, \hspace{0.5cm} (b) \hspace{0.2cm} \frac{n}{M}\sum_{k \ne \nu}{(t_k^T\mathcal{M}(\hat{l}_\nu I-\mathcal{M})^{-1}t_\nu)^2} \xrightarrow[]{p} 0.\]
\end{lemma}

\begin{lemma}\label{lemma2}
Under the assumptions of Proposition~\ref{prop2}, and \(<p_{A,\nu},u_{A,\nu}> \geq 0,\) for \(\beta_\nu=||\mathcal{R}_\nu \mathcal{D}_\nu e_\nu||,\)
as \(n \to \infty,\)
\[\frac{||a_\nu-e_\nu||}{\beta_\nu} \xrightarrow[]{a.s.} 1, \]
where if \(\beta_\nu=0\) for \(n\) sufficiently large, then \(||a_\nu-e_\nu||=0\) as well. Moreover, 
\[\sum_{n \in \mathbb{N}}{\mathbb{P}(||a_\nu-e_\nu|| \geq 2\beta_\nu)}<\infty.\]
\end{lemma}

\begin{lemma}\label{lemma3}
Under the assumptions of Proposition~\ref{prop2}, as \(n \to \infty,\)
\[\hat{l}_\nu R_\nu^2-\frac{N}{n} \xrightarrow[]{a.s.} 0.\]
\end{lemma}

\subsection{Proof of Lemma~\ref{lemma1}}\label{7.1}

\((a)\) Recall the identity
\[\frac{1}{(\hat{l}_\nu-x)^2}=\frac{1}{(l_\nu-x)^2} \cdot \frac{1}{(1-\frac{l_\nu-\hat{l}_\nu}{l_\nu-x})^2}=\sum_{m \geq 0}{\frac{(m+1)(l_\nu-\hat{l}_\nu)^m}{(l_\nu-x)^{m+2}}}\]
for \(x \in [0,\frac{l_\nu}{2}], |\hat{l}_\nu-l_\nu|<\frac{l_\nu}{2} \leq |x-l_\nu|,\) and denote by
\[s_k(i)=t_k^T\mathcal{M}(l_\nu I-\mathcal{M})^{-i-2}t_\nu\]
for \(k \ne \nu, i \geq 0.\)
Then as long as \(||\mathcal{M}|| \leq \frac{l_\nu}{2},\) the sum under consideration can be rewritten as
\[l_\nu^4\sum_{k \ne \nu}{(t_k^T\mathcal{M}(\hat{l}_\nu I-\mathcal{M})^{-2}t_\nu)^2}=l_\nu^4 \sum_{m \geq 0}{(l_\nu-\hat{l}_\nu)^m\sum_{k \ne \nu}{\sum_{0 \leq i \leq m}{(i+1)(m-i+1)s_k(i)s_k(m-i)}}};\]
hence, for any \(t>0,\) and \(n\) large enough such that \(l_\nu>2c_{K,\gamma},\)
\[\mathbb{P}(l_\nu^4\sum_{k \ne \nu}{(t_k^T\mathcal{M}(\hat{l}_\nu I-\mathcal{M})^{-2}t_\nu)^2} \geq t) \leq \mathbb{P}(|\hat{l}_\nu-l_\nu|>\frac{l_\nu}{32})+\mathbb{P}(||\mathcal{M}||>c_{K,\gamma})+\mathbb{P}(\frac{1}{n}||z_\nu||^2>2)+\]
\begin{equation}\label{808}
    +\sum_{m \geq 0, 0 \leq i \leq m}{\mathbb{P}(|\sum_{k \ne \nu}{s_k(i)s_k(m-i)}| \geq \frac{2^{m-1}t}{l_\nu^{m+4}} \hspace{0.05cm}| \hspace{0.05cm} ||\mathcal{M}|| \leq c_{K,\gamma}, \frac{1}{n}||z_\nu||^2 \leq 2)},
\end{equation}
because if none of the events on the right-hand side happens, then
\[l_\nu^4\sum_{k \ne \nu}{(t_k^T\mathcal{M}(\hat{l}_\nu I-\mathcal{M})^{-2}t_\nu)^2}=l_\nu^4 \sum_{m \geq 0}{(l_\nu-\hat{l}_\nu)^m\sum_{k \ne \nu}{\sum_{0 \leq i \leq m}{(i+1)(m-i+1)s_k(i)s_k(m-i)}}} \leq\]
\[\leq \sum_{m \geq 0}{\frac{l_\nu^{m+4}}{8^m}\sum_{0 \leq i \leq m}|{\sum_{k \ne \nu}{s_k(i)s_k(m-i)}}}| \leq \sum_{m \geq 0}{\frac{l_\nu^{m+4}}{8^m} \cdot (m+1)\frac{2^{m-1}t}{l_\nu^{m+4}}}=t\sum_{m \geq 0}{\frac{m+1}{2 \cdot 4^m}}=\frac{t}{2} \cdot \frac{1}{(1-\frac{1}{4})^2}<t.\]
To conclude the proof, it suffices to show the four terms on the right-hand side in (\ref{808}) are summable as functions of \(n \in \mathbb{N}\) from Borel-Cantelli lemma.
\par
For the first three components in (\ref{808}), their exponential decay ensues from (\ref{105}) and inequality (\ref{7}) as for \(n\) large enough, \(\frac{l_\nu}{64}>c_{K,\gamma},\) while for the third, Cauchy-Schwarz inequality yields
\[\sum_{m \geq 0, 0 \leq i \leq m}{\mathbb{P}(|\sum_{k \ne \nu}{s_k(i)s_k(m-i)}| \geq \frac{2^{m-1}t}{l_\nu^{m+4}}\hspace{0.05cm}| \hspace{0.05cm} ||\mathcal{M}|| \leq c_{K,\gamma}, \frac{1}{n}||z_\nu||^2 \leq 2)} \leq\]
\begin{equation}\label{809}
    \leq 2\sum_{m \geq 0, 0 \leq i \leq m}{\mathbb{P}(\sum_{k \ne \nu}{s^2_k(i)} \geq \frac{2^{m-1}t}{l_\nu^{2i+4}}\hspace{0.05cm}| \hspace{0.05cm} ||\mathcal{M}|| \leq c_{K,\gamma}, \frac{1}{n}||z_\nu||^2 \leq 2)}.
\end{equation}
We derive an upper bound for the terms in this new sum using Bernstein's inequality (theorem \(2.8.2\) in Vershynin~\cite{vershynin}): for \(t \geq 0\) and \(A \in \mathbb{R}^{n \times n}\) independent of \((z_k)_{k \ne \nu},\)
\[\mathbb{P}(|\sum_{k \ne \nu}{((z_k^tAz_\nu)^2-||Az_\nu||^2)}| \geq t|A,z_\nu) \leq 2\exp(-c\min{(\frac{t^2}{M \cdot ||Az_\nu||^4},\frac{t}{||Az_\nu||^2})}):\]
from proposition \(2.6.1\) in Vershynin~\cite{vershynin}, conditional on \(A,z_\nu,\) 
\[||z_k^tAz_\nu||_{\psi_2} \leq C||Az_\nu||,\]
which implies using lemma \(2.7.6\) that \((z_k^tAz_\nu)^2\) is subexponential with 
\[||(z_k^tAz_\nu)^2||_{\psi_1} \leq C^2||Az_\nu||^2,\]
where for a subexponential random variable \(x,\)
\[||x||_{\psi_1}=\inf{\{t>0: \mathbb{E}[\exp(|x|/t)] \leq 2\}}.\]
Since \(||Az_\nu|| \leq ||A|| \cdot ||z_\nu||,\) this inequality yields
\begin{equation}\label{707}
    \mathbb{P}(|\sum_{k \ne \nu}{((z_k^tAz_\nu)^2-||Az_\nu||^2)}| \geq t|A,z_\nu) \leq 2\exp(-c\min{(\frac{t^2}{M \cdot ||A||^4 \cdot ||z_\nu||^4},\frac{t}{||A||^2 \cdot ||z_\nu||^2})}).
\end{equation}
\par
Thus, for \(\mathcal{M}_{i,\nu}=\frac{l_\nu^{i+2}}{n} H\mathcal{M}(l_\nu I-\mathcal{M})^{-i-2}H^T,\)
\[\mathbb{P}(\sum_{k \ne \nu}{s^2_k(i)} \geq \frac{2^{m-1}t}{l_\nu^{2i+4}}|Z_B, z_\nu) \leq \mathbb{P}(|\sum_{k \ne \nu}{((z_k^T\mathcal{M}_{i,\nu}z_\nu)^2-||\mathcal{M}_{i,\nu}z_\nu||^2)}| \geq 2^{m-2}t|Z_B, z_\nu)+\]
\[+\mathbb{P}(M \cdot ||\mathcal{M}_{i,\nu}z_\nu||^2 \geq 2^{m-2}t|Z_B, z_\nu) \leq 2\exp(-c\min{(\frac{2^{2m-4}t^2}{M \cdot ||\mathcal{M}_{i,\nu}||^4 \cdot ||z_\nu||^4},\frac{2^{m-2}t}{||\mathcal{M}_{i,\nu}||^2 \cdot ||z_\nu||^2}}))+\]
\[+\mathbb{P}(M \cdot ||\mathcal{M}_{i,\nu}z_\nu||^2 \geq 2^{m-2}t|Z_B, z_\nu).\]
Furthermore, if \(||\mathcal{M}|| \leq c_{K,\gamma} \leq l_\nu(1-2^{-\frac{1}{4}}),\) then
\[M \cdot ||\mathcal{M}_{i,\nu}z_\nu||^2 \leq \frac{M}{n^2} \cdot  ||l_\nu^{i+2}\mathcal{M}(l_\nu I-\mathcal{M})^{-i-2}||^2 \cdot ||z_\nu||^2 \leq \frac{M}{n^2} \cdot c_{K,\gamma}^2 \cdot  2^{\frac{i+2}{2}} ||z_\nu||^2\]
as \(\frac{xl^{i+2}}{(l-x)^{i+2}} \leq c_{K,\gamma} \cdot 2^{\frac{i+2}{4}}\) for \(x \in [0,c_{K,\gamma}].\)
\par
Hence for \(n\) large enough such that \(c_{K,\gamma} \leq l_\nu(1-2^{-\frac{1}{4}}),\)
\[\sum_{m \geq 0, 0 \leq i \leq m}{\mathbb{P}(\sum_{k \ne \nu}{s^2_k(i)} \geq \frac{2^{m-1}t}{l_\nu^{2i+4}}|Z_B, z_\nu)} \leq \mathbb{P}(||\mathcal{M}||>c_{K,\gamma}|Z_B, z_\nu)+\mathbb{P}(\frac{M}{n^2} \cdot ||z_\nu||^2 \geq \frac{t}{8c_{K,\gamma}^2}|Z_B, z_\nu)+\]
\[+\sum_{m \geq 0, 0 \leq i \leq m}{2\exp(-c\min{(\frac{2^{m-6}t^2n^2}{Mc_{K,\gamma}^4 \cdot (\frac{1}{n}||z_\nu||^2)^2},\frac{2^{\frac{m-6}{2}}tn}{c_{K,\gamma}^2 \cdot \frac{1}{n}||z_\nu||^2})})},\]
from which for \(n\) sufficiently large so that \(\frac{n}{M} \cdot \frac{t}{8c^2_{K,\gamma}}>2,\) the sum in (\ref{809}) is at most
\[2\sum_{m \geq 0, 0 \leq i \leq m}{2\exp(-c\min{(\frac{2^{m-10}t^2n^2}{Mc_{K,\gamma}^4},\frac{2^{\frac{m-7}{2}}tn}{c_{K,\gamma}^2})})}.\]
Note that this is summable as a function of \(n \in \mathbb{N}\) from the following elementary inequality:
\[\sum_{m \geq 0, 0 \leq i \leq m}{\exp(-\min(2^{m}a,2^{\frac{m}{2}}b))} \leq \sum_{m \geq 0}{(m+1)(\exp(-2^{m}a)+\exp(-2^{\frac{m}{2}}b))} \leq \]
\[\leq \sum_{m \geq 0}{(m+1)(\exp(-(m+1)a)+\exp(-(m+1)\frac{b}{2} ))}=\frac{e^{-a}}{(1-e^{-a})^2}+\frac{e^{-b/2}}{(1-e^{-b/2})^2} \leq 4e^{-a}+4e^{-b/2},\]
if \(e^{-a},e^{-b/2} \leq \frac{1}{2}\) together with \(M \leq N \leq 2\gamma n\) for \(n\) sufficiently large.

\vspace{0.4cm}
\par
\((b)\) A similar approach to the one in \((a)\) can be used here:
\[\frac{1}{\hat{l}_\nu-x}=\frac{1}{l_\nu-x} \cdot \frac{1}{1-\frac{l_\nu-\hat{l}_\nu}{l_\nu-x}}=\sum_{m \geq 0}{\frac{(l_\nu-\hat{l}_\nu)^m}{(l_\nu-x)^{m+1}}}\]
for \(x \in [0,\frac{l_\nu}{2}]\) and \(|\hat{l}_\nu-l_\nu|<\frac{l_\nu}{2} \leq |x-l_\nu|;\) if \(||\mathcal{M}||<\frac{l_\nu}{2}, |\frac{\hat{l}_\nu}{l_\nu}-1|<\frac{1}{2},\) then
\[\frac{n}{M}\sum_{k \ne \nu}{(t_k^T\mathcal{M}(\hat{l}_\nu I-\mathcal{M})^{-1}t_\nu)^2}=\frac{1}{M}\sum_{m \geq 0}{(l_\nu-\hat{l}_\nu)^m\sum_{k \ne \nu}{\sum_{0 \leq i \leq m}{s'_k(i)s'_k(m-i)}}},\]
where \(s'_k(i)=\sqrt{n}t_k^T\mathcal{M}(l_\nu I - \mathcal{M})^{-i-1}t_\nu\) for \(k \ne \nu, i \geq 0;\) 
\[\mathbb{P}(\frac{n}{M}\sum_{k \ne \nu}{(t_k^T\mathcal{M}(\hat{l}_\nu I-\mathcal{M})^{-1}t_\nu)^2} \geq t) \leq \mathbb{P}(|\hat{l}_\nu-l_\nu|>\frac{l_\nu}{8})+\mathbb{P}(||\mathcal{M}||>c_{K,\gamma})+\mathbb{P}(\frac{1}{n}||z_\nu||^2>2)+\]
\[+\sum_{m \geq 0, 0 \leq i \leq m}{\mathbb{P}(|\sum_{k \ne \nu}{s'_k(i)s'_k(m-i)}| \geq \frac{2^{m-1}Mt}{l_\nu^{m}} \hspace{0.05cm}| \hspace{0.05cm} ||\mathcal{M}|| \leq c_{K,\gamma}, \frac{1}{n}||z_\nu||^2 \leq 2)};\]
for \(0 \leq i \leq m,\)
\[\mathbb{P}(\sum_{k \ne \nu}{(s'_k(i))^2} \geq \frac{2^{m-1}Mt}{l_\nu^{2i}}|Z_B, z_\nu) \leq \mathbb{P}(||\mathcal{M}'_{i,\nu}||^2 \cdot ||z_\nu||^2 \geq \frac{2^{m-2}t}{l_\nu^{2i}}|Z_B, z_\nu)+\]
\[+\mathbb{P}(\sum_{k \ne \nu}{((s'_k(i))^2-||\mathcal{M}'_{i,\nu}||^2)} \geq \frac{2^{m-2}Mt}{l_\nu^{2i}}|Z_B, z_\nu),\]
where \(\mathcal{M}'_{i,\nu}=\frac{1}{\sqrt{n}}H\mathcal{M}(l_\nu I-\mathcal{M})^{-i-1}H^T.\)
\par
For the first term in this last inequality, whenever \(||\mathcal{M}|| \leq c_{K,\gamma} \leq l_\nu(1-2^{-1/4}),\)
\[||\mathcal{M}'_{i,\nu}||^2 \leq \frac{1}{n} \cdot \frac{c^2_{K,\gamma}}{(\frac{l_\nu}{2^{1/4}})^{2(i+1)}},\]
from which 
\[\mathbb{P}(||\mathcal{M}'_{i,\nu}||^2 \cdot ||z_\nu||^2 \geq \frac{2^{m-2}t}{l_\nu^{2i}}) \leq \mathbb{P}(\frac{1}{n}||z_\nu||^2 \geq \frac{tl_\nu^2}{8c^2_{K,\gamma}}),\]
while for the second, (\ref{707}) yields the following upper bound
\[2\exp(-c\min{(\frac{2^{2m-4}M^2t^2}{Ml_\nu^{4i} \cdot ||\mathcal{M}'_{i,\nu}||^4 \cdot  ||z_\nu||^4},\frac{2^{m-2}Mt}{l_\nu^{2i} \cdot ||\mathcal{M}'_{i,\nu}||^2 \cdot ||z_\nu||^2)})}) \leq\]
\[\leq 2\exp(-c\min{(\frac{2^{m-5}Mt^2l_\nu^4}{c_{K,\gamma}^4 \cdot (\frac{1}{n}||z_\nu||^2)^2},\frac{2^{\frac{m}{2}-3}Mtl_\nu^2}{c_{K,\gamma}^2 \cdot \frac{1}{n}||z_\nu||^2})}).\]
\par
Hence for \(\frac{tl_\nu^2}{8c^2_{K,\gamma}}>2,\)
\[\mathbb{P}(\frac{n}{M}\sum_{k \ne \nu}{(t_k^T\mathcal{M}(\hat{l}_\nu I-\mathcal{M})^{-1}t_\nu)^2} \geq t) \leq \mathbb{P}(|\hat{l}_\nu-l_\nu|>\frac{l_\nu}{8})+\mathbb{P}(||\mathcal{M}||>c_{K,\gamma})+\mathbb{P}(\frac{1}{n}||z_\nu||^2>2)+\]
\begin{equation}\label{909}
    +c_1\exp(-c_2Mt^2l_v^4)+c_1\exp(-c_2Mtl_\nu^2)
\end{equation}
for some \(c_1>0,c_2=c_2(K,\gamma)>0,\) and the conclusion then ensues since \(l_\nu \to \infty.\)

\subsection{Proof of Lemma~\ref{lemma2}}\label{7.2}

Recall (\ref{13}) and (\ref{14}):
\[a_\nu-e_\nu=-\mathcal{R}_\nu \mathcal{D}_\nu e_\nu+r_\nu,\]
\[r_\nu=(<a_\nu,e_\nu>-1)e_\nu-\mathcal{R}_\nu \mathcal{D}_\nu (a_\nu-e_\nu)+(\hat{l}_\nu-l_\nu)\mathcal{R}_\nu(a_\nu-e_\nu),\]
(as in Proposition~\ref{prop2}, both (\ref{4}) and (\ref{5}) hold almost surely). (\ref{13}) provides
\[\beta_\nu-||r_\nu|| \leq ||a_\nu-e_\nu|| \leq \beta_\nu+||r_\nu||,\]
while from Proposition~\ref{prop1}, \(||\mathcal{R}_\nu D_\nu||+|\hat{l}_\nu-l_\nu| \cdot ||\mathcal{R}_\nu|| \xrightarrow[]{a.s.} 0,\) and \(||a_\nu-e_\nu|| \xrightarrow[]{a.s.} 0\) because 
\[<a_\nu,e_\nu>^2 \xrightarrow[]{a.s.} 1,\] 
and 
\[<a_\nu,e_\nu>=\frac{<p_{A,\nu},u_{A,\nu}>}{\sqrt{1-R_\nu^2}} \geq 0.\]
Since
\(2(1-<a_\nu,e_\nu>)=||a_\nu-e_\nu||^2,\) (\ref{14}) gives that for all \(\epsilon \in (0,\frac{1}{2}),\) almost surely for \(n\) large enough,
\[||r_\nu|| \leq 2||a_\nu-e_\nu||^2+\epsilon ||a_\nu-e_\nu|| \leq 2\epsilon ||a_\nu-e_\nu||,\] from which 
\[\frac{\beta_\nu}{1+2\epsilon} \leq  ||a_\nu-e_\nu|| \leq \frac{\beta_\nu}{1-2\epsilon},\]
yielding the first claim. Furthermore, these inequalities for \(\epsilon=\frac{1}{4}\) imply
\[\mathbb{P}(||a_\nu-u_\nu|| \geq 2||\beta_\nu||)  \leq \mathbb{P}(||a_\nu-e_\nu|| \geq \frac{1}{4})+\mathbb{P}(||\mathcal{R}_\nu \mathcal{D}_\nu||+|\hat{l}_\nu-l_\nu| \cdot ||\mathcal{R}_\nu|| \geq \frac{1}{4}),\]
and using Proposition~\ref{prop2},
\[||a_\nu-e_\nu||^2=2-2\sqrt{1-(\mathcal{P}^{\perp}_\nu a_\nu)^2}, \hspace{0.5cm} ||\mathcal{P}^{\perp}_\nu a_\nu|| \leq ||\mathcal{R}_\nu \mathcal{D}_\nu||+|\hat{l}_\nu-l_\nu| \cdot ||\mathcal{R}_\nu||,\]
from which
\[\mathbb{P}(||a_\nu-e_\nu|| \geq \frac{1}{4}) \leq \mathbb{P}(||\mathcal{R}_\nu \mathcal{D}_\nu||+|\hat{l}_\nu-l_\nu| \cdot ||\mathcal{R}_\nu|| \geq \frac{1}{8})\]
since \(2-2\sqrt{1-x^2} \leq 2x^2,\) for \(x^2 \leq 1.\)
\par
Moreover, Proposition~\ref{prop2} entails that for \(t>0\) 
\[\mathbb{P}(||\mathcal{R}_\nu \mathcal{D}_\nu||+|\hat{l}_\nu-l_\nu| \cdot ||\mathcal{R}_\nu|| \geq t) \leq \mathbb{P}(||\mathcal{R}_\nu \mathcal{D}_\nu|| \geq t/2)+\mathbb{P}(|\hat{l}_\nu-l_\nu| \cdot ||\mathcal{R}_\nu|| \geq t/2) \leq\]
\[\leq \mathbb{P}(||\frac{1}{n}Z_AZ_A^T-I+T^T\mathcal{M}(\hat{l}_\nu I -\mathcal{M})^{-1}T|| \geq tc(\epsilon_0,\gamma))+\mathbb{P}(|\frac{\hat{l}_\nu}{l_\nu}-1| \geq tc(\epsilon_0,\gamma)),\]
with the second terms are summable as functions of \(n \in \mathbb{N}\) from (\ref{105}), and
\[\mathbb{P}(||\frac{1}{n}Z_AZ_A^T-I+T^T\mathcal{M}(\hat{l}_\nu I -\mathcal{M})^{-1}T|| \geq t) \leq \mathbb{P}(||\frac{1}{n}Z_AZ_A^T-I|| \geq \frac{t}{2})+\mathbb{P}(||\frac{1}{n}Z_AZ_A^T|| \geq 2)+\]
\[+\mathbb{P}(||\mathcal{M}(\hat{l}_\nu I -\mathcal{M})^{-1}|| \geq \frac{t}{4})\]
because
\[||T^T\mathcal{M}(\hat{l}_\nu I -\mathcal{M})^{-1}T|| \leq ||\frac{1}{n}Z_AZ_A^T|| \cdot ||\mathcal{M}(\hat{l}_\nu I -\mathcal{M})^{-1}||,\]
which in conjunction with inequality (\ref{7}), Proposition~\ref{prop1}, and \(l_\nu \to \infty\) renders the first terms summable as functions of \(n \in \mathbb{N}\) as well.

\subsection{Proof of Lemma~\ref{lemma3}}\label{7.3}

As in Proposition~\ref{prop2}, equation (\ref{5}) holds almost surely:
\[a^T_\nu \Lambda^{1/2}T^T\mathcal{M}(\hat{l}_\nu I- \mathcal{M})^{-2}T\Lambda^{1/2} a_\nu=\frac{R_\nu^2}{1-R_\nu^2}.\]
\par
We show \(y_\nu-\frac{N}{n}=\frac{\hat{l}_\nu R_\nu^2}{1-R_\nu^2}-\frac{N}{n} \xrightarrow[]{a.s.} 0,\) which implies
\[\hat{l}_\nu R_\nu^2-\frac{N}{n}=y_\nu-R_\nu^2y_\nu-\frac{N}{n} \xrightarrow[]{a.s.} 0\]
because Proposition~\ref{prop2} entails \(R_\nu^2 \xrightarrow[]{a.s.}0.\) Denote by 
\[E_\nu=\hat{l}_\nu \Lambda^{1/2}T^T\mathcal{M}(\hat{l}_\nu I- \mathcal{M})^{-2}T\Lambda^{1/2}.\]
Since 
\[y_\nu = a_\nu^T E_\nu a_\nu,\] 
it suffices in turn to prove
\[e_\nu^TE_\nu e_\nu-\frac{N}{n} \xrightarrow[]{a.s.} 0, \hspace{0.5cm} e_\nu^T E_\nu (a_\nu-e_\nu) \xrightarrow[]{a.s.} 0, \hspace{0.5cm} (a_\nu-e_\nu)^TE_\nu (a_\nu-e_\nu) \xrightarrow[]{a.s.} 0.\]
\par
\(1.\) Consider the first term:
\[e_\nu^TE_\nu e_\nu=\hat{l}_\nu l_\nu \cdot t_\nu ^T \mathcal{M}(\hat{l}_\nu I- \mathcal{M})^{-2}t_\nu=\frac{l_\nu}{\hat{l}_\nu} \cdot t_\nu ^T \mathcal{M}(I-\frac{1}{\hat{l}_\nu}\mathcal{M})^{-2}t_\nu.\]
\par
Hanson-Wright inequality (\ref{9}) yields for any \(t \geq 0,\)
\[\mathbb{P}(|t_\nu ^T \mathcal{M}t_\nu-\frac{1}{n}tr(\mathcal{M})| \geq t|Z_B) \leq 2\exp(-c\min{(\frac{nt^2}{K^4 \cdot ||\mathcal{M}||^2},\frac{nt}{K^2 \cdot ||\mathcal{M}||})})\]
as \(t_\nu ^T \mathcal{M}t_\nu=\frac{1}{n}z_\nu ^T H \mathcal{M}H^Tz_\nu, tr(H\mathcal{M}H^T)=tr(\mathcal{M}),||H\mathcal{M}H^T||=||\mathcal{M}||\) from which
\[\mathbb{P}(|t_\nu ^T \mathcal{M}t_\nu-\frac{1}{n}tr(\mathcal{M})| \geq t) \leq \mathbb{P}(||\mathcal{M}||>c_{K,\gamma})+2\exp(-c\min{(\frac{nt^2}{K^4c_{K,\gamma}^2},\frac{nt}{K^2c_{K,\gamma}})}),\]
implying with Borel-Cantelli lemma that 
\[t_\nu ^T \mathcal{M}t_\nu-\frac{1}{n}tr(\mathcal{M}) \xrightarrow[]{a.s.} 0,\]
and consequently
\begin{equation}\label{913}
    t_\nu ^T \mathcal{M}t_\nu-\frac{N}{n} \xrightarrow[]{a.s.} 0
\end{equation}
because (\ref{726}) gives
\[\frac{1}{n}tr(\mathcal{M})-\frac{N}{n}=\frac{1}{n} tr(\frac{1}{n}Z_BZ_B^T)-\frac{N}{n}=-\frac{M}{n}+\frac{1}{n^2}\sum_{M+1 \leq i \leq N, 1 \leq j \leq n}{(z^2_{ij}-1)} \xrightarrow[]{a.s.} 0:\]
\[\mathbb{E}[(\frac{1}{n}tr(\mathcal{M})-\frac{N-M}{n})^2] \leq \sqrt{\mathbb{E}[(\frac{1}{n}tr(\mathcal{M})-\frac{N-M}{n})^4]} \leq \frac{\sqrt{c(K)}}{n(N-M)} \leq \frac{c(K,\gamma)}{n^2}.\]
Lastly, since almost surely
\[0 \leq t_\nu ^T \mathcal{M}(I-\frac{1}{\hat{l}_\nu}\mathcal{M})^{-2}t_\nu- t_\nu ^T \mathcal{M}t_\nu=t_\nu ^T \mathcal{M}((I-\frac{1}{\hat{l}_\nu}\mathcal{M})^{-2}-I)t_\nu \leq\]
\[\leq ||(I-\frac{1}{\hat{l}_\nu}\mathcal{M})^{-2}-I|| \cdot t_\nu ^T \mathcal{M}t_\nu \leq \frac{8c_{K,\gamma}}{\hat{l}_\nu} \cdot 2\gamma \to 0,\] from 
\[\frac{1}{(1-x)^2}-1=\frac{x(2-x)}{(1-x)^2} \leq \frac{2x}{(1/2)^2}=8x\] 
for \(0 \leq x \leq 1/2,\) and Proposition~\ref{prop2}, it follows using Proposition~\ref{prop1} that 
\[e_\nu^TE_\nu e_\nu-\frac{N}{n} = \frac{l_\nu}{\hat{l}_\nu} \cdot (t_\nu^T \mathcal{M}(I-\frac{1}{\hat{l}_\nu}\mathcal{M})^{-2}t_\nu-t_\nu ^T \mathcal{M}t_\nu+t_\nu ^T \mathcal{M}t_\nu-\frac{N}{n})+(\frac{l_\nu}{\hat{l}_\nu}-1) \cdot \frac{N}{n} \xrightarrow[]{a.s.} 0.\]

\vspace{0.5cm}
\par
\(2.\) Consider the second term:
\[2e_\nu^T E_\nu (a_\nu-e_\nu)=2\hat{l}_\nu \sqrt{l_\nu} \sum_{k}{t_k^T\mathcal{M}(\hat{l}_\nu I-\mathcal{M})^{-2}t_\nu \cdot \sqrt{l_k}(a_\nu-e_\nu)_k}.\]
\par
We have
\begin{equation}\label{911}
    \frac{1}{\hat{l}_\nu}a_\nu^T\Lambda a_\nu-1 \xrightarrow[]{a.s.} 0:
\end{equation}
(\ref{4}) implies
\[1-\frac{1}{\hat{l}_\nu}a_\nu^TS_{AA} a_\nu= \frac{1}{\hat{l}_\nu}a^T_\nu \Lambda^{1/2}T^T\mathcal{M}(\hat{l}_\nu I- \mathcal{M})^{-1}T\Lambda^{1/2} a_\nu,\]
from which 
\[1-\frac{1}{\hat{l}_\nu}a_\nu^TS_{AA} a_\nu \leq a^T_\nu \Lambda^{1/2}T^T\mathcal{M}(\hat{l}_\nu I- \mathcal{M})^{-2}T\Lambda^{1/2} a_\nu=\frac{R_\nu^2}{1-R_\nu^2}  \xrightarrow[]{a.s.} 0,\]
\[1-\frac{1}{\hat{l}_\nu}a_\nu^TS_{AA} a_\nu \geq \frac{1}{2} \cdot a^T_\nu \Lambda^{1/2}T^T\mathcal{M}(\hat{l}_\nu I- \mathcal{M})^{-2}T\Lambda^{1/2} a_\nu=\frac{R_\nu^2}{2(1-R_\nu^2)} \xrightarrow[]{a.s.} 0,\]
because 
\[\frac{x}{2(l-x)^2} \leq \frac{x}{l(l-x)} \leq \frac{x}{(l-x)^2},\] 
for \(l>0, x \in [0,\frac{l}{2}],\)
yielding
\[1-\frac{1}{\hat{l}_\nu}a_\nu^TS_{AA} a_\nu \xrightarrow[]{a.s.} 0;\]
lastly, notice that
\[\frac{1}{\lambda_{\max}(\frac{1}{n}Z_AZ_A^T)} \cdot \frac{1}{\hat{l}_\nu}a_\nu^TS_{AA}a_\nu \leq \frac{1}{\hat{l}_\nu}a_\nu^T\Lambda a_\nu \leq  \frac{1}{\lambda_{\min}(\frac{1}{n}Z_AZ_A^T)} \cdot \frac{1}{\hat{l}_\nu}a_\nu^TS_{AA} a_\nu,\]
and \(\lambda_{\min}(\frac{1}{n}Z_AZ_A^T),\lambda_{\max}(\frac{1}{n}Z_AZ_A^T) \xrightarrow[]{a.s.} 1\) from (\ref{7}).
\par
(\ref{911}) can be rewritten as
\begin{equation}\label{15}
    \frac{1}{\hat{l}_\nu}\sum_{k}{l_k((a_\nu-e_\nu)_k)^2} \xrightarrow[]{a.s.} 0
\end{equation}
because 
\[\frac{l_\nu}{\hat{l}_\nu}(<a_\nu,e_\nu>^2-1)-\frac{l_\nu}{\hat{l}_\nu}(<a_\nu,e_\nu>-1)^2=2 \cdot \frac{l_\nu}{\hat{l}_\nu}(<a_\nu,e_\nu>-1) \xrightarrow[]{a.s.} 0\]
from Lemma~\ref{lemma2} and Proposition~\ref{prop1}. Returning to our term, Cauchy-Schwarz inequality entails
\[2|e_\nu^T E_\nu (a_\nu-e_\nu)| \leq 2\hat{l}_\nu \sqrt{l_\nu} \cdot \sqrt{\sum_{k}{(t_k^T\mathcal{M}(\hat{l}_\nu I-\mathcal{M})^{-2}t_\nu)^2} \cdot \sum_{k}{l_k((a_\nu-e_\nu)_k)^2}},\]
and in light of (\ref{15}) and Proposition~\ref{prop1}, it suffices to show that 
\[l_\nu^4\sum_{k}{(t_k^T\mathcal{M}(\hat{l}_\nu I-\mathcal{M})^{-2}t_\nu)^2}\]
is almost surely bounded to conclude that \(e_\nu^T E_\nu (a_\nu-e_\nu) \xrightarrow[]{a.s.} 0:\) from part \((a)\) of Lemma~\ref{lemma1}, 
\[l_\nu^4 \sum_{k \ne \nu}{(t_k^T\mathcal{M}(\hat{l}_\nu I-\mathcal{M})^{-2}t_\nu)^2} \xrightarrow[]{a.s.} 0,\]
and with probability one for \(n\) sufficiently large, from Proposition~\ref{prop1} and (\ref{913}),
\[0 \leq l_\nu^2t_\nu^T\mathcal{M}(\hat{l}_\nu I-\mathcal{M})^{-2}t_\nu=(\frac{l_\nu}{\hat{l}_\nu})^2 \cdot \hat{l}_\nu^2t_\nu^T\mathcal{M}(\hat{l}_\nu I-\mathcal{M})^{-2}t_\nu \leq 2 \cdot 4c_{K,\gamma} \cdot t_\nu^T\mathcal{M}t_\nu \leq 8c_{K,\gamma} \cdot 2\gamma,\]
as \(\frac{xl^2}{(l-x)^2} \leq 4c_{K,\gamma} \cdot x\) for \(x \in [0,\frac{l}{2}], x \leq c_{K,\gamma},\) and \(l_\nu \to \infty.\)

\vspace{0.3cm}
\par
\(3.\) Consider the third term:
\[0 \leq (a_\nu-e_\nu)^TE_\nu(a_\nu-e_\nu) \leq \hat{l}_\nu ||\mathcal{M}(\hat{l}_\nu I-\mathcal{M})^{-2}|| \cdot ||\Lambda^{1/2}(a_\nu-e_\nu)||^2 \leq \frac{4c_{K,\gamma} ||\Lambda^{1/2}(a_\nu-e_\nu)||^2}{\hat{l}_\nu}\]
almost surely as \(\frac{x}{(l-x)^2} \leq \frac{4c_{K,\gamma}}{l^2}\) for \(x \in [0,\frac{l}{2}], x \leq c_{K,\gamma}.\) Hence this term tends to zero almost surely from (\ref{15}) and Proposition~\ref{prop1}.

\begin{appendix}
The following concentration inequalities will be used repeatedly throughout the proofs (\(c>0\) is a universal constant independent of any of the parameters):

\vspace{0.3cm}
\textit{Two-Sided Bound on Subgaussian Matrices} (Vershynin~\cite{vershynin}, theorem \(4.6.1\))
\par
For any random matrix \(A \in \mathbb{R}^{p \times q}\) whose entries are independent, of mean zero, subgaussian with \(\max_{1 \leq i \leq p,1 \leq j \leq q}{||A_{ij}||_{\psi_2}} \leq K,\) and all \(t \geq 0,\)
\begin{equation}\label{7}\tag{SM}
    \sqrt{p}-cK^2(\sqrt{q}+t) \leq s_q(A) \leq s_1(A) \leq \sqrt{p}+cK^2(\sqrt{q}+t),
\end{equation}
with probability at least \(1-2\exp(-t^2),\) where \(s_i(\cdot)\) is the \(i^{th}\) largest singular value.

\vspace{0.3cm}
\textit{Hanson-Wright Inequality} (Vershynin~\cite{vershynin}, theorem \(6.2.1\)): 
\par
For \(y \in \mathbb{R}^{p}\) a random vector with independent entries of mean zero, \(\max_{1 \leq i \leq p}{||y_i||_{\psi_2}} \leq K, C \in \mathbb{R}^{p \times p},\) and all \(t \geq 0,\)
\begin{equation}\label{9}\tag{HW}
    \mathbb{P}(|y^T C y - \mathbb{E}[y^T C y]| \geq t) \leq 2\exp(-c\min{(\frac{t^2}{pK^4||C||^2},\frac{t}{K^2||C||})}).
\end{equation}

\vspace{0.3cm}
\textit{Analogue of Hanson-Wright Inequality:} 
\par
For \(y,y' \in \mathbb{R}^{p}\) random vectors with entries \(y_1, y_2,\hspace{0.05cm} ... \hspace{0.05cm}, y_p, y'_1, y'_2, \hspace{0.05cm} ... \hspace{0.05cm}, y'_p\) independent, of mean zero, \(\newline \max_{1 \leq i \leq p}{||y_i||_{\psi_2}}, \max_{1 \leq i \leq p}{||y'_i||_{\psi_2}} \leq K,\) and all \(t \geq 0,\)
\begin{equation}\label{10}\tag{aHW}
    \mathbb{P}(|y^T C y'| \geq t) \leq 2\exp(-c\min{(\frac{t^2}{pK^4||C||^2},\frac{t}{K^2||C||})}).
\end{equation}

\par
\vspace{0.2cm}
\textit{Proof:} This result can be derived from the proof of its counterpart as follows: lemma \(6.2.3\) from Vershynin~\cite{vershynin} gives that for \(y,y'\) as above, \(g,g' \in \mathbb{R}^p\) independent Gaussian random vectors of mean zero, covariance the identity matrix, and all \(\lambda \in \mathbb{R},\)
\[\mathbb{E}[\exp(\lambda y^TCy')] \leq \mathbb{E}[\exp(c_1K^2\lambda g^T Cg')].\]
\par
Take \(\sigma_1,\sigma_2, \hspace{0.05cm} ... \hspace{0.05cm}, \sigma_p\) to be the eigenvalues of \(C.\) Then by rotational invariance,
\[\mathbb{E}[\exp(c_1K^2\lambda g^T Cg')]=\mathbb{E}[\exp(c_1K^2\lambda g^T diag(\sigma_1,\sigma_2,\hspace{0.05cm} ...\hspace{0.05cm}, \sigma_p) g')]=\]
\[=\prod_{1 \leq i \leq p}{\mathbb{E}[\exp(c_1K^2\lambda \sigma_i g_i g_i')]} =\prod_{1 \leq i \leq p}{\mathbb{E}[\mathbb{E}[\exp(c_1K^2\lambda \sigma_i g_i g_i')|g_i]]}=\]
\[=\prod_{1 \leq i \leq p}{\mathbb{E}[\exp(c_1^2K^4\lambda^2 \sigma_i^2 g_i^2)]}  \leq \exp(pc_2||C||^2 K^4 \lambda^2)\] 
for \(|\lambda| \leq \frac{c_3}{c_1K^2||C||},\) from \(\mathbb{E}[\exp(\lambda g_i)]=\exp(\lambda^2), \mathbb{E}[\exp(\tau^2X^2)] \leq \exp(c_4\tau^2)\) for \(||X||_{\psi_2} \leq ||g_1||_{\psi_2}\) and \(|\tau| \leq c_3\) (this is property \((iii)\) of proposition \(2.5.2\) in Vershynin~\cite{vershynin}). 
\par
Hence, for \(|\lambda| \leq \frac{c_3}{c_1K^2||C||},\)
\[\mathbb{E}[\exp(\lambda y^TCy')] \leq \exp(pc_2||C||^2 K^4 \lambda^2),\]
which implies (arguing as for Bernstein's inequality) that for some universal constant \(c>0\) and all \(t \geq 0,\)
\[\mathbb{P}(|y^T C y'| \geq t) \leq 2\exp(-c\min{(\frac{t^2}{pK^4||C||^2},\frac{t}{K^2||C||})}).\]
\end{appendix}

\end{document}